\documentclass[11pt]{article}

\usepackage[all]{xy}

\newtheorem{tm}{Theorem}[subsection]

\newtheorem{pr}[tm]{Proposition}
\newtheorem{rmk}[tm]{Remark}
\newtheorem{cor}[tm]{Corollary}
\newtheorem{ex}[tm]{Example}

\newtheorem{??}[tm]{Question}
\newtheorem{defi}[tm]{Definition}

\font\tenmsb=msbm10
\font\sevenmsb=msbm7
\font\fivemsb=msbm5
    
\newfam\msbfam
\textfont\msbfam=\tenmsb
\scriptfont\msbfam=\sevenmsb
\scriptscriptfont\msbfam=\fivemsb
\def\Bbb#1{{\fam\msbfam #1}}

\font\teneufm=eufm10
\font\seveneufm=eufm7
\font\fiveeufm=eufm5
\newfam\eufmfam
\textfont\eufmfam=\teneufm
\scriptfont\eufmfam=\seveneufm
\scriptscriptfont\eufmfam=\fiveeufm
\def\frak#1{{\fam\eufmfam\relax#1}}

\def\lorw{\longrightarrow}
\newcommand\n{\noindent}

\newcommand\ci{\cite}

\newcommand\rat{{\Bbb Q}}
\newcommand\comp{{\Bbb C}}
\newcommand\real{{\Bbb R}}
\newcommand\zed{{\Bbb Z}}
\newcommand\nat{{\Bbb N}}
\newcommand\pn[1]{{\Bbb P}^{#1}}

\newcommand\blacksquare{{\hspace*{\fill} $\fbox{}$}}

\newcommand\e{\epsilon}

\newcommand{\coke}{ \hbox{\rm Coker} }

\newcommand{\im}{ \hbox{\rm Im} }

\newcommand{\ke}{ \hbox{\rm Ker} }

\newcommand{\ptd}[1]{ \,^{\frak p}\!\tau_{ \leq {#1} } }

\newcommand{\ptu}[1]{ \,^{\frak p}\!\tau_{ \geq {#1} } }

\newcommand{\td}[1]{ \tau_{ \leq {#1} } }

\newcommand{\pe}{ {\cal P }  }

\newcommand{\pc}[2]{ \,^{\frak p}\!{\cal H}^{#1}({#2})   }
\newcommand{\tu}[1]{ \tau_{ \geq {#1} } }

\newcommand{\ql}{\rat_{l}}
\newcommand{\qlb}{ {\overline{\rat}}_{l} }

\newcommand{\fb}{\Bbb F}
\newcommand{\fq}{{\Bbb F}_q}

\title{The decomposition theorem, \\ perverse sheaves \\
and the topology of algebraic maps\\
} 
\author{
Mark Andrea A.  de Cataldo
\, 
and Luca Migliorini\thanks{Partially supported by GNSAGA and PRIN 2007 project ``Spazi di moduli e teoria di Lie''}
}
\date{}

\begin{document}\maketitle

\begin{abstract}
We give a motivated introduction to the theory of perverse sheaves, 
culminating in the decomposition theorem of  Beilinson, Bernstein, Deligne and Gabber.
A goal of this survey is to show how the theory develops 
naturally from classical constructions used in the study of topological properties 
of algebraic varieties. 
While most proofs are omitted, we discuss several approaches to the decomposition 
theorem, indicate some important applications and examples.  
\end{abstract}

\tableofcontents

\section{Overview}
\label{secintro}
The theory of perverse sheaves and one of its crowning achievements,
the decomposition theorem,
are at the heart of a revolution  which has taken  place
over the last thirty years in algebra, representation
theory and algebraic geometry.

The decomposition theorem is a powerful tool for  investigating the topological properties
 of proper maps between algebraic varieties and is
 the deepest  known fact  relating their homological, Hodge-theoretic  and 
  arithmetic properties.  
 
 In  this $\S$\ref{secintro},  we try to motivate the statement
of this  theorem as a natural
 outgrowth of the  investigations on the topological properties of algebraic
 varieties begun with  Lefschetz and culminated
 in the spectacular results obtained with the development of Hodge theory and
 \'etale cohomology.  
  We  gloss over many crucial technical details  in favor
 of rendering a more panoramic picture; the  appendices in $\S$\ref{cemetto}
 offer a partial remedy to these omissions.
 We state the  classical Lefschetz and Hodge  theorems for projective manifolds in \S\ref{tcav}
 and Deligne's results on families of projective manifolds
 in \S\ref{rel2}.  In \S\ref{sav}, we
briefly discuss singular varieties and the appearance  and role of 
mixed Hodge structures and intersection cohomology.
 In  $\S$\ref{rel3}, we state the
 decomposition theorem in terms of intersection
cohomology  without any reference
 to perverse sheaves. 
The known proofs, however, use in an essential way
the theory of perverse sheaves which, in turn,   is deeply rooted
in the formalism of sheaves and derived categories. 
We offer a   ``crash course" on sheaves  in $\S$\ref{subsec-crash-course}.
With these notions and ideas in hand,  in \S\ref{statem} we  state
the decomposition theorem in terms of intersection complexes
(rather than in terms of intersection cohomology groups).   We also
state two important related results: the  
 relative hard Lefschetz and semisimplicity theorems. 
\S\ref{ictse} reviews the generalization to singular maps  of the
now classical  properties of the monodromy representation in cohomology
for  a family of projective manifolds. 
\S\ref{afex} discusses surface and threefold
examples of the statement of the decomposition theorem.
\S\ref{dtmHs} overviews the mixed Hodge structures 
arising from the decomposition theorem.
  We provide a timeline for the main results mentioned
 in this overview in $\S$\ref{hire}.

We have tried, and have surely failed in some ways, to write this survey so that most of it can be read  by non experts
and so that each chapter can be read independently of the others. For example, a reader
interested in the decomposition theorem and in its applications
could read $\S$\ref{secintro},  the first half of $\S$\ref{dtappl}
and skim through the second half on geometrization, while
 a reader interested in the proofs
could read $\S$\ref{secintro} and $\S$\ref{3stooges}. Perhaps, at that point,
the reader may be motivated to read more about perverse sheaves.
 
   $\S$\ref{secptsdy} is an introduction to  perverse sheaves. In this survey,
   we deal only with middle perversity, i.e. with a special case of perverse sheaves.
   It seemed natural to us to start this section with a discussion of intersection cohomology.
  In \S\ref{subsecoervsh}, we  define perverse sheaves, discuss their first properties,
   as well as  their natural categorical framework, i.e.
$t$-structures. In  \S\ref{tpfmc},  we introduce the perverse filtration in cohomology
   and its geometric description via the Lefschetz hyperplane theorem.
   \S\ref{subsecpervco} reviews the basic properties
 of the cohomology functors associated with  the perverse $t$-structure.
 \S\ref{subsectex} is about the Lefschetz hyperplane theorem for intersection
 cohomology.  In \S\ref{subsectie}, we  review the properties of the  intermediate extension
 functor, of which
 intersection complexes are a key example.

 In $\S$\ref{3stooges}, we discuss the  three  known  approaches to
 the decomposition theorem: the original one, due to A. Beilinson, J. Bernstein,
 P. Deligne and O. Gabber, via the arithmetic
 properties of varieties over finite fields, the one of M. Saito, via mixed Hodge
 modules, and ours, via classical Hodge theory. Each approach
 highlights  different aspects of this important theorem.

$\S$\ref{dtappl} contains 
  a sampling of   applications of the theory of perverse sheaves and,
 in particular, of the decomposition theorem.
 The applications
  range from algebraic geometry to 
representation theory and  to combinatorics. 
While the first  half of \S\ref{dtappl}, on toric and on semismall maps,  is  targeted to a general audience,  the second half,  on the  geometrization of Hecke algebras  and
of the Satake isomorphism,  is technically  more demanding. 
Due to the fact that the  recent and  exciting development \ci{ngo}  in the Langlands
program makes use of a result  that deals with the decomposition
theorem with ``large fibers," we have included a brief discussion
of B.C. Ng\^o's support theorem in \S\ref{nst}.

The appendix \S\ref{cemetto} contains a brief definition of quasi projective varieties
(\S\ref{algvar}),
of  pure and mixed Hodge structures, the statement of the hard Lefschetz theorem 
and of the Hodge-Riemann relations (\S\ref{pam}), a description of the formalism
of derived categories (\S\ref{tfidy}),  a discussion of how the more classical objects
in algebraic topology relate to this formalism (\S\ref{famfa}), a discussion
of the nearby and vanishing cycle functors (\S\ref{psiphi}), as well as their
unipotent counterparts (\S\ref{unipotpsiphi}), two descriptions
of the category of perverse sheaves (\S\ref{structpervshvs}) and, finally,
a formulary for the derived category (\S\ref{formulary}).


Unless otherwise stated,   a variety is an irreducible
complex algebraic variety and a map is a map of varieties.
We work  with  sheaves of rational vector spaces,
 so that the cohomology groups are rational vector spaces.

\smallskip
{\bf Acknowledgments.}
 We thank Pierre Deligne, Mark Goresky, Tom Haines,
Andrea Maffei and
 Laurentiu Maxim for many conversations.
We thank I.A.S. Princeton for the great working conditions
in the a.y. 2006/2007. 
During the preparation of this paper the second-named author has
also been guest 
of the following institutions: I.C.T.P. in Trieste, Centro di Ricerca 
Matematica E. de Giorgi in Pisa.

\subsection{The  topology of complex projective manifolds: Lefschetz
and Hodge theorems}
\label{tcav}
Complex algebraic varieties provided an important motivation
for the development of algebraic topology from its earliest
days. On the other hand, algebraic varieties
and algebraic maps enjoy many truly remarkable topological properties
that are not shared by other classes of spaces and maps. These special features
were first exploited by Lefschetz (\ci{lef})
(who claimed to have ``planted the harpoon of algebraic topology into the body
of the whale of algebraic geometry" \ci{Lefschetz}, p.13)
and they are almost completely summed up in the statement of the decomposition theorem and 
of its embellishments. 

The classical precursors to the decomposition theorem include
the theorems of Lefschetz, Hodge, Deligne, and the invariant cycle theorems. 
In the next few paragraphs, we discuss the Lefschetz and Hodge  theorems
and 
the Hodge-Riemann relations.
Together with  Deligne's Theorem \ref{delteo},
these  precursors  are in fact  
essential tools in  the three  known proofs ($\S$\ref{3stooges})
of the decomposition theorem.

 Let $X$ be a   nonsingular complex $n$-dimensional   projective 
 variety embedded in some projective space $X \subseteq 
 \pn{N}$, and let $D= {\cal H} \cap X$ be the intersection of $X$ with a generic
 hyperplane ${\cal H}\subseteq \pn{N}$. Recall that   we use  cohomology with
 rational coefficients.  A standard textbook reference for what follows
 is  \ci{gh}; see also \ci{voisinbook,decbook}.

The {\em  Lefschetz hyperplane  theorem}
states that the restriction map $H^i(X) \to H^i(D)$ is an isomorphism
for $i <n-1$ and is injective for $i=n-1$.

The cup product with the first Chern class of the hyperplane
bundle gives a mapping $\cup c_1({\cal H}): H^i(X) \to H^{i+2}(X)$
which can be identified with the composition
$H^i(X) \to H^i(D) \to  H^{i+2}(X)$, the latter being a ``Gysin" homomorphism.  

The
{\em  hard Lefschetz theorem} states that for all $0 \leq i \leq n$ the $i$-fold
iteration of the cup product operation  is an isomorphism
$$
(\cup c_1({\cal H}))^i: H^{n-i}(X) \stackrel{\simeq}\lorw H^{n+i}(X).
$$

The {\em Hodge decomposition}
is a canonical decomposition 
\[H^i(X,\comp) = \bigoplus_{p+q=i} H^{p,q}(X).\]
The summand $H^{p,q}(X)$ consists of  cohomology
classes  on $X$ which can be represented by a closed    differential form
on $X$ of type $(p,q)$
(i.e. one whose local expression involves $p$ $dz$'s and $q$ $d \overline{z}$'s).

For every fixed index $0 \leq i \leq n$, 
define a bilinear form $S^{\cal H}$ on
$H^{n-i}(X)$ by 
$$(a,b)  \longmapsto S^{\cal  H}(a,b):=  \int_X{(c_1({\cal H}))^i\wedge a \wedge b} = 
\deg{([X] \cap ( (c_1({\cal H}))^i \cup a \cup b))},$$
where $[X]$ denotes the fundamental homology class of the naturally oriented $X$.

The  hard Lefschetz theorem is equivalent to the nondegeneracy of the forms $S^{\cal H}$.

The {\em Hodge-Riemann bilinear relations} ($\S$\ref{pam}, (\ref{hrbr}))  establish their signature properties.

\subsection{Families of smooth projective varieties}
\label{rel2}

 If $f:X \to Y$ is a  $C^{\infty}$ fiber bundle with smooth compact fiber $F$, let
$\underline{H}^j(F)$ denote the local  system on $Y$ whose fiber 
at the point $y\in Y$ is $H^j(f^{-1}(y))$.    There are the associated {\em Leray spectral sequence}
\begin{equation}
\label{eqn-Leray} 
E_2^{i,j} =H^i(Y; \underline{H}^j(F)) \Longrightarrow H^{i+j}(X)
\end{equation}
and   the {\em monodromy representation}
\begin{equation}\label{eqn-rep}
 \rho_i:\pi_1(Y,y_0) \to GL(H^i(F)).\end{equation}
 
Even if $Y$ is simply connected, the Leray spectral sequence can be
nontrivial, for example, the Hopf fibration $f:S^3 \to S^2$. 

We define a {\em family of projective manifolds} to be
  a  proper holomorphic submersion  $f: X \to Y$ of nonsingular varieties
   that 
 factors through some product $Y \times \pn{N}$ and for which the fibers are connected
 projective manifolds. The nonsingular hypersurfaces of a fixed degree
 in some projective space give an interesting example.
By a  classical result of Ehresmann, 
such a map is also a  $C^{\infty}$ fiber bundle.

The results that follow  are due to
Deligne \ci{dess, ho2}. Recall that a representation is said to be {\em
irreducible} if it does not admit a non trivial   invariant subspace, i.e. if it is simple
in the category of representations. 
\begin{tm}  
\label{delteo}
{\bf (Decomposition and semisimplicity for families of projective
manifolds)}
Suppose $f:X \to Y$ is a  family of projective manifolds.
Then
\begin{enumerate}
\item The Leray spectral sequence \mbox{\rm(\ref{eqn-Leray})} degenerates at the $E_2$-page and induces an isomorphism \[H^i(X) \cong \bigoplus_{a+b=i}H^a(Y;
\underline{H}^b(F)).\]
\item The representation \mbox{\rm (\ref{eqn-rep})} is semisimple:  it is a direct sum of
irreducible representations.
\end{enumerate}\end{tm}

Part 1. gives a rather complete description of the cohomology of $X$.    Part 2. is remarkable 
because the fundamental group of $Y$ can be infinite.

\begin{rmk}
\label{notcanonical}
{\rm
Theorem \ref{delteo}, part 1 is stated using cohomology. Deligne proved
a stronger, sheaf-theoretic statement; see Theorem \ref{dss}.
}
\end{rmk}

\begin{rmk}
\label{ex17}
{ \rm 
For singular maps, the Leray spectral sequence is very seldom degenerate.
If $f: X \to Y$ is a  resolution of the singularities
of a projective variety $Y$ whose cohomology
has a mixed Hodge structure which is not pure, then
$f^*$ cannot be injective, and this prohibits degeneration in view of the
edge-sequence.}
\end{rmk}

The following is the global invariant cycle theorem.
We shall come back to this later in $\S$\ref{ictse}, where we give
some generalizations, and in $\S$\ref{hire}, where we give
some references.

\begin{tm}
\label{gict00}
Suppose $f: X \to Y$ is a family of projective manifolds.
Then 
$$
H^i(F_{y_0})^{\pi_1(Y, y_0)} = \im \, \{ H^i(X) \lorw 
H^i(F_{y_0}) \},
$$
i.e. the monodromy invariants are precisely the classes obtained
by  restriction
from the total space of the family.
\end{tm}

Although the classical Lefschetz-Hodge theorems described in
$\S$\ref{tcav} and the results described in this section
appear to be very different from each other, the decomposition
theorem forms a 
beautiful common generalization which holds also in the presence
of singularities.

\subsection{Singular algebraic varieties}
\label{sav}
The Lefschetz  and Hodge theorems fail if $X$ is singular. There
are two 
somewhat complementary approaches to generalize these statements to singular projective
varieties. They  involve mixed Hodge theory \ci{ho2, ho3} and intersection
cohomology \ci{goma1, goma2}  (see  also \ci{borel}). 
 
In {\em mixed Hodge theory} the topological invariant studied is the same
 investigated for nonsingular varieties, namely, the cohomology groups
 of the variety, whereas the structure with which it is endowed changes.
 See \ci{durfee} for an elementary and nice introduction.
 The $(p,q)$-decomposition of classical Hodge theory is replaced
 by a more complicated structure: 
the rational cohomology groups $H^i(X)$ are endowed with an increasing filtration
$W$ (the weight filtration)
$W_0 \subseteq W_1 \subseteq \ldots \subseteq W_{2i} = H^{i}(X)$,
 and the  complexifications of the graded pieces $W_k/W_{k-1}$ have a $(p,q)$-decomposition 
 of weight $k$, that is $p+q=k$.
Such a structure, called a {\em mixed Hodge structure}, exists canonically
 on any algebraic variety and
 satisfies several fundamental restrictions on the weights, such as:

\begin{enumerate}
\item
if  $X$  is nonsingular, but  possibly non-compact,  then the weight filtration 
on $H^i(X)$ {\em starts} at $W_i$, that is
$W_rH^i(X)=0$ for $r<i;$

\item
 if $X$ is compact, but possibly   singular, then the  weight filtration on $H^i(X)$ {\em ends}
 at $W_i$, that is
$W_rH^i(X)=W_iH^i(X)=H^i(X)$ for $r\geq i$.
\end{enumerate}
\begin{ex}
\label{exmHs}
{\rm Let $X=\comp^*$; then $H^1(X)\simeq \rat$ has weight $2$
and  the classes in  $H^1(X)$
are of type $(1,1)$.
Let $X$ be a rational irreducible curve with a node
(topologically, this is a pinched torus, or also the  two-sphere
with the north and south poles identified);
then $H^1(X)
\simeq \rat$ has weight $0$ and
the classes in  $H^1(X)$ are of type $(0,0)$.}
\end{ex}

In  {\em intersection cohomology theory}, by contrast,  it is the topological 
invariant which is changed, 
whereas the $(p,q)$-structure turns out to be the same.
The {\em intersection cohomology groups} $I\!H^i(X)$ ($\S$\ref{subsecintcoh})
can be described using geometric ``cycles" on the possibly singular variety $X$,
and this gives a concrete way to compute simple examples. There is a natural
homomorphism $H^i(X) \to I\!H^i(X)$ which is an isomorphism
when $X$ is nonsingular. The groups $I\!H^i(X)$
are finite dimensional; they satisfy 
the Mayer-Vietoris theorem and    the K\"unneth Formula.
These groups are not homotopy invariant but, in compensation,
they have the following additional features:
they satisfy
 Poincar\'e duality,    the Lefschetz  theorems and, if $X$ is projective, they admit a 
 pure  Hodge structure.

\begin{ex}
\label{nodal}
{\rm Let $X$ be the nodal curve of Example \ref{exmHs}.
Then $I\!H^1(X)=0$.}
\end{ex}
\begin{ex}
\label{iccone}
{\rm   Let $E \subseteq \pn{N}_\comp$ be a nonsingular projective variety of dimension $n-1$,
and let $Y\subseteq \comp^{N+1}$ be its affine cone
with vertex  $o$. 
The intersection cohomology groups can be easily computed (see
\ci{borel} and also Example \ref{exicd}): 
$$
I\!H^i(Y)=0 \hbox{ for } i\geq n \qquad I\!H^i(Y)  \simeq H^i(Y \setminus \{o\}) \hbox{ for }i<n.
$$

}
\end{ex}

There is a  twisted version of intersection (co)homology with values in a 
local system $L$ defined on a Zariski dense nonsingular  open subset of the variety
$X$.  
Intersection cohomology with twisted coefficients is denoted 
 $I\!H^*(X,L)$ and it 
appears in the  statement of the decomposition theorem.

\subsection{Decomposition   and  hard Lefschetz 
in intersection cohomology
}
\label{rel3}
The decomposition theorem is a result
about certain  complexes of sheaves on varieties.
In this section, we state a provisional, yet suggestive form 
that involves only  intersection cohomology groups.

\begin{tm}
\label{dtfirstform} 
{\rm ({\bf Decomposition theorem for intersection cohomology groups})}
Let $f:X \to Y$ be a proper map of  varieties. 
There exist finitely many  triples $(Y_a, L_a, d_a)$
made of locally closed,  smooth and  irreducible algebraic  subvarieties  $Y_a\subseteq Y$, 
semisimple local systems $L_a$ on $Y_a$
and integer numbers $d_a$, such that for every 
 open set $U \subseteq
Y$ there is an isomorphism
\begin{equation}
\label {dtnaive}
I\!H^r(f^{-1}U)\simeq \bigoplus_a I\!H^{r-d_a}( U \cap\overline{Y}_a,L_a).
\end{equation}
\end{tm}

The triples $(Y_a, L_a, d_a)$ are essentially unique,
independent of $U$,  and they are described
in \ci{herdlef, decmightam}.
Setting $U= Y$  we get a formula for $I\!H^*(X)$  and
therefore, if $X$ is nonsingular, a formula for $H^*(X)$.
If  $f:X\to Y$ is a family of projective manifolds, 
then (\ref{dtnaive}) coincides with  the decomposition
in Theorem \ref{delteo}, part 1. On the opposite side of the spectrum,
if $f: X \to Y$ is a resolution of the singularities of $Y$, i.e.
$X$ is nonsingular and 
$f$ is an isomorphism outside a closed subvariety of $Y$, 
then we can  deduce that the intersection cohomology groups
$I\! H^*(Y)$ are  direct summands of $H^*(X)$.

If $X$ is singular, then
there is no analogous  direct sum decomposition formula for $H^*(X)$. 
Intersection cohomology turns out to be precisely the topological invariant
apt to deal with singular varieties and maps.
The  notion of intersection cohomology is needed 
even when  $X$ and $Y$  are nonsingular, but the map $f$  is not a submersion.

\begin{rmk}
\label{itisnotcanonical}
{\rm
{\bf (The splitting is not canonical)}
The decomposition map (\ref{dtnaive}) is not uniquely defined. This is analogous to the 
elementary fact that a  filtration on a vector space can always be given
in terms of a direct sum decomposition, but the filtration does not determine
in a natural way
the summands as subspaces of the given vector space. In the case when $X$ is quasi projective,
one can make distinguished choices  which realize the summands as mixed
Hodge substructures of a canonical mixed Hodge structure on 
$I\!H^*(X)$ 
 (see \ci{decmigseattle,decII} and $\S$\ref{dtmHs}, 5).
 }
 \end{rmk}

If ${\cal L}$ is a hyperplane line bundle on  a  projective variety  $Y$,
then   the hard Lefschetz theorem 
for the intersection cohomology groups
of $Y$ holds, i.e.
for every  integer $k \geq 0$, the  $i$-th iterated cup product
\begin{equation}
\label{00hlic}
c_1({\cal L})^i: I\!H^{\dim{Y}-i}({Y}) 
\stackrel{\simeq}\lorw
I\!H^{\dim{Y}+i}(Y)
\end{equation} is an isomorphism.
Recall that intersection cohomology is not  a ring, however, the cup
product with a cohomology class is well-defined and  intersection cohomology
is a module over cohomology.

The 
analogue of Theorem \ref{chl}.(3) (hard Lefschetz,  
Lefschetz decomposition and   Hodge-Riemann relations)  holds
for the intersection cohomology groups  $I\!H^*(Y)$
of a  singular projective variety $Y$.

\subsection{Crash course on sheaves and   derived categories}
 \label{subsec-crash-course}
The statement of   Theorem \ref{dtfirstform} involves only  the notion of intersection cohomology.
We do not know of a  general method for  
proving the decomposition (\ref{dtnaive})
without first proving   the analogous
decomposition,  Theorem \ref{stdteo},
 at the level of complexes of sheaves.

The language and theory of sheaves and 
 homological algebra,  specifically derived categories
and
perverse sheaves,
plays an essential role in all  the known proofs of the decomposition theorem, as well as in its numerous applications.

In this section, we collect the few facts about sheaves and 
 derived categories needed in order to
 understand the statement of the decomposition Theorem \ref{stdteo}. We amplify
 and complement this crash course in the appendices in  $\S$\ref{cemetto} and  in 
 section $\S$\ref{secptsdy} on perverse sheaves. 
Standard references are \ci{borel, gel-man, goma2,  k-s, iv}.

{\bf 1. Complexes of sheaves.} 
Most of the constructions in homological algebra
involve complexes. For example, if $Z$ is a ${\cal C}^{\infty}$ manifold, in order to 
compute the cohomology of the constant sheaf $\real_Z$, we replace it
by the complex of sheaves  of  differential forms, and  then  we take the complex 
of  global sections, i.e. the de Rham complex.
More generally, to define the cohomology of a sheaf $A$ on a 
topological space $Z$, we choose an injective, or flabby, resolution, for instance the one defined by Godement, 
$$
\xymatrix{
0 \ar[r] & A  \ar[r]  &  I^{0} \ar[r]^{d^{0}} & \ldots   \ar[r]^{d^{i-1}}  &   I^i  \ar[r]^{d^i} &
I^{i+1} \ar[r] &  \ldots} 
$$
then consider the complex of abelian groups
$$
\xymatrix{
 0 \ar[r] & \Gamma( I^{0}) \ar[r]^{d^{0}} & \ldots  \ar[r]^{d^{i-1}}  &  \Gamma( I^i)  \ar[r]^{d^i} &
\Gamma(I^{i+1}) \ar[r] &  \ldots} 
$$
and finally take its cohomology.
The derived category is a formalism developed in order  
to work systematically with complexes 
of sheaves with a notion of morphism
which is far more flexible than that of morphism of complexes; 
for instance, 
two different resolutions of the same sheaf are isomorphic 
in the derived category.
Let $Z$ be a topological space.  
We consider sheaves
of $\rat$-vector spaces on $Z$.   
A  bounded {\em complex of sheaves}
$K$ is a diagram
$$
\xymatrix{
\ldots \ar[r]  &  K^{i-1} \ar[r]^{d^{i-1}} &   K^i  \ar[r]^{d^i} &
K^{i+1} \ar[r] &  \ldots} 
$$
with $K^i=0$ for $|i| \gg 0$ and
satisfying $d^i\circ d^{i-1} = 0$ for every $i$.   The
{\em shifted} complex $K[n]$ is the complex with $K[n]^i=K^{n+i}$
and differentials $d_{K[n]}= (-1)^n d_K$.
Complexes of sheaves form an Abelian category and we may form 
the {\em cohomology sheaf}
${\cal H}^i(K) = \ke (d^i)/ \im (d^{i-1})$ which is a sheaf whose stalk
at a point  $x\in Z$ is the   cohomology of the complex of stalks at $x$.

{\bf 2. Quasi-isomorphisms and resolutions. }
A morphism $K \to L$ of complexes of sheaves is a {\em 
quasi-isomorphism} if it induces isomorphisms ${\cal H}^i(K)\cong
{\cal H}^i(L)$ of cohomology sheaves, i.e.  if the induced
map at the level of the  stalks of the cohomology sheaves
is an isomorphism at each point $z \in Z$.
An {\em injective (flabby, fine) resolution} of a complex $K$
is a quasi-isomorphism $K \to I$, where $I$ is  a complex with injective (flabby, fine)
components. Such a resolution always exists
for a bounded below  complex. 
The cohomology  groups $H^*(Z,K)$ of $K$
are  defined to be the cohomology groups of the complex of global sections $\Gamma(I)$ 
of $I$. 
As soon as one identifies sheaves with the complexes of sheaves concentrated in degree $0$, this definition of the groups  $H^*(Z,K)$  extends 
the definition of the  cohomology groups  of a single sheaf
given above to  the case of 
bounded (below) complexes.

A quasi-isomorphism $K \to L$ induces isomorphisms on the 
cohomology, $H^i(U,K) \cong H^i(U,L)$ of any open set $U\subset
Z$ and these isomorphisms are compatible with the maps 
induced by inclusions and with Mayer-Vietoris sequences.

{\bf 3. The derived category.}
The  {\em derived category} $D(Z)$  is a category whose objects
are  the  complexes of sheaves, but whose morphisms 
have been cooked up in such a way that every quasi-isomorphism 
$S \to T$ becomes an isomorphism in $D(Z)$ (i.e. it has a unique 
inverse morphism). In this way, quite different complexes of sheaves
that realize the same cohomology theory (such as the complex 
of singular cochains and the complex of differential forms on a 
${\cal C}^{\infty}$ manifold)  become isomorphic in $D(Z)$.
The definition of the 
morphisms in the derived category is done by first identifying
morphisms of complexes which are 
homotopic to each other, and then  by
formally adding inverses to quasi-isomorphisms.
The second step   is strongly 
reminiscent of the construction of the rational numbers as the field of fractions
of the  ring of integers, and the   necessary  calculus of fractions is made possible
in view of the first step.
There is the analogous notion of {\em bounded derived category}
$D^b(Z)$, where the objects are the  bounded complexes of sheaves.
The bounded derived category sits inside the derived category and the embedding
$D^b(Z) \subseteq D(Z)$ is full. Similarly, for complexes bounded below
(i.e. ${\cal H}^i (K) =0, \forall i \ll 0)$ and the corresponding 
category
$D^+(Z)\subseteq D(Z)$, etc.

{\bf 4. Derived functors.} The main feature of the  derived category is
the possibility of defining  derived functors. We discuss the case of cohomology
and the case of the push-forward via a continuous map.
If $I$ is a  bounded below complex of injective (flabby, or even fine) sheaves
on $Z$, the cohomology $H^i(Z,I)$ is the cohomology of 
the complex of abelian groups
\label{cplxsec}
\[ 
\xymatrix{\ldots \ar[r] &  \Gamma (Z,  I^{i-1}) \ar[r] & \Gamma(Z,  I^{i}) \ar[r] &
\Gamma(Z, 
 I^{i+1}) \ar[r] & \cdots}  \]
which can be considered as an object, denoted $R\Gamma(Z,I)$  of the 
bounded below derived category of a point $D^+(pt)$.
However, if the complex is not injective, as the example of the constant sheaf on a 
${\cal C}^{\infty}$
manifold shows, this procedure gives the wrong answer, as the complexes of global sections of two quasi-isomorphic complexes are not necessarily quasi-isomorphic.
Every bounded  below complex $K$ admits  a bounded   below injective resolution $K\to I$, unique up to a unique isomorphism in $D^+(Z)$.
The complex  of global sections $R\Gamma(Z, K):= \Gamma(Z, I)$ (a flabby resolution can be used  as well and, if there is one,
also a fine one) is well-defined up to unique isomorphism in the derived category
$D^+(pt)\subseteq D(pt)$. For our limited purposes, note that
we always work with bounded complexes whose resolutions can be chosen to be 
bounded, i.e. we can and do  work within  $D^b(Z)$, etc.

A similar construction  arises when $f:W \to Z$ is a continuous mapping:  
if $I$ is a bounded below complex of injective sheaves on $W$, then the push forward 
complex  $f_*(I)$ is a complex of sheaves on $Z$ that satisfies 
\begin{equation}\label{eqn-pushforward}
 H^i(U, f_*(I)) \cong H^i(f^{-1}(U), I)
\end{equation}
for any open set $U \subseteq Z$.  
However if a  bounded below complex $C$  on $W$ is not injective, then
(\ref{eqn-pushforward}) may fail, and $C \in D^+(W)$ should first be replaced 
by an injective resolution before pushing forward.  
The resulting complex of sheaves on $Z$  is well defined up to
canonical isomorphism in $D^+(Z)$, is denoted $Rf_*C$
and is called the (derived) {\em direct image} of $C$. Its cohomology sheaves
are sheaves on $Z$,
are denoted $R^if_* C$ and are called the {\em $i$-th direct image sheaves}.
Note that if $f$ maps $W$ to a point,  then $Rf_*C=R\Gamma(W,C)$ and 
$R^if_* C=H^i(W,C)$. 

When $f: W \to Z$ is a continuous map of  locally compact
spaces, a similar process, that starts
with the functor {\em direct image with proper supports $f_!$},   yields
the  functor  {\em derived direct image with proper supports}   $Rf_!:
D^+(W) \to D^+(Z)$.  There is a map of functors 
$Rf_! \to Rf_* $ which is an isomorphism if $f$ is proper.
Under quite general hypotheses, always satisfied by algebraic maps of
algebraic varieties,
given a map $f: W \to Z$, there are  the {\em inverse image} and 
{\em extraordinary
inverse image}  functors 
$f^*,f^!: D^b(Z) \to D^b(W)$.  See $\S$\ref{tfidy} for a list of
the   properties of these four functors $Rf_*,  Rf_!, f^*$  and $f^!$,
as well as 
for their relation to Verdier duality.

{\bf 5. Constructible sheaves.}  (See \ci{goma2}.) From now on, suppose $Z$ is a complex algebraic variety.
A subset $V \subset Z$ is {\em constructible} if it is
obtained from a finite sequence of unions, intersections,
or complements of algebraic subvarieties of $Z$.
A {\em local system}  on $Z$ is a locally constant sheaf on $Z$ with finite dimensional
stalks. A local system  on $Z$ corresponds to a finite dimensional
representation of the fundamental group
of  $Z$.
A complex of sheaves $K$ has 
 {\em constructible cohomology
sheaves}  if
there exists a decomposition $Z = \coprod_{\alpha} Z_{\alpha}$
into finitely many constructible subsets such that
each of the cohomology sheaves ${\cal H}^i(K)$ is locally
constant along each $Z_{\alpha}$ with finite dimensional
stalks.  This implies that the
limit 
\begin{equation}\label{eqn-stalk}
{\cal H}^i_x(K) :=\lim_{\to}H^i(U_x,K)\end{equation} 
is attained by any ``regular''
neighborhood $U_x$ of the point $x$ (for example, one may
embed (locally) $Z$ into a manifold and
take $U_x = Z \cap B_{\epsilon}(x)$ to be the intersection
of $Z$ with a sufficiently small ball centered at $x$). It
also implies  that $H^i(Z,K)$ is finite
dimensional.  Constructibility prevents the cohomology 
sheaves from exhibiting Cantor-set-like behavior. 

Most of the complexes of sheaves arising naturally from
geometric constructions on varieties  are  bounded and
have constructible cohomology sheaves. 

>From now on,  in this survey, unless otherwise stated,
bounded complexes with constructible  cohomology sheaves
are simply called {\em constructible complexes}.

The  {\em constructible 
bounded derived category} ${\cal D}_Z$
is  defined to be the   full subcategory 
of the  bounded derived category $D^b(Z)$ whose objects are the constructible complexes.
This subcategory is stable under the Verdier
duality functor, i.e. the dual of a  constructible 
complex is a bounded constructible complex,
it is stable under Hom, tensor products, vanishing and nearby
cycles functors,  and it is  well-behaved with respect to 
the functors $Rf_*,Rf_!,f^*,f^!$ associated with an algebraic map
$f: W \to Z$, i.e. $Rf_*, Rf_!: {\cal D}_W \to {\cal D}_Z$ and 
$f^*, f^!: {\cal D}_Z \to {\cal D}_W$.

{\bf 6.  Perverse sheaves,  intersection complexes.}
A {\em perverse sheaf} is a constructible complex with certain restrictions
(see $\S$\ref{subsecoervsh}) on the dimension of the support of its stalk cohomology and
of its stalk cohomology with compact supports (i.e. the analogue
with compact supports of (\ref{eqn-stalk})). These restrictions
are called the {\em support} and {\em co-support conditions}, respectively.

Let $U\subset Z$ be a nonsingular Zariski open subset and let $L$ be a
local system on $U$.  The {\em intersection complex} (\ci{goma2})
$IC_Z(L)$ is a complex of sheaves on $Z$, which extends the complex $L[\dim Z]$
on $U$ and is determined, up
to unique isomorphism in ${\cal D}_Z$,  by   support
and co-support
conditions that are slightly stronger  than the ones used to define perverse sheaves; see  equations (\ref{sprt}) and (\ref{csprt}) in \S\ref{subsecintcoh}. In particular, intersection 
complexes are perverse sheaves.
Up to a dimensional shift, the cohomology groups of the intersection complex $IC_Z(L)$
are the   the intersection cohomology groups of $Z$ twisted by
the system of local coefficients  $L$:
$H^i(Z, IC_Z(L))
= I\!H^{\dim{Z}+i}(Z,L)$. 

The category of perverse sheaves is
Abelian  and Artinian (see $\S$\ref{tfidy}):  every perverse sheaf
is an iterated extension of finitely many simple perverse sheaves. 
The simple perverse sheaves on $Z$
are  the intersection complexes
$IC_Y(L)$ of  irreducible subvarieties $Y\subset Z$ and irreducible
local systems $L$ defined on a nonsingular Zariski open subset of $Y$.

{\bf 7.  Perverse  cohomology sheaves,  perverse spectral sequence.}
The (ordinary) constructible sheaves, thought of as
the constructible  complexes  which are concentrated in degree $0$, form an Abelian full subcategory of the
constructible derived category ${\cal D}_Z$. An object $K$ of 
${\cal D}_Z$ is isomorphic to an object of this subcategory
if and only if ${\cal H}^i(K)=0$  for every $i \neq 0$.
There is a similar characterization of the category of perverse sheaves:
every  constructible complex $K \in {\cal D}_Z$ comes  equipped with a 
canonical collection of perverse sheaves on $Z$,   the {\em perverse cohomology
sheaves} $\pc{i}{K}$, $i \in \zed$.  
The perverse sheaves are 
characterized, among  the constructible complexes, by the property that $\pc{i}{K}=0$ for every  $i \neq 0$.

Just as there is the {\em Grothendieck spectral sequence}
\[
E^{l,m}_2= H^{l}(Z, {\cal H}^m(K))  \Longrightarrow H^{l+m}(Z,K),
\]
abutting to the  {\em standard} (or Grothendieck) filtration, 
there is the  {\em perverse spectral sequence} 
\[
E^{l,m}_2= H^{l}(Z, \pc{m}{K})  \Longrightarrow H^{l+m}(Z,K),
\]
abutting to the {perverse} filtration.
Similarly, for  the cohomology groups with compact supports $H_c^*(Z,K)$.

Let $f: W \to Z$ be a map of varieties and $C \in {\cal D}_W$. We have $
H^*(W, C) = H^*(Z, Rf_* C)$ and $H^*_c (W, C) = H^*_c(Z, Rf_! C)$.
The {\em perverse Leray spectral sequence and filtration} for
$H^*(W,C)$ and $H^*_c(W,C)$ are defined to be the perverse
spectral sequence and filtrations for $H^*(Z, Rf_*C) $ and $H^*_c(Z, Rf_! C)$,
respectively.

\begin{rmk}
\label{ptrestr}
{\rm 
If  $U$
is a nonempty,   nonsingular  and pure dimensional
open subset of $Z$ on which all the cohomology sheaves ${\cal H}^i(K)$ are local systems, then the restriction  to $U$ of
$\pc{m}{K}$    and ${\cal H}^{m- \dim Z}(K)[\dim{Z}]$ coincide. 
In general, the two differ:
in Example \ref{sfcefibr}, we have $\pc{0}{Rf_* \rat_X[2]} \simeq 
IC_Y (R^1) \oplus T_{\Sigma}$. This   illustrates
the non triviality of the notion of perverse cohomology sheaf.
}
\end{rmk}

\subsection{Decomposition, semisimplicity and  relative hard Lefschetz
theorems}
\label{statem}
Having dealt with some preliminaries on sheaves and derived categories,
we now state 
\begin{tm}
\label{stdteo}
{\rm ({\bf Decomposition and semisimplicity theorems})}
Let $f: X  \to Y$ be a proper map of complex algebraic varieties. There is an  isomorphism in the constructible bounded  derived category ${\cal D}_Y$:
\begin{equation}
\label{404}
Rf_* IC_X \; \simeq \; \bigoplus_{i\in \zed}    \pc{i}{Rf_*IC_X}  [-i].
  \end{equation}
  Furthermore,  the perverse sheaves
 $ \pc{i}{Rf_*IC_X}$ are semisimple, i.e.
  there is a decomposition into finitely
  many 
  disjoint locally closed and nonsingular subvarieties
  $Y= \coprod S_\beta$ and a canonical decomposition 
  into a direct sum of intersection complexes
  of   semisimple local systems
  \begin{equation}
  \label{4004}
  \pc{i}{ Rf_* IC_X } \; \simeq \;
   \bigoplus_{\beta}{
IC_{\overline{S_{\beta}}} (L_{\beta}) }.
\end{equation}
\end{tm}

The decomposition theorem is usually understood to be the combination
of (\ref{404}) and (\ref{4004}), i.e. the existence of a finite collection of triples
$(Y_a, L_a, d_a)$ as in theorem \ref{dtfirstform} such that
we have a direct sum decomposition 
\begin{equation}
\label{00dt}
Rf_* IC_X\; \simeq\; \bigoplus_a IC_{\overline{Y_a}}(L_a) [\dim{X} -\dim{Y_a} -d_a].
\end{equation}
Recalling that $I\!H^* (X) = H^{* - \dim{X}}(X, IC_X)$,
the cohomological shifts in the formula above
are chosen so that they match 
the ones of Theorem \ref{dtfirstform}, which is in fact a consequence
of (\ref{00dt}). The 
local systems $L_a$ are semisimple and the collection of triples
$(Y_a, L_a, d_a)$
is essentially unique.

The direct sum decomposition
(\ref{404})  is finite and $i$ ranges in the interval $[-r(f), r(f)]$,
where
$r(f)$ is the defect of semismallness of the map $f$ (see \S\ref{outlinedecmig}, part 2, and
\ci{decmightam}). In view of the properness of $f$ and of the fact that
$IC_X$ is a self-dual complex (i.e. it coincides with its own dual), 
Poincar\'e Verdier duality (cf. \S\ref{formulary}, duality exchanges), implies 
the existence of a canonical isomorphism
\begin{equation}
\label{pvddteq}\pc{-i}{f_*IC_X} \,\simeq \, \pc{i}{f_*IC_X}^{\vee}.\end{equation}
This important symmetry between
the summands in (\ref{404})  should not be confused with the somewhat deeper relative hard Lefschetz theorem, which is discussed below.

\begin{rmk}
\label{itisnotcanonical2}
{\rm
{\bf (The splitting is not canonical)}
The splittings (\ref{404}) and (\ref{00dt}) are not uniquely determined.
See Remark \ref{itisnotcanonical}.
}
\end{rmk}

It seems worthwhile to list {\em some} important and
 immediate consequences of Theorem \ref{stdteo}.

\begin{enumerate}
\item
The isomorphism (\ref{404}) implies immediately that the
perverse Leray spectral sequence 
$$
 E_2^{l,m}:= H^{l}(Y,{\pc{m}{Rf_*IC_X}})  \Longrightarrow I\!H^{\dim{X}+l+m}(X,\rat)
$$
is $E_{2}$-degenerate. 

\item
If $f: X \to Y$ is a resolution of the singularities of a variety $Y$, i.e. $X$ is 
nonsingular and $f$
is proper and an  isomorphism away from a proper closed subset of $Y$, then
one of the summands in (\ref{404}) is $IC_Y$ and we deduce that the intersection cohomology
of $Y$ is  (noncanonically) a  direct summand of the cohomology of any of its resolutions. Such resolutions exist, by a fundamental result
of H. Hironaka.

\item
If $f: X \to Y$ is a proper submersion of nonsingular varieties, then, in view of 
Remark \ref{ptrestr}, the decomposition (\ref{00dt}) can be re-written as
\[ Rf_* \rat_X \simeq \bigoplus R^if_* \rat_X [-i]\]
and one recovers Deligne's theorem (\ci{dess}) for families of projective manifolds
(a weaker form of which is  the $E_2$-degeneration
of the Leray spectral sequence for such maps stated in Theorem \ref{delteo},
part 1).
The semisimplicity statement of Theorem \ref{stdteo} corresponds then to
Theorem \ref{delteo}, part 2.

\end{enumerate}

As the name suggests, the {\em relative hard Lefschetz theorem} stated below
is the relative version of the classical hard Lefschetz theorem
seen in  $\S$\ref{tcav}, i.e. it is a statement that occurs in connection
with a map of varieties which, when applied to the special case of the
map of a projective manifold to a point, yields the classical
hard Lefschetz theorem. The relative version is closely linked to the decomposition theorem
as it expresses a  symmetry among the summands in (\ref{404}).

The symmetry in question  arises when considering the operation of cupping with the first Chern class
of  {\em a hyperplane line  bundle} on the domain of the map $f: X \to Y$.
The hyperplane bundle on projective space is the holomorphic
line bundle whose sections vanish precisely on  linear hyperplanes.
A hyperplane bundle on a quasi projective variety $X$
is the restriction to $X$ of the hyperplane line bundle
for some embedding $X \subseteq \pn{N}$.

 The first Chern class  of a line bundle  $\eta$ on $X$   
 yields, for every $i\geq 0$,  maps $\eta^i: Rf_* IC_X \to Rf_* IC_X [2i]$
and, by taking the perverse cohomology sheaves,  we obtain maps
of perverse sheaves $\eta^i: \pc{-i}{ Rf_* IC_X} \lorw \pc{i}{Rf_* IC_X}$.
 
\begin{tm}
\label{strhleo} 
{\rm ({\bf Relative hard Lefschetz theorem})}
Let $f: X \to Y$ be a proper map of varieties with 
$X$ quasi  projective  and let  $\eta$  be the first Chern class
of a hyperplane
 line bundle on $X$. 
Then we have isomorphisms
\begin{equation}
    \label{303}
  \eta^{i}: \pc{-i}{Rf_*IC_X} \stackrel{\simeq}\lorw 
\pc{i}{Rf_*IC_X}.
\end{equation}
\end{tm}

If $f$ is also a proper submersion, then we simply recover the classical hard Lefschetz
on the fibers of the map.
As mentioned above, if we apply this result to the special case
$f: X \to pt$, where $X$ is a projective manifold, then we obtain the classical hard Lefschetz.
If $X$ is a possibly singular  projective variety, then we obtain the hard Lefschetz theorem
in intersection cohomology ($\S$\ref{rel3}).

\begin{rmk}
\label{equiva}
{\rm 
Theorems \ref{stdteo} and \ref{strhleo}  also apply to $Rf_* IC_X(L)$
for   certain classes of  local systems $L$ (see \ci{bbd, samhp}).
}
\end{rmk}

 \begin{ex}
 \label{nonalgnodec}
 {\rm Let
$X=\pn{1}_{\comp} \times \comp$ and $Y$  be the space obtained collapsing the set $\pn{1}_{\comp} \times \{o\}$
to a point. This is  not a  complex algebraic  map 
and  (\ref{4004})
does not hold.
}
\end{ex}

\begin{ex}
\label{exhopf}
{\rm 
Let $f: ( \comp^2 \setminus \{0 \})/\zed =:X \to \pn{1}$ be the fibration 
in elliptic curves associated with
a Hopf surface.  Hopf surfaces are compact complex manifolds. Since $\pi_1(X) \simeq \zed$, we have $b_1(X)=1$ so that $X$ is not 
algebraic. In particular,  
though the map $f$ is  a proper
holomorphic submersion,
it  is not an algebraic map and Deligne's theorem,
and hence the decomposition theorem,  does not apply.
In fact, $Rf_{*} \rat_X$ does not split, for if it did, then $b_1(X) =2$.}
\end{ex}

\subsection{Invariant Cycle theorems}
\label{ictse}
The following theorem, in its local and global form, follows quite directly from the decomposition theorem.
It generalizes previous results, which assume that $X$ is smooth. 
For references, see the end of $\S$\ref{hire}.

In a nutshell, the global invariant cycle Theorem \ref{gict00} can be re--stated
as asserting   that if 
$f: X \to Y$ is a family of projective manifolds,
then the monodromy invariants
$H^*(F_y)^{\pi_1 (Y, y)}$ on the cohomology of a fiber
are precisely the image of the restriction map
$H^*(X) \to H^*(F_y)$ from the total space of the family. (Clearly, the image of the restriction map
is made of invariant classes, and the deep
assertion is that every invariant class is global, i.e. it comes from $X$.) In view of the generalization
given in Theorem  \ref{lict} below, we conveniently re-state this as the fact 
that the natural ``edge" map
\[
H^i(X) \lorw H^0 (Y, R^if_* \rat_X)
\qquad
is  \;\;surjective.
\]

\begin{tm}
\label{lict}
{\rm ({\bf  Global and local  invariant cycle theorems})}
Let $f: X \to Y$ be a proper map. 
Let $U \subseteq Y$ be  a Zariski open subset on which  the sheaf $R^i f_* (IC_X)$
is locally constant.
Then the following  assertions hold.
\begin{enumerate}
\item
{\rm
  ({\bf Global})}  The natural restriction map
$$
I\!H^i (X) \lorw H^0 (U, R^if_* IC_X) \qquad \mbox{is surjective.}
$$
\item
{\rm ({\bf Local})}
 Let $u\in U$  and $B_u \subseteq U$ be  the intersection
 with a sufficiently small Euclidean ball (chosen with respect to any local embedding
 of $(Y,u)$ into a manifold)
centered at $u$. Then the natural restriction/retraction map
$$
H^i(f^{-1}(u), IC_X) \simeq  H^i ( f^{-1} (B_u), IC_X)
\lorw H^0 ( B_u, R^i f_* IC_X) \qquad 
\mbox{is surjective}.
$$
\end{enumerate}
\end{tm}

\subsection{A few examples}
\label{afex}
In this section we discuss the statement of the decomposition theorem
in the following three examples: the resolution of singularities of 
a singular surface, the resolution of the affine cone over a projective nonsingular
surface
and a fibration of a surface onto a curve. More details can be found 
in \ci{decmigleiden}.

\begin{ex}
\label{exrsu}
{\rm 
Let $f: X \to Y$ be a resolution of the singularities of a singular
 surface $Y$. Assume that we  have a single singular point
 $y\in Y$ with $f^{-1} (y)= E$ a finite union of curves on $X$.
 Since $X$ is nonsingular, $IC_X = \rat_X [2]$ and
 we have  an isomorphism 
 $$
 Rf_* \rat_X[2] \simeq IC_Y \oplus T,
 $$
 where $T$ is a skyscraper sheaf at $y$ with stalk $T= H_2(E)$.
 }
\end{ex}

\begin{ex}
\label{cone}
{\rm 
Let  $S\subseteq {\Bbb P}^N_{\comp} $ be an embedded  projective
nonsingular  surface
and $Y\subseteq {\Bbb A}^{N+1}$ 
be the corresponding  threefold affine cone over $S$. Let $f: X \to Y$ be the blowing up of $Y$ at the vertex $y$. This is a resolution of the singularities
of $Y$, it is an isomorphism outside the vertex of the cone
and the fiber over the vertex is a copy of $S$. We have an   isomorphism
$$
Rf_* \rat_X [3] \simeq T_{-1}[1] \oplus  (IC_Y \oplus T_0)    \oplus
T_1[-1],
$$ 
where the $T_j$ are skyscraper sheaves at $y$ with stalks
 $T_{1}\simeq T_{-1}\simeq  H_4(S)$ and $T_0 \simeq H_3(S)$.
}
\end{ex}

\begin{ex}
\label{smallres}
{\rm
Let $S \subseteq \pn{3}$ be the nonsingular  quadric. The affine cone $Y$ over $S$
admits a resolution as in Example \ref{cone}. It also admits  resolutions
$f: X' \to Y$,  obtained  by blowing up a plane passing through
the vertex. In this case the exceptional fiber is isomorphic to $\pn{1}$ and
we have $Rf_* \rat_{X'}[3] = IC_Y$.
}
\end{ex}

\begin{ex}
\label{sfcefibr}
{\rm
Let $f: X \to Y$ be a projective map
with connected fibers from a smooth surface $X$ onto
a smooth curve $Y$. Let $\Sigma \subseteq Y$ be the finite set of critical values and let
$U = Y \setminus \Sigma$ be its complement. The map $f$ is a ${\cal C}^{\infty}$
fiber bundle over $U$ with typical fiber a compact oriented surface of some
fixed genus $g$.
Let $R^1 = (R^1 f_* \rat_X)_{|U}$ be the rank $2g$ local system
on $U$ with stalk the first cohomology
of the typical fiber. We have an   isomorphism
$$
Rf_* \rat_X [2]  \simeq \rat_Y [2] \oplus (  IC_Y(R^1) \oplus T_{\Sigma}   ) \oplus \rat_Y, 
$$
where $T_{\Sigma}$ is a skyscraper sheaf over $\Sigma$ with stalks
$T_s \simeq H_2(f^{-1} (s))/ \langle [f^{-1}(s) ]
\rangle$ at  $s \in \Sigma$.

}
\end{ex}

In all three examples the target space is a union
$Y= U \coprod \Sigma$ and we have two corresponding types of summands.
The summands of type $T$   
consists of classes which can be represented by 
cycles supported  over the exceptional set $\Sigma$.
This is precisely the kind of statement which lies at the heart of the decomposition theorem.
There are classes which can be represented by intersection cohomology classes
of local systems on  $Y$ and classes which 
can be represented by intersection cohomology classes of local systems supported over  smaller strata, 
and the cohomology of $X$ is the direct sum of these two subspaces. 
Suggestively speaking,
{\em it is as if the intersection cohomology relative to a stratum singled out precisely the classes which 
cannot be squeezed in the inverse image by $f$ of a smaller stratum}.

\subsection{The decomposition theorem and mixed Hodge structures}
\label{dtmHs}
The proof of the hard Lefschetz theorem in intersection cohomology
appears in \ci{bbd}. Therefore, at that point in time,  intersection cohomology
was known to enjoy the two Lefschetz theorems and Poincar\'e duality
(\ci{bbd, goma1,goma2}).
The question concerning a possible Hodge structure in intersection cohomology,
as well as other Hodge-theoretic questions, was very natural at that juncture
(cf. \ci{bbd}, p.165).

The work of M. Saito \ci{samhp, samhm}  settled these issues completely
with the use of mixed Hodge modules.
The reader  interested in the precise statements and generalizations is referred to Saito's papers (for  brief summaries, see \ci{durfeesaito} and  $\S$\ref{ms}).

In this section, we  summarize some of the mixed-Hodge-theoretic properties
of  the intersection cohomology  of complex  quasi projective varieties
that we have re-proved using  classical Hodge theory
(see $\S$\ref{dmapp}). 

The proofs can be found in
\ci{decmightam, decmigseattle, decmigso3,  decII}. More precisely,
the results for  projective varieties and the maps between them (in this case,
all Hodge structures are pure) are found in
\ci{decmightam, decmigseattle} and the extension to 
quasi projective varieties and the proper maps 
between them is found in \ci{decII}, which builds heavily on
\ci{decmigso3}.

Let us fix the set-up.
Let   $f: X \to Y$ be a  proper map of quasi projective varieties.
The intersection cohomology groups $I\!H^*(X)$  and $I\!H^*_c(X)$ are naturally filtered by the 
{\em perverse Leray  filtration} $P_*$, where $P_p I\!H^*(X)  \subseteq I\!H^*(X)$
and $P_p I\!H^*_c(X)$ are the images in cohomology and in cohomology with compact supports
 of the direct sum 
of terms $i'$ with $i' \leq p$ in the decomposition theorem (\ref{404}).
Up to re-numbering, this is the  filtration abutment of
the perverse Leray spectral sequence met in the crash course $\S$\ref{subsec-crash-course}  and it  can   be defined and described geometrically regardless
of the decomposition theorem (\ref{404}); see $\S$\ref{tpfmc}. We abbreviate  mixed Hodge structures  as
mHs.

\begin{enumerate}
\item
The intersection cohomology groups $I\!H^*(Y)$ and $I\!H^*_c(Y)$
carry  natural mHs. If $f:X \to Y $ is a resolution of the singularities
of $Y$, then  these mHs are    canonical subquotients
of the mHs on $H^*(X)$ and on $H^*_c(X)$, respectively.  
If $Y$ is a projective manifold, then the mHs is pure and it coincides
with the classical one (Hodge decomposition). If $Y$ is nonsingular, then 
the mHs coincide with Deligne's mHs on cohomology (see $\S$\ref{pam}).
The intersection bilinear  pairing in intersection cohomology
is compatible with  the mHs, i.e. the resulting map
$I\!H^{n-j}(Y) \lorw (I\!H_c^{n+j})^{\vee} (-n)$ is an isomorphism
of mHs.
The natural map $H^j(Y) \lorw I\!H^j(Y)$ is a map of mHs; if $Y$ is projective, then the  kernel is the subspace
$W_{j-1}$ of classes of Deligne weight $\leq j-1$.

\item
If $Y$ is a projective variety and $\eta$ is an hyperplane line bundle on $Y$, then
the hard Lefschetz theorem in intersection cohomology  
of $\S$\ref{rel3} holds.  In fact, the obvious transpositions from cohomology 
to intersection cohomology of the statements   in $\S$\ref{pam},   Theorem \ref{chl}
hold. 

\item
The subspaces $P_p$    of the  perverse Leray filtrations
in $I\!H^*(X)$ and in $I\!H^*_c(X)$
are mixed Hodge substructures of the mHs  mentioned in 1. The 
graded  spaces of these filtrations (i.e. $P_p/P_{p+1}$)  for $I\!H^*(X)$ and for
$I\!H^*_c(X)$ inherit the natural quotient  mHs
and they coincide (up to a shift in  cohomological degree) with the cohomology
and cohomology with compact supports   of the perverse cohomology sheaves
$\pc{p}{Rf_* IC_X}$. We call these spaces the {\em perverse cohomology
groups}. 

\item
The splitting of the perverse cohomology groups associated with the canonical splitting
(\ref{4004}) of the decomposition theorem takes place in the category of mHs.

\item
There exist splittings (\ref{404}) for the decomposition theorem
which induce  isomorphisms of mHs in cohomology and in cohomology with compact supports.
(Note that this  statement is stronger than the one above:  while these splittings take place
in $I\!H^*(X)$ and in  $I\!H^*_c(X)$, the previous ones take place
 in the perverse cohomology groups
which are subquotients of $I\!H^*(X)$ and of $I\!H^*_c(X)$.)

\item
The mHs we introduce coincide with the ones obtained by M. Saito
using mixed Hodge modules.
\end{enumerate}

\subsection{Historical  and other remarks}
\label{hire}
In this section we offer
few remarks that describe  the timeline for some of the results
mentioned in this survey.
We make no pretense to  historical completeness. 
For an account
 of the development of intersection cohomology,
 see the historical remarks in \ci{gomacsmt} and the
 survey  \ci{kle}.

By the late 1920's S. Lefschetz had ``proofs" of the  Lefschetz hyperplane    and hard Lefschetz theorems 
in singular cohomology (see \ci{lamotke} for an interesting
discussion of Lefschetz's proofs). Lefschetz's proof of the hard Lefschetz theorem is incomplete.

The Hodge decomposition theorem of cohomology classes
 into $(p,q)$-harmonic parts appears in W. Hodge's book \ci{hodge}. This is   where one also finds
 the first complete proof of the hard Lefschetz theorem (see also \ci{aweilherm}). The proof  of the $(p,q)$ decomposition in \ci{hodge} is not complete, and the missing analytical step was supplied by H. Weyl  (\ci{weyl}).

 S.S.  Chern gave a proof of the hard Lefschetz in the 1950's (see \ci{gh}) which  
 still relies on Hodge theory and  exploits the action of 
$sl_2 (\comp)$ on the differential forms on a K\"ahler manifold.

 In the 1950's R. Thom outlined a Morse-theoretic approach 
 to the hyperplane theorem which was worked out in detail by 
 A. Andreotti and Frankel
\ci{af}   (see \ci{milnor}) and  by R.  Bott \ci{bot}.

 The Hodge  decomposition
 is the blueprint for the definition of  pure  and mixed Hodge structures
 given by P. Griffiths and by  P. Deligne, respectively.
 The subject of how this decomposition varies in  a  family of projective manifolds
 and eventually degenerates has been studied, starting in the late 1960's, by 
 P. Griffiths and his school. The degeneration of the Leray spectral sequence for families of projective
 manifolds 
 was proved by P. Deligne in 1968.
 
 In 1980, Deligne  \ci{weil2} gave a new proof and a vast generalization
 of the hard Lefschetz theorem by  proving this result for varieties over finite fields
 and then inferring from this fact  the result over the complex numbers.
 (One usually says that one ``lifts the result from positive characteristic to characteristic zero,"
 see below.)
In particular, the hard Lefschetz theorem is
 proved for varieties defined over an algebraically closed field. By a result of M. Artin,
 the Lefschetz hyperplane theorem also holds in this generality.  
 
   Poincar\'e duality for intersection cohomology
is proved in \ci{goma1}. 
 The Lefschetz hyperplane theorem in intersection cohomology
is  proved
in \ci{goma2} and amplified in \ci{gomacsmt}. The hard Lefschetz theorem
 for the  intersection cohomology of projective varieties  is proved in \ci{bbd}.
 In the 1980's,  M. Saito  (\ci{samhp, samhm})    proved  that in the projective case
  these groups  admit a pure Hodge structure (i.e. a  $(p,q)$-Hodge decomposition),
  re-proved that
 they satisfy the hard Lefschetz theorem and proved
  the Hodge-Riemann bilinear relations.  In the 2000's, we re-proved 
these results in \ci{herdlef, decmightam}.

The decomposition theorem (\ref{dtnaive})  for the intersection cohomology groups had been conjectured in 1980 by S. Gelfand and
R. MacPherson. Note that they did not mention perverse sheaves. In fact,
the decomposition theorem (\ref{00dt}) only needs
the notion of intersection cohomology in order to be
 formulated.

The decomposition, semisimplicity and relative hard Lefschetz theorems
in $\S$\ref{statem} 
were proved by A. Beilinson, J.  Bernstein, P. Deligne and O. Gabber in 1982 (\ci{bbd}).
They first proved it  for proper maps of varieties defined   the algebraic closure
of finite fields, and then they lifted
 the result to  characteristic zero, i.e. for proper maps of complex algebraic
varieties. In fact, they prove the result for the proper direct image of 
complexes of geometric origin (see Definition  $\S$\ref{ssgops} in $\S$\ref{dafac})
and the intersection complex $IC_X$ is a special and important  example of a complex of geometric origin.
They also proved the invariant cycle results summarized in Theorem \ref{lict}.
Finally, they proved the hard Lefschetz theorem 
(\ref{00hlic}) for intersection cohomology as a special case of their relative hard Lefschetz theorem.
The equivariant version of these results  are proved in   \ci{bl}. 

At that juncture, it was natural to ask: 1)  for  a proof of the decomposition
theorem, semisimplicity and relative hard Lefschetz theorems for
complex varieties that 
uses transcendental methods; about the existence of 
Hodge structures  in intersection cohomology
(pure in the compact case, mixed in the general case), 2) 
about Hodge-Riemann relations in intersection cohomology
(in analogy with the ones for the  singular cohomology
of projective, or K\"ahler, manifolds; see Theorem \ref{chl} in $\S$\ref{pam}),
3) 
about possible extensions of the decomposition theorem etc.  to intersection complexes with twisted coefficients
underlying a polarized variation of pure Hodge structures,
4) about  suitable extensions to quasi projective varieties and mixed Hodge structures,
and  finally 5)
about generalizations of all these results to the K\"ahler case (e.g.
for proper  holomorphic maps $f: X \to Y$, where $X$ is a complex analytic space
which admits a  proper surjective  and generically finite map
onto it, e.g. a resolution of singularities, from  a complex K\"ahler manifold).

All these questions have been answered in the work of M. Saito \ci{samhp,samhm}
in the  1980's.

 The case of $IC_X$  (i.e. untwisted coefficients)  and of quasi projective varieties
 has been re-proved 
by us   using classical Hodge theory (see $\S$\ref{dtmHs}).

Finally, let us discuss the invariant cycle theorems. For families of
projective manifolds,
the global case was proved by  P. Deligne, in \ci{ho2}, 4.1.1.
The local case, conjectured and shown to hold for families of curves  by P. Griffiths in \ci{griff}, Conjecture. 8.1,  was proved by P.  Deligne in \ci{weil2}. 
For  Hodge-theoretic approaches to the local case,  see \ci{clem,steen,elz, nava}.
The ``singular" case, i.e. Theorem \ref{lict}, is proved in
 \ci{bbd}, p.164; see also \ci{samhp}.

 \section{Perverse sheaves}
 \label{secptsdy}
 Perverse sheaves have become an important tool in the study of singular
 spaces as they enjoy many of  the local and global properties
 of   the constant sheaf  that hold    on  nonsingular  spaces, but  that   fail
 on singular ones. They are fundamental mathematical objects
 whose importance goes beyond their role in the proof of the  decomposition theorem.

Here are some of the highlights of the theory of perverse sheaves.
 The reader can consult \ci{bbd,k-s,dimca}. Recall that we are dealing
 with $\rat$-coefficients and with middle-perversity only.
 We refer to  $\S$\ref{subsec-crash-course} and $\S$\ref{cemetto}
for more details and amplifications.

Historically,  perverse sheaves arose naturally from the theory of D-modules,
i.e. the sheaf-theoretic re-formulation of
linear systems of  partial differential equations:
 The "solution sheaf" of a holonomic D-module with regular
singularities is a perverse sheaf, and
this (Riemann-Hilbert correspondence) defines a functor from 
the category of holonomic D-modules with regular singularities 
to perverse sheaves.

 Even though the D-modules side of the story   is a necessary
complement to the more topological-oriented approach presented here,
for lack of competence, we  do not  treat it in this paper.
A partial list of references is  \ci{boreldm, bernst, ka1, ka2, ka3, meb1, meb2,
bjork}. 

Let $Y$ be a complex algebraic variety.
 Like the category of constructible sheaves,
the category $\pe_Y$ of perverse sheaves is a full
 subcategory of the 
constructible derived
category ${\cal D}_Y$. The category $\pe_Y$  is Abelian, Noetherian  and
Artinian (i.e. every perverse sheaf is a finite iterated extensions
of simple perverse sheaves). 
The simple perverse sheaves on $Y$
are the intersection
complexes $IC_W(L)$  associated with an irreducible  and closed subvariety $W \subseteq Y$ and an
irreducible local system $L$ (on a Zariski-dense open nonsingular subvariety
of $W$). 
Since $\pe_Y$ is an Abelian category, any morphism in $\pe_Y$
admits a  (``perverse'')
kernel and (``perverse'') cokernel.
 Given a complex $K\in {\cal D}_Y$, there are the 
(``perverse'') cohomology sheaves $\pc{i}{K} \in  \pe_Y$. 
A theorem of A. Beilinson's states  that the  bounded derived category of $\pe_Y$
is again ${\cal D}_Y$.  Many operations work better in the category of perverse 
sheaves than in the category of sheaves, e.g.  the duality  and vanishing cycles
functors preserve perverse sheaves.
The Lefschetz hyperplane theorem holds 
for perverse sheaves.  Specialization over a curve takes perverse sheaves
to perverse sheaves.  The intersection cohomology
of a projective variety satisfies the Hodge-Lefschetz 
theorems and Poincar\'e duality.

\subsection{Intersection cohomology}
\label{subsecintcoh}
The intersection cohomology complex of a complex
algebraic variety $Y$ is a special case of a perverse sheaf
and every perverse sheaf is a finite iterated extension of intersection complexes.
It seems appropriate to start a discussion of perverse sheaves
with this most important example.

Given a complex $n$-dimensional algebraic
variety $Y$ and a local system $L$ on a nonsingular Zariski-dense open subvariety
$U\subseteq Y$,  there exists a constructible complex of sheaves  $IC_Y (L) \in {\cal D}_Y$, 
unique  up to canonical isomorphism in ${\cal D}_Y$, such that
$IC(L)_{\mid U} \cong L$ and:
\begin{equation}
\label{sprt}
 \dim \left\{ y \in Y|\ {\cal H}_y^i(IC(L)) \ne 0\right\} 
<  -i, \; \mbox{if $i> -n$,}
\end{equation}

\begin{equation}
\label{csprt}
 \dim \left\{ y \in Y |\ {\cal H}^{i}_{c,y}(IC(L)) \ne 0 \right\} < 
 i, \; \mbox{if $i<n$,}
\end{equation}
where, for any complex $S$ of sheaves, 
\[{\cal H}^i_{c,y}(S) = \lim_{\leftarrow} H^i_c(U_y,S)\]
is the local {\em compactly supported} cohomology at $x$. (As 
explained in the ``crash course'' $\S$\ref{subsec-crash-course}, 
if $S$ is constructible, then the above limit is attained by
any regular neighborhood $U_y$ of $y$.) The intersection complex $IC_Y(L)$ 
is sometimes called the {\em intermediate extension} of $L$. Its (shifted) cohomology
is the intersection cohomology
of $Y$ with coefficients in $L$, i.e.  $I\!H^{n+*}(Y,L):= H^*(Y, IC_Y(L))$.   
The reader can consult  \ci{goma1, goma2} and 
\ci{borel, dimca}.

Even though 
intersection cohomology lacks functoriality with respect to algebraic maps
(however, see \ci{bbfgk}),
the intersection cohomology groups of projective
varieties enjoy the same properties of Hodge-Lefschetz-Poincar\'e-type
as the singular cohomology of  projective manifolds. 
Poincar\'e duality takes the form $I\!H^k(Y) \simeq I\!H^{2n-k}(Y)^{\vee}$
and follows formally from the canonical isomorphism
$IC_Y \simeq  IC_Y^{\vee}$ stemming from Poincar\'e-Verdier duality;
in particular, there is a non degenerate geometric intersection pairing
$$I\!H^i (Y) \times I\!H^{2n-i}(Y) \lorw  \rat, \qquad
(a,b) \longmapsto a\cdot b;$$
on the other hand there is no cup-product. As to the other properties,
i.e. the two Lefschetz theorems, the Hodge decomposition and the Hodge-Riemann
bilinear relations, see
$\S$\ref{dtmHs} and \S\ref{dmapp}.

\subsection{Examples of intersection cohomology}
\label{6010}
\begin{ex}
\label{exicd}
{\rm 
Let $E^{n-1} \subseteq
\pn{N}$ be  a projective manifold, $Y^n \subseteq {\Bbb A}^{N+1}$
be the associated affine cone. The   link  $\cal L$ of $Y$ at the vertex
$o$ of the cone, i.e.  the intersection of $Y$ with a sufficiently small
Euclidean sphere centered at $o$,  is an oriented
compact smooth manifold of real dimension $2n-1$ and is an $S^1$-fibration over $E$.
The cohomology groups of $\cal L$ are
$$
H_{2n-1 -j}({\cal L}) = H^j({\cal L}) = P^j(E),  \;\;  0 \leq
j \leq n-1, \qquad H^{n-1+j} ({\cal L}) = P^{n-j}(E), \; \; 
0 \leq j \leq n.
$$
where $P^j(E) \subseteq H^j(E)$ is the subspace of primitive vectors
for the given embedding of $E$, i.e. the kernel of 
cupping with the appropriate power of the first Chern class
 of ${\cal O}_E (-E)$.
The Poincar\'e intersection form  on $H^*({\cal L})$ is non degenerate,
as usual,  and also because of the Hodge-Riemann bilinear relations 
(\ref{hrbr})
on $E$.

The intersection cohomology groups of $Y$ are
$$
I\!H^j(Y) \simeq P^j(E)=H^j({\cal L}), \;\; 0 \leq j \leq n-1,
\qquad
I\!H^j(Y) =0, \;\;  n \leq j \leq  2n.
$$
The intersection cohomology with compact supports of $Y$
are
$$
I\!H^{2n-j}_c(Y) \simeq  H_j({\cal L}), \;\;
0 \leq j \leq n-1, \qquad
I\!H^{2n-j}_c(Y) =0, \;\;   n \leq j\leq 2n.
$$ 
We thus see 
that, in this case,   the Poincar\'e duality isomorphism $I\!H^j(Y) \simeq I\!H_c^{2n-j}(Y)^{\vee}$
 stems from  the classical Poincar\'e duality on the link. 
 }
 \end{ex}

 In  the remaining part of this section, 
 we complement some examples  of intersection complexes and groups with 
 some further information expressed using the language   of  perverse sheaves 
 which we discuss in the next few sections.

\begin{ex}
\label{exica}{\rm
 Let $Y$ be the projective cone  over a nonsingular curve
 $C \subseteq  \pn{N}$ of genus $g$.
 The  cohomology groups are
 $$
 H^{0}(Y) = \rat,  \quad H^1(Y) =0,  \quad    H^2(Y)=\rat , \quad
H^3(Y) = \rat^{2g} , \quad
 H^4(Y) =\rat.
 $$  
 The intersection cohomology groups are:
 $$
 I\!H^{0}(Y) = \rat , \quad
 I\!H^{1}(Y) =\rat^{2g}, \quad
  I\!H^{2}(Y)  = \rat,  \quad
  I\!H^{3}(Y) =\rat^{2g}, \quad
  I\!H^{4}(Y) = \rat.
  $$ 
  Note the  failure of Poincar\'e duality in cohomology
and its restoration via intersection cohomology.
   There is a canonical resolution $f: X \to Y$ of the singularities of $Y$ obtained by blowing up
the vertex of $Y$. The decomposition theorem
yields  
a splitting  exact sequence of perverse sheaves on $Y$:
$$
\xymatrix{
0 \ar[r] &  IC_Y \ar[r] &  f_{*} \rat_X[2] \ar[r]&    H^2(C) [0] \ar[r]&  0.}
$$}
\end{ex}

\begin{ex}
\label{exicb}
{\rm We now re-visit Example  \ref{nonalgnodec}.
Let $f: X \to Y$ be the  space obtained by contracting 
to a point $v \in Y$,
 the zero section $C  \subseteq \pn{1} \times C=:X$. This example is analogous to
 the one in Example \ref{exica}, except that  $Y$ is not a complex algebraic
 variety.
  The  cohomology groups are 
  $$H^{0}(Y) = \rat, \quad  H^1(Y) =0, \quad H^2(Y)=\rat, \quad
H_3(Y) = \rat^{2g}, \quad   H_4(Y) = \rat.
$$
The stratified space $Y$ has strata of even codimension and we can define
its intersection complex etc.
The intersection  cohomology groups are:
$$I\!H^{0}(Y) = \rat, \quad  I\!H^{1}(Y) =\rat^{2g}, \quad
I\!H_{2}(Y) =0, \quad I\!H^{3}(Y)= \rat^{2g}, \quad 
I\!H_{4}(Y)  = \rat.
$$  
Note the  failure of Poincar\'e duality in cohomology
and its restoration via intersection cohomology.
 There is a natural epimorphism of
perverse sheaves $\tau :  f_{*} \rat_X[2]  \lorw H^2(C) [0]$.
There are 
non splitting  exact sequences in ${\cal P}_Y:$
$$
0 \lorw  \ke{\,\tau}  \lorw f_{*} \rat_X[2]  \lorw H^2(C) [0]  \lorw 0, \qquad 
0 \lorw IC_Y \lorw \ke{\, \tau}  \lorw \rat_v[0] \lorw 0.
$$
The complex  $f_{*} \rat_X[2]$ is a perverse sheaf on $Y$  obtained by  two-step-extension procedure
involving intersection complexes (two of which are skyscraper sheaves).
The intersection cohomology complexes $IC_Y$ and $\rat_v$ of $Y$ and $v \in Y$
appear in this process, but not as direct summands.
The conclusion of the decomposition theorem does not hold for this map $f$.
}\end{ex}
 
\begin{ex}
\label{exicc}
{\rm
 Let $Y$ be the projective cone over the quadric $\pn{1}\times \pn{1} \simeq 
Q \subseteq \pn{3}$. The odd cohomology is trivial.  The even cohomology is as follows:
$$H^0 (Y) =0, \quad H^2(Y) =\rat, \quad H^4(Y) = \rat^2, \quad
  H^6(Y) = \rat.
  $$
The intersection cohomology groups are the same as the cohomology
groups, except  that  $I\!H^2(Y) =\rat^2$.  Note the  failure of Poincar\'e duality in homology
and its restoration via intersection homology. There are at least two different
and interesting resolutions of the singularities of $Y$: the ordinary blow up of the vertex
$o \in Y$
$f: X \to Y$ which has fiber $f^{-1}(o)\simeq Q$, and the blow up of any line
on the cone  through the origin $f': X' \to Y$ which has fiber
${f'}^{-1} (o) \simeq \pn{1}$.  The decomposition theorem
yields (cf. Example \ref{cone})
$$
f_* \rat_X[3] = IC_Y \oplus \rat_o [1] \oplus \rat_o [-1], \qquad
f'_* \rat_{X'} [3] = IC_Y.
$$
}
\end{ex}

\begin{ex}
\label{buff}
{\rm
Let $E$ be the rank two local system on the punctured complex line $\comp^*$
defined by the automorphism of $e_1 \mapsto e_1, $ $ e_2 \mapsto e_1+e_2$.
It fits into the non trivial extension
$$
\xymatrix{
0 \ar[r] &  \rat_{\comp^*} \ar[r] &  E \ar[r]^\phi & \rat_{\comp^*} \ar[r] & 0.
}
$$
Note that $E$ is self-dual. If we shift this extension by $[1]$,
then we get a non split  exact sequence of perverse sheaves
in ${\cal P}_{\comp^*}$. Let $j: \comp^* \to \comp$ be the open immersion.
The complex $IC_{\comp}(E) = R^0j_* E[1]$ is a single sheaf in
cohomological degree $-1$ with generic stalk $\rat^2$ and stalk $\rat$ at the origin $0 \in
\comp$.
In fact, this stalk is given by the space of invariants which is spanned
by the single  vector $e_1$. We remark, in passing, that given any local 
system $L$ on $\comp^*$, we have that $IC_{\comp}(L) = R^0j_* L[1]$.
There is the monic map $\rat_{\comp}[1] \to IC_{\comp}(E)$. The cokernel $K'$
is the nontrivial extension, unique since $\mbox{Hom} (\rat_{\comp}, \rat_{ \{ 0 \} }) =
\rat$ is one dimensional,
$$
\xymatrix{
0 \ar[r] & \rat_{ \{ 0 \} } \ar[r] &  K' \ar[r] &  \rat_{\comp} [1] \ar[r] & 0.
}$$
Note that while  the perverse sheaf $IC_{\comp}(E)$, being an intermediate extension ($\S$\ref{subsectie}), has no subobjects and no quotients
supported at $\{ 0 \}$, it has a subquotient supported 
at $\{ 0 \}$, namely   the perverse sheaf $\rat_{ \{ 0 \} }$.
We shall meet this example again later in Example \ref{buffo}, in the context of
the non exactness of the intermediate extension functor.
}
\end{ex}

\begin{ex}
\label{exnnormcross}
{\rm Let $\Delta \subseteq \comp^n$ be the subset 
$\Delta=\{(x_1, \ldots,x_n)\in \comp^n \,: \prod x_i=0    \}$.
The datum of $n$ commuting endomorphisms $T_1, \ldots, T_n$ of a $\rat$-vector space $V$ defines a local system $L$ on $(\comp ^*)^n=\comp^n \setminus \Delta$
whose stalk at some base point $p$ is identified with $V$, and $T_i$ is the monodromy along the path ``turning around the divisor $x_i=0$." The vector space $V$ has a natural structure of a $\zed^n= \pi_1((\comp ^*)^n,p)$-module. The complex which computes the group cohomology $H^{\bullet}(\zed^n,V)$ of $V$ can be described as follows: Let $e_1, \ldots, e_n$ be the canonical basis of $\rat^n$, and, for 
$I=(i_0, \ldots , i_k)$, set $e_I=e_{i_0}\wedge \ldots \wedge e_{i_k}$.
Define a complex   $(C,d)$ by setting
\[C^k=\bigoplus_{0<i_0<\ldots <i_k<n}V \otimes e_I, \qquad 
d(v \otimes e_I)= \sum N_i(v)\otimes e_i\wedge e_I,\]
with $N_i:=T_i- I$. 
Since $(\comp ^*)^n$ has no higher homotopy groups,   we have the quasi isomorphism
$(j_*L)_0 \stackrel{qis}\simeq (C, d)$.
Let $$
\widetilde{C^k}=\bigoplus_{0<i_0<\ldots <i_k<n}N_I V \otimes e_I, 
$$
where $N_I:=N_{i_0}\circ \ldots \circ N_{i_k}$.
It is clear that $(\widetilde{C}, d)$ 
is a subcomplex of $(C, d)$. There is   a  natural  isomorphism
$(IC_{\comp^n}(L))_0  \simeq (\widetilde{C}, d)$.
The particularly important case in which $L$ underlies a polarized variation 
of Hodge structures has been investigated in depth in \ci{cks} and \ci{kk}.
}
\end{ex}

\subsection{Definition  and first properties of perverse sheaves}
\label{subsecoervsh}

Let $K \in {\cal D}_Y$ be a constructible complex on the variety $Y$.
Recall that the support of a sheaf is the closure of the set of points where the sheaf has non trivial
stalks.
We say that $K$ satisfies the  {\em support condition} if
$$
\dim{ \{\mbox{Supp}} \, {\cal H}^{-i}(K) \} \leq i,  \quad \forall \, i \in \zed.
$$
We  say that $K$ satisfies the  {\em co-support condition}
if the Verdier dual $K^{\vee}$ ($\S$\ref{tfidy}) satisfies the conditions of support.

By Verdier duality, we have
${\cal H}^i_y (K^{\vee})\simeq {\cal H}^{-i}_{c,y}(K)^{\vee}$, so that
we may write the support and co-support conditions as follows:
\begin{equation}
\label{pssprt}
 \dim \left\{ y \in Y|\ {\cal H}_y^i(K) \ne 0\right\} 
\leq  -i, \;  \forall \mbox{ $i\in \zed$,}
\end{equation}

\begin{equation}
\label{pscsprt}
 \dim  \left\{ y \in Y |\ {\cal H}^{i}_{c,y}(K) \ne 0 \right\}  \leq 
 i, \;  \forall \mbox{$i \in \zed$.}
\end{equation}

\begin{defi}
\label{defps}
{\rm A {\em perverse sheaf} on $Y$ is a constructible complex $K$ in ${\cal D}_Y$ that
 satisfies the conditions of support and co-support. The category $\pe_Y$ of perverse
 sheaves is the full subcategory of ${\cal D}_Y$ whose objects are the  perverse
 sheaves.
}
\end{defi}

A complex $K$ is perverse iff $K^{\vee}$ is perverse. The defining conditions of intersection complexes in 
$\S$\ref{subsecintcoh} are a stricter versions of the support and co-support conditions given above
It follows that  intersection complexes are  special perverse sheaves.

Figure 1  below
\begin{figure}[!ht]
\begin{picture}(600,200)(-90,-90)
\linethickness{.7pt}\setlength{\unitlength}{1.2pt}
\multiput(0,40)(10,0){5}{c}  \multiput(0,-40)(10,0){5}{x}
\multiput(20,30)(10,0){3}{c} \multiput(20,-30)(10,0){3}{x}
\multiput(30,20)(10,0){2}{c} \multiput(30,-20)(10,0){2}{x}
\multiput(40,10)(10,0){1}{c} \multiput(40,-10)(10,0){1}{x}

\multiput(100,60)(10,0){7}{c}  \multiput(100,-60)(10,0){7}{x}
\multiput(110,50)(10,0){6}{c}  \multiput(110,-50)(10,0){6}{x}
\multiput(120,40)(10,0){5}{c} \multiput(120,-40)(10,0){5}{x}
\multiput(130,30)(10,0){4}{c} \multiput(130,-30)(10,0){4}{x}
\multiput(140,20)(10,0){3}{c} \multiput(140,-20)(10,0){3}{x}
\multiput(150,10)(10,0){2}{c} \multiput(150,-10)(10,0){2}{x}
\multiput(160,0)(10,0){1}{$\bullet$}

\put(-10,-70){\vector(0,1){140}}
\put(-10,-70){\vector(1,0){190}}
\put(-22,75){degree}\put(190,-72){codim $Y_{\alpha}$}

\put(-20,0){0}\put(-20,10){1}\put(-20,20){2}\put(-20,30){3}
\put(-20,40){4}\put(-20,50){5}\put(-20,60){6}

\put(-23,-10){-1}\put(-23,-20){-2}\put(-23,-30){-3}
\put(-23,-40){-4}\put(-23,-50){-5}\put(-23,-60){-6}

\put(0,-80){0}\put(10,-80){1}\put(20,-80){2}\put(30,-80){3}\put(40,-80){4}
\put(100,-80){0} \put(110,-80){1}\put(120,-80){2}\put(130,-80){3}
\put(140,-80){4}\put(150,-80){5}\put(160,-80){6}

\end{picture}\caption{support conditions for $IC$ (left) and for a 
perverse sheaf (right)}
\end{figure}
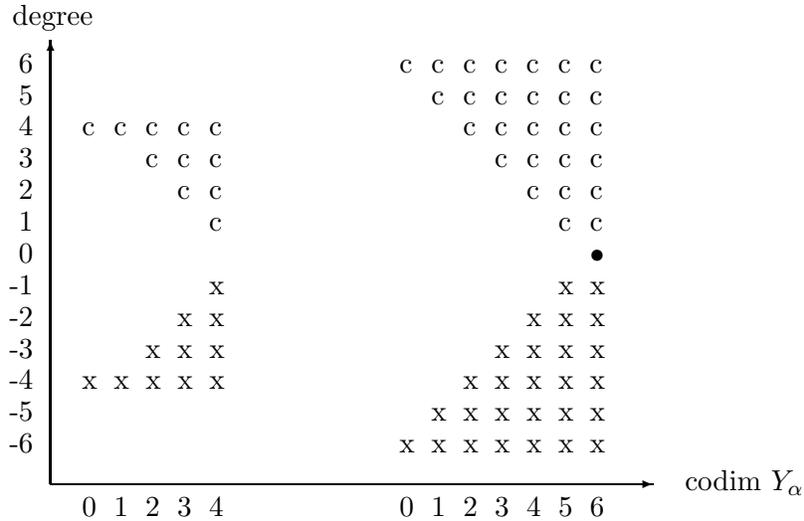
illustrates the support and cosupport conditions 
for intersection cohomology on a variety of dimension 4 (left)
and a perverse sheaf on a variety of dimension 6 (right). 
The symbol ``c'' means that compactly supported stalk cohomology can be non-zero at that place,
while the symbol ``x'' means that stalk cohomology can be 
non-zero at that place.
Note that the $\bullet$ symbol shows that, for a perverse sheaf, there is a place
at which both compactly supported and ordinary cohomology can be non-zero.
As explained in  $\S$\ref{macvilon}, the natural map $
{\cal H}_{c,y}^i(-) \to {\cal H}_y^i(-)$ governs the splitting behaviour of the perverse sheaf.

 Denote by ${\cal P}_Y$
the full subcategory of ${\cal D}_Y$ whose objects are perverse sheaves.
Denote by  $^{\frak p}\!{\cal D}_Y^{\leq 0}$ ($^{\frak p}\!{\cal D}_Y^{\geq 0}$,   resp.) 
the full subcategory of ${\cal D}_Y$
with objects the complexes satisfying the conditions of support (co-support, resp.).
Clearly, $^{\frak p}\!{\cal D}_Y^{\leq 0} \cap \, ^{\frak p}\!{\cal D}_Y^{\geq 0} = {\cal P}_Y$.  These data give rise to the {\em middle perversity $t$-structure}
on ${\cal D}_Y$ (see $\S$\ref{tfidy}).

\begin{tm}
\label{itistst}
The datum of the conditions of (co)support
together with the associated full subcategories  $(^{\frak p}\!{\cal D}_Y^{\leq 0}, 
\, ^{\frak p}\!{\cal D}_Y^{\geq 0})$ 
yields a $t$-structure on ${\cal D}_Y$,
called the middle perversity  $t$-structure,  with heart
$^{\frak p}\!{\cal D}_Y^{\leq 0} \cap \, ^{\frak p}\!{\cal D}_Y^{\geq 0}$
the category of perverse sheaves ${\cal P}_Y$.
\end{tm}

The resulting truncation  and cohomology functors   are denoted, for every $i\in \zed$:
  $$ \ptd{i} : {\cal D}_Y \lorw \, ^{\frak p}\!{\cal D}_Y^{\leq i}, \quad 
    \ptu{i} : {\cal D}_Y \lorw \, ^{\frak p}\!{\cal D}_Y^{\geq i}, 
    \quad
    ^{\frak p}\!{\cal H}^0 = \ptu{0}\ptd{0}, \; 
    ^{\frak p}\!{\cal H}^i = \, ^{\frak p}\!{\cal H}^0 \circ [i] \;  : 
{\cal D}_Y \lorw {\cal P}_Y.
    $$
    In particular, any complex $K \in {\cal D}_Y$ has ``perverse cohomology sheaves"
    $\pc{i}{K} \in \pe_Y$.
    
The key point in the proof  is to show the existence of $\ptu{0}$ and $\ptd{0}$.
The construction of these perverse truncation functors
involves only the four functors $f^*,f_*, f_!,f^!$ for open and closed immersions
and standard truncation.
See \ci{bbd},  or \ci{k-s}. Complete and brief summaries 
 can be found in \ci{decmightamv2,decmigleiden}.

Middle-perversity is   well-behaved with respect to Verdier  duality: 
the Verdier duality functor  $D: \pe_Y \to \pe_Y$ is an equivalence and 
we have canonical isomorphisms
$$
^{\frak p}\!{\cal H}^i \circ   D \simeq   D \circ \, ^{\frak p}\!{\cal H}^{-i}, 
\qquad
\ptd{i} \circ   D \simeq   D \circ  \ptu{-i},  \qquad
\ptu{i} \circ   D \simeq  D \circ  \ptd{-i}.
$$

 It is not difficult to show, by using the perverse cohomology functors (see \S\ref{subsecpervco}),
that $\pe_Y$ is an Abelian category. As it is customary when dealing with 
Abelian categories, when we say that  $A \subseteq B$ ($A$ is included in $B$),
we mean that there is a monomorphism $A \to B$. 
The Abelian category $\pe_Y$  is  Noetherian (i.e. every increasing sequence of perverse subsheaves of a perverse sheaf
must stabilize) so that,  by Verdier duality,
it is also Artinian (i.e. every decreasing sequence stabilizes). The category of constructible sheaves is Abelian and  
Noetherian, but not Artinian.

A. Beilinson  \ci{beili} has proved that, remarkably,  the bounded
derived category of perverse sheaves $D^b ({\cal P}_Y)$ is equivalent to ${\cal D}_Y$.
There is a second, also remarkable, equivalence due to M. Nori.
Let $D^b (CS_Y)$ be 
the  bounded derived category of the 
category of constructible sheaves on $Y$ (the objects
are  bounded complexes of constructible sheaves).
There is a natural inclusion of categories $D^b(CS_Y) \subseteq
{\cal D}_Y$ (recall that the objects of ${\cal D}_Y$ are bounded 
complexes of sheaves 
whose cohomology sheaves are constructible).
M. Nori \ci{nori} has proved that the inclusion $D^b (CS_Y) \subseteq
{\cal D}_Y$ is an equivalence of categories.
This is a striking  instance of the
  phenomenon that a category  can arise as a derived category in
fundamentally  different ways: ${\cal D}_Y \simeq D^b (\pe_Y) \simeq
D^b (CS_Y)$.

 Perverse sheaves, just like ordinary sheaves, form a stack 
(\ci{bbd}, 3.2), i.e. suitably compatible systems of perverse sheaves 
can be glued to form a single perverse sheaves, and similarly for
compatible systems of  morphisms of perverse sheaves.
 This is not the case for the objects and morphisms
 of ${\cal D}_Y$; e.g. a non trivial extension of vector bundles
 yields a non zero morphism in the derived category that
 restricts to zero on the open sets of a suitable open covering,
 i.e. where the extension
 restricts to trivial extensions.

\begin{ex}
\label{exofpervp}
{\rm 
Let $Y$ be a point.  The standard and perverse $t$-structure coincide.
A complex $K \in {\cal D}_{pt}$ is perverse iff 
it is isomorphic in ${\cal D}_{pt}$  to a complex
concentrated in degree zero iff  ${\cal H}^{j}(K)=0$ for every $j\neq 0$.
}
\end{ex}

\begin{ex}
\label{isosing}
{\rm 
If $Y$ is a variety of dimension $n$, then the complex  $\rat_Y[n]$ trivially
satisfies  the conditions of support.
If $n=\dim{Y}=0,1$, then $\rat_Y[n]$ is perverse.
On a surface $Y$ with isolated singularities,  $\rat_Y[2]$ is perverse iff the singularity is unibranch, 
e.g. if the surface is normal.
If $(Y,y)$ is a germ of a threefold isolated singularity, then  $\rat_Y[3]$  is perverse iff the singularity is unibranch
and  $H^1(Y\setminus y) =0$. 
}
\end{ex}

\begin{ex}
\label{exviatex}
{\rm 
The direct image $f_*\rat_X[n]$ via a proper semismall map  $f:X \to Y$, where $X$ is a nonsingular $n$-dimensional nonsingular variety, is perverse   (see Proposition
\ref{ss}); e.g. a generically finite map of surfaces is semismall.  For an interesting, non semisimple, perverse sheaf arising from a non algebraic semismall  map,
see Example \ref{exicb}.}
\end{ex}

Perverse sheaves are stable under the following functors:
intermediate extension, nearby and vanishing cycle (see  $\S$\ref{psiphi}).

Let $i: Z \to Y$ be the closed immersion of a subvariety of $Y$.
One has the functor $i_{*} : {\cal P}_Z \to {\cal P}_Y$. 
This functor is fully faithful, i.e. it induces a bijection
on the Hom-sets.  It is customary, e.g. in the statement of the decomposition
theorem, to drop the symbol $i_*$.

Let  $Z$ be an irreducible   closed subvariety of $Y$ and $ L$
be  an irreducible
(i.e. without trivial local subsystems, i.e. simple in the category of local systems)  local system on a non-empty Zariski open
subvariety  of  the regular part $Z_{reg}$ of $Z$. 
Recall that a simple object in an abelian category is one without
trivial subobjects.
The complex $IC_Z ( L)$ is a simple object
of the category ${\cal P}_Y$. 
Conversely, every simple object of ${\cal P}_Y$ has this form. This follows 
from the following  proposition \ci{bbd}, which yields a
direct proof of the fact that $\pe_Y$ is Artinian. 

Recall
that by an inclusion $A \subseteq B$, we mean the  existence
of a monomorphism $A \to B$, so that by a chain of inclusions,
we mean a chain of monomorphisms. 
The following caveat may be useful, as it points out that
the usual set-theoretic intuition about injectivity and surjectivity
may be misleading when dealing with perverse sheaves, or with Abelian categories in general. Let $j : \comp^* \to \comp$ be the open immersion.
We have a natural  injection of sheaves  $j_! \rat_{{\comp}^*}\to \rat_{\comp}$. On the other hand,
one can see that the  induced map of perverse sheaves
$j_! \rat_{{\comp}^*}[1]  \to \rat_{\comp}[1]$ is not a monomorphism; in fact, it is an epimorphism.

\begin{pr}
\label{pyart} 
{\rm ({\bf Composition series})}
Let $P \in {\cal P}_Y$. There is a  finite decreasing filtration
$$
 P = Q_1 \supseteq Q_{2}  \supseteq \ldots   \supseteq
  Q_{\lambda} = 0,
 $$
where the quotients $Q_i/Q_{i-1}$ are  simple perverse sheaves on $Y$.
Every simple perverse sheaf if of the form $IC_{\overline{Z}}(L)$, where 
$Z \subseteq Y$ is an irreducible and nonsingular subvariety and $L$
is an irreducible local system on $Z$.
\end{pr}

As  usual, in this kind of situation, e.g. the Jordan-H\"older
theorem for finite groups, the filtration is not unique, but the
constituents of $P$, i.e. the non trivial  simple quotients,
and  their multiplicities are uniquely determined.

\subsection{The perverse filtration}
\label{tpfmc}
The theory of $t$-structures  coupled with Verdier's formalism
of spectral objects
(cf. \ci{shockwave}, Appendix), endows the cohomology groups
$H^*(Y,K)$ with the  canonical {\em perverse filtration $P$}
defined by $P^pH^*(Y,K):= \im\, \{H^*(Y, \ptd{-p} K) \to H^*(Y,K)\}$, which   is,
up to re-numbering
 the abutment of the {\em perverse spectral sequence}
 $H^{p}(Y, \pc{q}{K}) \Longrightarrow H^*(Y,K)$. 
 See $\S$\ref{subsec-crash-course}.(7).
 Similarly, for cohomology with compact supports. 
  
  In \ci{decmigso3},  we give a geometric description of the  perverse filtration on
the cohomology and on  the cohomology with compact supports 
of a  constructible complex on a quasi projective variety. 
The paper \ci{decII} gives an alternative proof
with the applications to mixed Hodge theory mentioned in
$\S$\ref{dtmHs}; the paper  \ci{deabday} proves 
similar results for the standard filtration on cohomology with compact supports.  

 The description is in terms of restriction to generic hyperplane sections and it is somewhat unexpected, especially if one views
the constructions leading to perverse sheaves  as 
transcendental and hyperplane sections as more algebro-geometric. 
If $f: X \to Y$ is a map of quasi projective varieties and $C \in {\cal D}_X$,
then our results   yield a similar  geometric description of the
perverse Leray filtration on $H^*(X,C)$ and on $H^*_c(X,C)$ induced by the map $f$.

We now describe the perverse filtration on the cohomology
groups $H^*(Y,K)$  when $Y$ is {\em affine}.
Let $Y_*= \{Y \supseteq Y_{-1} \supseteq  \ldots \supseteq 
Y_{-n}\}$ be a sequence of closed subvarieties; we call this data
an $n$-flag. Basic sheaf theory  
  endows $H^*(Y,K)$ with the 
flag filtration $F$, abutment of the  spectral sequence associated with
 the filtration by closed subsets $Y_* \subseteq Y$:
$E_1^{p,q}= H^{p+q}(Y_{p}, Y_{p-1}, K_{|Y_{p}}) \Longrightarrow
H^*(Y,K)$. We have $F^p H^*(Y,K)= \ke \,\{ H^*(Y,K) \to H^*(Y_{p-1}, K_{|Y_{p-1}}) \}$.
For an arbitrary $n$-flag, the perverse and flag filtrations
are unrelated. If $Y$ is {\em affine} of dimension $n$  and 
the $n$-flag is obtained using  $n$   hyperplane sections
in sufficiently general position, then
\begin{equation}
\label{1}
P^p H^j(Y,K) = F^{p+j} H^j(Y,K).
\end{equation}

\subsection{Perverse cohomology}
\label{subsecpervco}
The functor $^{\frak p}\!{\cal H}^0: {\cal D}_Y \to {\cal P}_Y$ sends a complex $K$
 to its iterated truncation
$\ptd{0} \ptu{0} K$. This functor is cohomological. In particular, given a
distinguished  triangle
$K' \to K \to K'' \to K'[1]$, one has a long exact sequence
\[\xymatrix{
 \ldots  \ar[r]  &  \pc{j}{K'}   \ar[r] &  \pc{j}{K}  \ar[r] &   \pc{j}{K''} 
 \ar[r]  & \pc{j+1}{K'}  \ar[r] & \ldots
}\]

Kernels and cokernels in ${\cal P}_Y$ can be seen  via perverse cohomology.
Let $f: K \to K'$ be an arrow in ${\cal P}_Y$. View it in ${\cal D}_Y$, cone it and
 obtain a distinguished
triangle
\[ \xymatrix{
K \ar[r]^f &  K' \ar[r]  & \mbox{Cone}(f) \ar[r] & K[1].
} \]
Take the associated long exact sequence of perverse cohomology
$$
\xymatrix{
0 \ar[r]    & 
\pc{-1}{ \mbox{Cone} (f) } \ar[r] & 
 K \ar[r]^{f} & K' \ar[r] &  \pc{0}{ \mbox{Cone} (f) } 
\ar[r] & 0.}
$$
One verifies that ${\cal P}_Y$ is abelian by setting 
$$ 
\ke{\,f} := \, \pc{-1}{\mbox{Cone} (f) }, \qquad 
\coke{\,f}: =  \,  \pc{0}{ \mbox{Cone} (f) } .
$$

\begin{ex}
\label{exdicosest}
{\rm
Consider the natural map $a: \rat_Y[ n] \to  IC_Y$. 
Since $\rat_Y[n] \in \, ^p\!{\cal D}_Y^{\leq 0}$, and $IC_Y$ does not admit 
non trivial subquotients, the long exact sequence of perverse cohomology sheaves
yields the following short exact sequences
$$
\pc{l<0}{ \mbox{Cone}(a) } \simeq \, \pc{l <0}{\rat_Y[n]}, \qquad
0 \to \, \pc{0}{ \mbox{Cone}(a) } \lorw
\, \pc{0}{ \rat_Y [n] } \lorw IC_Y \to 0.
$$
If $Y$ is a normal surface, then $\rat_Y[2]$ is perverse and we are left with
the short exact sequences in ${\cal P}_Y$
$$
\xymatrix{
0 \ar[r] &  \pc{0}{ \mbox{Cone}(a) } \ar[r] &
 \rat_Y [2]  \ar[r]^a & IC_Y \ar[r] & 0.}
$$
By taking the long exact sequence associated with ${\cal H}^j$, one sees that
$\pc{0}{ \mbox{Cone}(a) }$ reduces to a skyscraper sheaf
supported at the singular points of $Y$  in cohomological
degree zero and stalk computed by the cohomology of
the link  of $Y$ at $y:$ ${\cal H}^{-1} (IC_Y)_y = H^1({\cal L}_y)$.
Note that, in general,  the short exact sequence does not split, i.e. $\rat_X [2]$ is not necessarily
a  semisimple perverse sheaf.
}
\end{ex}

\begin{ex}
\label{exscopp}
{\rm ({\bf Blowing up with smooth centers})
Let $X \to Y$ be the blowing up of  a manifold $Y$ along a codimension
$r+1$ submanifold $Z\subseteq Y$.  One has an isomorphism in ${\cal D}_Y:$
$$
f_{*} \rat_X  \; \simeq \; \rat_Y[0]  \oplus  \bigoplus_{j=1}^r \rat_Z [-2j].
$$
If $r+1$ is odd (the even case is analogous and left to the reader), then 
$$\pc{0}{ f_{*} \rat_X [n] }  = \rat_Y [n],\qquad \pc{j}{ (f_{*} \rat_X [n]  }  = \rat_Z [\dim{Z}],
\;\;0 <  | j | \leq r/2. $$
We have three sets  of summands,  i.e.  $(j>0, j=0, j<0)$.
Poincar\'e-Verdier duality exchanges the first and third  sets and fixes the second.
The relative hard Lefschetz theorem  identifies  the first set with the third.
}
\end{ex}

\begin{ex}
\label{smoothmorpervdec}
{\rm ({\bf Families of projective manifolds})
Let $f: X \to Y$ be a family of $d$-dimensional
projective manifolds and let $n := \dim{X}$.
Theorem \ref{delteo} is the cohomological consequence
of a stronger sheaf-theoretic  result (cf. \ci{dess}): there is a direct sum decomposition in ${\cal D}_Y$
$$
f_{*} \rat_X  \; \simeq \;  \bigoplus_{j =0}^{2d} R^j f_* \rat_X [-j].
$$
We have 
\[ \pc{j}{ f_{*} \rat_X [n]  } = R^{d+j}f_* \rat_X [\dim{Y}],  \quad j \in \zed.\]
If we apply Poincar\'e duality and hard Lefschetz to the fibers of $f$
we obtain the following isomorphisms (where the second one is obtained by cupping
with  $c_1(H)^j$, where $H$ is a hyperplane bundle on $X$):
\[\pc{j}{ f_{*} \rat_X [n]  }  \simeq \pc{-j}{ f_{*} \rat_X [n]  }^{\vee}, \;\;\;\forall j \in \zed,
\qquad 
\pc{-j}{ f_{*} \rat_X [n]  } \simeq \pc{j}{ f_{*} \rat_X [n]  }, \;\;\; \forall j \geq 0.\]
}
\end{ex}

\subsection{$t$-exactness and the Lefschetz hyperplane theorem}
\label{subsectex}
The following prototypical Lefschetz-type  result is a consequence
of  the left $t$-exactness  of affine maps (cf. \S\ref{tfidy}).
\begin{pr}
\label{10wl10}
Let $f: X \to Y$ be a proper map, $C \in  \, ^{\frak p}\!{\cal D}^{\geq 0}_{X}$.
Let $Z \subseteq X$ be a closed subvariety,  $U: = X \setminus Z$. 
There is the commutative diagram of maps
$$
\xymatrix{
U \ar[r]^{j} \ar[dr]_{h}& X \ar[d]^{f} & \ar[l]_{i} Z \ar[dl]^{g} \\
&Y. &
}
$$
Assume that $h$ is affine.
Then 
$$\xymatrix{
\pc{j}{f_* C} \ar[r] & \pc{j}{g_{*} i^{*}C}}
$$
is an  isomorphism for $j \leq -2$ and  is a monomorphism  for $j=-1$.
\end{pr}
{\em Proof.} By applying $f_{!}=f_{*}$ to the  distinguished
triangle $j_{!}j^{!} C \to C \to i_{*}i^{*}C \stackrel{[1]}\to$
we get the distinguished  triangle
$$
\xymatrix{
h_{!} j^{*}C \ar[r] &  f_{*} C \ar[r] & 
 g_{*} i^{*} C \ar[r]^{[1]} & .}
$$
Since $h$ is affine, $h_{!}$ is left $t$-exact, so that
$$
\pc{j}{h_{!} j^{*}C} = 0 \qquad \forall j <0.
$$
The result follows by taking the long exact sequence  of perverse cohomology.
\blacksquare

\bigskip
Taking $C= IC_Y$ and $f$ to be the map to a point, and observing that
$i^* IC_Y [-1] = IC_Z$, gives the following Lefschetz hyperplane theorem 
(\ci{goma2} Theorem 7.1) in intersection cohomology.

\begin{tm}
\label{tmwlt}
{\rm ({\bf   Lefschetz hyperplane theorem for intersection cohomology})}
Let $Y$ be an irreducible  projective variety of dimension $n$ and $Z\subseteq Y$ 
be a general hyperplane section.
The restriction 
$$
I\!H^l(Y) \lorw  I\!H^l(Z)
\qquad \mbox{is  an isomorphism for $l\leq n-2$ and  monic for $l=n-1$.}
$$
\end{tm}

\begin{rmk}
\label{8010}
{\rm One has the dual result for the Gysin map in the positive cohomological  degree range.
Similar conclusions 
hold for the cohomology groups  of any  perverse sheaf on $Y$ 
(see \ci{beili}, Lemma 3.3).
}
\end{rmk}

Another related  special case of  Proposition \ref{10wl10},
 used in  \ci{bbd} and in \ci{decmightam}
as  one step towards the   proof of  the relative hard Lefschetz
theorem, 
arises as follows. Let ${\Bbb P} \supseteq X' \to Y'$  be a proper map,
let ${\Bbb P}^{\vee}$ be the dual projective space to ${\Bbb P}$, 
whose points parametrize the hyperplanes in ${\Bbb P}$, 
let  $Z := \{ (x',H)\, | \; x' \in H\} \subseteq  X:= X' \times {\Bbb P}^{\vee}$ be the universal hyperplane
section, $Y:= Y '\times {\Bbb P}^{\vee}$. In this case, 
$IC_Z= i^* IC_X [-1] $. We have

\begin{tm}
\label{8020}
{\rm ({\bf Relative Lefschetz hyperplane theorem})}
The natural map 
$$
\pc{j}{f_* IC_X} \lorw \pc{j+1}{g_* IC_Z}  \quad
\mbox{is an isomorphism for $j \leq -2$ and monic for $j=-1$.}
$$
\end{tm}

\subsection{Intermediate extensions}
\label{subsectie}
A standard reference is \ci{bbd}.
Let $j: U\to  Y$ be  a locally closed embedding 
$Y$ and $i:  \overline{U} \setminus  U =: Z\to Y$.
Given a perverse sheaf  $Q$ on $U$, the {\em intermediate extension}
(often called the ``middle extension") 
$j_{!*}: \pe_U \to \pe_{\overline{U}}$ is a simple operation
that produces  distinguished perverse extensions to  $\overline{U} $ and hence to $Y$.

Intersection complexes are intermediate extensions: 
let $L$ be a local system on  a nonsingular open and dense subvariety
$U$ of an irreducible $d$-dimensional  variety  $Y$; then
$IC_Y(L) = j_{!*} L[d]$.

Let  $Q \in {\cal P}_U$. The natural map
 $
j_! Q \lorw j_{*} Q
$ induces the natural map in perverse cohomology
$a: \pc{0}{j_! Q} \to \, \pc{0}{j_{*} Q}$. 
The {\em intermediate extension} of $Q \in {\cal P}_U$
is the perverse sheaf  
$$j_{!*} Q: = \im{\, (a)} \in  {\cal P}_{\overline U} \subseteq 
{\cal P}_Y. $$
There is the canonical  factorization in the abelian categories
${\cal P}_{\overline U} \subseteq {\cal P}_Y$  
$$
\xymatrix{
\pc{0}{j_! Q} \ar[r]^(.58){epic} &  j_{!*} Q \ar[r]^(.4){monic} & \pc{0}{j_{*} Q}.}
$$
The intermediate extension $j_{!*} Q$
admits several useful  characterizations. For example:
\begin{enumerate}
\item it is   the unique  extension
of $Q \in {\cal P}_U$ to $ {\cal P}_{\overline{U}} \subseteq {\cal P}_Y$ 
with neither  subobjects, nor quotients  supported on $Z;$
\item
 it is the unique  extension $\widetilde{Q}$  of $Q \in  
{\cal P}_U$ to $ {\cal P}_{\overline{U}} \subseteq {\cal P}_Y$ 
such that $i^* \widetilde{Q} \in \,^{\frak p}\!{\cal D}_Z^{\leq -1}$ and 
$i^! \widetilde{Q} \in \,^{\frak p}\!{\cal D}_Z^{\geq 1}$. 
\end{enumerate}

There are an  additional characterization of   and a precise formula 
involving standard truncation and derived push-forwards for 
the intermediate extension functor (cf. \ci{bbd}, 2.1.9 and  2.1.11)
both of which involve stratifications. 
This formula  implies that:
i) if  $j:U\to Y$ is an open immersion of    irreducible varieties and $U$
is  nonsingular of  dimension $d$ then $j_{!*} L[d]$ is canonically isomorphic to $IC_Y(L)$,
ii)
if $Y$ is a nonsingular curve, then 
$j_{!*} L[1] = IC_Y(L) = R^0j_* L [1]$.  

We leave to the reader the task to formulate in precise terms  and verify that
the intermediate extension of an intermediate extension is an intermediate extension.

An intersection cohomology complex, being an  intermediate extension,
does not admit  subobjects or  quotients  supported on proper subvarieties of its support.

The intermediate extension functor $j_{!*}: {\cal P}_U \to {\cal P}_Y$ is not exact
in a funny way. 
Let $0 \to P \stackrel{a}\to Q \stackrel{b}\to R \to 0$ be exact in ${\cal P}_U$.
Recall that $j_!$ is right $t$-exact and that $j_*$ is left $t$-exact.
There is the display
with  exact rows: 
$$
\xymatrix{
\ldots \ar[r] & j_! P  \ar[r] \ar[d]^{epic} & j_! Q \ar[r]^{epic} \ar[d]^{epic} 
& j_! R  \ar[r] \ar[d]^{epic} & 0 \\
              & j_{!*} P \ar[d]^{monic}     & j_{!*}Q \ar[d]^{monic}      & j_{!*} R
\ar[d]^{monic}           \\
0 \ar[r]      & j_* P   \ar[r]^{monic}&  j_* Q \ar[r] & j_* R \ar[r] & \ldots
}
$$   
It is a  simple diagram-chasing exercise  to complete the middle
row functorially with a necessarily monic $j_{!*} (a)$ and 
a necessarily epic $j_{!*} (b)$.
It follows that the intermediate extension functor preserves
monic and epic maps.
What fails is the exactness ``in the middle:''
in general $\ke{ \,  j_{!*}(b) \, } / \im{ \,  j_{!*}(a)  \, } \neq 0$.

\begin{ex}
\label{buffo}
{\rm
Let $E[1]$ be the perverse sheaf on $\comp^*$
discussed in Example \ref{buff}; recall that it fits
in the non split short exact sequence of perverse sheaves:
$$
0 \lorw \rat [1] \stackrel{a}\lorw E[1] \stackrel{b}\lorw \rat [1] \lorw 0.
$$
Let $j: \comp^* \to \comp$ be the open immersion.
We have the commutative diagram of perverse sheaves with exact top and bottom rows:
$$
\xymatrix{
0 \ar[r] & j_! \rat [1] \ar[r] \ar[d]^{epic}  &      j_! E [1] \ar[r]
\ar[d]^{epic}   & j_! \rat [1] \ar[r] \ar[d]^{epic} &  0 \\
         & \rat_{\comp} [1] \ar[r]^(.4){monic}_(.4){j_{!*} (a)} \ar[d]^{monic} &   R^0j_* E [1]  
\ar[r]^(.57){epic}_(.57){j_{!*} (b)}
\ar[d]^{monic}  &  \rat_{\comp} [1] \ar[d]^{monic} &  \\
0 \ar[r] &j_* \rat [1] \ar[r]  &      j_* E [1] \ar[r]   &j_* \rat [1] \ar[r] &  0.
}
$$
The middle row, i.e. the   one of middle extensions,
 is not exact in the middle. In fact,   inspection of the stalks at the origin
yields the non exact sequence 
$$
\xymatrix{
0 \ar[r] &  \rat \ar[r]^(.45){\simeq} &  \rat \ar[r]^{0}&  \rat \ar[r] & 0.}
$$
This failure prohibits exactness in the middle.
The inclusion $\im j_{!*}(a) \subseteq \ke j_{!*}(b)$ is strict: 
$K:= \ke{ \, j_{!*} (b) \, }  $ is the unique non trivial extension,
$\mbox{Hom} ( \rat_{ \{ 0 \} }, \rat_{\comp}[2]) = \rat$,
$$
\xymatrix{
0 \ar[r] & \rat_{\comp}[1] \ar[r] &  K \ar[r] & \rat_{ \{0 \} } \ar[r] & 0.}
$$
The reader can check,  e.g. using the  self-duality of $E$, that  $K^{\vee} = K'$
($K'$ as in Ex. \ref{buff}).
}
\end{ex}

\medskip

Property 1, characterizing intermediate extensions, is  used in the construction
of composition series for perverse sheaves in Proposition \ref{pyart}.
If follows that $j_{!*} Q $ is  a simple perverse sheaf on $Y$  iff $Q$ is
a  simple perverse sheaf on $U$.

\begin{ex}
\label{nomap}
{\rm ({\bf Intersection cohomology complexes
with different supports})
Let $Z_1, Z_2 \subseteq Y$ be irreducible closed subvarieties
with $Z_1 \neq Z_2$ 
(note that  we are allowing $Z_1 \cap Z_2 \neq \emptyset$).
Let $IC_{Z_i}(L_i)$, $i=1,2$ be intersection cohomology complexes. Then
(cf. \ci{goma2}, Theorem  3.5)
\[\mbox{Hom} (IC_{Z_1}(L_1), IC_{Z_2} (L_2) )=0.\]
In fact, the kernel (cokernel, resp.) of any such map would have to be either
zero, or supported on $Z_1$ ($Z_2$, resp.), in which case, it is easy to 
conclude by  virtue of characterization 1 given above.
}
\end{ex}

Here is a nice application of what above. Let $f: X \to Y$ be a proper and semismall map
of irreducible proper varieties; see $\S$\ref{semismall}. The decomposition theorem
yields a (canonical in this case) splitting
$$
f_* IC_X = \bigoplus IC_{Z_a} (L_a).
$$
Poincar\'e duality on $IC_X$ yields  a canonical isomorphism 
$e: f_* IC_X \simeq (f_* IC_X)^{\vee}$ which,  by Example \ref{nomap}, 
 is a direct sum map. It follows that the direct 
summands $I\!H^*(Z_a, L_a) \subseteq I\!H^*(X)$ are mutually orthogonal
with respect to the Poincar\'e pairing.

\section{Three approaches to the decomposition theorem}
\label{3stooges}
\subsection{The proof of Beilinson, Bernstein, Deligne and Gabber}
\label{8000}
The original proof \ci{bbd}
of the decomposition theorem for proper maps of complex algebraic varieties uses 
in an essential way
the language of the 
\'etale cohomology  of $l$-adic sheaves and the arithmetic properties of varieties defined over finite fields.
  
   In this section we try to introduce the reader
to some of the main ideas in \ci{bbd}. 

Let us first give a very brief and rough
summary of these ideas. 
The theory of weights, i.e. of the eigenvalues of
the Frobenius automorphisms on the stalks of $l$-adic sheaves
on varieties defined over finite fields, leads to the  notion of    pure complexes.
There are many pure complexes:
O. Gabber  proved that
the intersection cohomology complex of a variety is a pure perverse sheaf.
 The push-forward via a proper map of algebraic varieties
 defined over a finite field of a pure complex is a pure complex. 
 After passing to an   algebraic closure of the finite field,
a pure complex splits as direct sum
of shifted intersection 
complexes with   coefficients in  lisse irreducible sheaves
(a lisse $l$-adic sheaf is the $l$-adic analogue
of a local system). We thus obtain the decomposition
theorem for the proper push-forward of a pure complex,
e.g. the intersection complex of a variety, at least after
passing  to the algebraic closure of the finite field.

Associated with a map of  complex algebraic varieties  there are
companion  maps of  varieties defined over   finite fields.   
There is the  class of constructible complexes of geometric origin
over complex varieties. A complex of geometric origin
over a complex variety admits  $l$-adic counterparts
on the companion varieties defined over the finite fields.
The intersection complex is of geometric origin.

The decomposition result  over the algebraic closures of the finite fields
 is   shown to imply the analogous result in  (i.e. it lifts to)
 the   complex algebraic setting and we  finally obtain
the decomposition theorem in the complex setting.

The  idea that results over finite fields can be used to prove
results over the complex numbers  is rooted in  the classical result that
a  system of rational polynomial equations has a solution
over an algebraic number field if   it has a solution
modulo an infinite number of prime numbers. 
 
There are several  appearances of this idea  in the literature, often
in connection with a beautiful discovery.
Here are few: P. Deligne-D. Mumford's proof  \ci{de-mumf} that 
the moduli space
of curves of a given genus is irreducible
in any characteristic,  
S. Mori's proof \ci{mori} of  Hartshorne's conjecture, P. Deligne and L. Illusie's
algebraic proof  \ci{de-illusie} of the Kodaira vanishing theorem and of the degeneration
of  the Hodge to de Rham spectral sequence
 (see  the nice survey \ci{illusie}). 
 
 A precursor of the techniques used in lifting
 the decomposition theorem from finite fields to the complex numbers
 is P. Deligne's  proof  (\ci{weil2}) of 
 the hard Lefschetz
 theorem. 
 
 We do not discuss further the ``lifting" technique and we refer  the reader 
 to \ci{bbd}, $\S$6.

\medskip
The goal of the remaining part of this section is
to introduce the reader to constructible
$\qlb$-sheaves ($\S$\ref{sscqls}), weights,  pure complexes 
and their structure ($\S$\ref{sswei}, $\S$\ref{sssopuco}),
to discuss the decomposition, semisimplicity and hard Lefschetz
 theorems 
in the context of pure complexes over finite fields and over their algebraic closures
($\S$\ref{tdecthof}), and  to   state the decomposition theorem etc. for complexes of geometric origin on complex algebraic varieties ($\S$\ref{dafac}).

We hope that
our stating separately the results over finite fields,  over their algebraic closures and over the complex numbers
may help the reader better understand the whole picture and perhaps  justifies the
tediousness of these repetitions.

\medskip
Let us fix some notation.
A variety over a field is a separated scheme of finite type over that field.
For a quick summary on quasi projective varieties (which is all we need here) see $\S$\ref{algvar}.
Let $\fq$ be a finite field, let $\fb$ be a fixed algebraic closure of $\fq$ and
let $\mbox{Gal} (\fb/\fq)$ be the Galois group. 
This group is profinite,  isomorphic to
the profinite completion of $\zed$, and it admits as topological generator
the geometric Frobenius $Fr:= \varphi^{-1}$, where $\varphi: \fb \to \fb$,
$t \mapsto t^q$ is the arithmetic Frobenius (see Remark \ref{didk}).
Let $l\neq 
\mbox{char} \, \fq$ be a fixed prime number, let  $\zed_l$ be the
ring of $l$-adic integers, i.e. the projective limit
of the system $\zed/ l^n \zed$ (abbreviated by $\zed/l^n$),  let  $\ql$ be the $l$-adic numbers,
i.e. the quotient field of $\zed_l$, and let  $\qlb$ 
be a fixed algebraic closure
of $\ql$. Recall that $\zed_l$ is uncountable and that $\qlb \simeq \comp$, non canonically.

\subsubsection{Constructible $\qlb$-sheaves}
\label{sscqls}
Let   $X_0$ be an algebraic  variety  defined over a finite field $\fq$. 
We refer to \ci{bbd, weil2}, and to the introductory \ci{freitag}, $\S$12,  for the definitions of the
category and $D^b_c(X_0, \qlb)$ of constructible complexes
of  $\qlb$-sheaves. These categories are stable under the usual operations
$f^*, f_*, f_!, f^!$, derived Hom and tensor product, duality and vanishing and nearby cycles.
With some homological restrictions on Tor groups, the standard and  the 
middle perverse $t$-structure  are also defined, and one obtains the category
${\cal P}(X_0,\qlb)$ of perverse sheaves on $X_0$.
If $X$ is the ${\Bbb F}$-variety obtained
from $X_0$ by extending the scalars to ${\Bbb F}$, then we obtain
in the same way the categories $D^b_c (X, \qlb)$ and ${\cal P}(X, \qlb)$
which are also stable under the usual operations mentioned above.

The construction of these categories and functors and the verification of their
fundamental properties 
 requires a massive background
(a large part of Grothendieck et al. S.G.A. seminars is devoted to this task)
and
has lead P. Deligne to complete
the proof of the Weil Conjectures (\ci{weil1}), one of the crowning achievements
of $20^{\rm th}$ century mathematics. 

For the purpose of this survey, let us just say that
we will  mostly think of   $D^b_c(X_0, \qlb)$ etc., 
by   analogy  with the perhaps more geometric   constructible derived categories
${\cal D}_{X}$ 
associated with complex varieties. There is one important difference:
the  action of the Frobenius automorphism.

\subsubsection{Weights and purity}
\label{sswei}
In positive characteristic, the \'etale cohomology of algebraic varieties
presents a feature that is absent in characteristic zero: the eigenvalues
of Frobenius, i.e. {\em weights}.

Let $X_0$ be a variety over the finite field $\fq$.
Suppression of the index $-_0$
denotes extension of scalars from $\fq$ to $ \fb = \overline{\Bbb F}_q$.
For example,  if $F_0$ is a $\qlb$-sheaf on $X_0$,
then we denote its pull-back to $X$ by $F$.

To give a $\qlb$-sheaf $F_0$ on the one-point variety $\mbox{Spec} \, \fq$
is equivalent to giving a 
finite dimensional continuous $\qlb$-representation
of the Galois group $\mbox{Gal} (\fb/\fq)$. The pull-back
$F$
to $\mbox{Spec} \, \fb$ is the sheaf given by the
underlying $\qlb$-vector space of the representation (i.e. we ``forget"
the representation; this is because the Galois group
$\mbox{Gal}(\fb/\fb)$ is trivial).
This is called the {\em stalk} of $F_0$ at the point. 

It is important
to keep in mind that the sheaf $F_0$ on ${\rm Spec}\, \fq$ must be thought of as
the pair given by the vector space {\em and} the representation, while its
pull-back $F$ to ${\rm Spec}\, \fb$ is just the datum
of the vector space. This partially explains
why  the decomposition theorem
holds over the algebraic closure $\fb$, but not necessarily over the finite field
$\fq$, where the splittings  have to be compatible
with the Frobenius action. 

There are restrictions
on the  representations arising in this context:
e.g. in the case of a $\ql$-sheaf of rank one on $\mbox{Spec} \, \fq$,
keeping in mind that the Galois group is compact,
continuity implies  that $Fr \in \mbox{Gal} (\fb/\fq)$ must act by units
in $\zed_l \subseteq \qlb$. 

\begin{rmk}
\label{didk}
{\rm
It is often useful to keep in mind the following
roughly approximated  picture when thinking about
the extension $\fb/\fq$: think of the one point variety   $\mbox{Spec} \,\fq$ as being
a  circle $S^1$;   think of the extension
$\fb/\fq$ as being the universal covering space  $\real \to S^1$
with deck group  given by translations by integers;  think of the Galois group
as being the deck group; given an $l$-adic sheaf $F_0$
on the one point variety  $\mbox{Spec}\, \fq$, think of 
the action of Frobenius  on the stalk of this sheaf as the
$\zed$-action on a sheaf  on $\real$  pull-back of a sheaf  on $S^1$. 
}
\end{rmk}

 For every $n \geq 1$, the finite set $X_0 ( {\Bbb F}_{q^n})$  of closed points in $X_0$ which
are defined over
the  degree $n$ extension 
$\fq \subseteq {\Bbb F}_{q^n}$ 
is precisely the set of closed points (we are using the Zariski topology) of $X$
which are fixed under the action of the $n$-th iterate, $Fr^n: X \to X$,
 of the geometric Frobenius
$Fr: X \to X$. 
Recall that if, for example,  $X_0$ is defined by a system of  polynomials
$\{P_i(T)\}$
in $\fq[T_1, \ldots, T_N]$, then 
a  closed point of $X_0$ defined over ${\Bbb F}_{q^n}$ can be identified with
an $N$-tuple $(a_1, \ldots, a_N) \in {\Bbb F}_{q^n}^N$ which is a   solution
of the system of polynomial equations $P_i(T)=0$.

Let $x \in X_0({\Bbb F}_{q^n})$ be such a
$Fr^n$-fixed point. The $\qlb$-sheaf $F_0$
restricted to $x$ has stalk the $\qlb$-vector space $F_x$
on which $Fr^n$ acts as an automorphism.

\begin{defi}
\label{pupu}
{\rm ({\bf Punctually pure})
The $\qlb$-sheaf
 $ F_{0}$  on $X_0$ is {\em punctually pure of weight} $w$ ($w \in \zed$)
if,  for every $n \geq 1$ and every  $x \in X_0 ({\Bbb F}_{q^n})$,
the eigenvalues of the action of   $Fr^n$  on $F_{x}$ are 
algebraic numbers
such  that all of  their complex algebraic conjugates 
 have absolute value $q^{n\, w/2}$.
}
\end{defi}

For example, on $\mbox{Spec}  \, \fq$,  the sheaf $\qlb$ is has weight $0$, while
the Tate-twisted 
$\qlb(1)$ has weight $-2$. If $X_0$ is a nonsingular projective curve
of genus $g$,
then  the \'etale cohomology group $H^1_{et}(X_0, \qlb)$ can be viewed as
an $l$-adic sheaf on  $ {\rm Spec}\, \fq$ with weight $1$.

The eigenvalues of Frobenius are naturally  elements of $\qlb$.
While  $\qlb\simeq \comp$, there is no natural
isomorphism between them. 
However, since $\rat \subseteq \qlb$, it makes sense to
 request that the eigenvalues are algebraic numbers (i.e. their being
 algebraic is independent of the choice of an isomorphism
 $\qlb\simeq \comp$). Once a number is algebraic,
 the set of its algebraic conjugates  is well-defined independently of
 a choice of an 
 isomorphism $\qlb \simeq \comp$, and this  renders meaningful the request 
on the absolute values.
This is a strong request:
$1+\sqrt{2}$ and 
$1-\sqrt{2}$ are  algebraic conjugates, however, they 
 have different  absolute values.

\begin{defi}
\label{mixeds}
{\rm ({\bf Mixed sheaf, weights})
A $\qlb$-sheaf $ F_{0}$ on $X_{0}$
is {\em mixed} if it admits a finite filtration with  
punctually pure successive quotients.
The {\em weights} of  a mixed  $F_{0}$ are  the weights of 
the non-zero quotients.
}
\end{defi}

\begin{defi}
\label{mixdbc}
{\rm ({\bf Mixed  and pure complexes})
The category $D^b_m(X_0, \qlb)$ of 
{\em mixed complexes} is the full subcategory of
$D^b_c(X_0, \qlb)$ given by those complexes whose
cohomology sheaves are mixed.
A complex $K_0 \in D^b_m(X_0, \qlb)$ is {\em pure
of weight $w$} if the cohomology sheaves ${\cal H}^i (K_0)$ are punctually pure of weights
$\leq w+i$ and the same is true for its Verdier dual $K_0^{\vee}$.
}
\end{defi}

The following theorem is proved in \ci{bbd} (see $\S$3.3.1 and $\S$6.2.3) and is a key step 
towards the proof of  the decomposition theorem given in
\ci{bbd}.  Note that
the special case when $X_0$ is nonsingular and projective and $Y_0 =
\mbox{Spec} \, \fq$ yields a proof of  the main result in \ci{weil1}, i.e.
the completion of the proof of the Weil conjectures. 
\begin{tm}
\label{dipp}
{\rm ({\bf Purity for  proper maps or relative Weil conjectures})}
Let $K_0$ be pure of weight $w$ and $f_0:X_0 \to Y_0$ be a  proper map
of $\fq$-varieties.
Then
${f_0}_* K_0$ is pure of weight $w$.
\end{tm}

\subsubsection{The structure of pure complexes}
\label{sssopuco}
In this section we state O. Gabber purity theorem and 
discuss the special splitting  features of pure complexes.

The following result of O. Gabber \ci{gabber} was never published. A proof 
appears 
in \ci{bbd}, Corollaire 5.4.3
and it is  summarized in \ci{brili}. This result makes it 
clear that the class of pure complexes contains  many
geometrically relevant objects.

Recall that lisse $\qlb$-sheaves are the $\qlb$-analogues of local systems
in the classical topology.

\begin{tm}
\label{pic}
{\rm ({\bf Gabber purity theorem})}
The intersection cohomology complex $IC_{X_0}$ of a  connected 
pure $d$-dimensional variety ${X_0}$ is pure
of weight $d$.
More generally, if $L$ is a pure lisse $\qlb$-sheaf
of weight $w$
on a connected,  pure $d$-dimensional subvariety $j:Z_0 \to X_0$, 
then $IC_{\overline{Z}_0}(L):= j_{!*} L[d]$ is a pure perverse sheaf
of weight $ w +d$.
\end{tm}

\medskip
The following result (\ci{bbd}, Corollaire 5.3.4) generalizes 
Gabber's Purity theorem and  is  another    key  step in the
proof  in \ci{bbd} of the  decomposition, semisimplicity and relative hard Lefschetz theorems over the complex numbers.

\begin{tm}
\label{mmm}
{\rm ({\bf Mixed and simple is pure})}
Let $P_0 \in {\cal P}_m (X_0, \qlb)$ be a simple  mixed
perverse $\qlb$-sheaf. Then $P_0$ is pure.
\end{tm}

The following  theorem summarizes the basic splitting properties of pure complexes.
The proofs can be found in
\ci{bbd},   Th\'eor\`emes 5.4.1, 5.4.5 and  5.4.6, and Corollaire  5.3.8.

\begin{tm}
\label{decviaf} {\rm  ({\bf Purity and decompositions})}

\begin{enumerate}
\item
Let $K_0 \in D^b_m(X_0, \qlb)$ be pure of weight $w$. 
Each $\pc{i}{K_0}$ is a pure  perverse sheaf of weight $w+i$.
There is an isomorphism
in $D^b_c(X, \qlb)$
$$K \simeq  \bigoplus_i \pc{i}{K}[-i]. $$

\item
Let $P_0 \in {\cal P}_m(X_0, \qlb)$  be a pure perverse 
$\qlb$-sheaf on $X_0$.
The  pull-back $P$ to $X$ splits
in ${\cal P} (X, \qlb)$
 as a direct sum
of intersection cohomology complexes associated with
lisse irreducible   sheaves on subvarieties of  $X$. 
\end{enumerate}
\end{tm}

\begin{rmk}
\label{tt00}
{\rm The splittings above do not necessarily hold over $X_0$.
 }
\end{rmk}

If $K_0 \in D^b_c (X_0,\qlb)$, then the cohomology groups
$H^*(X, K)$ on $X$  are finite dimensional
$\qlb$-vector spaces with  a continuous $\mbox{Gal} (\fb/ \fq)$-action
and one can speak about the weights of $H^*(X, K)$, so that
 the notions of weights and  purity
extend to this context. In particular, this applies to the Ext-groups below.
 
We would like to give the reader a feeling of why  weigths
are related to splitting behaviors. These behaviors are governed
by the Ext groups. Let $K_0, L_0 \in D^b_m (X_0, \qlb)$.
The natural map $\mbox{Ext}^1(K_0, L_0) \to \mbox{Ext}^1
(K,L)$ factors through the  space of Frobenius invariants $\mbox{Ext}^1(K,L)^{Fr}$
which is of pure of  weight zero. If $K_0$ has weights $\leq w$ and $L_0$
has weights $\geq w'$, then $\mbox{Ext}^1 (K,L)$ has weights
$\geq 1 +w' -w$. If $w'=w$,  then  $\mbox{Ext}^1 (K,L)$ has weights
$\geq 1 $, so that $\mbox{Ext}^1(K,L)^{Fr}$ is trivial. 
The upshot is that given the right weights, 
 a nontrivial extension class over $\fq$ must become
 trivial over $\fb$  and splittings may  ensue (but only
over the algebraic closure).

\subsubsection{The decomposition  over $\fb$}
\label{tdecthof}
With Theorem \ref{decviaf} in hand, it is immediate to prove
the following theorem, which is one of the main results in \ci{bbd}.

\begin{tm}
\label{puredt} {\rm  ({\bf Decomposition theorem and semisimplicity over $\fb$})}
Let $f_0: X_0 \to Y_0$ be a proper morphism of $\fq$-varieties,
$K_0 \in D^b_c(X_0, \qlb)$ be pure and $f:X\to Y$ and $K$ 
be the corresponding data over $\fb$. There is an isomorphism in $D^b_c(Y, \qlb)$
\begin{equation}
\label{splf}
f_*K \;  \simeq \; \bigoplus_i \pc{i}{f_*K} [-i],
\end{equation}
where each $\pc{i}{f_*K}$ splits as a direct sum of intersection
cohomology complexes associated with lisse irreducible sheaves on
subvarieties of  $Y$.
In particular, $f_*K$ is semisimple, i.e. the unshifted
summands $\pc{i}{K}$ are semisimple perverse sheaves on $Y$.
\end{tm}

We now turn our attention to the relative hard Lefschetz,
also proved in \ci{bbd}.
Let  $f_0: X_0 \to Y_0$ be a morphism of $\fq$-varieties, 
$\eta_0$ be the first Chern class of a line bundle $\eta_0$
on $X_0$. This defines a natural transformation $\eta_0: {f_0}_* \to {f_0}_* [2](1)$.
Here $(1)$ is the Tate twist, lowering the weigths by two; the reader unfamiliar
with this notion, may ignore the twist and still get a good idea of the meaning of the
statements. By iterating, we obtain maps
 $\eta_0^i: {f_0}_* \to {f_0}_* [2i](i)$, $i \geq 0$. In particular, it defines
natural transformations $\eta_0^i:  \pc{-i}{{f_0}_* (-)} \to \pc{i}{{f_0}_* (-)}(i)$.

\begin{tm}
\label{rhlet}
{\rm ({\bf Relative hard Lefschetz over $\fq$ and $\fb$})}
Let $P_0$ be a pure perverse sheaf on $X_0$. Assume that $X_0$ is quasi projective
and that $\eta_0$ is a hyperplane bundle. Then the iterated cup product operation
induces isomorphisms
$$
\eta_0^i :  \pc{-i}{{f_0}_* P_0}   \stackrel{\simeq}\lorw  \pc{i}{{f_0}_* P_0}(i), \qquad \forall i\geq 0.
$$
The same holds over ${\fb}$ (with the understanding that $P$ should come from a $P_0$).
\end{tm}

\begin{rmk}
\label{equello}
{\rm
The case $Y_0=pt$,  $P_0 = IC_{X_0}$, yields   the hard Lefschetz theorem
for intersection cohomology (over $\fb_0$ and over $\fb$). Using the same technique
``from ${\fb}$ to $\comp$" in  \ci{bbd}, $\S$6,  one sees that theorem \ref{rhlet} implies the 
 hard Lefschetz theorem for the intersection cohomology
of complex projective varieties. 
An important precursor of the relative Hard Lefschetz
theorem is P. Deligne's algebraic proof in 
\ci{weil2} of the classical hard Lefschetz theorem.}
\end{rmk}

\subsubsection{The decomposition theorem for complex varieties}
\label{dafac}
The technique ``from $\fb$ to $\comp$'' is used in
\ci{bbd}, $\S$6 to deduce the results
of this section on complex algebraic varieties,
from
the results of  the previous  $\S$\ref{tdecthof} on varieties defined
over finite fields.

 Let $X$ be a complex variety. Consider the categories ${\cal D}_X$
of bounded constructible complexes of sheaves of {\em complex}  vector spaces 
and its full sub-category of complex perverse sheaves ${\cal P}_X$.
Recall that every perverse sheaf admits a finite filtration with simple quotients
called the constituents of the perverse sheaf.

\begin{defi} 
{\rm
({\bf Perverse sheaves of geometric origin})
 A perverse sheaf $P \in {\cal P}_X$  is said to be of {\em geometric origin}
if it belongs to the smallest set 
such that

\smallskip
\n
(a) it contains the constant sheaf $\comp_{pt}$ on a point,

\smallskip

and that is stable under the following operations

\smallskip
\n
(b) for every map  $f$, take the
simple constituents of 
$\pc{i}{T(-)}$, where $T= f^{*},f_{*}, f_!, f^!$,

\n
(c) take the simple constituents  of
$\pc{i}{ -\otimes -}$, $\pc{i}{R{\cal H}om (-,-)}$.
}
\end{defi}

As a  first example on a variety $Z$ one may
start with the map $g:Z \to pt$,
take $g^{*} \comp_{pt} = \comp_Z$, and
set $P$ to be any  simple
constituent
of  one of the perverse complexes $\pc{i}{\comp_Z}$. 
If $f: Z \to W$ is a map, one can take a simple constituent
of $\pc{j}{f_{*} P}$ as an example on $W$.
Another example consists of taking a simple local system
of geometric origin $L$ on a connected and smooth 
Zariski open subvariety $j: U \to X$
and setting $P := j_{!*} L [\dim{U}]$. 
Using either construction, we verify immediately
 that   the intersection cohomology complex  
 of a variety is  of geometric origin. 

\begin{defi}
\label{ssgops}
{\rm
({\bf Semisimple complexes  of geometric origin})
A perverse sheaf $P$ on $X$ is said to be {\em semisimple 
of geometric origin} if it is a direct sum of simple perverse sheaves
of geometric origin.
A  constructible complex  $K \in {\cal D}_X$ is said to be {\em semisimple of geometric origin}
if there is an isomorphism $K \simeq \oplus \pc{i}{K}[-i]$
in ${\cal D}_X$
and each perverse cohomology complex $\pc{i}{K}$ 
is semisimple of geometric origin.
}
\end{defi}

We can now state the decomposition theorem and
the relative hard Lefschetz theorems as they are   stated and proved in
\ci{bbd}. If $X$ is irreducible, then 
 $IC_X$ is  simple of geometric origin so that   the two theorems below apply
 to $K= IC_X$. The proofs can be found in \ci{bbd}, Th\'eor\`emes 6.2.5, 6.2.10.
 Note that while the results  proved there are for sheaves of
$\comp$-vector spaces, one can deduce easily the variant for
sheaves of  $\rat$-vector spaces.

\begin{tm}
\label{dtgo}
{\rm ({\bf Decomposition theorem over $\comp$})}
Let $f: X \to  Y$ be a proper morphism of complex varieties.
If  $K \in {\cal D}_X$ is  semisimple of geometric origin, then so is $f_*K$.
\end{tm}

\begin{tm}
\label{rrhlco}
{\rm ({\bf Relative hard Lefschetz theorem over $\comp$})}
 Let $ f : X \to Y$
be a projective morphism, $P$ a perverse sheaf on $X$ which is semisimple
of  geometric origin,
$\eta$ the first Chern class of an $f$-ample line bundle on $X$. Then
the iterated cup product operation induces isomorphism
$$
\eta^i \;:\;  \pc{-i}{f_* P}\; \stackrel{\simeq}\lorw \; \pc{i}{f_* P}, \qquad \forall \, i \geq 0.
$$
\end{tm}

\subsection{M. Saito's approach via mixed Hodge modules}
\label{ms}
The authors of \ci{bbd} (cf. p.165) left open   two questions:
 whether the decomposition theorem holds for the push forward of
the intersection cohomology complex of a local system underlying a polarizable variation of 
pure Hodge structures
and whether it holds in the K\"ahler context. 
(Not all local systems as above are of geometric origin.)

In his remarkable work on the subject, M. Saito  answered the first question in 
the affirmative  in \ci{samhp} and the second question in the affirmative
in the case of $IC_X$  in \ci{satohoku};
we refer the reader to M. Saito's paper for the precise formulations in the K\"ahler context. In fact, he  developed in \ci{samhm} a  general theory
of compatibility of mixed Hodge structures with the various functors, and in the process
he completed the extension of the Hodge-Lefschetz theorems from
the cohomology of projective manifolds, to the intersection cohomology
of projective varieties.

There are at least two important new ideas in his work.
The former is that the Hodge filtration is to be obtained 
by a filtration at level of   $D$-modules.  A precursor of this idea
is Griffiths' filtration by the order of the pole.
The latter is that the properties of his  mixed Hodge modules are defined and tested
using the vanishing cycle functor.

Saito's  approach is  deeply rooted in the theory of $D$-modules
and, due to our  ignorance on the subject, it will not be
explained here.  We refer to Saito's papers \ci{samhp,samhm,satohoku}.
For a more detailed overview, see \ci{bry-zu}.
The papers \ci{saitointro} and  \ci{durfeesaito} contain brief summaries of the results of the
theory.  See also \ci{ms}.

Due to the  importance of these results,
we would like to discuss   very informally Saito's achievements in the hope
that even a  very rough outline can be helpful to some.
For simplicity only, we restrict ourselves to complex algebraic varieties
(some results hold  for complex analytic spaces).

Saito has constructed, for every variety $Y$, an  abelian category 
MHM($Y$) of mixed Hodge modules on $Y$.
The construction  is a tour-de-force which
uses  induction on  dimension via a systematic 
use of  the vanishing cycle functors
associated  with germs of holomorphic maps. 
It is in the derived category $D^b (\mbox{MHM}(Y))$ that
Saito's results on mixed Hodge structures can be stated and proved.
If one is interested only  in the decomposition and relative hard Lefschetz theorems,
then it will suffice to work with  the categories $MH(Y,w)$ below.

One starts with the abelian and semisimple category of polarizable Hodge modules
of some weigth $MH(Y,w)$.
Philosophically they correspond to perverse pure complexes in
$\qlb$-adic theory.
Recall that, on a smooth variety, the Riemann-Hilbert correspondence
assigns to a regular holonomic $D$-module a perverse sheaf with complex coefficients.
Roughly speaking, the simple objects are certain filtered
regular holonomic $D$-modules
$({\cal M}, F)$.
The $D$-module ${\cal M}$ corresponds, via an easy extension  of the
Riemann-Hilbert correspondence to singular varieties, to the intersection cohomology complex
of the complexification of a rational local system underlying a polarizable simple  variation of
pure Hodge structures of some weight (we omit the bookkeeping of weights).

Mixed Hodge modules correspond philosophically
to perverse mixed complexes and  are, roughly speaking,
certain  bifiltered regular holonomic $D$-modules $({\cal M}, W, F)$
with the property that the graded objects $Gr^W_i {\cal M}$ are polarizable
Hodge modules of weight $i$. The resulting abelian category
$\mbox{MHM}(Y)$ is not semisimple. However, the extensions
are not arbitrary, as they are controlled by the vanishing cycle functor.
The extended Riemann-Hilbert correspondence
assigns to the pair $({\cal M}, W)$ a filtered perverse sheaf
$(P, W)$ and this data extends to a functor of $t$-categories
$$
\frak{r} \; : \;  D^b (\mbox{MHM}(Y)) \lorw {\cal D}_Y,
$$
with the standard $t$-structure on  
$D^b (\mbox{MHM}(Y))$
and the perverse
$t$-structure on ${\cal D}_Y$.
Beilinson's equivalence theorem \ci{beili}, i.e. ${\cal D}_Y \simeq D^b (\pe_Y)$,
  is used here, and in the rest of
this theory, in an essential way.
In fact, there is a second $t$-structure, say $\tau'$, on  $D^b (\mbox{MHM}(Y))$ corresponding
to the standard one on  ${\cal D}_Y;$ see \ci{samhm}, Remarks 4.6.

The usual operations on $D$-modules induce a collection of operations 
on $D^b(\mbox{MHM}(Y))$ that
correspond to the usual operations on the categories
 ${\cal D}_Y$,  i.e. $f^*, f_*, f_!, f^!$, tensor products,
Hom, Verdier duality, nearby and vanishing cycle functors
(cf. \ci{samhm}, Th. 0.1).

In the case when $Y$ is a  point, the category
$\mbox{MHM}(pt)$ is naturally equivalent to the category of  graded polarizable
rational  mixed Hodge structures (cf. \ci{samhm}, p.319); here ``graded" means
that one has  polarizations on the graded pieces of the weight filtration.
At the end of the day, the $W$ and $F$ filtrations produce
two filtrations on the cohomology and on the cohomology with compact supports
of a complex in the  image of $\frak r$
and give rise to mixed Hodge structures compatible with the usual
operations. 
The functor $\frak r$ is exact and faithful,  but not fully faithful (the map on Hom sets is
injective, but  not surjective), not even
over a point:
in fact, a  pure Hodge structure of weigth $1$ and rank $2$,
e.g. $H^1$ of an elliptic curve, is irreducible as a Hodge structure, but not as
a vector space.

The constant sheaf $\rat_Y$ is  in the image
of the functor $\frak r$ and 
Saito's theory recovers Deligne's functorial mixed Hodge theory
of complex varieties  \ci{ho2,ho3}.
 See \ci{samhm}, p. 328  and \ci{saitomathann}, Corollary 4.3.

As mentioned above, mixed Hodge modules are a Hodge-theoretic analogue
of the arithmetic mixed perverse sheaves discussed in $\S$\ref{8000}.
A mixed Hodge module  $({\cal M}, W, F) \in MHM(Y)$ is said to be pure of weight
$k$ if $Gr^W_i {\cal M} =0$, for all $i \neq k$. In this case it is, by definition, 
a polarizable Hodge module so that
a   mixed Hodge module which is of some pure
weight is analogous  to an arithmetic  pure perverse sheaf.

Saito proves the  analogue of the arithmetic
Corollary  \ref{dipp}, i.e.  if $f$ is proper, then $f_*$ preserves weights.
Though the context and the details are vastly different, 
the rest of the story unfolds by analogy with the arithmetic case 
discussed in $\S$\ref{8000}. 
A complex in $D^b ( \mbox{MHM} (Y)$)
is said to be semisimple 
if it is a direct sum
of shifted mixed   Hodge modules which are simple and pure
of some weight (= polarizable Hodge modules, i.e.
associated with a  simple variation of polarizable   pure Hodge structures).

In what follows, note that the  faithful functor $\frak r$
commutes, up to natural equivalence,  with the usual operations, e.g.
${\frak r} ({\cal H}^j (M)) =
\pc{j}{  {\frak r} (M) }$, $ f_* ( {\frak r} (M) ) = {\frak r} ( f_* (M))$.

\begin{tm}
\label{sa1}
{\rm ({\bf Decomposition theorem for polarizable  Hodge modules})}
Let $f: X \to Y$ be proper and $M \in D^b ( \mbox{MHM} (X) )$
be semisimple.
Then
the direct image  $f_* M\in D^b ( \mbox{MHM} (Y) )$ is semisimple. More precisely, 
if $M \in \mbox{MHM}(X)$ is semisimple and pure, then
$$
f_* M  \; \simeq  \; \bigoplus_{j \in \zed} {\cal H}^j (f_* M) [-j]
$$
where 
the ${\cal H}^j (f_* M) \in \mbox{MHM} (Y)$ are semisimple and  pure.
\end{tm}
\begin{tm}
{\rm ({\bf Relative hard Lefschetz for  polarizable Hodge modules})}

\n
Let $f: X \to Y$ be projective, $M\in  \mbox{MHM} (X) $
be semisimple and pure and $\eta \in H^2(X, \rat)$ be the first  Chern class
of an $f$-ample line bundle on $X$. 
Then the  iterated cup product map is an isomorphism
$$
\eta^j \, : \, {\cal H}^{-j} (f_* M ) \stackrel{\simeq}\lorw 
{\cal H}^j ( f_* M )
$$
of semisimple and pure mixed Hodge modules. 
\end{tm}

The proof relies on  an  inductive use, via Lefschetz pencils,  of 
S. Zucker's \ci{zucker}
results on Hodge theory for degenerating coefficients
in one variable.

The  intersection cohomology complex 
of a polarizable variation of pure Hodge structures is
the perverse sheaf associated with a pure mixed Hodge module (= polarizable Hodge module).
This fact is not as automatic as in the case of the constant sheaf, for it requires
the verification of the conditions of vanishing-cycle-functor-type involved in the definition
of the category of polarizable Hodge modules. One may view this
fact as the analogue of Gabber's purity theorem \ref{pic}.

M. Saito thus  establishes   the decomposition and the relative
hard Lefschetz theorems  for  coefficients in the intersection cohomology
complex $IC_X (L)$
of a polarizable variation of pure Hodge structures,
 with the additional fact
that one has mixed Hodge structures on the cohomology of the  summands on $Y$
and that the (non-canonical) splittings on the intersection cohomology group $I\!H(X, L)$
are  compatible with the mixed Hodge structures of the summands. 
He has also established the hard Lefschetz theorem and the Hodge-Riemann
bilinear relations for the intersection cohomology groups of projective varieties.

Saito's results complete the verification
of the Hodge-Lefschetz package for the intersection cohomology
groups of a variety $Y$, thus yielding the wanted generalization
of the classical results in $\S$\ref{tcav} to singular varieties.

The perverse and the standard truncations in ${\cal D}_Y$
correspond to the standard and to the  above-mentioned $\tau'$ truncations
in $D^b ( \mbox{MHM} (Y))$, respectively. See \ci{samhm}, p. 224 and Remarks 4.6.
It follows  that  the following spectral sequences
 associated with complexes  $K \in
{\frak r}(D^b  ( \mbox{MHM} (Y) )) \subseteq {\cal D}_Y$
are spectral sequences of mixed Hodge structures:

\begin{enumerate}
\item
the perverse spectral sequence;
\item  the Grothendieck spectral  sequence; 
\item the perverse Leray spectral sequence  associated with a map $f: X \to Y$
\item
the Leray spectral sequence associated
with  a map $f: X \to Y$.
\end{enumerate}

\begin{rmk}
\label{020}
{\rm 
C. Sabbah, \ci{sabbah} and T. Mochizuki \ci{mochizuki}
have extended the range of applicability of the decomposition theorem
to the case of intersection cohomology complexes associated 
with semisimple local systems
on quasi-projective varieties. They use, among other ideas,   M. Saito's $D$-modules approach.
}
\end{rmk}

\subsection{A proof via classical Hodge theory}
\label{dmapp}
Let us summarize some of our joint work on the subject of the decomposition theorem.

\begin{itemize}
\item
Our paper \ci{decmightam} gives a geometric proof of the decomposition  theorem
for the push forward  $f_* IC_X$ of the intersection cohomology complex
via a proper map  $f: X \to Y$ of complex algebraic varieties,  and complements it
with a series of Hodge-theoretic results in the case when $Y$ is projective.
In particular, we endow the intersection cohomology groups
of a projective variety with a pure Hodge structure.
These results are stated  in the case when $X$ is nonsingular and projective
as Theorem \ref{mtmdecmig} below. The statements in the case when $X$ is projective,
but possibly singular, are essentially identical to the ones
in Theorem \ref{mtmdecmig}, except that one is required to
replace $\rat_X[n]$ with $IC_X$ (see \ci{decmightam}). 

\item
In the paper \ci{decmigseattle},   we show how  to choose,
when $X$ and $Y$  are  projective,  splitting isomorphisms
in the decomposition theorem  so that they are compatible
with the various
Hodge structures found in \ci{decmightam}. 

\item
The extension to the  quasi projective context of the results in \ci{decmightam, decmigseattle} is contained 
in  \ci{decII}, which  builds on
\ci{decmigso3}.    Since  these papers deal with non compact varieties,
the statements involve mixed Hodge structures.
These results are listed
in $\S$\ref{dtmHs}. 
  \end{itemize}

Most  of  results mentioned above have been obtained  earlier and in greater generality
by  M. Saito  in \ci{samhp, samhm} by the use of mixed Hodge modules.  
While our  approach uses  heavily the theory of perverse sheaves,
it  ultimately rests  on classical and mixed  Hodge theory.

The  proof  of the decomposition theorem in \ci{decmightam} is geometric in the sense that: 

\begin{itemize}
\item
it identifies  the refined intersection
forms on the fibers of the map $f$ as  the agent responsible
for  the splitting behavior of $f_* IC_X$ and
\item it provides a geometric 
interpretation of the perverse  Leray filtration on $I\!H^*(X)$.
\end{itemize}

Since the mixed-Hodge-theoretic results are surveyed in $\S$\ref{dtmHs},
  in this section we mostly concentrate on outlining the  proof of the decomposition theorem
  given in \ci{decmightam}.
  
In the following two sections  $\S$\ref{decmigresults} and $\S$\ref{outlinedecmig},
 we list the results contained in \ci{decmightam}
and give an outline of the proofs in the key special case of a projective map 
$f: X \to Y$ of irreducible projective varieties with $X$ {\em nonsingular} of dimension $n$.

 We fix embeddings $X \subseteq {\Bbb P}$ and $Y\subseteq {\Bbb P}'$
into some projective spaces.
We denote by ${\Bbb P}^{\vee}$
the projective space ``dual" to ${\Bbb P}$, i.e. the projective space
of hyperplanes in ${\Bbb P}$.
 Let
$\eta$ and $L$ be the corresponding hyperplane line bundles 
on $X$ and $Y$, respectively and let $L':=f^*L$. We  denote with the same symbol a line bundle, its first Chern class and the operation of cupping with it.

\subsubsection{The results when $X$ is projective and nonsingular}
\label{decmigresults}
The following theorem summarizes some of the main results in \ci{decmightam}
when $X$ is projective and nonsingular. The results  hold
in the singular case as well, provided 
we replace $\rat_X[n]$ with $IC_X$. However, since  the proof of the singular
case relies on the proof of the nonsingular case, and this latter presents
 all the essential difficulties
(see \ci{decmightamv2}), we prefer to discuss the nonsingular case
only. 
Most of the results that follow hold in the case when $X$ and $Y$ are quasi 
projective (see \S\ref{dtmHs} and \ci{decII}).
  Recall that since $X$ is nonsingular
 of dimension $n$, then $IC_X \simeq \rat_X[n]$.

\begin{tm}
\label{mtmdecmig}
Let $f: X \to Y$ be a proper map of  projective varieties, with $X$ nonsingular
of dimension $n$. The following statements hold.

\begin{enumerate}
\item
 {\rm  ({\bf Decomposition theorem})} $f_*\rat_X[n]$ splits as a direct sum of shifted intersection cohomology
complexes with twisted coefficients on subvarieties of $Y$ (cf.  
{\rm $\S$\ref{statem}.(\ref{404}).(\ref{4004}))}. 

\item
{\rm  ({\bf Semisimplicity theorem})} The summands are semisimple, i.e.
 the local systems {\rm (\ref{4004})}   giving the twisted coefficients are semisimple.
They
are described  below, following   the refined intersection form theorem.
\item
{\rm ({\bf Relative hard Lefschetz theorem})}

Cupping with $\eta$ yields isomorphisms 
\[
\eta^i : \pc{-i}{f_* \rat_X [n] } \simeq \pc{i}{f_* \rat_X [n]}, \qquad  \forall \; i \geq 0.
\]

\item
{\rm ({\bf  Hodge structure theorem})}
The perverse $t$-structure yields  the perverse filtration   
$$P^pH(X)= \im\, \{H(Y, \ptd{-p} f_* \rat_X [n]) \to H(Y,f_* \rat_X [n])\}$$
on the cohomology groups
$H(X)$.  This filtration is  by Hodge substructures and 
the perverse cohomology groups
$$
H^{a-n}(Y, \pc{b}{f_* \rat_X [n]}\simeq P^{-b}H^{a}(X)/P^{-b+1}H^a(X)=
H^a_b(X)
$$
i.e. the
graded groups of the perverse filtration, inherit a pure Hodge structure.

\item
{\rm ({\bf Hard Lefschetz theorems
for perverse cohomology groups})}  The collection
of  perverse cohomology groups $H^*(Y, \pc{*}{f_* \rat_X [n]}$ satisfy
 the conclusion of the hard Lefschetz theorem
with respect to  cupping  with  $\eta$ on $X$ and with respect
to  cupping with an $L$ on $Y$, 
namely:

The cup product with $\eta^i: H^*(Y,\pc{-i}{f_* \rat_X [n]} ) \to H^{*+2i}(Y,\pc{i}{f_* \rat_X [n]} )$ is an isomorphism for all $i\geq 0$.

The cup product with $L^l: H^{-l}(Y,\pc{i}{f_* \rat_X [n]} ) \to H^{l}(Y,\pc{i}{f_* \rat_X [n]} )$ is an isomorphism for all $l\geq 0$ and all $i$.

\item
{\rm ({\bf The perverse filtration on $H^*(X)$})}
The
perverse  filtration  on the groups
$H^{r}(X)$
is given by the following equation
(where it is understood that a linear map with a non-positive exponent
is defined to be the identity and that kernels and images are inside of
$H^r(X)$): 

$$
P^pH^r(X)\; = \;  \sum_{a+b=n-(p+r)}\ke \, L'^{a+1}\cap \im \, L'^{-b}.
$$

\item
{\rm ({\bf Generalized  Lefschetz decomposition and Hodge-Riemann bilinear relations}) }
Let $i,j \in \zed$ and consider the perverse cohomology groups of 4.
 Define $P^{-j}_{-i}: = \ke{\,\eta^{i+1}}\cap \ke{\,
L^{j+1}} \subseteq H^{n-i-j}_{-i}(X)$ if $i,j \geq 0$ and $P^{-j}_{-i}:= 0$, otherwise.
 There is a Lefschetz-type
direct sum decomposition (the $(\eta, L)$-decomposition)
 into pure Hodge substructures
\[
H^{n-i-j}_{-i} (X) = \bigoplus_{l,m\, \in \zed} 
\eta^{-i+l} L^{-j+m} P^{j-2m}_{i-2l},\]
Define, for $i,j \geq 0$,  bilinear forms  on $H^{n-i-j}_{-j} (X)$
\[
S^{\eta L}_{ij} (\alpha, \beta) := \int_X \eta^i \wedge L^j \wedge \alpha
\wedge \beta.\] These forms are well-defined and, 
using the hard Lefschetz theorems 5., they can be suitably defined
for every $i,j \in \zed$.
The bilinear forms $S^{\eta,L}_{ij}$ are non degenerate and orthogonal
with respect to the $(\eta, L)$-decomposition.
Up to the   sign $(-1)^{i+j-m-l+1}$, these forms are a polarization
(see {\rm  $\S$\ref{pam}}, especially {\rm (\ref{hrbr}))} of each $(\eta, L)$-direct summand.

\item
{\rm  ({\bf Generalized Grauert contractibility criterion})}
Fix $y\in Y$ and $ j \in \zed$.
The natural class map, obtained by composing
push forward in homology with Poincar\'e duality   $$H_{n -j} (f^{-1}(y)) \lorw H^{n +j}(X)$$ is
naturally  filtered in view of the decomposition theorem. 
The   resulting graded class map
$$H_{n -j,j} (f^{-1}(y)) \lorw H^{n +j}_j(X)$$
is an injection of pure Hodge structures polarized 
in view of the generalized Hodge-Riemann relations  7.

\item
{\rm ({\bf Refined intersection form theorem})}
The refined intersection form 
\[ H_{n-j}(f^{-1}(y) ) \lorw H^{n+j}(f^{-1}(y))\]
(see $\S$\ref{famfa}, Refined intersection forms) is  naturally filtered in view of the decomposition theorem,
and the resulting graded refined intersection form
\[H_{n-j, k} (f^{-1}(y) ) \lorw H^{n+j}_k (f^{-1}(y))
\quad
\mbox{ is zero for
$j\neq k$ and an isomorphism for $j=k$.}
\]
\end{enumerate}
\end{tm}

\subsubsection{An outline of the proof of Theorem \ref{mtmdecmig}}
\label{outlinedecmig}

We start by sketching the proof in  the non-trivial toy model
of a semismall map (\ci{decmigsemi}),  as many important steps 
appear already in this case.
 We refer to $\S$\ref{semismall}
for basic definitions and facts concerning this remarkable class of  maps.

\medskip
{\bf 1. The case of semismall maps.}

\smallskip
There is no loss of generality in assuming that
the map $f$ is surjective. Since a semismall map is generically finite, we  have
$n =\dim{X}= \dim{Y}$.
We proceed by induction on $n=\dim Y$ and  prove all the results 
of Theorem \ref{mtmdecmig}. 

By the semismallness assumption,
we have that    $\pc{j}{f_* \rat_X [n] }=0$ for every $j \neq 0$, so that
the relative hard Lefschetz is trivial and so is the perverse filtration. 
The first point to show is that,
{\em from the point of view of the Hodge-Lefschetz package,
$L'=f^{*}L$
  behaves  as if it were  a hyperplane line bundle,}  even though it is not
 (it is trivial along the fibers of the map $f$): 
{\em all the theorems in {\rm $\S$\ref{chl}} hold with  $L'$ replacing  $\eta$.} 

{\em The hard Lefschetz theorem for $L'$.}
By induction, we assume that the statements in Theorem
\ref{mtmdecmig} hold for all semismall 
maps between varieties of dimension less than $n$. 
Let $D\subseteq Y$ be a generic hyperplane section. The map 
$f^{-1}(D) \to D$ is still semismall.
Since $f_* \rat_X [n] $ is perverse, in the range $i\geq 2$ 
the Lefschetz theorem on hyperplane sections for perverse sheaves
  (see $\S$\ref{subsectex}) reduces the hard Lefschetz 
for $L'^i$ on $X$ to that for $L'^{i-1}$  on $f^{-1}(D)$. 
In  the critical case $i=1$, the cup product with $L'$
factors as
$H^{n-1}(X) \to H^{n-1}(f^{-1}(D)) \to H^{n+1}(X)$,
where 
the first map is injective and the second is surjective.
As explained in the ``inductive  approach to hard Lefschetz'' paragraph
of \S\ref{pam}, the inductive Hodge-Riemann relations for the restriction of $L'$ to 
$f^{-1}(D)$ give the hard Lefschetz theorem for the cup product with 
$L'$.  

{\em The approximation trick.} 
We must  prove the Hodge-Riemann relations for  
the space of primitives $P^{n}_{L'}=\,\ke \, L': H^{n} (X)\to  H^{n+2}(X)$
(for use in the case when $\dim{X}=n+1$).
The hard Lefschetz theorem discussed above implies  that
$\dim P^{n}_{L'}= b_{n}- b_{n-2}$  and that 
the decomposition 
$H^{n} (X)=P^{n}_{L'}\oplus  L' H^{n-2} (X)$
is orthogonal with respect to the Poincar\'e pairing, just as if $L'$ were a hyperplane
bundle. 
In particular, the restriction of the Poincar\'e pairing $S(\alpha, \beta)= \int_X \alpha \wedge \beta$
to $P_{L'}^{n}$ is nondegenerate.
The bilinear form  $\widetilde{S}(\alpha, \beta ):=S( \alpha,C \beta)$
($C$ is the Weil operator; see
\S\ref{pam}) is still nondegenerate. 
The class $L'$ is on the boundary of the ample cone: for any positive integer $r$, the class $L'+\frac{1}{r}\eta$ is  ample, and we have the classical Hodge-Riemann relations on the subspace $P^{n}_{L_{r}}:= \ke \, (L'+\frac{1}{r}\eta)\subseteq  H^{n} (X)$:  the remark
made above on the dimension of $P^{n}_{L'}$ implies  that any class $\alpha \in P^{n}_{L'}$ is the limit of classes $\alpha_r \in P^{n}_{L_{r}}$, 
so that the restriction of $\widetilde{S}$ to $P^{n}_{L'}$ is semidefinite; since
it is  also nondegenerate,  the Hodge-Riemann bilinear relations follow.

{\em Decomposition and semisimplicity.}
To prove the decomposition and semisimplicity theorems,
we proceed one stratum at the time; higher dimensional strata are dealt with
inductively by cutting transversally with a generic 
hyperplane section $D$ on $Y$, so that one is re-conduced to the semismall map $f^{-1}(D)\to D$ where the dimension of a positive dimensional stratum
on $Y$ has decreased by one unit on $D$.
The really significant case left is that of a zero-dimensional relevant stratum $S$. 
As explained in $\S$ \ref{interform}, the splitting of the perverse sheaf $f_* \rat_X [n]$ into a direct sum of intersection cohomology complexes with twisted coefficients on subvarieties of $Y$
is equivalent to the nondegeneracy of the refined intersection form 
$I:H_{n}(f^{-1}(y)) \times H_{n}(f^{-1}(y))\lorw \rat$, for 
$y \in S$. 

In order to establish the nondegeneracy of the refined intersection  forms $I$,
we turn to mixed Hodge theory ($\S$\ref{pam}) and  use   the following result of P. Deligne
(cf. \ci{ho3}, Proposition 8.2.6):

 {\bf (weight miracle)}:
 {\em 
if $Z\subseteq U \subseteq X$ are inclusions with $X$ 
a nonsingular compact variety, $U\subseteq X$ a Zariski dense open subvariety and
$Z\subseteq U$ a closed subvariety of $X$,
then the images in $H^j(Z, \rat)$ of the restriction maps
from $X$ and from $U$ coincide.}

Thanks to  the 
weight miracle, $H_{n}(f^{-1}(y))$ injects in $H^{n}(X)$ 
 as a Hodge substructure.
Since, for a general section $D$, we have  
$f^{-1}(y)\cap f^{-1}(D)=\emptyset$, 
we see that $H_{n}(f^{-1}(y))$ is contained  in $P^{n}_{L'}$.
The restriction of the Poincar\'e pairing to $H_{n}(f^{-1}(y))$  
is thus a polarization and is hence  nondegenerate. The same is thus true for
the refined intersection form  $I$. 

As  noted already in \ref{dtss}, the local systems involved have finite monodromy, hence they are obviously semisimple.
 This concludes our discussion of the semismall case.

\medskip
{\bf 2. The general case: extracting the semismall ``soul" of a map.}

\smallskip
The  proof is by induction  on the  the pair of indices  $(\dim{Y}, r(f))$,
where $r(f)= \dim{X \times_Y X} -\dim{X} $ is the {\em defect of semismallness} of the map $f$.
To give an idea of the role played by $r(f)$ let us say that in the 
decomposition theorem $\S$\ref{statem}.(\ref{404}), the direct sum ranges precisely
in the interval $[-r(f), r(f)]$.  The inductive hypothesis
takes the following form:
all the statements of Theorem \ref{mtmdecmig} hold
for all proper maps $g:X' \to Y'$ with either $r(g)<r(f)$, or  with $r(g)=r(f)$ and
$\dim Y'< \dim Y$.   Let $n:=\dim X$. The induction starts with the verification
of Theorem \ref{mtmdecmig} in the case when $Y$ is a point, in which case
the results boil down to the classical result of Hodge-Lefschetz theory outlined 
in $\S$\ref{tcav} and listed more succinctly   in  Theorem \ref{chl}.

\medskip
{\bf 2a. The universal hyperplane section and relative hard Lefschetz.}

Let  $g: X'\subseteq X \times {\Bbb P}^{\vee}  \to Y'=Y \times {\Bbb P}^{\vee} $  
be the universal hyperplane section. 
If $r(f) >0$, then $r(g)<r(f)$ and, by induction, 
Theorem  \ref{mtmdecmig} holds for  $g$.
As in the classical case (cf. \S\ref{chl}), the relative Lefschetz hyperplane Theorem \ref{8020}
implies the relative hard Lefschetz theorem for $f$ except for $i=1$,
where we have
 the  factorization of the cup product map with  $\eta$
$$
\xymatrix{
\pc{-1}{f_* \rat_X[n]} \ar[r]^{\rho} &
\pc{0}{g_* \rat_{X'}[n-1] } 
\ar[r]^{\gamma} &  \pc{1}{f_* \rat_X[n]}.}
$$
The first map is a monomorphism and the second is an epimorphism.
We argue as in the proof of the hard Lefschetz theorem via the semisimplicity of monodromy: we use 
an argument similar to the identification
of the monodromy invariants  of a Lefschetz pencil with  the image 
of the cohomology of a variety into  the cohomology of a hyperplane section, 
and we couple it
with  the semisimplicity (inductive assumption!)
of $\pc{0}{g_* \rat_{X'}[n-1] }$ to show that: 
\begin{pr}
\label{splitsummand}
The image of 
$\pc{-1}{f_* \rat_X[n]}$ in  $\pc{0}{g_* \rat_{X'}[n+1] }$ is a split summand applied isomorphically
onto $\pc{1}{f_* \rat_X[n]}$ by $\gamma$.
\end{pr} 

The relative hard Lefschetz for $f$  follows and,
by applying   Deligne's Lefschetz splitting criterion,
Theorem \ref{delsplit}, we conclude that 
$f_* \rat_X [n] \simeq \oplus_i \pc{i}{f_* \rat_X[n]}[-i]$.

>From the statements known for $g$ by induction,
we  get that $\pc{i}{f_* \rat_X[n]}$ is a direct sum of intersection 
cohomology complexes of semisimple local systems for all $i \neq 0$.
Moreover,  for all $i \neq 0$,
the associated perverse  cohomology groups verify the hard Lefschetz theorem and the Hodge-Riemann relations 
with respect to cupping with $L$.

What is left to investigate is the zero perversity complex $\pc{0}{f_* \rat_X[n]}$. 
Again in analogy with the classical case, 
we can ``shave off" another piece which comes from the hyperplane section  and dispose of by
  using the inductive hypothesis.
In fact, 
the analogue of the  primitive Lefschetz decomposition theorem \ref{chl}.2. holds: 
by setting, for every
 $i \geq 0$, 
${\cal P}^{-i}:= \ke \, 
\{\eta^{i+1}:\pc{-i}{f_* \rat_X[n]} \to \pc{i+2}{f_* \rat_X[n]}\}$ 
we have canonical direct sum decompositions:
\begin{equation}
\label{decompo}
\pc{-i}{f_* \rat_X[n]} = \bigoplus_{r\geq 0} \eta^r  {\cal P}^{-i-2r},
\qquad
 \pc{i}{f_* \rat_X[n]} = \bigoplus_{r\geq 0} \eta^{i+r}  {\cal P}^{-i-2r}.
\end{equation}
The only remaining pieces for which we have to prove the statements of Theorem \ref{mtmdecmig}
are  the perverse sheaf ${\cal P}^{0}$ and its cohomology $H^*(Y,{\cal P}^{0})$
which, in view of the primitive decomposition, is a summand of
the perverse cohomology group
$H^{*+n}_0(X)$.
(The analogy with the classical study of algebraic varieties by means of hyperplane sections is as follows:
the new  cohomology classes, i.e. the ones not coming from a hyperplane section, 
appear only in the middle dimension $P^{n}=\ke\, \{ \eta: H^{n}(X) \to H^{n+2}(X)\}$. In this game, ``middle dimension" is re-centered at zero.)
We are left with proving:
\begin{enumerate}

\item The Hodge package of  \S\ref{chl} holds for $H^*(Y, {\cal P}^{0})$ with respect to cupping with $L$.
\item ${\cal P}^{0}$ is a direct sum of
twisted  intersection cohomology complexes. 
\item The twisting  local systems are semisimple.
\end{enumerate}

{\bf 2b. The Hodge package for ${\cal P}^0$.}

The main intuition behind the proof of  the statements 1. and 2. above,  
which was  inspired also by the 
illuminating discussion of the decomposition theorem contained in \ci{macicm},
is that $H^*(Y, {\cal P}^{0})$ is the ``semismall soul of the map $f$,"
that is {\em it behaves as the cohomology  of a (virtual) nonsingular 
projective variety with a semismall map to $Y$}.
In order to handle the group $H^*(Y, {\cal P}^{0})$, we
mimic  the proof of the decomposition theorem for semismall maps.

One of the main difficulties in   \ci{decmightam}
is  that, in order to use classical Hodge theory,
we have to prove at the outset   that the perverse Leray filtration is Hodge-theoretic, 
i.e. that the subspaces  $P^pH^*(X) \subseteq H^*(X)$
(cf. $\S$\ref{decmigresults}.(4)) are Hodge substructures of the natural
Hodge structure on $H^*(X)$.
The geometric description of the perverse filtration
in \ci{decmigso3} (see \S\ref{tpfmc}) implies
that this fact holds   for every  algebraic map,
proper or not,  to a quasi projective variety, 
and the  proof in \ci{decmigso3} is independent of the decomposition theorem. 
It follows that the geometric description of the perverse filtration in 
\ci{decmigso3}  can therefore be used to  yield a  considerable simplification
of the line of  reasoning   in \ci{decmightam} for it 
 endows, at the outset,   the perverse cohomology groups
$H^a_b(X)$ 
with a natural Hodge structure, compatible with the primitive
Lefschetz decompositions stemming from (\ref{decompo}), and 
 with respect to which the cup product maps
$L:H^*(Y, {\cal P}^i) \to
H^{*+2}(Y,{\cal P}^i) $ and 
$\eta:P^kH^*(X) \to
P^{k-2}H^{*+2}(X)$ are Hodge maps of type $(1,1)$. 

We start by proving 1., i.e. the Hodge package for $H^*(Y, {\cal P}^0)$.
The argument for  the hard Lefschetz 
isomorphism  $L^i: H^{-i}(Y, {\cal P}^0) \simeq H^{i}(Y, {\cal P}^0)$ is completely analogous to the one used for a semismall map:
the  Lefschetz theorem on hyperplane sections for the perverse sheaf ${\cal P}^0$ and the inductive hypothesis
(for a generic hyperplane
 section $D\subseteq Y$,  we have
 $f': f^{-1}(D) \to D$ and 
  ${\cal P}^0$ restricts, up to a shift, 
to the analogous complex ${\cal P}'^0$  for $f'$) yield immediately the theorem in the range $i\geq 2$ and also yield  a factorization of $L: H^{-1}(Y, {\cal P}^0) 
\to H^{1}(Y, {\cal P}^0)$ as the composition of the injective restriction to $D$
and the surjective  Gysin map. Again by the inductive hypotheses, the Poincar\'e pairing polarizes $\ke \, L:H^0(D,{\cal P}'^0) \to H^2(D,{\cal P}'^0)$, and,
as in the classical case,  this proves the remaining case $i=1$.

The most delicate point is to prove that the Riemann Hodge relations hold for 
$P^{00}:= \ke \, \{L:H^0(Y,{\cal P}^0) \to H^2(Y,{\cal P}^0)\}$. 
The Poincar\'e pairing  induces a bilinear form $S$ 
on $H^{n}(X)=H^0(f_* \rat_X [n])$  
and on its subquotient  $H^0(Y,{\cal P}^0)$.
This is because we have the following orthogonality
relation $P^1H(X) \subseteq P^0H(X)^{\perp}$.
More  is true:  $S$ is nondegenerate  
on $P^0H^n(X)/P^1H^n(X) = H^n_0(X)$ and the $(\eta,L)$ decomposition is orthogonal
so that  the restriction of $S$ to the summand
$P^{00}$ is nondegenerate.
The Hodge-Riemann relations are then proved with an ``approximation trick''
similar, although more involved, to the one used in the semismall case.
We consider the subspace
$\Lambda=\lim_r  \ke (L'+\frac{1}{r}\eta)\subseteq  H^{n} (X)$.
Clearly, we have  $\Lambda \subseteq \ke \, L'$ and the hard Lefschetz theorem
implies that $\ke\,  L' \subseteq P^0H^n(X)$.  The 
 nondegenerate form $\widetilde{S}$  is semidefinite
on $\Lambda/\Lambda\cap P^1H^n(X)$.  It follows that it is a 
 a polarization.  A polarization restricted to a Hodge
 substructure is  still a polarization. The Hodge-Riemann relations for $P^{00}$ 
follow from  the inclusion of Hodge structures
$P^{00} \subseteq  \Lambda/\Lambda\cap P^1H^n(X)$.

\medskip
{\bf 2c. Semisimplicity.}

We need to prove that ${\cal P}^0$ splits as a direct sum of intersection cohomology complexes of semisimple local systems.
As in the case of semismall maps, higher dimensional strata are disposed-of
by induction on the dimension
of $Y$ and by cutting with generic hyperplane sections of $Y$. 
One is left to prove the critical  case of a zero-dimensional stratum.
Again by the splitting criterion of Remark
\ref{splitdecmig}, we have to prove that, for any point $y$ in the zero dimensional stratum, denoting by $i:y \to Y$ the closed imbedding, 
$\iota: {\cal H}^{0}(i^!{\cal P}^0)\to {\cal H}^{0}(i^{*}{\cal P}^0)$
is an isomorphism. Given the decomposition (\ref{decompo}), 
${\cal H}^{0}(i^!{\cal P}^0)$ is a direct summand of $H_n(f^{-1}(0))$
and ${\cal H}^{0}(i^{*}{\cal P}^0)$ is a direct summand
of $H^n(f^{-1}(0))$, so that
the map $\iota$   is the restriction to these summands of the refined intersection
form ($\S$\ref{famfa}) on $f^{-1}(0)$.
Although in general,  the  map $H_n(f^{-1}(0))\to H^n(X)$ is not injective,
the   weight miracle  is used to prove that
the map ${\cal H}^{0}(i^!{\cal P}^0) \to H^{n}_0 (X)$
is an injection
with image a pure Hodge substructure  of the Hodge structure
we  have on $H^n_0(X)$ (by virtue of the 
geometric description of the perverse filtration \ci{decmigso3}
mentioned above).
Since this  image lands automatically in the $L'$-primitive part,
we conclude that the descended intersection form  polarizes this image,
hence
 $\iota$ is an isomorphism and we have the desired splitting
into intersection cohomology complexes.

We still have to establish the semisimplicity
of the local systems in (\ref{4004}) (and hence of the ones
appearing in ${\cal P}^0$).
 This is accomplished by exhibiting them
as quotients of local systems associated with smooth proper maps
and are hence semisimple by the semisimplicity for smooth proper maps 
Theorem \ref{dss}.

This concludes the outline of the proof of Theorem \ref{mtmdecmig}.

\section{Applications of perverse sheaves and of the decomposition theorem } 
\label{dtappl}
In this section, 
we   give, without any pretense of completeness, 
a sample of remarkable applications of the theory of perverse sheaves and
of the decomposition theorem.

We focus mostly  on the complex case, although most of
the discussion goes through over a field of positive characteristic,
with constructible $\rat$-sheaves replaced by $l$-adic ones.

In this chapter, we use the machinery of derived categories 
and functors and some results on perverse sheaves. The
notions introduced in our crash course may not be sufficient
to follow the (few) proofs included.
We refer to $\S$\ref{tfidy}, to the references quoted there,
and to \S\ref{secptsdy}. In particular,
we adopt the simplified notation $f_*,f_!$ for the derived functors $Rf_*,Rf_!$.

\subsection{Toric varieties and combinatorics of polytopes}
\label{toric}
 The purpose of this section is to state and explain the content of Theorem
\ref{hvectic} on how the combinatorics of rational polytopes in Euclidean space
 relates
to the intersection cohomology groups of the associated toric varieties. 
Theorem \ref{hvectic} is stated in $\S$\ref{hpolysec} and we work out two
examples in $\S$\ref{twoex}, where the decomposition theorem
is seen in action in situations where, we hope, 
the minimal  background we provide in this section  is sufficient to follow the arguments.

For the basic definitions concerning  toric varieties, we refer to  \ci{fulton,oda}. 
The recent survey  \ci{braden} contains many 
historical details,  motivation, a discussion of open problems and recent results, and an 
extensive bibliography. 

We will adopt the point of view of polytopes, which we find more appealing to intuition.

Recall that a $d$-dimensional normal projective complex variety $X$ is
a toric variety if it has an action of the complex torus
$T=(\comp^*)^d$ with finitely many orbits. In this case,  there is a
moment map $\mu:X \to \real^d$ whose image is a $d$-dimensional convex
polyhedron $P$, whose vertices have rational coordinates, and which
determines the toric variety $X$ up to isomorphism. The mapping $\mu$
determines an order-preserving one to one correspondence between the orbits
of $T$ and the faces of $P$ as follows. For each orbit ${\cal O}
\subseteq X$ the image $\mu({\cal O})\subseteq P$ is the interior $F^0$
of a face $F \subseteq P$. Moreover, $\dim_{\comp}({\cal O})
=\dim_{\real}(F)$
and the fibers of $\mu:{\cal O} \to F^0$ are diffeomorphic to the
compact torus $(S^1)^{\dim F}$. For $i=0, \ldots, d-1$, let $f_i$
be the number of $i$-dimensional faces of $P$.
We denote by $X_P$ the 
projective toric variety associated with $P$.
A $d$-dimensional simplex $\Sigma_d$ is the convex envelope of $d+1$ affinely independent points $v_0, \ldots, v_d$
in $\real ^d$.   $X_{\Sigma_d}$ is a possibly weighted  $d$-dimensional projective space.
A polytope is said to be simplicial if its faces are simplices,
We  say that a toric variety is $\rat $-smooth when it has only finite quotient singularities.
A map of varieties $f:\widetilde{X} \to X$, both of which are toric,
 is called a toric resolution if it is birational, equivariant 
with respect to the torus action, and $\widetilde{X}$ is $\rat $-smooth.
\medskip

The following is well known:

\begin{pr}
A toric variety $X_P$ is $\rat$-smooth if and only if $P$ is simplicial.
\end{pr}

\subsubsection{The $h$-polynomial}
\label{hpolysec}

Let  $P$ be  a simplicial $d$-dimensional polytope with 
number of faces encoded by the ``face vector'' $(f_0, \ldots f_{d-1})$. Define the associated 
 ``$h $-polynomial''
\begin{equation}
\label{hpoly}
h(P,t)=(t-1)^d+f_0(t-1)^{d-1}+ \ldots +f_{d-1}.
\end{equation}
The simplicial  toric variety $X_P$ has a decomposition as a disjoint union 
of locally closed subsets, each isomorphic to the quotient of an affine space 
by a finite commutative group. 
This decomposition can be used to compute the rational cohomology groups $H^*(X_P,\rat)$, and
we have the following proposition, 
see \ci{fulton}, Section 5.2 for a detailed proof:

\begin{pr}
\label{hvect}
Let $P$ be a simplicial rational polytope, with ``$h $-polynomial'' $h(P,t)= \sum_0^d h_k(P)t^k$.
Then
$$
H^{2k+1}(X_P, \rat)=0 \quad  \hbox{ and } \quad   {\rm dim \,} H^{2k}(X_P, \rat) = h_{k}(P).
$$
\end{pr}

Poincar\'e duality and the hard Lefschetz theorem imply the  following

\begin{cor}
\label{prophevct}
We have the following relations: 
$$
h_k(P)=h_{d-k}(P)\; \;  \hbox{ for } 0\leq k \leq d , \qquad h_{k-1}(P) \leq h_k(P) \;\;\hbox{ for } 0\leq k \leq \lfloor d/2 \rfloor.
$$
\end{cor}

Corollary \ref{prophevct} amounts to a set of non trivial relations among the face numbers $f_i$, 
and gives necessary conditions for a sequence $(a_0, \ldots a_{d-1}) \in \nat^d$ to be the face vector of 
a simplicial polytope.
Exploiting more fully the content of the hard Lefschetz theorem, it is possible to characterize completely  the sequences  in $\nat^d$ which occur as
the  face vectors of some simplicial polytope;  
see \ci{braden}, Theorem 1.1.

\medskip

The polynomial 
\begin{equation}
\label{gpoly}
g(P,t)=h_0+(h_1-h_0)t+ \ldots +(h_{[d/2]}-h_{[d/2]-1})t^{[d/2]}
\end{equation}
has, by  Corollary  \ref{prophevct}, positive coefficients and  uniquely determines $h$.
The coefficient $g_l=h_l-h_{l-1}$ is the dimension of the primitive cohomology
($\S$\ref{pam}) of $X_P$ in degree $l$.

\begin{ex}
\label{ex1}
{\rm Let $\Sigma_d$ be the $d$-dimensional simplex. We have 
$f_0=d+1={d+1\choose1}, \ldots, f_i={d+1 \choose i+1}$ and 
$$
h(\Sigma_d, t)=
(t-1)^d+{d+1\choose1}  + \ldots +{d+1 \choose i+1}(t-1)^{d-i-1}+ \ldots {d+1\choose d}=1+t+\ldots +t^d,
$$
so that  $h_i=1$ and $g(\Sigma_d, t)=1$, consistently with the fact that $X_{\Sigma_d}=\pn{d}$.  

\n
Let $C_2$ be the square, convex envelope of the four points $(\pm 1,0), (0,\pm 1)$. 
We have  $f_0=4, f_1=4$,
$h(C_2,t)=(t-1)^2+4(t-1)+4=t^2+2t+1$, and $g(C_2,t)=1+t$. 
In fact, $X_{C_2} =\pn{1} \times \pn{1}$.

\n
Similarly, for the octahedron $O_3$, convex envelope of $(\pm 1, 0 , 0), (0, \pm 1, 0 ), (0,0,\pm 1)$, 
we have $f_0=6, f_1=12, f_2=8$,  $h(O_3,t)=t^3+3t^2+3t+1$ and $g(O_3,t)=2t+1$. This is in accordance
with the Betti numbers of  $X_{O_3}=(\pn{1})^3$.
}
\end{ex}

If the polytope is not simplicial, so that the
toric variety is not $\rat$-smooth, neither  Poincar\'e duality, 
nor the hard Lefschetz theorem
necessarily hold for the cohomology groups. 
Furthermore, as shown in  \ci{mcconnel},  the ordinary  cohomology groups of a singular toric variety 
is not a  purely combinatorial invariant, but depends also on some geometric data of the polytope, e.g. 
the measures of the angles between the faces of the polytope.
The situation drastically simplifies when considering 
intersection cohomology groups. In fact, Poincar\'e duality
and the hard Lefschetz theorem hold for intersection cohomology, 
so that  the ``generalized'' $h$-polynomial
  $h(P,t)= \sum_0^d
h_k(P)t^k$, where $h_k(P):= \dim \, I\!H^{2k}(X_P, \rat)$
satisfies the conclusions of  Corollary  \ref{prophevct}. 
Furthermore, it turns out that the polynomial 
$h(P,t)$ is a combinatorial invariant, i.e.  it can be defined only in terms of the 
partially ordered set of faces of the polytope $P$.
 Note that when the polytope $P$ is simplicial,
so that the toric variety $X_P$
is $\rat$-smooth, then $H^*(X_P,\rat)= I\!H^*(X_P, \rat)$. Hence,
in this case,  by Proposition \ref{hvect}, the generalized $h$-polynomial defined below coincides with the one defined earlier and we can denote the two  in the same way.

We now give the combinatorial definitions of 
the $h$ and $g$ polynomials for a not necessarily simplicial polytope. 

\begin{defi}
\label{generh}
{\rm Suppose $P$ is a polytope of dimension $d$ and that the polynomials $g(Q,t)$ and $h(P,t)$
have been defined for all convex polytopes $Q$ of dimension less than $d$. 
We set 
$$
h(P,t)=\sum_{F<P}g(F,t)(t-1)^{d-1-{\dim} F},
$$
where the sum is extended to all proper faces $F$ of $P$ including the empty face $\emptyset$, for which 
$g(\emptyset,t)=h(\emptyset,t)=1 $ and ${\rm dim \,} \emptyset=-1$. 
The polynomial $g(P,t)$ is defined from $h(P,t)$ as in (\ref{gpoly}).}
\end{defi}

We note that these definitions coincide with the previous ones given in (\ref{hpoly}) and (\ref{gpoly}) 
if $P$ is simplicial, since $g(\Sigma,t)=1;$ see Example \ref{ex1}. In fact, we have the following 

\begin{tm}
\label{hvectic}  {\rm(\ci{fieseler})}
Let $P$ be a rational polytope. Then
$$h(P,t)=\sum_{F<P}g(F,t)(t-1)^{d-1-{\rm dim F}}= 
\sum {\rm dim \,}I\!H^{2k}(X_P, \rat)t^k. $$
\end{tm}

\smallskip

Given a subdivision $\widetilde{P}$ of the polytope $P$, there is a corresponding map $X_{\widetilde{P}} \to X_P$.
The toric orbits of $X_P$ provide a stratification for $f$.
The fibers over toric orbits, 
which properties can be read-off from the combinatorics of the subdivision,
are unions of toric varieties glued along toric subvarieties; for a discussion, see
 \ci{yauetal}. It is well known (cf.  \ci{fulton}, Section 2.6) that any polytope 
becomes simplicial after a sequence of  subdivisions. 

Theorem \ref{hvectic} on the dimension of the intersection cohomology groups 
of a toric variety can be proved  by exploiting the decomposition
theorem for a resolution defined by a subdivision of the polytope $P$. A sketch of a proof along these lines has
been given by R. MacPherson in several talks in 1982. J. Bernstein and A. Khovanskii also 
developed proofs which have not been  published.

\subsubsection{Two worked out examples of toric resolutions}
\label{twoex}

We describe  MacPherson's approach to Theorem \ref{hvectic} via the decomposition
theorem  in the special cases  of   subdivision of the cube of dimension 3 and
4. The general case can be proved along these lines.

\smallskip
\n
Let $C_i$ be the $i$-dimensional cube. It is not simplicial
if  $i>2$, and the $k$-dimensional faces of $C_i$
are $k$-dimensional cubes $C_k$. 
The three dimensional cube $C_3$. has $8$ faces of dimension $0$ and $12$ faces of dimension $1$ which are of course simplicial; there are
  $6$ faces of dimension $2$, for which we have  already computed 
  $g(C_2,t)=1+t$.  It follows that
\begin{equation}
\label{3cube}
   h(C_3,t)=(t-1)^3+8(t-1)^2+12(t-1)+6(1+t)=1+5t+5t^2+t^3 \hbox{ and }
   g(C_3,t)=1+4t.
\end{equation}

Similarly, the four dimensional cube   $C_4$ has $16$ faces of dimension $0$,  $32$ faces of dimension $1$, which are  
all simplicial, $24$ faces of dimension $2$, which are equal to $C_2$, and finally $8$ faces of dimension $3$,
which are  equal to $C_3$. Thus
\begin{equation}
\label{4cube}
h(C_4,t)=(t-1)^4+16(t-1)^3+32(t-1)^2+24(1+t)(t-1)+8(1+4t)=t^4+12t^3+14t^2+12t+1.
\end{equation}

The $3$-dimensional cube $C_3$ has a simplicial subdivision $C_3'$ which does not add any vertex,
and divides every two-dimensional face into two simplices by adding its diagonal, 
see the picture in \ci{fulton}, p.50. The resulting map $f:X_{C_3'} \to X_{C_3}$ is an isomorphism
outside the six singular points of  $X_{C_3}$, and the fibers over this points are isomorphic to $\pn{1}$.  
The $f$-vector of $C_3'$ has $f_0=8$, $f_1=18$ and $f_2=12$ and  $h$-polynomial 
$h(C'_3,t)=t^3+5t^2+5t+1$ which equals the $h$-polynomial $h(C_3,t)$ computed above. 
This equality reflects the fact that $f$ is a small resolution in the sense of  \ref{small}, 
so that $H^i(X_{C_3'})=I\!H^i(X_{C_3})$.

We discuss the decomposition theorem for the map $f:X_{\widetilde{C_3}} \to X_{C_3}$
where $\widetilde{C_3}$ is obtained by the following decomposition of $C_3$:
 for each of the six two-dimensional faces $F_i$,
we add its barycenter $P_{F_i}$ as a new vertex, and we join $P_{F_i}$ with each vertex of $F_i$.
We obtain in this way a simplicial polytope $\widetilde{C_3}$  with $14$ vertices, $36$ edges and $24$ two-dimensional simplices. Its $h$-polynomial is $h(\widetilde{C_3}, t)=t^3+11t^2+11t+1$. The map $f$ is an isomorphism away from the six points $p_1, \ldots, p_6$ corresponding to the two-dimensional faces of $C_3$. The fibers $D_i$ over each point $p_i$ 
is the toric variety corresponding to $C_2$, i.e. $\pn{1} \times\pn{1}$, in particular $H^4(D_i)=\rat$, and 
$^{\frak p}\!{\cal H}^{\pm1}(f_*\rat_{X_{\widetilde{C_3}}}[3])\simeq \oplus \rat_{p_i} $. 
The decomposition theorem for $f$ reads as follows:
$$
f_*\rat_{X_{\widetilde{C_3}}}[3] \; \simeq  \;
 IC_{{C_3}}\,  \oplus  \,   
(  \oplus_i      \rat_{p_i}[1]) \, \oplus \,
 ( \oplus_i \rat_{p_i}[-1])      
$$
and 
$$
H^l(X_{\widetilde{C_3}}) \simeq I\!H^l(X_{C_3}) \hbox { for } l\neq 2,4, 
 \qquad \dim H^l(X_{\widetilde{C_3}}) 
= \dim \, I\!H^l(X_{C_3})+6     \hbox{ for } l=2,4.
$$
It follows that 
$\sum {\rm dim \,}I\!H^{2k}(X_{C_3})t^k=\sum {\rm dim \,}H^{2k}(X_{\widetilde{C_3}})t^k -6t-6t^2=h(\widetilde{C_3},t)-6t-6t^2=t^3+5t^2+5t+1=
h(C_3,t)$, as already computed.

Finally, as a more challenging example, we  consider the four-dimensional cube $C_4$. We subdivide it by adding as new 
vertices the barycenters of the $8$ three-dimensional faces and of the $24$ two-dimensional faces. It is not hard to see that
the $f$-vector of the resulting simplicial polytope $\widetilde{C_4}$ is $(48,240,384,192)$ and  
$h(\widetilde{C_4}, t)=t^4+44t^3+102t^2+44t+1$. The geometry of the map $f:X_{\widetilde{C_4}}
 \to X_{C_4}$ which is relevant to the decomposition theorem is the following.
The $24$ two-dimensional faces correspond to  rational curves $\overline{O}_i$, closures of one-dimensional orbits
$O_i$, along which the map $f$ is locally trivial and looks, on a 
normal slice, just as the map $X_{\widetilde{C_3}} \to X_{C_3}$ examined in the example above.
The fiber over each of the $8$ points $p_i$ corresponding to the three-dimensional faces is 
isomorphic to $X_{\widetilde{C_3}}$. Each  point $p_i$ is the intersection of the six 
rational curves $\overline{O}_{i_j}$ corresponding to the six faces of the three-dimensional cube associated with $p_i$.
The last crucial piece of information is that {\it the local systems arising in the 
decomposition theorem are in fact trivial}.
Roughly speaking, this follows from the fact that the fibers of the map $f$ along a fixed orbit 
depend only on the combinatorics of the subdivision of the  corresponding face. 
We thus have   $^{\frak p}\!{\cal H}^{\pm1}(f_*\rat_{X_{\widetilde{C_4}}}[4])_{|O_i} \simeq \oplus_i \rat_{O_i}[1]$
and  $^{\frak p}\!{\cal H}^{\pm2}(f_*\rat_{X_{\widetilde{C_4}}}[4])\simeq \oplus_i H^6(f^{-1}(p_i)) 
\simeq \oplus_i H^6(\widetilde{C_3})_{p_i}\simeq \oplus_i \rat_{p_i} $.
The decomposition theorem  reads:
$$
f_*\rat_{X_{\widetilde{C_4}}}[4]\; \simeq\;  IC_{{C_4}}\, 
 \oplus\, (\oplus_i V_{p_i})  \,
\oplus\,  (\oplus_i(  IC_{\overline{O}_i}[1]
\oplus IC_{\overline{O}_i}[-1]))\,
\oplus \, (\oplus_i( \rat_{p_i}[2] \oplus \rat_{p_i}[-2])).
$$
The vector spaces $V_{p_i}$ are subspaces of $H^4(f^{-1}(p_i))$,
and contribute to the zero perversity term     
 $^{\frak p}\!{\cal H}^{0}(f_*\rat_{X_{\widetilde{C_4}}}[4])$. In order to determine their dimension, 
we compute the stalk 
$${\cal H}^0(f_*\rat_{X_{\widetilde{C_4}}}[4] )_{p_i}=H^4(f^{-1}(p_i))=H^4(\widetilde{C_3}).$$
As we saw above, ${\rm dim \,}H^4(\widetilde{C_3})=11$.  
By the support condition ${\cal H}^0({\cal IC}_{C_4})=0$ and, since 
${\cal IC}_{\overline{O}_i}=\rat_{\overline{O}_i}[1]$,
we get 
$$11={\rm dim \,}{\cal H}^0(f_*\rat_{X_{\widetilde{C_4}}}[4] )_{p_i}=
{\rm dim \,}V_{p_i} \oplus (\oplus_{\overline{O}_j \ni p_i}{\cal H}^{-1}({\cal IC}_{\overline{O}_j}))=
{\rm dim \, }V_{p_i}+6,$$ 
since only six curves $\overline{O}_j$ pass through $p_i$.
Hence ${\rm dim }V_{p_i}=5$ and finally 
$$
f_*\rat_{X_{\widetilde{C_4}}}[4]\simeq {\cal IC}_{C_4} 
\oplus(\oplus_{i=1}^8 (\rat^{\oplus 5}_{p_i} \oplus \rat_{p_i}[2] \oplus \rat_{p_i}[-2]) 
\oplus (\oplus_{i=1}^{24} (\rat_{\overline{O}_i} \oplus \rat_{\overline{O}_i}[2])).
$$
By taking the cohomology we get:

\noindent
$ \sum \hbox{ dim }I\!H^{2k}(X_{C_4})t^k=\sum \hbox{ dim } H^{2k}(X_{\widetilde{C_4}})t^k -8(t+5t^2+t^3)-24(t+2t^2+t^3)=
t^4+44t^3+102t^2+44t+1-8(t+5t^2+t^3)-24(t+2t^2+t^3)= t^4+12t^3+14t^2+12t+1=
h(C_4,t)$, as computed in (\ref{4cube}).

\subsection{Semismall maps}
\label{semismall}

Semismall maps occupy a very special place
in the applications of the theory of perverse sheaves to geometric 
representation theory. Surprisingly,  many maps which
arise naturally from Lie-theoretic objects are semismall. 
In a sense which we will try to illustrate
in the discussion of the examples below, the semismallness 
of a map is related to the  semisimplicity 
of the algebraic object under consideration. 
We  limit ourselves to proper and surjective semismall maps with a nonsingular domain.

In the case of semismall maps, the decomposition theorem takes the  particularly simple form of Theorem \ref{dtss}. 
Corollary  \ref{tseoendal}, on the semisimplicity
of the algebra of endomorphisms of the direct image, is a simple consequence.

As we have showed in  \ci{herdlef} the proof of Theorem 
\ref{dtss} is reduced to the proof of 
  the non degeneration
of certain bilinear forms defined on the homology groups of the fibers via intersection theory. We discuss
this point of view in $\S$\ref{interform}.

We discuss two examples of semismall maps: the resolution of the nilpotent cone
($\S$\ref{exssm1}) and the resolution of the $n$-th symmetric product
of a  nonsingular surfaces via Hilbert scheme of  $n$ points on the surface
($\S$\ref{exssm2}). In the first case, the decomposition theorem
leads to a simplified description of the Springer correspondence;
this correspondence (see  
 Theorem \ref{geomweylgroup})
gives  a geometric realization of the Weyl group of a semisimple linear algebraic group and its representations.
In the second case, we recall the basic geometric facts about 
Hilbert schemes that lead to 
the remarkably explicit Theorem  \ref{dthilbtm}.

A {\em stratification } for $f$ is  a decomposition of $Y$ into finitely many 
locally closed nonsingular subsets such that $f^{-1}(S_k) \to S_k$ is a topologically locally trivial fibration.
The subsets $S_k$ are called {\em strata}.

The following easy observation makes perverse sheaves enter this picture.
\begin{pr}
\label{ss}
Let $X$ be a connected nonsingular  $n$-dimensional variety, 
and $f:X \to Y$ be a proper surjective map of varieties. 
Let $Y = \coprod_{k=0}^n S_k$ be a stratification for $f$.
Let $y_k \in S_k$ and set $d_k:={\dim}f^{-1}(y_k)={\dim}f^{-1}(S_k)-{\dim}S_k$.
The following are equivalent:

\begin{enumerate}
\item
  $f_* \rat_X[n]$ is a perverse sheaf on $Y$;

\item
 $\dim {X \times_Y X }\leq n$; 

\item
 $\dim{S_k}+2d_k \leq \dim{X}$,   for every  $k=0, \ldots ,n$. 
\end{enumerate}
\end{pr}

\begin{defi}
{\rm A proper and surjective map $f$ satisfying one of the equivalent properties in Proposition \ref{ss}
is said to be  {\em semismall}.}
\end{defi}

\begin{defi}
{\rm
Let $X, Y,S_k$ and $ d_k$ be as in Proposition \ref{ss}. 
A stratum $S_k$ is said to be {\em relevant } if $\dim{S_k}+2d_k = \dim{X}$.
}  
\end{defi}

A semismall map $f: X \to Y$   must be finite over an open dense
stratum in $Y$  in view of property 3. Hence,  semismall maps are generically finite.
The converse is not true, e.g. the blowing-up of a point in $\comp^3$.
Note that, since ${\dim Y}=\dim X$, a relevant stratum has  even codimension.

\begin{rmk}
\label{small}
{\rm If the stronger inequalities  
$\dim{ S_k}+2d_k < \dim{X}$ is required to hold for every non-dense 
stratum, then the  map is said to be {\em small}. 
In this case, $f_* \rat_X[n]$ satisfies the support and co-support conditions
for intersection cohomology (\ref{sprt},\ref{csprt} of $\S$\ref{subsecintcoh}). 
Hence, if $Y_{o} \subseteq Y$ denotes a  nonsingular dense open subset 
over which $f$ is a covering,  then we have that  $f_* \rat_X[n]=IC_{Y}(L)$, where 
$L$ is the local system $f_*\rat_{X_{|Y_{o}}}$. 
}
\end{rmk}

\begin{ex}
{\rm  Surjective maps between surfaces are always semismall.
A surjective map  of threefolds is semismall iff no divisor $D \subseteq X$
is contracted to a point on $Y$.
}
\end{ex}

A great wealth of examples of semismall maps is furnished by contractions of 
(holomorphic) symplectic varieties,
which we now describe. A nonsingular quasi-projective complex variety
is called {\em holomorphic symplectic} if there is a holomorphic 
$2$-form $\omega \in \Gamma(X, \Omega_X^2)$ which is closed and nondegenerate, that is $d\omega=0$, and
$\omega^{\frac{{\dim} X}{2}}$ does not vanish at any point. The following is proved in \ci{kaledin}:

\begin{tm}
Let $X$ be a quasi-projective holomorphic symplectic variety, and $f:X \to Y$ a
projective birational map. Then $f$ is semismall. 
\end{tm}

Some important examples of semismall maps which are  contractions of  holomorphic symplectic varieties will be considered in \S\ref{exssm1}
and \S\ref{exssm2}.

The decomposition theorem for  a semismall map takes a particularly
simple form: the only contributions
come from the  relevant strata  $S$  and they consist
of non trivial summands $IC_{\overline{S}}(L)$ where  the local systems $L$ turn out to have finite monodromy.
 
Let $S$ be a relevant stratum, $y \in S$ and let
$E_1, \ldots, E_{{l}}$ be the irreducible $\dim{S}$-dimensional components of $f^{-1}(y)$.
The monodromy of the $E_i$'s defines a group homomorphism $\rho_S$ of the
fundamental group
$ \pi_1(S,y)$ to the group of permutations of the $E_i$'s, 
and, correspondingly, a local system of $\rat$-vector spaces $L_S$.
The semisimplicity of $L_S$ then follows immediately from the fact that the monodromy factors through a finite group. 
With this notation, let us give the statement of 
the decomposition theorem in the case of semismall maps:

\begin{tm}
\label{dtss}
 {\rm \bf (Decomposition theorem for semismall maps) } Let $I_{rel}$ be the set of relevant strata, and, for each $S \in I_{rel}$, let
$L_S$ be the corresponding local system with finite monodromy defined above.
There is a canonical isomorphism  in ${\cal P}_Y:$
\begin{equation}
\label{dtssmall1}
f_* \rat_{X}[n] \;\simeq \;
\bigoplus_{S \in I_{rel}} IC_{\overline{S}}(L_S).
\end{equation}
\end{tm}

Let $Irr(\pi_1(S))$ be the set of irreducible representations
of $\pi_1(S,y).$ 
For $\chi \in Irr(\pi_1(S))$, we denote by ${\cal L}_{\chi}$
the corresponding local system on $S$.
We have an isotypical decomposition in the category $\pi_1(S)$-${\rm Mod}$
of representations of $\pi_1(S)$
$$\rho_S\;  \simeq \; \bigoplus_{\chi \in Irr(\pi_1(S))} \chi \otimes V_S^{\chi},$$
where $V_S^{\chi}$ is a vector space whose dimension is 
the multiplicity of the representation $\chi$ in $\rho_S$.
Correspondingly, we have a decomposition of local systems
$L_S=\bigoplus_{\chi \in Irr(\pi_1(S))}{\cal L}_{\chi} \otimes V_S^{\chi}$, 
and, for each term $IC_{\overline{S}}(L_S)$ in (\ref{dtssmall1}), an isotypical decomposition
\begin{equation}
\label{decomssls}
IC_{\overline{S}}(L_S)\; \simeq  \; \bigoplus_{\chi \in Irr(\pi_1(S))} IC_{\overline{S}}({\cal L}_{\chi})\otimes V_S^{\chi}.
\end{equation}

\medskip
The second special feature of semismall maps concerns the endomorphism
algebra ${\rm End}_{{\cal D}_Y}(f_* \rat_X[n])$;
see  \ci{chrissginzburg,cortihana}  for details.

By Schur's lemma,
for $\chi$ and ${\cal L}_{\chi}$ as above, we have that
$$
{\rm End}_{{\cal D}_Y}(IC({\cal L}_{\chi}))=
{\rm End}_{{\cal D}_S}({\cal L}_{\chi})=
{\rm End}_{ \pi_1(S)-{\rm Mod}} (\chi)
$$
is a division ring $R_{\chi}$.
The intersection cohomology sheaves $IC({\cal L}_{\chi})$
are simple objects in the category of perverse sheaves and  Theorem \ref{dtss}
can be restated by saying that $f_* \rat_X[n]$ is a semisimple perverse sheaf.
We thus have the following 
\begin{cor}
\label{tseoendal}
{\rm \bf (Semisimplicity of the endomorphism algebra) }Let $f:X \to Y$ be a semismall map. Then the endomorphism algebra 
${\rm End}_{{\cal D}_Y}(f_* \rat_X[n])$ is semisimple, that is 
isomorphic to a direct sum of matrix algebras over division rings. In fact, 
we have:
\begin{equation}
\label{semisimplicity}
{\rm End}_{{\cal D}_Y}(f_* \rat_X[n]) \;\simeq  \;
\bigoplus_{S \in I_{rel}} {\rm End}_{{\cal D}_Y}(IC_{\overline{S}}(L_S))
\;\simeq\;
\bigoplus_{{S \in I_{rel}}\atop {\chi \in Irr(\pi_1(S))}}R_{\chi}\otimes{\rm End}(V_S^{\chi})
\end{equation}
\end{cor}

Furthermore, if $H^{BM}_{2n}(X \times_Y X)$ is given the structure of algebra 
coming from the composition of correspondences, then
there is an isomorphism of algebras see \ci{cortihana},
Lemma 2.23
\begin{equation}
\label{endo}
{\rm End}_{{\cal D}_Y}(f_* \rat_X[n])\simeq H^{BM}_{2n}(X \times_Y X),
\end{equation}

The endomorphism algebra contains in particular the idempotents giving the projection of
$f_* \rat_X[n]$ on the irreducible summand of the 
canonical decomposition (\ref{dtssmall1}).
Since, again by semismallness,  $H^{BM}_{2n}(X \times_Y X)$ is the top
dimensional Borel Moore homology, it is generated by the irreducible components of $X \times_Y X$.
The projectors are therefore  realized by algebraic correspondences. 

This has been pursued in \ci{decmigsemi}, where we prove, 
 in accordance with the general philosophy of \ci{cortihana},
 a ``motivic'' refinement of the decomposition theorem in the case of semismall maps. 
In particular, it is possible to construct a (relative) 
Chow motive corresponding to the intersection cohomology groups 
of singular varieties which admit a semismall resolution.

\subsubsection{Semismall maps and intersection forms}
\label{interform}
Let $f : X \to Y$ be a semismall map. Every stratum yields a  bilinear form
on a certain homology group which has a neat geometric interpretation in terms of basic  intersection theory on $X$. 
Theorem \ref{pasolini} below, states that 
the decomposition theorem for the semismall map $f$ turns out to be equivalent to  
the non degeneracy of all these intersection forms.  

Let us describe these intersection forms.
If a stratum is
not relevant, then, as noted below, the construction  that follows yields a trivial homology group. 
Let $S\subseteq Y$ be a relevant stratum, and $y \in S$.
Let $\Sigma$ be a local  trasversal slice to $S$ at $y$, given for example by 
intersecting a small contractible  Euclidean neighborhood of $y$ with the  complete intersection of 
$\dim{S}$ general  hyperplane sections in $Y$ passing through $y$. 
The restriction $f_|:f^{-1}(\Sigma) \to \Sigma$ is still semismall and 
$d=\dim{ f^{-1}(y)}=(1/2) \dim{f^{-1}( \Sigma)}$.
By composing the chain of maps:
$$
H_{2d}(f^{-1}(y))=H_{2d}^{BM}(f^{-1}(y))\lorw H_{2d}^{BM}(f^{-1}(\Sigma))
\simeq H^{2d}(f^{-1}(\Sigma)) \lorw H^{2d}(f^{-1}(y)),
$$
where the first map is the push-forward with respect to a closed inclusion 
and the second is the restriction,
we obtain  the {\em intersection pairing} (cf. $\S$\ref{famfa})
{\em associated with the relevant stratum $S$}
$$
I_S:H_{2d}(f^{-1}(y)) \times H_{2d}(f^{-1}(y))\lorw \rat.
$$
Of course, we have used the usual identification
$\mbox{Bil}(U,U) \simeq  \mbox{Hom}(U, U^{\vee})$.
If the stratum is not relevant then $\dim{f^{-1}(y)} <d$ and  $H_{2d}(f^{-1}(y))=0$,
and the intersection form is defined, it is trivial and also nondegenerate,
in the sense that the corresponding  linear map is an isomorphism of trivial vector spaces.

A basis of $H_{2d}(f^{-1}(y))$ is given by the classes of the 
$d$-dimensional irreducible components $E_1, \ldots, E_l$ of 
$f^{-1}(y).$ The intersection pairing $I_S$ is then represented
by the intersection matrix $||E_i \cdot \,E_j||$ of these components, computed  in
the, possibly disconnected, manifold 
$f^{-1}(\Sigma)$.

In what follows, for simplicity only, let us assume that 
 $S=S_k$ is a connected stratum of dimension $k$, relevant or not.
Let $U= \coprod_{k' >k} S_{k'}$ be the union of the strata of dimension strictly bigger
than $k$
and  $U'=U \coprod S$.
Denote by $i:S \to U' \longleftarrow U:j$ 
the corresponding imbeddings. 
The intersection map
$
H_{2d}(f^{-1}(y)) \to H^{2d}(f^{-1}(y))
$
is then  identified with the natural map of stalks
$$
\xymatrix{{\cal H}^{-d}(i^!f_* \rat_{U'}[n])_y \ar[r] &
          {\cal H}^{-d}(i^{*}f_* \rat_{U'}[n])_y.}$$

By Remark  \ref{splitdecmig}, the non-degeneracy of $I_S$ is equivalent to the 
existence of a canonical isomorphism:
\begin{equation}
\label{spli}
f_* \rat_{U'}[n]\, \simeq\, j_{!*}f_* \rat_{U'}[n]  \,
\bigoplus\, {\cal H}^{-{\rm dim }S}(i^{!}f_* \rat_{U'}[n])[{\rm dim }S].
\end{equation}

It follows that the splitting behavior of $f_* \rat_X[n]$ is governed precisely
by the nondegeneracy of the forms $I_S$.  

In our paper
cf. \ci{herdlef}, we proved, using classical Hodge-Lefschetz theory, that for every relevant stratum $S$ with typical
fiber of dimension $d$, 
the form $I_S$ has a precise sign. In particular, all forms $I_S$ are non degenerate. We summarize
these results in the following

\begin{tm}
\label{pasolini}
Let $f: X \to Y$ be a semismall map with $X$ nonsingular. 
Then the  statement of the decomposition theorem is equivalent to
the non degeneracy of the intersection forms $I_S$. These forms
are nondegenerate and if a connected component of a   stratum $S$  is relevant with 
typical
fiber of dimension $d$, then the form $(-1)^d I_S$ is positive definite.
\end{tm}

\subsubsection{Examples of semismall maps I:   Springer theory }
\label{exssm1}
References for what follows are  \ci{chrissginzburg,springerbourbaki}.
Let $G$ be a semisimple connected linear algebraic group with Lie algebra $\frak g$, 
let $T \subseteq G$
be  a maximal torus,
let $B$  be a Borel subgroup containing $T$ and  let $W$ the Weyl group.
The {\em flag variety} $G/B$ is complete and parametrizes the Borel subalgebras 
of $\frak g$.
We recall that an element $x\in \frak g$ is {\em nilpotent} (resp. {\em semisimple}) if the endomorphism $[x, -]:\frak g \to \frak g$ is nilpotent (resp. diagonalizable). If $\dim{ \ke\,{ [x, \,-\,]}}$ equals the dimension of $T$, then $x$ is said to be {\em regular}.

Let ${\cal N}\subseteq  \frak g$ the cone of nilpotent elements of $ \frak g$. 
It can  be easily  shown (cf. \ci{chrissginzburg}) that
$$
\widetilde{\cal N}=\{ (x, \frak c)\in {\cal N}\times G/B: \frak c \hbox{ is a Borel subalgebra of  } \frak g \hbox{ and } x \in {\cal N} \cap \frak c \}
$$
is isomorphic to the cotangent bundle $T^* \,G/B$ of the flag variety $G/B$, and is therefore endowed with a natural (exact) holomorphic symplectic form. 
The map $p:\widetilde{\cal N}  \to {\cal N}\subseteq  \frak g$, defined as $p(x, \frak c)=x$, is surjective, since every nilpotent element is contained in a 
Borel subalgebra, generically one-to-one, since a generic nilpotent element is contained in exactly one Borel subalgebra, proper, since 
$G/B$ is complete, and semismall, since $\widetilde{\cal N}$ 
is holomorphic symplectic. The  map $p$ is called the Springer resolution. 

\begin{ex}
\label{sl2}
{\rm If $G=SL_2$,  then the flag variety $G/B=\pn{1}$ and
the cotangent space is the total space of the line bundle ${\cal O}_{\pn{1}}(-2)$. The  variety
obtained by contracting   the zero-section to a point is isomorphic to the cone
with equation  $z^2=xy$ in $\comp^3$. If $\{H,X,Y\}$ denotes the usual basis of ${\frak s \frak l}_2$,
 the matrix $zH+xX-yY$ is nilpotent precisely when $z^2=xy$. 
}
\end{ex}

The aim of  the  {\em Springer correspondence} is to get an algebra isomorphism
between the  rational group algebra of the Weyl group $W$ of $\frak{g}$ and the algebra of correspondences of $\widetilde{\cal N}$
$$ 
\xymatrix{
\rat [W] \ar[r]^(.32){\simeq} & H^{BM}_{2{\dim \widetilde{\cal N}}}( \widetilde{\cal N}\times_{\cal N}\widetilde{\cal N})}$$
 so that the elements of the Weyl group will correspond to certain correspondences
in the fiber product above.

The Springer correspondence is realized as follows.
One constructs  an action of the Weyl group $W$ on
$p_* \rat_{\widetilde{\cal N}}[{\dim} \widetilde{\cal N} ]$.  This action extends
to an algebra homomorphism $\rat [W] \to {\rm End}_{{\cal D}_{\cal N}}(p_* \rat_{\widetilde{\cal N}}[{\dim \widetilde{\cal N}} ])$ which is verified to be an isomorphism. Finally, one uses (\ref{endo}).

We  now sketch,  following \ci{luszunipo} (see also \ci{bomac1,bomac2}),
the construction of the desired $W$-action.  
By a theorem of Chevalley, there is a map $q:{\frak g}\to {\frak t}/W$ defined as follows:
any $x \in {\frak g}$ has a unique expression $x=x_{ss}+x_{n}$ where
$x_{ss} $ is semisimple, $x_n$ is nilpotent and commutes with
  $x_{ss}$. Then $x_{ss}$ is conjugate to an element of 
$ {\frak t}$, well defined up to the action of $W$. The quotient  ${\frak t}/W$ is an affine space.
Let us denote by ${\frak t}^{rs}={\frak t} \setminus \{ \hbox{\rm root hyperplanes }\}$,  
the set of regular elements in ${\frak t}$,
and by  ${\frak g}^{rs}= q^{-1}({\frak t}^{rs}/W)$ the set of regular semisimple elements in ${\frak g}$.
The set  ${\frak t}^{rs}/W$
is the complement of a divisor  $ \Delta \subseteq {\frak t}/W$.
The map $q:{\frak g}^{rs} \to {\frak t}^{rs}/W $ 
is a fibration with fiber $G/T$.  There is the   monodromy representation 
$\rho: \pi_1({\frak t}^{rs}/W) \to {\rm Aut}(H^*(G/T))$.

\begin{ex}
{\rm Let $G=SL_n$. The map $q$ sends  a traceless matrix to  the coefficients 
of its characteristic polynomial. The set ${\frak t}^{rs}/W={\frak t}/W \setminus \Delta$ is the set of polynomials 
with  distinct roots. The statement that the map $q:{\frak g}^{rs} \to {\frak t}^{rs}/W $ is a fibration
 boils down to 
the fact that a matrix commuting with a diagonal matrix with distinct eigenvalues must be diagonal, and that the 
adjoint orbit of such a matrix is closed in ${\frak s \frak l}_n$.}
\end{ex}

Let us define  
$$\widetilde{\frak g}=\{
(x, \frak c)\in \frak g \times G/B : \frak c \hbox{ is a Borel subalgebra of  } \frak g \hbox{ and } x \in \frak c
\}.$$
Let $p:\widetilde{\frak g} \to {\frak g}$ be the projection to the
first factor. This  map ``contains" the Springer resolution in the sense that
$\widetilde{\cal N} = p^{-1}({\cal N})\subseteq \widetilde{\frak g}$.

 On the other hand, the $G$-orbits in  $\widetilde{\frak g}$ of regular semisimple
elements (i.e. for which the corresponding $x$ is regular semisimple) are affine varieties
isomorphic to  $G/T$ and diffeomorphic to $\widetilde{\cal N}$.

\begin{ex}
\label{doublepoint}
{\rm In the case discussed in Example  \ref{sl2},
one  considers the family of affine quadrics $Y_t \subseteq \comp^3$ of equation
$z^2=xy+t$ for $t \in \comp$.
For $ t\neq 0$, $Y_t$ is diffeomorphic, but not isomorphic, to 
$T^*{\pn{1}}$, while, for $t=0$, $Y_0$ is the nilpotent cone of  ${\frak s \frak l}_2$. Pulling back this family by the map $t \to t^2$, we get the family $z^2=xy+t^2$ which admits a simultaneous (small) resolution, whose fiber at $t=0$ is the map $T^*{\pn{1}} \to Y_0$. 
}
\end{ex}

The Weyl group acts simply transitively on the set of Borel subgroups 
containing a regular semisimple element. 
 Setting $\widetilde{\frak g}^{rs}=p^{-1}({\frak g}^{rs})$, this observation 
 leads to the following:

\begin{pr}
\label{simul} 
The restriction $p':\widetilde{\frak g}^{rs} \to {\frak g}^{rs}$
is an (unramified) covering   map with Galois group $W$. 
The map $p:\widetilde{\frak g} \to \frak g$ is small.
\end{pr}

We summarize what we have discussed so far in the following diagram (the map $r$
is  defined below):
$$
\xymatrix{
         &                  &                        & \widetilde{\cal N}  \ar@{^{(}->}[dl] \ar[dd]_(.6){p}|\hole           &                            \\
\widetilde{\frak g}^{rs} \ar[dd]_(.6){p'}\ar@{^{(}->}[rr] & & \widetilde{\frak g} \ar[dd]_(.6)p\ar[rr]    &                   &  {\frak t}\ar[dd]      \\
                      &   G/T   \ar@{^{(}->}[dl]\ar[rr]^(.7){r}|\hole   &      &  {\cal N}\ar@{^{(}->}[dl]^i         &                             \\  
\widetilde{\frak g}^{rs}/W= {\frak g}^{rs}  \ar@{^{(}->}[rr]    & &   {\frak g}  \ar[rr]^q           &                   &  {\frak t}/W 
}
$$ 
Let $L=p'_*\rat_{\widetilde{\frak g}^{rs}}$ be the  local system 
associated with the $W$-covering,
which, by its very definition, 
is endowed with  an action of the Weyl group $W$. Since any map between local
systems extends uniquely to a map between the associated intersection cohomology complexes 
(see \ci{goma2}, Theorem 3.5),  we have an action of $W$ on
$IC_{\frak g}(L)$. 
Since $p$ is  small, by Remark \ref{small},
$IC_{\frak g}(L)=p_* \rat_{\widetilde{\frak g}}[\dim{\frak g}]$. 

In particular, there is an action of $W$ on 
$ i^*p_* \rat_{\widetilde{\frak g}}[{\dim} \frak g]=p_* \rat_{\widetilde{\cal N}}[{\dim} \widetilde{\cal N}]$,
and this is the sought-for $W$-action.

A perhaps more intuitive way to realize this action is the following.
We have  ${\cal N}=q^{-1}(0)$.
There is a continuous retraction map $r:G/T \to {\cal N}$.
Since the affine variety $G/T$ is diffeomorphic to $\widetilde{\cal N}$,
we have an isomorphism:
$$ r_* \rat_{G/T}[{\dim} \widetilde{\cal N} ] \simeq p_* \rat_{\widetilde{\cal N}}[{\dim} \widetilde{\cal N} ]. $$
As we  have already observed, the monodromy  of the fibration   
$q:{\frak g}^{rs} \to {\frak t}^{rs}/W $ gives an action of
$\pi_1({\frak t}^{rs}/W )$ on  $r_* \rat_{G/T}[{\dim} \widetilde{\cal N}]$.
There is an exact sequence of groups:
$$
0 \to \pi_1({\frak t}^{rs}) \to \pi_1({\frak t}^{rs}/W )\to W  \to 0
$$
and the existence of the simultaneous resolution $\widetilde{\frak g}$
shows that the monodromy factors through  an action of $W$,  and this yields the desired 
alternative description of the $W$ action on $p_* $ etc.

As mentioned earlier, the $W$-action extends to an algebra homomorphism
$$
\rat[W] \lorw {\rm End}_{{\cal D}_{\cal N}}(p_* \rat_{\widetilde{\cal N}}[{\dim} \widetilde{\cal N} ])
\;\; = \;\; H^{BM}_{2{\dim}\, \widetilde{\cal N}}(\widetilde{\cal N}\times_{\cal N}\widetilde{\cal N} )
$$
and we have
\begin{tm} {\rm (\ci{bomac1})}
\label{geomweylgroup}
The map 
\begin{equation}
\label{isoweyl}
\rat[W] \lorw H^{BM}_{2{\dim}\, \widetilde{\cal N}}(\widetilde{\cal N}\times_{\cal N}\widetilde{\cal N} )
\end{equation}
constructed above is an isomorphism of algebras.
\end{tm}

We thus have a geometric construction of the group ring of the Weyl group $W$
as an algebra of (relative) correspondences on $\widetilde{\cal N}$,
and a natural basis given by the irreducible components of 
$\widetilde{\cal N}\times_{\cal N}\widetilde{\cal N}$.

A deeper investigation of the isomorphism (\ref{isoweyl})
sheds light on the irreducible representations  of 
$\rat[W]$, or, equivalently, of $W$, by giving a natural geometric construction of these
representations: 
the nilpotent cone ${\cal N}$ has a natural $G$-invariant stratification, 
given by the orbits of the adjoint action contained in ${\cal N}$, i.e. by
the conjugacy classes of nilpotent elements.
Let ${\rm Conj}({\cal N})$
be the set of conjugacy classes of nilpotent elements in ${\frak g}$.
For $[x] \in {\rm Conj}({\cal N})$,  let $x$ be a representative, and denote by
${\frak B}_x:=p^{-1}(x)$ the fiber over $x$ and  by $S_x=Gx$ the stratum 
 of  ${\cal N}$ containing $x$.

\begin{ex}
{\rm Let $G=SL_n$. Each conjugacy class contains exactly one matrix which is a sum of Jordan matrices,
so that
the $G$-orbits are parameterized by the partitions of the integer $n$. The open dense stratum of ${\cal N}$
corresponds to the Jordan block of length $n$.}
\end{ex}

It can be proved (cf.  \ci{steinberg, spa}),  that every stratum $S_x$ is relevant and that all the components
of ${\frak B}_x$ have the same dimension $d_x$.
The vector space $H_{2d_x}({\frak B}_x)$, 
generated by the irreducible components of ${\frak B}_x$,
is, by construction, a representation of $W$. This representation is not necessarily irreducible, as the  finite group of connected components of the stabilizer $G_x$ of $x$ acts. This action commutes with the action of $W$, and splits $H_{2d_x}({\frak B}_x)$. It can be showed that 
every  irreducible representation of $W$ is realized as a direct summand
of some $H_{2d_x}({\frak B}_x)$. At this point, we refer the reader
to the original papers \ci{spri1, bomac1, bomac2},
and  to the book \ci{chrissginzburg}.

 \subsubsection{Examples of semismall maps II:  Hilbert schemes of  points}
\label{exssm2}
A reference for what follows is  \ci{nakajilectures}. 
The $n$-th symmetric product  $(\comp^2)^{(n)}=(\comp^2)^n/{\cal S}_n$,
parametrizing $0$-cycles  $Z=\sum_k n_k p_k$  of $\comp^2$ of length
$n$, is singular. Singularities appear when some points come together, 
that is at cycles $\sum_k n_k p_k$ where some multiplicity is bigger than one.
The Hilbert scheme  $X= (\comp^2)^{[n]}$ is a certain a resolution of singularities of  $n$-th  $(\comp^2)^{(n)}=(\comp^2)^n/{\cal S}_n$ which keeps tracks of the ``tangent'' information when two or more points collapse.
For instance,  $X= (\comp^2)^{[2]}$ is the blow-up of  
$X=(\comp^2)^{(2)}$ along the diagonal, consisting of cycles of type $2p_1$.

When $n$-points come together at a point $p_0$ of coordinates $(x_0,y_0)$
this tangent information is encoded  as a  {\em scheme} structure supported
on the point parametrizing the cycle $n p_0$. This  scheme structure is given
in terms of an  ideal $I$ of the ring of polynomials in two indeterminates $\comp[X,Y]$
with radical $\sqrt I=(X-x_0,Y-y_0)$, the maximal ideal of the point $p_0$,
such that 
${\dim}_{\comp}{\comp[X,Y]/I}=n$.
In general, the points of the variety $X= (\comp^2)^{[n]}$ parametrize ideals $I \subseteq \comp[X,Y]$ 
such that ${\dim}_{\comp}{ \comp[X,Y]/I}=n$. Every such ideal is the product 
$\prod I_k$ 
of ideals supported at points $p_k \in \comp^2$ and we can  
associate with it the $0$-cycle $Z(I):= \sum_k n_k p_k$, where 
$n_k={\dim}_{\comp}{\comp[X,Y]/I_k}$, called the support of this ideal. 
Then $n= \sum_k n_k$ and $Z(I)$ is a point in the symmetric product $(\comp^2)^{(n)}=(\comp^2)^n/{\cal S}_n$. The   {\em Hilbert-Chow map} $\pi:(\comp^2)^{[n]} \to (\comp^2)^{(n)}$,
sending with $I$ to  its support $Z(I)$, is well-defined and proper. It is an isomorphism precisely on the set $(\comp^2)^{(n)}_{reg}$
corresponding to cycles $p_1+ \ldots + p_n$ consisting of $n$ distinct points, as in this case there is only one possible scheme structure. 
Let $(x_1,y_1), \ldots , (x_n,y_n)$ be coordinates on $(\comp^2)^n$.
The form $ \sum_k dx_k \wedge dy_k$ on $(\comp^2)^n$ is ${\cal S}_n$-invariant and descends to a closed and nondegenerate form on  $(\comp^2)^{(n)}_{reg}$. A local computation shows that its pullback by $\pi$ extends to a symplectic form on $(\comp^2)^{[n]}$.  
In particular $\pi$ is semismall (this can be also verified directly).
The subvariety $ (\comp^2)^{[n]}_0$ of subschemes supported at $0$ 
is called the punctual Hilbert scheme of length $n$. 
Its points parametrize the $n$-dimensional quotient rings of $\comp[X,Y]/(X,Y)^{n+1}$. 
These punctual Hilbert schemes   have been studied in depth,
see \ci{iarrobino, briancon}, for example.  They are  irreducible, of dimension $n-1$, and admit 
a disjoint-union-decomposition
into affine spaces. Clearly, $ (\comp^2)^{[n]}_0\simeq (\pi^{-1}(np))_{red}$, for every $p \in \comp^2$. Similarly, if $Z:= \sum_k n_k p_k$ with 
$p_i \neq p_j$  for all $i\neq j$, then 
$(\pi^{-1}(Z))_{red}\simeq \prod_i  (\comp^2)^{[n_i]}_0$.
The construction can be globalized, in the sense that, for any nonsingular surface $S$, the Hilbert scheme 
$S^{[n]}$ is nonsingular and there is a map $\pi:S^{[n]} \to S^{(n)}$ which is semismall, 
and locally, in the analytic topology, isomorphic to $\pi:(\comp^2)^{[n]} \to (\comp^2)^{(n)}$.
There also exists a  version of $S^{[n]}$ for  a symplectic manifold 
$S$ of real dimension four, which was defined and 
investigated by C.Voisin in  \ci{voisin}.

To describe the strata of the map $\pi$, we denote by ${\frak P}_n$ the set of  partitions     of the natural number $n$.          . 
Let $\nu=(\nu_1, \ldots,\nu_{l(\nu)}) \in {\frak P}_n$, 
so that  $\nu_1 \geq \nu_2\geq \ldots \geq \nu_{l(\nu)}$ and $\sum_i \nu_i=n$.
We will also write $\nu= 1^{a_1}2^{a_2}\ldots n^{a_n}$, with $\sum ka_k=n$, where  $a_i$ is the number of times that
the number  $i$ appears in the partition $\nu$. Clearly $l(\nu)= \sum a_i$.
We consider the following stratification of $S^{(n)}$: for $\nu \in {\frak P}_n $
we set 
$$
S_{(\nu)}=\{ \mbox{$0$-cycles}  \subseteq S^{(n)}\hbox{ of type } \nu_1p_1+ \ldots + 
\nu_{l(\nu)}p_{l(\nu)} \hbox{ with } p_i \neq p_j, \; \forall i\neq j \}.
$$
Set $S_{[\nu]}=\pi^{-1}(S_{(\nu)})$ (with the reduced structure). The variety $S_{(\nu)}$ is nonsingular of dimension $2l(\nu)$.
It can be shown that $\pi: S_{[\nu]} \to S_{(\nu)}$ is locally trivial with fiber
isomorphic to the product $ \prod_i (\comp^2)^{\nu_i}_0$ of punctual Hilbert schemes.
In particular, {\em the fibers of $\pi$ are irreducible,} hence the local systems occurring in (\ref{dtss}) are constant of rank one.
Furthermore, the closures $\overline{S}_{(\nu)}$ and their desingularization can be explicitly determined.
If $\nu$ and $\mu$ are two partitions, we say that $\mu \leq \nu$ if there exists a decomposition
$I_1, \ldots, I_{l(\mu)}$ of the set $\{1, \ldots ,l(\nu)\}$ such that
$\mu_1=\sum_{i\in I_1} \nu_i, \ldots ,\mu_{l(\mu)}=\sum_{i\in I_{l(\mu)}} \nu_i. $
Then
$$
\overline{S}_{(\nu)}=\coprod _{\mu \leq \nu}S_{(\mu)}.
$$
This reflects just the fact that a cycle $\sum \nu_ip_i \in S_{(\nu)}$ can degenerate to a cycle 
in which some of the $p_i'$s come together.
If $\nu= 1^{a_1}2^{a_2}\ldots n^{a_n}$, we set $S^{(\nu)}= \prod_i S^{(a_i)}$
 (product of symmetric products). The variety $S^{(\nu)}$
has dimension $2l(\nu)$, and there is a natural {\em finite} map $\nu:S^{(\nu)} \to \overline{S}_{(\nu)}$, which is an isomorphism
when restricted to $\nu^{-1}(S_{(\nu)})$. Since $S^{(\nu)}$ has only quotient singularities, it is normal, so that
$\nu:S^{(\nu)} \to \overline{S}_{(\nu)}$ is the normalization map, and $IC_{\overline{S}_{(\nu)}}=\nu_*\rat_{S^{(\nu)}}[2l(\nu)]$.
\begin{tm}
\label{dthilbtm}
The decomposition theorem {\rm \ref{dtss}} for $\pi: S^{[n]} \to S^{(n)}$ gives a canonical isomorphism:

\begin{equation}
\label{dthilbert}
\pi_*\rat_{S^{[n]}}[2n] \; \; \simeq \bigoplus_{\nu \in {\frak P}_n}\nu_*\rat_{S^{(\nu)}}[2l(\nu)].
\end{equation}
Taking cohomology, {\rm (\ref{dthilbert})} gives
\begin{equation}
\label{vfcl}
H^i(S^{[n]},\rat)= \bigoplus_{\nu \in {\frak P}_n}H^{i+2l(\nu)-2n}(S^{(\nu)},\rat).
\end{equation}

\end{tm}

This explicit form was given by L. G\"ottsche and W. Soergel in \ci{gottso} as an application
of M. Saito's  theorem \ci{samhp}. 
Since $S^{(n)}$ is the quotient of the nonsingular variety $S^n$ by the finite group ${\cal S}_n$,
its {\em rational} cohomology $H^i(S^{(n)},\rat)$ is just the ${\cal S}_n$-invariant part of $H^i(S^n,\rat)$.
In \ci{macdona}, MacDonald determines  the dimension of such invariant subspace. 
His result is more easily stated in terms of 
generating functions:
$$
\sum {\dim}\,  H^i(S^{(n)},\rat)t^iq^n=
\frac{(1+tq)^{b_1(S)}(1+t^{3}q)^{b_3(S)}}
{(1-q)^{b_0(S)}(1-t^{2}q)^{b_2(S)}(1-t^{4}q)^{b_4(S)}}.
$$ 
With the help of this formula and (\ref{vfcl}), 
we find ``G\"ottsche's Formula'' 
for the generating 
function of the Betti numbers of the Hilbert scheme:
\begin{equation}
\label{gottfor}
\sum_{i,n} {\dim} \, H^i(S^{[n]},\rat)t^iq^n= 
\prod_{m=1}^{\infty}\frac{(1+t^{2m-1}q^m)^{b_1(S)}(1+t^{2m+1}q^m)^{b_3(S)}}
{(1-t^{2m-2}q^m)^{b_0(S)}(1-t^{2m}q^m)^{b_2(S)}(1-t^{2m+2}q^m)^{b_4(S)}     }.
\end{equation}
\begin{rmk}
{\rm  Setting $t$$=$$-1$, we get the following simple formula for the 
generating function for the Euler characteristic:
$$
\sum_{n=0}^{\infty} \chi(S^{[n]})q^n= 
 \prod_{m=1}^{\infty}\frac{1}{(1-q^m)^{\chi(S)}}.
$$
See \ci{decaeuler}, for a rather elementary  derivation of this formula.
}
\end{rmk}

G\"ottsche's Formula appeared first in \ci{gotts}, following some preliminary work in the case $S=\comp^2$
by Ellingsrud and Stromme, (\ci{ellstr1, ellstr2}). The original proof relies on the Weil conjectures, 
and on a delicate counting of points over a finite field with the help of the cellular structure of 
the punctual Hilbert scheme following from Ellingsrud and Stromme's results.

\subsection{The functions-sheaves dictionary and geometrization}
\label{foncfasc}
 In $\S$\ref{hecke} and $\S$\ref{satake} we discuss two 
rather deep applications of the decomposition theorem to geometric representation
theory.  
Even though the applications can be stated and proved within the realm of complex geometry, they have been inspired by the Grothendieck's 
philosophy (\ci{grbourbaki}) (see also \ci{lau}, \S 1.1) of the 
{\em dictionnaire fonctions-faisceaux} for varieties defined over finite fields
This section is devoted to a brief explanation of  this philosophy.
The reader  who is unfamiliar with algebraic geometry over finite fields, 
may look at $\S$\ref{8000}.

Suppose that  $X_0$ is a variety defined over a finite field $\fq$ of cardinality $q$.
Associated with  any complex  of $l$-adic sheaves $K_0$ on $X_0$,  there is the function
 $t_{K_0}:X_0(\fq) \to \qlb$:
$$
t_{K_0}(x)= \sum_i (-1)^i {\rm Trace}(Fr:{\cal H}_x^i(K) \to {\cal H}_x^i(K)), 
$$
where $Fr$ is the Frobenius endomorphism of ${\cal H}^i_x(K)$. This function is additive with respect to distinguished triangles in $D^b_c(X_0, \qlb)$, 
multiplicative with respect to tensor
products of complexes, compatible with pullbacks, and satisfies the 
Grothendieck trace formula: If $f:X_0 \to Y_0$ is a proper map of 
$\fq$-schemes and $K_0 \in D^b_c(X_0, \qlb)$, then, for $y \in Y(\fq)$, 
$$
t_{f_*K_0}(y)=\sum_{x \in f^{-1}(y) \cap X_0(\fq)} t_{K_0}(x).
$$

Since $l$-adic sheaves on $X_0$ yield these trace-like functions, one may
think of 
replacing certain class of 
 {\em functions} on $X_0$ with (complexes of) $l$-adic {\em sheaves} on 
$X_0$.

The philosophy of {\em geometrization} is rooted in the fact that
quite often functions arising  from representation theory or combinatorics 
can be interpreted as associated with  sheaves -often perverse sheaves-  on  
algebraic varieties and   theorems about such functions become consequences of 
theorems about the corresponding sheaves.

If the  cohomology sheaves of $K_0$ are zero in odd degree and,
for every $i$, the 
eigenvalues of Frobenius (not just their absolute values!) on 
${\cal H}^{2i}_x(K_0)$ are equal to 
$q^i$,
then the function $t_{K_0}$  satisfies  the relation
$t_{K_0}(x)=\sum_i \dim {\cal H}_x^{2i}(K_0)q^i$.

 We can modify this  formula so that it makes sense  for a constructible
complex of sheaves $K$ 
on  a complex algebraic variety $X$. We do so by considering $q$ to be a 
free variable. If 
 the stalk cohomology vanishes in all odd degrees, then we obtain
 a Poincar\'e-like polynomial for $K$.
 
This is the case in the two examples we discuss in \S\ref{hecke} and $\S$ \ref{satake}.

In $\S$\ref{hecke} we show that the 
Kazhdan-Lusztig polynomials, which are
   associated   in a purely combinatorial way
(see (\ref{klp}))  with  the Weyl group $W$ of an algebraic 
group $G$, may be interpreted, via the functions-sheaves dictionary, 
as the Poincar\'e-like polynomials   of   the intersection  complexes of 
Schubert varieties  in the flag variety $G/B$ of $G$.
This fact allows a geometric interpretation of the Hecke algebra of $W$ as an algebra of equivariant perverse sheaves on the flag varieties.

 Similarly, in \S\ref{satake} we treat the case of   a certain class of functions arising from the classical Satake isomorphism and which 
are associated,  via the  functions-sheaves dictionary, 
with the  intersection  complexes of certain subvarieties of the
(infinite dimensional) affine Grassmannian.
This leads to a  geometrization of the classical  Satake isomorphism. 

In both of these situations, the strategy
towards geometrization is similar. We start with an algebra of functions on a group $G$ with some invariance property. For instance, in the case treated in \S\ref{hecke}, the group is a Chevalley group and the functions are the  left and right  invariant
with respect to  a fixed  Borel subgroup,  and  in  the case
treated in $\S$\ref{satake}, they are the functions on an algebraic group over a local field which are left and right invariant with respect to the maximal compact subgroup of points over its ring of integers.
This algebra has a natural basis, consisting of characteristic functions of double cosets which correspond, via the functions-sheaves dictionary, to the constant sheaves concentrated on some subvarieties of   a variety 
associated with $G$, i.e   the flag variety in the case of \S\ref{hecke} and the affine Grassmannian in the case of \S\ref{satake}.
In each  of the two situations, there is   another basis
which is more significant 
from the group theoretic point of view, as  it
carries representation theoretic information:
it affords a description of representations of Hecke algebra via 
the $W$-graph in the first case, and describes the weight decomposition of the representations of the Langlands dual group in the second.
 The matrix relating the natural and the group-theoretic
 bases  singles out a set of functions, the Kazhdan-Lusztig polynomials in \S\ref{hecke} and the functions of formula (\ref{sake}) in  \S\ref{satake}.
 
 The upshot is that
 in both cases it turns out that these functions are those associated,  via the 
functions-sheaves dictionary, with  the  {\em intersection complexes} of the
aforementioned subvarieties.

 In $\S$\ref{hecke} and $\S$\ref{satake}, the main role is played by
certain $G$-equivariant perverse sheaves and by the notion
of  $l$-adic purity.
The decomposition theorem allows to greatly simplify the arguments
and to clarify the overall picture.

\subsection{Schubert varieties and Kazhdan-Lusztig polynomials}
\label{hecke}
The connection between the Kazhdan-Lusztig polynomials
associated with the Weyl group of a semisimple 
linear algebraic group and the intersection cohomology groups
of the Schubert varieties of the associated flag variety played 
an important role in  
the development of the theory of perverse sheaves.
This connection was worked out by D. Kazhdan 
and G. Lusztig, (\ci{kl1},\ci{kl2}),
following discussions with R. MacPherson and P. Deligne. 

We quickly review the basic definitions in the more general framework of Coxeter groups, 
see \ci{humph} for more details on this beautiful subject.
Let $(W,S)$ be a Coxeter group, that is, a group $W$ with a set of generators $S$
which satisfy relations $(ss')^{m(ss')}=1$ with $m(s,s)=1$ and $m(s,s')\geq 2$ if $s \neq s'$. 
Any element $w \in W$ has an expression $w=s_1 \cdots s_n$ with $s_i \in S$, and the {\em length} $l(w)$ of $w \in W$ is the minimal number of $s_i$'s appearing in a such expression. 
For the definition of the Bruhat order, a partial order on $W$ compatible with lengths, see \ci{humph} 5.9.

\begin{ex}
\label{simm}
{\rm Let $W={\cal S}_{n+1}$, the symmetric group. The   the set of transpositions $s_i=(i,i+1)$
yields a set of generators 
$S=\{s_1, \ldots, s_n\}$.
}
\end{ex}

A basic object associated with $(W, S)$ is the Hecke algebra  
$\frak{H}$.  
It is a free module over the ring $ {\Bbb Z} [q^{1/2}, q^{-1/2}]$ with basis $\{T_w\}_{w \in W}$
and ring structure 
$$
T_wT_{w'}=T_{ww'}, \;  \hbox{ if } \,  l(ww')=l(w)+l(w'), 
 \qquad T_sT_w=(q-1)T_{w}+qT_{sw},\;  \hbox{ if }\, l(sw)<l(w).
$$ 

As the  following two  examples show,
Hecke algebras often arise as convolution algebras in Lie theory.
Recall that, given a locally compact topological group $G$ with Haar measure $dg$, the convolution product of two compactly supported measurable functions 
$f_1,f_2:G \lorw \zed$ is defined as 

\begin{equation}
\label{convolution}
f_1 * f_2(h)=  \int_{G}
f_1(g)f_2(g^{-1}h)dg.
\end{equation}
In the case of a finite group, the ring of $\zed$-valued functions with respect to the convolution product is thus canonically isomorphic to the group ring
$\zed[G]$.

\begin{ex}
\label{cheva}
{\rm Let $G_q$ be a Chevalley group over the finite field with $q$
  elements $\fq$, e.g. the general linear group $GL_n(\fq)$, the
  symplectic group $Sp_{2n}(\fq)$ or the orthogonal group $O_n(\fq)$.
Let $B_q \subseteq G_q$ be
  a Borel subgroup, and $W$  be the Weyl
  group. We consider functions $f: G \lorw \zed$ which are left and right $B_q$-invariant, that is, such that $f(b_1gb_2)=f(g)$ for all $b_1,b_2 \in B_q$ and $g \in G_q$. The convolution of two such functions is still left and right 
$B_q$-invariant  and the corresponding algebra  $\zed [B_q \backslash G_q /B_q]$ is generated by the  characteristic functions of the
  double  $B_q$-cosets.
In \ci{iwa}, Iwahori proved that the Bruhat decomposition
$$ G_q = \coprod_{w \in W} B_q w B_q$$ determines an algebra isomorphism 
between 
  $\zed [B_q \backslash G_q /B_q]$, and the 
  Hecke algebra $\frak{H}$ of $W$, where  
the indeterminate $q$ is 
specialized to the cardinality of the field. 
The survey \ci{curtis} gives a useful summary
  of the properties of this algebra and its relevance to the 
  representation theory of groups of Lie-type.
 }
\end{ex}

\begin{ex}
\label{heckeaffine}
{\rm Let ${\cal K}$ be a local field and ${\cal O}$ be its ring of integers,
and denote by $\pi$ a generator of the unique maximal ideal ${\frak p}$ of ${\cal O}$. 
E.g. ${\cal K}= \rat_p$,  ${\cal O}= \zed _p$, and  $\pi=p \in \zed _p$, or 
${\cal K}= \fq((T))$, the field of formal Laurent series with coefficients in a finite field,  ${\cal O}= \fq[[T]]$, and  $\pi=T$. Denote by
$q$ the cardinality of the residue field $k={\cal O}/\pi$.
Let $G$ be a simply connected reductive  group {\em split} over  ${\cal K}$, 
that is, $G$ contains a maximal torus $T$ whose set of ${\cal K}$-points is $T({\cal K})=({\cal K}^*)^r$.  
Let $W^{\rm aff}$ be its affine Weyl group, i.e. the semidirect product of $W$ with the coroot lattice of $G$, see \ci{humph}, \S 4.2. There is the ``reduction  mod-$\frak p$'' map 
$\pi:G({\cal O}) \to G(k)$. Let $B':=\pi^{-1}(B)$ be the inverse image of a Borel subgroup of $G(k)$.
For instance, if $G=SL_2$ with the usual choice of  a positive root,
and ${\cal K}=\rat_p$, then the ``Iwahori subgroup'' $B'$ consists of matrices 
in $SL_2(\zed _p)$ whose entry on the upper right corner is a multiple of $p$.
Iwahori and Matsumoto, \ci{iwamatsu} proved that the algebra 
$\zed [B' \backslash G({\cal K}) /B']$, generated by the  characteristic functions of the
  double  $B'$-cosets, endowed with the
  convolution product,
 is isomorphic to the Hecke algebra for $W^{\rm aff}$.
As in  Example \ref{cheva}, the double $B'$-cosets are parameterized, via  a  Bruhat-type decomposition, by  $W^{\rm aff}$, and the basis $T_w$ of
their characteristic functions satisfies the two defining relations of the
Hecke algebra of  $W^{\rm aff}$.  
The closely related ``spherical Hecke algebra'' 
will be discussed in $\S$\ref{satake}, in connection with the geometric Satake isomorphism.
}
\end{ex}

It follows from the second defining relation of the Hecke algebra that $T_s$ is invertible for $s \in S:$ $T_s^{-1}=q^{-1}(T_s-(q-1)T_e)$.
This implies that $T_w$ is invertible for all $w$. 

The algebra ${\frak H}$ admits two commuting involutions $\iota$ and $\sigma$, defined by 
$$
\iota(q^{1/2})=q^{-1/2}, \quad \iota(T_w)=T_{w^{-1}}^{-1},
\qquad
\hbox{ and }
\qquad
\sigma(q^{1/2})=q^{-1/2}, \quad \sigma(T_w)=(-1/q)^{l(w)}T_w.
$$

The following is proved in \ci{kl1}:
\begin{tm}
There exists a unique 
$\zed [q^{1/2}, q^{-1/2}]$-basis $\{C_w \}$ of ${\frak H}$   with 
the following properties:
\begin{equation}
\label{klp}
\iota(C_w)=C_w \qquad C_w=(-1)^{l(w)}q^{l(w)/2} \sum_{v\leq w}(-q)^{-l(v)}P_{v,w}(q^{-1})T_v
\end{equation}
with $P_{v,w} \in \zed [q]$ of degree at most $1/2(l(w)-l(v)-1)$, if $v<w$, and $P_{w,w}=1$.
\end{tm}

The polynomials $P_{v,w}$ are called the {\em Kazhdan-Lusztig polynomials of $(W,S)$}.

\begin{rmk}
\label{sn}
{\rm For $s \in S$, we have that  $C_s=q^{-1/2}(T_s-qT_e)$ satisfies (\ref{klp}), hence $P_{s,s}=P_{e,s}=1$.
A direct computation shows that
if $W={\cal S}_3$, then
$P_{v,w}=1$ for all $v,w$. In contrast, 
if $W= {\cal S}_4$, then  $P_{s_1s_3,\,s_1s_3s_2s_3s_1}=P_{s_2,\,s_2s_1s_3s_2}=1+q$.
}
\end{rmk}

\medskip

Let $G$ be a semisimple linear algebraic group, $T$ be a maximal torus,
$W=N(T)/T$ be  the Weyl group. Choose a system of simple roots.
Each simple root yields a reflection about the hyperplane associated with the root.
The set of these reflections is known  to generate $W$. Let  $B$ be  the Borel subgroup
containing $T$  and associated with the choice of the simple roots:
this means that the Lie algebra of $B$ is spanned
by the Lie algebra of $T$ and the positive roots spaces.
If $w \in W$, then  we  denote a representative of $w$  in $N(T)$  by the same letter.

The flag variety $X=G/B$ parametrizes the Borel subgroups via the map $gB \to gBg^{-1}$.
The $B$-action on $X$ gives the ``Bruhat decomposition''  by $B$-orbits $X=\coprod_{w \in W} X_w$.
The Schubert cell $X_w$ is the $B$-orbit of $wB$.
It is well known, see \ci{borelalggrp}, that  $X_w \simeq \comp^{l(w)}$ and 
$\overline{X}_w= \coprod_{v\leq w}X_v$, where $\leq$ is the Bruhat ordering.
Hence the Schubert variety $\overline{X}_w$ is endowed with a natural $B$-invariant cell decomposition. 
  
\begin{ex}
\label{aenne}
{\rm Let $G=SL_{n+1}$, $B$ be the subgroup of upper triangular matrices, $T$ be the subgroup of diagonal matrices.
Then $W\simeq {\cal S}_{n+1}$, and the choice of $B$ correspond to $S=\{s_1, \ldots s_n\}$ as in 
Example \ref{simm}. 
Clearly $X_e=\overline{X}_e$ is the  point $B$, and $\overline{X}_{w_0}=X$, if $w_0$ denotes 
the longest element of $W$. If $s \in S$ then $\overline{X}_s\simeq \pn{1}$.
If 
$\{o\}\subseteq \comp^0 \subseteq \comp^1 \subseteq \ldots \subseteq \comp^n$ is the flag determined by 
the canonical basis of $\comp^n$, then 
$\overline{X}_{s_i}$ parameterizes the flags $\{o\}\subseteq V_1\subseteq \ldots \subseteq V_{n-1} \subseteq \comp^n$
such that $V_k=\comp^k$ for all $k\neq i$. One such flag is determined by the line $V_i/V_{i-1} \subseteq 
V_{i+1}/V_{i-1}$. If $l(w)\geq 2$, then  the Schubert variety $\overline{X}_w$ is, in general, singular.
The flags ${\Bbb V} =\{o\}\subseteq V_1\subseteq \ldots \subseteq V_{n-1} \subseteq \comp^n $
in a Schubert cell $X_w$ can be described in terms of dimension of the intersections
$V_i \cap \comp^j$ as follows:
$$
X_w =\{ {\Bbb V}:  \dim{V_i \cap \comp^j}= w_{ij}\}, \;\;  \hbox{ where }
w_{ij}=\sharp\{k \leq i \hbox{ such that }  w(k)\leq j \} \}.
$$
}
\end{ex}

Since $B$ acts transitively 
on any Schubert cell, it follows that 
$ \dim {\cal H}^i(IC_{\overline{X}_w})_x$ depends only 
on the cell $X_v$ containing the point $x$.

Set, for $v\leq w$, $h^i(\overline{X}_w)_v:=\dim{\cal H}^i(IC_{\overline{X}_w})_x$
for $x$ any point in $X_v$.
Define, for $v\leq w$, 
the Poincar\'e polynomial $\widetilde{P}_{v,w}(q)= \sum_i h^{i-l(w)}(\overline{X}_w)_vq^{i/2}$.

We have the following remarkable and surprising fact, which yields a geometric
interpretation of the Kazhdan-Lusztig polynomials in terms of dimensions
of stalks of cohomology sheaves of intersection complexes of Schubert varieties.

\begin{tm}
\label{PtildequalsP}
{\rm (\ci{kl2})}
We have $P_{v,w}(q)=\widetilde{P}_{v,w}(q)$.
In particular, if $i+l(w)$ is odd, then
${\cal H}^i(IC_{\overline{X}_w})=0$, 
and the coefficients of the Kazhdan-Lusztig polynomials $P_{v,w}(q)$ are non negative.
\end{tm}

\begin{rmk}
{\rm 
Theorem \ref{PtildequalsP} implies   that $P_{v,w}=1$ for all $v\leq w$
iff $IC_{\overline{X}_w}=\rat_{\overline{X}_w}[l(w)]$. This happens, for instance, for $SL_3$
(cf. \ref{sn}). The Schubert varieties of $SL_3$ are in fact smooth.
}
\end{rmk}

\begin{rmk}
{\rm To our knowledge, there is  no purely combinatorial proof of the non negativity of the coefficients of the Kazhdan-Lusztig polynomials. This fact illustrates the power of the geometric interpretation. }
\end{rmk}

\begin{rmk}
{\rm In the same paper \ci{kl1} in which the polynomials $P_{v,w}$
are introduced, Kazhdan and Lusztig conjecture a formula, 
involving the values $P_{v,w}(1)$, for the multiplicity of a representation
in the Jordan-H\"older sequences  of Verma modules. 
The proofs of this conjecture, due independently to Beilinson-Bernstein 
(\ci{beibe1})
and Brylinski-Kashiwara (\ci{{brylikashi}}), make essential use of the geometric interpretation
given by Theorem \ref{PtildequalsP} of the Kazhdan-Lusztig polynomials. See \ci{springerbourbaki}, \S 3, for  the necessary definitions and a sketch of the proof.
}
\end{rmk}

\begin{rmk}
\label{degree}
{\rm Since ${\dim} X_v=l(v)$, the support conditions (\ref{sprt}) of $\S$ \ref{subsecintcoh} 
for intersection cohomology imply that  
 if $v<w$, then
 ${\cal H}^{i-l(w)}(IC_{\overline{X}_w})_v=0$ for $i-l(w)\geq -l(v)$. It follows that 
the degree of $\widetilde{P}_{v,w}(q)$ is a most $1/2(l(w)-l(v)-1)$, 
as required by the definition of the Kazhdan-Lusztig polynomials. Furthermore, 
as $(IC_{\overline{X}_w})_{|{X_w}}=\rat_{X_w}[l(w)]$,
we have $P_{w,w}=1$
}
\end{rmk}

The original proof of  Theorem \ref{PtildequalsP},
given in \ci{kl2}, 
is inspired to the ``functions-sheaves dictionary'' 
briefly discussed in $\S$\ref{foncfasc}, and
does not use the decomposition theorem, 
but, rather,  the purity  of the  intersection cohomology  complex
in the $l$-adic context (see $\S$\ref{8000})  and the Lefschetz Trace Formula, \ci{grbourbaki}.
As seen in  Remark \ref{degree}, the polynomials $\widetilde{P}_{v,w}$ 
satisfy 
the first property (\ref{klp}) on the degree. 
It thus remains to show the invariance under the involution  $\iota$.
Kazhdan and Lusztig directly
show that ${\cal H}^i(IC_{\overline{X}_w})=0$ if  
$i+l(w)$ is odd, and that the Frobenius map acts on 
${\cal H}^{2i-l(w)}(IC_{\overline{X}_w} )$ with eigenvalues equal to $q^i$,
so that, up to a shift,
$\widetilde{P}_{v,w}(q)=t_{IC_{\overline{X}_w}}(x)$, if $x \in X_v(\fq)$.
Once this is shown, the invariance under the involution  $\iota$ turns out 
to be equivalent to  the Poincar\'e duality theorem 
for intersection cohomology, \S\ref{subsecintcoh}. 

For another approach, again based on
the  purity of $l$-adic
intersection cohomology complex, see \ci{luszvogan}.

\medskip 
An approach to Theorem \ref{PtildequalsP}  
due to MacPherson,  gives also a topological 
description of the Hecke algebra.
It is based on exploiting the decomposition theorem for 
the Bott-Samelson variety (see \ci{bottsa,dema}), which is
a $G$-equivariant resolution of 
a variety closely related to the Schubert cell $\overline{X}_w$.

Another proof, which still relies on
applying the decomposition theorem to the 
resolutions of the Schubert varieties mentioned above,
was later worked out by T. Haines, in \ci{haines}. 
It exploits the fact that the fibers 
of the resolution have a decomposition 
as a disjoint union of affine spaces. This latter 
approach works  
with the flag variety as well as  with the 
(infinite dimensional) affine flag variety.

\subsection{The Geometric Satake isomorphism}
\label{satake}
We now discuss an  analogue of the constructions 
described in $\S$\ref{hecke}, culminating in a geometrization 
of the spherical Hecke algebra and the Satake isomorphism. 
In this case, the Schubert subvarieties 
will be replaced by certain subvarieties 
$\overline{\rm Orb}_{\lambda}$ 
of the affine Grassmannian ${\cal G}{\cal R}_{G}$.

Let us first recall, following  the clear exposition \ci{gross},
the basic statement of the classical Satake isomorphism (\ci{satake}). 

Let ${\cal K}$, ${\cal O}$, ${\frak p}$, $\pi$ and $q$ be as in 
\S\ref{heckeaffine}. 
We  let  $G$ be  a reductive linear
algebraic group split over ${\cal K}$. We denote by $G({\cal K})$  the set of ${\cal K}$-points  and  by  $K=G({\cal O})$,  the
set of ${\cal O}$-points, a compact subgroup of $G({\cal K})$. 
Similarly to Examples \ref{cheva}, \ref{heckeaffine},
the {\em spherical} Hecke algebra ${\cal H}(G({\cal K}),G({\cal O}))$ is  defined to be the set of
$K$$-$$K$-invariant locally constant $\zed $-valued functions on $G({\cal K})$ endowed
with the convolution product (\ref{convolution}) where the Haar measure is 
normalized so that the volume of $K$ is $1$.

 The group $X_{\bullet}(T):=\mbox{Hom}({\cal K}^*, T({\cal K}))$
 of co-characters of $T$ is free abelian, and carries  a natural
action of the Weyl group $W$.   
The choice of a set of positive roots
singles out a  system of positive coroots in $X_{\bullet}(T)$ as well 
as  the
positive chamber 
$$
X_{\bullet}(T)^+=\{ \lambda \in X_{\bullet}(T) \, s.t. \, 
\lambda(\alpha)\geq 0 \, if\, \alpha>0\},
$$
which is a fundamental domain for the action of $W$. 
Given $\lambda,\mu \in  X_{\bullet}(T)$, we say that 
$\lambda \geq \mu$ if $\lambda - \mu$ is a sum of positive coroots.

Every $\lambda \in X_{\bullet}(T)$ defines an element $\lambda(\pi)
\in K$, and one has the following  Cartan-type decomposition:
$$
G= \coprod_{\lambda \in X_{\bullet}(T)^+ }K\lambda(\pi)K.
$$
The characteristic functions $C_{\lambda}$ of the double
cosets $ K\lambda(\pi)K, $ for $\lambda \in
X_{\bullet}(T)^+$, give a $\zed$-basis of   ${\cal H}(G,K)$.
The spherical Hecke algebra is commutative.

\begin{rmk}
\label{torus1}
For the torus $T$, we have ${\cal H}(T({\cal K}),T({\cal O}))\simeq \zed [X_{\bullet}(T)]$.
\end{rmk}

\begin{ex}
\label{glp}
{\rm
Let $G=GL_n$. With the usual choice of positive roots, an element $\lambda \in
X_{\bullet}(T)^+$ is of the form $diag( t^{a_1}, \ldots t^{a_n})$, with
$a_1 \geq a_2 \geq \ldots \geq a_n$.
The above decomposition boils down to the fact that
a  matrix can be reduced to  diagonal form
by multiplying it on the left and on the right by elementary matrices.
}
\end{ex}

The {\em Langlands dual} $^L G$ of $G$, is the reductive group whose 
root datum is the co-root datum of $G$
and whose co-root datum is the root datum of $G$.
For a very nice description of these notions see \ci{spri0}.
The representation ring  $R(^L G)$  of  $^L G$ is isomorphic to 
$\zed [X_{\bullet}(T)]^W$. 

\begin{tm}
\label{classsat}
{\bf (The classical Satake isomorphism)}
There is an isomorphism of algebras:
\begin{equation}
\label{classatiso}
{\cal S}: {\cal H}(G({\cal K}),G({\cal O})) \otimes \zed
[q^{1/2},q^{-1/2}]  \stackrel{\simeq }{\lorw} 
R(^LG) \otimes \zed [q^{1/2},q^{-1/2}] .
\end{equation}

\end{tm}

\begin{rmk}
\label{nontrivial}
{\rm The $\zed$-module $R(^LG)$ has a basis $[V_{\lambda}]$
  parameterized by  $\lambda \in
X_{\bullet}(T)^+$, where $V_{\lambda}$ is the irreducible
representation with highest weight $\lambda$. It may be tempting
to think that the inverse ${\cal S}^{-1}$ 
 sends   $[V_{\lambda}]$   to  the characteristic function  $C_{\lambda}$ of the double
coset $ K\lambda(\pi)K$.   However, this does  not work: there exist integers
$d_{\lambda}(\mu)$, defined for $\mu \in X_{\bullet}(T)^+$, with 
$\mu < \lambda$ such that the  more complicated formula 
\begin{equation}
\label{sake}
{\cal S}^{-1}([V_{\lambda}])= q^{-\rho(\lambda)} (C_{\lambda} +
\sum_{\stackrel{\mu \in X_{\bullet}(T)^+}{\mu < \lambda }} d_{\lambda}(\mu)
C_{\mu}), 
\end{equation}
where  $\rho= (1/2) \sum_{\alpha >0} \alpha$, 
holds instead.}
\end{rmk}

The Satake isomorphism is remarkable in the sense that it relates $G$ and $^LG$.
A priori, it is unclear that the two should be related at all, beyond
the defining exchanging property. The isomorphism gives, in principle,
a recipe to construct the  (representation ring of the) Langlands dual $^LG$ of $G$
from the datum of the ring of functions on the double coset space $K\backslash G/K$.  

A striking application of the theory of perverse sheaves is the
``geometrization'' of this isomorphism. The whole subject was started
by the important work of Lusztig \ci{luszqmult, luszunipo}. 
In this work, it is shown that the Kazhdan-Lusztig polynomials
associated with a group closely related to $W^{\rm aff}$
are the Poincar\'e polynomials of the  intersection cohomology sheaves of the  singular varieties   $\overline{\rm Orb}_{\lambda}$, for $\lambda \in X_{\bullet}(T)$,  inside the affine Grassmannian ${\cal G}{\cal R}_{G}$ defined below, and coincide with the weight multiplicities $d_{\lambda}(\mu)$  of the representation $V_{\lambda}$ appearing in  formula (\ref{sake}). 
As a consequence, he showed that, if we set    $I\!H^*(\overline{\rm Orb}_{\lambda})= \oplus_l I\!H^l(\overline{\rm Orb}_{\lambda})$,  then we have 
$\dim I\!H^*(\overline{\rm Orb}_{\lambda})=\dim V_{\lambda}$ and that the tensor product operation  $V_{\lambda}\otimes V_{\nu}$ correspond to a ``convolution'' operation $IC_{\overline{\rm Orb}_{\lambda}}\star IC_{\overline{\rm Orb}_{\nu}}$.

The geometric significance of Lusztig's result was clarified by the work of 
Ginzburg \ci{ginzburg} and  Mirkovi\'c-Vilonen \ci{mivi}.

We  quickly review  the geometry involved, according to the paper \ci{mivi}.
We work over the field of complex numbers. The analogue of the coset space $G({\cal K})/G({\cal O}))$
of $\S$\ref{hecke} is the  affine Grassmannian, which we now introduce; see \ci{bd} for a thorough treatment.
Let $G$ be a reductive algebraic group over $\comp$, let $\comp[[t]]$ be the ring of formal power series and 
$\comp((t))$ its fraction field of Laurent series. 
The quotient ${\cal G}{\cal R}_{G}=G(\comp((t)))/ G(\comp[[t]])$
is called the {\em affine Grassmannian}: it is an ind-variety, i.e.
a countable  increasing union of projective varieties.

\begin{rmk}
\label{lattices}
{\rm Let  $G=SL_n(\comp)$. The points of ${\cal G}{\cal R}_{SL_n(\comp)}$
parametrize special lattices in the $\comp((t))$-vector space $V=\comp((t))^n$. 
A special lattice is  a $\comp[[t]]$-module 
$M \subseteq V$ such that $t^N\comp[[t]]^n \subseteq M \subseteq t^{-N}\comp[[t]]^n$ for some
$N$, and $\bigwedge^n M=\comp[[t]]$. The action of $SL_n(\comp((t)))$ on the set of special lattices
is transitive, and $SL_n(\comp[[t]])$ is the stabilizer
of the lattice $M=\comp[[t]]^n$. }
\end{rmk}

\begin{rmk}
\label{grtorus}
{\rm The set of points of the affine Grassmannian ${\cal G}{\cal R}_{T}$
of a torus $T$ is easily seen to be  $X_{\bullet}(T)$ (see Remark \ref{torus1}). 
The scheme structure is somewhat subtler, as it turns out to be non reduced.}
\end{rmk}

Set ${\cal K}=\comp((t))$, and ${\cal O}=\comp[[t]]$.
The imbedding $T \subseteq G$ of the maximal torus gives a map
${\cal G}{\cal R}_{T} \to {\cal G}{\cal R}_{G}$; thus,
by Remark \ref{grtorus}, we can identify $X_{\bullet}(T)$ with a subset of  
${\cal G}{\cal R}_{G}$. 
 It turns out that the group $G({\cal O})$ acts on ${\cal G}{\cal R}_{G}$ 
with finite dimensional orbits.
We  still denote by  $\lambda$ the point of the affine Grassmannian corresponding to  $\lambda \in  X_{\bullet}(T)$, and denote
its $G({\cal O})$-orbit
by   ${\rm Orb}_{\lambda} \subseteq {\cal G}{\cal R}_{G}$  
(cf.  \ci{bd}).

\begin{pr}
\label{afforbits}
{\rm (\ci{bd}, 5.3)}
There is a decomposition 
${\cal G}{\cal R}_{G}= \coprod_{\lambda \in X_{\bullet}(T)^+ }{\rm Orb}_{\lambda}$.
Furthermore, every orbit ${\rm Orb}_{\lambda}$ has the structure of a 
vector bundle over a rational homogeneous variety,
it is connected 
and simply connected,  
$$
\dim{\rm Orb_{\lambda}}=2\rho(\lambda) \qquad \hbox{ and } \qquad 
\overline{\rm Orb}_{\lambda}=\coprod_{\mu \leq \lambda}{\rm Orb}_{\mu}. 
$$
\end{pr}

Proposition  \ref{afforbits} implies that the category 
${\cal P}_{G({\cal O})}$ of perverse sheaves which are constructible
with respect to the decomposition in $G({\cal O})$-orbits is generated
by the intersection cohomology complexes $IC_{\overline{\rm Orb}_{\lambda}}$. 
Lusztig has proved in  \ci{luszqmult} that the cohomology sheaves 
${\cal H}^i(IC_{\overline{\rm Orb}_{\lambda}})$ are different from zero only in one parity. 
Together with the fact that the dimensions of all $G({\cal O})$-orbits
in the same connected component of   
${\cal G}{\cal R}_{G}$ have the same parity, this implies  that ${\cal P}_{G({\cal O})}$ is a semisimple category. Its objects are automatically 
$G({\cal O})$-equivariant perverse sheaves. 

The Tannakian formalism, see \ci{demi}, singles out the categories which are 
equivalent to categories of representations of affine groups schemes 
and it gives a precise prescription
for reconstructing the group scheme from its category of representations.  
The geometrization of the Satake isomorphism essentially states that the category 
${\cal P}_{G({\cal O})}$ is equivalent to the category of representations ${\rm Repr}(^LG)$ of 
the Langlands dual group $^LG$, so it
yields a recipe to re-construct this dual group.
More precisely, it is necessary   to endow ${\cal P}_{G({\cal O})} $ with the structure of rigid tensor category with a ``fiber functor.'' 
Essentially, this means that there must be
1) a 
bilinear functor $
\star: {\cal P}_{G({\cal O})}\times {\cal P}_{G({\cal O})} \to {\cal P}_{G({\cal O})}$
with compatible associativity and commutativity constraints,  i.e.  functorial isomorphisms
$ A_1 \star (A_2 \star A_3) \stackrel{\simeq}{\to} (A_1 \star A_2) \star A_3$
and    $ A_1 \star A_2 \stackrel{\simeq}{\to} A_2 \star A_1$, 
and  2) an exact functor,  called the fiber functor, $F :{\cal P}_{G({\cal O})}  \to {\rm Vect}_{\rat}$ which is a tensor functor,
i.e. there is a functorial isomorphism
$F(A_1 \star A_2) \stackrel{\simeq}{\to}F(A_1)\otimes F(A_2)$.

\begin{rmk}
\label{tannaka}
{\rm For the category of representations of a group, the product is given by the tensor product of representations, while the fiber functor $\omega$ associates with a representation its underlying vector space.}
\end{rmk}

In fact, there exists a geometrically defined ``convolution product''
$$
\star\; : \;\, {\cal P}_{G({\cal O})}\times {\cal P}_{G({\cal O})} \lorw {\cal P}_{G({\cal O})}
$$
with ``associativity and commutativity constraints,'' such that
the  cohomology functor  $H: {\cal P}_{G({\cal O})}  \to {\rm Vect}_{\rat}$
is a tensor functor. The construction  of this geometric convolution product 
is reviewed below, see \ref{convo}.

 We  state the geometric Satake isomorphism  as follows:

\begin{tm}
\label{gsi}

{\rm {\bf (Geometric Satake isomorphism)}}There is an equivalence of 
tensor categories
\begin{equation}
\label{gsif}
{\cal S}_{geom}:({\cal P}_{G({\cal O})}, \star, H) \stackrel{\simeq}{\lorw} ({\rm Rep}(^LG), \otimes, \omega).
\end{equation}

\end{tm}

\begin{rmk}
\label{real}{\rm
Nadler investigated (\ci{nadler}) a subcategory of perverse sheaves on the affine Grassmannian of a real form $G_{\real}$ of $G$ and proved that it is equivalent to
the category of representations
of a reductive subgroup $^LH$ of $^LG$. This establishes a real
version of the Geometric Satake isomorphism and, as a corollary,
the decomposition theorem is shown to hold for several {\em real} algebraic maps arising in Lie theory. 
}
\end{rmk}

We discuss only two main points of the construction of \ci{mivi}, 
the definition of the convolution product 
and the use of the `` semi-infinite'' orbits to construct the weight functors. 
We omit all technical details and  refer the reader to \ci{mivi}. 

\medskip
{\bf The convolution product.}
In the following description of the convolution 
product we  treat the spaces involved  as if they were honest varieties. See 
\ci{gaits}  for a detailed account. 
Let us consider the diagram:
$$
\xymatrix{
{\cal G}{\cal R}_{G} & G({\cal K}) \ar[r]^{\pi} & {\cal G}{\cal R}_{G} &      \\
G({\cal K}) \times_{G({\cal O})} {\cal G}{\cal R}_{G} \ar[u]^p&  G({\cal K}) \times {\cal G}{\cal R}_{G} \ar[u]\ar[l]_(.4)q \ar[r]^{\pi \times {\rm Id}} & {\cal G}{\cal R}_{G}\times{\cal G}{\cal R}_{G}\ar[u]_{p_1} \ar[r]^(.6){p_2} & {\cal G}{\cal R}_{G}.
}
$$
The map $q: G({\cal K}) \times {\cal G}{\cal R}_{G} \to G({\cal K}) \times_{G({\cal O})} {\cal G}{\cal R}_{G}$ is the quotient map by the action of $G({\cal O})$, the map $p:G({\cal K}) \times_{G({\cal O})} {\cal G}{\cal R}_{G} \to  {\cal G}{\cal R}_{G}$ is the ``action'' map, 
$p(g, hG({\cal O}))=ghG({\cal O})$.
If $A_1,A_2 \in {\cal P}_{G({\cal O})}$, then $(\pi \times {\rm Id})^*(p_1^*(A_1) \otimes p_2^*(A_2))$ on $G({\cal K}) \times {\cal G}{\cal R}_{G}$
descends to  $G({\cal K}) \times_{G({\cal O})} {\cal G}{\cal R}_{G}$, that is, there exists a unique complex of sheaves $A_1 \widetilde{\otimes}A_2$ on $G({\cal K}) \times_{G({\cal O})} {\cal G}{\cal R}_{G}$ with the property that 
$(\pi \times {\rm Id})^*(p_1^*(A_1) \otimes p_2^*(A_2))=q^*(A_1 \widetilde{\otimes}A_2)$, and we set
 $A_1 \star A_2 := p_*(A_1 \widetilde{\otimes}A_2)$.

The following fact is referred to as ``Miraculous theorem'' in \ci{bd}:
\begin{tm}
\label{convo}
If $A_1,A_2 \in {\cal P}_{G({\cal O})}$, then $A_1 \star A_2 \in {\cal P}_{G({\cal O})}$.
\end{tm}

The key reason why this theorem holds is that the map $p $ enjoys a strong form of semismallness.

First of all the  complex $A_1 \widetilde{\otimes}A_2$ is constructible with respect to the decomposition 
$$G({\cal K}) \times_{G({\cal O})} {\cal G}{\cal R}_{G}= \coprod {\cal S}_{\lambda, \mu}
\qquad  \hbox{ with }
\qquad  {\cal S}_{\lambda, \mu}=\pi^{-1}({\rm Orb}_{\lambda})\times_{G({\cal O})}{\rm Orb}_{\mu}.
$$

\begin{pr}
\label{strss}
The map $p:G({\cal K}) \times_{G({\cal O})} {\cal G}{\cal R}_{G} \to 
{\cal G}{\cal R}_{G}$ is stratified semismall, in the sense, that for any 
${\cal S}_{\lambda, \mu}$, the map $p_{|\overline{\cal S}_{\lambda, \mu}}:
\overline{\cal S}_{\lambda, \mu} \to p(\overline{\cal S}_{\lambda, \mu})$
is semismall. As a consequence $p_*$ sends perverse sheaves 
constructible with respect to the decomposition $\{{\cal S}_{\lambda, \mu}\}$,
to perverse sheaves on ${\cal G}{\cal R}_{G}$ constructible with respect to 
the decomposition $\{ {\rm Orb}_{\lambda}  \}$.
\end{pr}

\begin{rmk}
\label{costraint}
{\rm 
While the ``associativity constraints'' of the convolution product are almost immediate from its definition, 
the commutativity constraints are far subtler (see  \ci{mivi}, and  also \ci{gaits}).
}
\end{rmk}

{\bf The weight functor.} The cohomology functor 
$H(-):=\oplus_l H^l(-)$  is a fiber functor for the category
${\cal P}_{G({\cal O})}$. In particular,  it is a tensor functor:
$H(A_1 \star A_2 )\simeq H(A_1)\otimes H(A_2)$. In order to verify this, Mirkovi\'c and Vilonen decompose this functor
as a direct sum of functors $H_{\mu}$ parameterized by $\mu \in X_{\bullet}(T)$. This decomposition is meant to mirror 
the weight decomposition of a representation of $^LG$. It is realized  by introducing certain 
ind-subvarieties $N_{\mu}$ which have a ``cellular'' property with respect to any $A \in {\cal P}_{G({\cal O})}$, 
in the sense that at most one compactly supported cohomology 
group does not vanish.
Let $U$ be the unipotent radical of the Borel group $B$,
and $U({\cal K})$ be  the corresponding subgroup of $G({\cal K})$.
The $U({\cal K})$-orbits in the affine Grassmannian are
neither of finite dimension nor of finite codimension. 
It can be shown that they are parameterized by $X_{\bullet}(T)$. If, as before, 
we still denote by  $\nu$ the point of the affine Grassmannian corresponding to a cocharacter $\nu \in  X_{\bullet}(T)$, and  set  $S_{\nu}:=U({\cal K})\nu$, then we have ${\cal G}{\cal R}_{G}=
\coprod_{\nu \in X_{\bullet}(T)}S_{\nu}$.

\begin{pr}
\label{cells}
For any $A \in {\cal P}_{G({\cal O})}$, we have 
$
H^l_c(S_{\nu}, A)=0 \hbox{ for } l \neq 2\rho(\nu)
$
and
$H^{2\rho(\nu)}_c(S_{\nu}, IC_{\overline{\rm Orb}_{\lambda}})$
is canonically isomorphic to the vector space generated by the 
irreducible components of $\overline{\rm Orb}_{\lambda}\cap S_{\nu}$.
In particular, the functor  
$ H^{2\rho(\nu)}_c(S_{\nu}, \, - \, ):{\cal P}_{G({\cal O})}  \to {\rm Vect}_{\rat}$ sending     $A \in {\cal P}_{G({\cal O})}$ to $H^{2\rho(\nu)}_c(S_{\nu}, A )$ is exact, and 
$$H( A   ) \;\; :=\bigoplus_{l\in \zed}H( {\cal G}{\cal R}_{G},  A   )=
\bigoplus_{\nu \in X_{\bullet}(T)}H^{2\rho(\nu)}_c(S_{\nu}, A ).$$
\end{pr}

\begin{rmk}
\label{cmnts}
{\rm 
Let $A \in {\cal P}_{G({\cal O})}$.
Since in the equivalence of categories of Theorem \ref{gsi},
the fiber functors $H$ corresponds to $\omega$, the decomposition 
\[H( A ) \;\; =
\bigoplus_{\nu \in X_{\bullet}(T)}H^{2\rho(\nu)}_c(S_{\nu},A)\]
of Proposition  \ref{cells} of the cohomology of $A$
must reflect a decomposition of the underlying vector 
space of the representation ${\cal S}_{geom}(A)$. 
In fact, this is  the weight decomposition 
of the corresponding representation of $^LG$.

An aspect of  the Geometric Satake correspondence which we find particularly beautiful is that, up to a re-normalization, the intersection cohomology complex
$IC_{\overline{\rm Orb}_{\lambda}}$ correspond, via the Geometric 
Satake isomorphism, to the irreducible representation $V(\lambda)$ of $^LG$ with highest weight $\lambda$. This explains  (see Remark \ref{nontrivial})
why the class of $V(\lambda)$ is not easily expressed in terms of the characteristic function $C_{\lambda}$ of the double coset $ K\lambda(\pi)K$, which 
corresponds, in the function-sheaves dictionary of \ref{foncfasc}, to the constant sheaf on $\overline{\rm Orb}_{\lambda}$,
and once again emphasizes the fundamental nature of intersection cohomology. 
Furthermore, in view of   Proposition \ref{cells},
the irreducible components of $\overline{\rm Orb}_{\lambda}\cap S_{\nu}$
as $\nu$ varies in  $X_{\bullet}(T)$, give a canonical basis for $V(\lambda)$.
These components are now called {\em Mirkovi\'c-Vilonen cycles}. 

}
\end{rmk}

\medskip
The classical Satake isomorphism \ref{classsat} for ${\cal K}= \fq((T))$
may be recovered from the geometric Satake isomorphism \ref{gsi} by considering the Grothendieck group of the two tensor categories. In fact,
the Grothendieck ring  of the  category ${\rm Repr}(^LG)$ 
is the representation ring  $R(^L G)$, while the functions-sheaves dictionary identifies the Grothendieck ring of ${\cal P}_{G({\cal O})}$ with the spherical Hecke algebra ${\cal H}(G({\cal K}),G({\cal O}))$.

\subsection{Ng\^o's support theorem}
\label{nst}

We thank G. Laumon and B.C.  Ng\^o for very useful conversations.
The paper \ci{ngo}  is devoted to the proof
of the fundamental lemma in the Langland's program,  a long-standing 
and deep conjecture
concerning Lie groups. 
For its complexity, depth and wealth of   applications to 
representation theory, this  paper
deserves a separate treatment, which  we do not provide here. 
In this section,  instead, we give a brief and rough discussion of
 B.C. Ng\^o  {\em support
theorem}  (\ci{ngo}, Theorem 7.1.13). This result,
which we state in a slightly weaker form in Theorem \ref{ngostm},
 can be stated and proved without any reference 
to the context of the Langlands program,  and  is of great independent
geometrical interest.
Under the favourable assumptions which are explained in the sections that follow, 
it gives
a precise characterization of the supports of the perverse sheaves  
which enter the decomposition theorem for a map $f: M \to S$
acted upon in a fiber-preserving manner by a family of
commutative algebraic  groups $g: P\to S$. 
This seems to be 
one of the first cases in which the decomposition theorem is studied
in depth in the context of a non generically finite map, i.e.  of a map
with large fibers.
For expository reasons, we state these results over the complex numbers, even though the main use in \ci{ngo} is  in the $l$-adic context over a finite field. 

The determination of  the simple summands 
$IC_{\overline{Y_a}}(L_{a})$ 
 appearing in the decomposition theorem   (\ref{00dt})  is
a difficult problem. The determination of the supports
$\overline{Y_{a}}$ does not seem to be  easier. 
In fact, consider Examples \ref{cone}, \ref{smallres}: the vertex of the cone
is certainly   a special locus in both cases,   however, it appears
as a support of a summand in the decomposition theorem only
in Example \ref{cone}.

One important ingredient of Ngo's proof of the support theorem
is the following result of Goresky and MacPherson which,
in the case of equidimensional maps,  yields
an a priori constraint on the codimension of subvarieties
supporting simple summands in the decomposition theorem.
The proof  is a simple and elegant application of 
the symmetry (\ref{pvddteq}) arising from Poincar\'e-Verdier duality,
  and can be found in \ci{ngo}, Appendice A, Th\'eor\`emes 2 and 3.

\begin{tm}
\label{spgoma}
Let $f: X \to Y$ be a proper map of algebraic varieties, with $X$ nonsingular. 
Assume that   all the fibers of $f$  have   the same dimension $d$.
Assume that $Z\subseteq Y$ is an irreducible subvariety 
which is the support of a non zero summand appearing in the decompositon
theorem  {\rm (\ref{00dt})}. Then
\[ {\rm codim} \,({Z}) \leq d.\]
If, in addition, the fibers are irreducible, then one has strict inequality
in what above.
\end{tm}

The basic idea in the proof is that a bigger codimension, coupled
with duality, would force the corresponding summand to contribute
a nontrivial summand to the direct image sheaf
$R^jf_* \rat_X$ for $j>2d$, contradicting the fact that the fibers have dimension
$d$.

 In order to  state Ng\^o's support
theorem  let us fix some notation. 

Let $f: M \to S$   be a proper and flat map
of relative dimension $d$ with reduced fibers and where  $M$ and $S$  are smooth
irreducible  varieties. 
The map $f$ is assumed to be endowed with an action of a 
commutative group scheme $g:P \to S$ of relative dimension $d$.
A group scheme is a map $g:P \to S$ together with $S$-maps
 $e: S \to P$, $m: P\times_S P \to P$  and $\iota: P \to P$ that
satisfy the usual axioms of a group. Each fiber $g^{-1}(s)$ is 
an algebraic group,
and a group scheme can be seen as a family of  groups.
In this context, an action is an  $S$-map $a: P \times_S M \to M$  
 commuting with the projections
to $S$ that satisfies the usual requirements  of an action, suitably modified
to the ``relative to $S$" situation.

Let $g:P \to S$  be as above and   {\em with connected fibres}, and let $s \in S$. By a classical result of Chevalley, 
there is a canonical exact sequence $1 \to R_{s} \to P_{s} \to A_{s} \to 1$
of algebraic groups, where $R_s$ is affine (thus a product  $\comp^{\alpha_s} \times
{\comp^*}^{\mu_s}$ of  additive and multiplicative groups), 
and   $A_s$ is an abelian variety.
There is the  well-defined function $\delta: S \to {\Bbb N}$,  $s \mapsto \delta_s:= \dim R_{s}= \alpha_s + \mu_s$. 
This function is upper-semicontinuous, i.e. it jumps up on Zariski-closed subsets. In particular,
there is a  partition of $S= \coprod_{\delta \geq 0} S_{\delta}$ into locally closed subvarieties where the 
invariant $\delta$ is constant. We assume furthermore that {\em $P$ acts with affine stabilizers}: for any $m \in M$, the isotropy subgroup of $m$ is an affine subgroup of $P_{f(m)}$.

We need the two notions of $\delta$-regularity and of polarizability.
An  $S$-group scheme $g:P\to S$ as above  is {\em $\delta$-regular}
if \[\mbox{codim}_S (S_{\delta}) \geq \delta.\] 
The Tate sheaf is the sheaf 
\[ T(P) := R^{2d-1}g_! \rat\]
whose stalk  at $s \in S$ is, by base change,  the homology group $H_1(P_s)$.
We say that $T_{\rat}(P)$ is {\em polarizable} if there is an alternating bilinear pairing
\[ T(P) \otimes T(P) \lorw \rat_S\]
that factors through $H_1(A_s)$ at every point  and induces on it a perfect pairing.

We can now state the following
\begin{tm}
\label{ngostm} {\rm ({\bf Ng\^o support theorem})}
Let $f: M\to S$, $g: P \to S$ be as above. Assume that $P\to S$ is $\delta$-regular
and that $T(P)$ is polarizable.
A closed irreducible subvariety $Z \subseteq S$ is the support
of a nontrivial simple summand appearing in the decomposition theorem
for $f: M \to S$ if and only if there is a Zariski dense open subvariety
$Z^0 \subseteq Z$ such that the sheaf $R^{2d} f_*\rat$  is locally constant on $Z^0$
and $Z$ is maximal with respect to this property.
\end{tm}  

\begin{rmk}
\label{qep}
{\rm
Theorem \ref{ngostm} is applied to the case when $f: M \to S$
is a suitable open subset of  the Hitchin fibration associated with a Lie group.
The hypothesis of $\delta$-regularity is verified with the aid of Riemann-Roch
and of the deformation theory of Higgs bundles. 
Over the complex numbers an infinitesimal argument shows that a group scheme (variety) associated with a completely integrable algebraic system is always $\delta$-regular. The hypothesis of polarizability
is verified using the classical Weil pairing. See \ci{ngo}. 
}
\end{rmk}

The statement of the support theorem,  is remarkable because it tells us where to look
for the supports of the summands of the decomposition theorem: they are those varieties
closures of  (maximal) parts of $S$ over which the single sheaf $R^{2d}f_* \rat$ is a locally constant. On the other hand, since the fibers are assumed to be reduced,
the sheaf $R^{2d}f_*\rat$ is the linearization of the sheaf of finite sets given by the irreducible components of the fibers of $f$. This fact makes the determination
of these supports an approachable problem. 
For example, suppose that the map $f$ has irreducible fibres; denote by 
$j:S_{reg} \to S$ the imbedding of the open set of regular values 
of $f$ and by $R^i$ the local systems $R^if_*\rat$ on $S_{reg}$. 
Then  $f_*\rat[\dim M]=\oplus_i IC(R^i)[d-i]$, that is, there are no summands beyond those which are determined by the ``fibration part''.

\subsection{Decomposition up to  homological cobordism and  signature}
\label{csmax}
We want to mention, without any detail, a purely topological counterpart
of the decomposition theorem. Recall that this result holds only  in the
algebraic context, 
e.g. it fails for proper holomorphic maps of complex manifolds.

In the topological context, Cappell and Shaneson \cite{CS1} 
introduce a notion of                        
cobordism for complexes of sheaves and prove  a general topological result
for maps between  Whitney                               
stratified space with only even  codimension                                   
strata   that    in  the case of a proper algebraic map $f: X \to Y$,  
identifies, up to cobordism, $f_*IC_X$ with $\pc{0}{f_* IC_X}$ and 
its splitting as in the decomposition
theorem.  For a related question, see 
 \ci{gomaprob}, D., Problem 6.
                                      
The decomposition up to cobordism is sufficient                      
to provide exact formulae                                                                                              
for many topological invariants, such as
Goresky-MacPherson $L$-classes and                         
signature thus generalizing the classical Chern-Hirzebruch-Serre multiplicativity property                       
of  the signature for  smooth fiber bundles with no monodromy to the case
of stratified maps (see \cite{CS2,CS3,S}).                                     
                                                                                                    
In the case of complex algebraic varieties,   one may also  look at the MacPherson Chern                       
classes \cite{M}, the Baum-Fulton-MacPherson Todd classes \cite{BFM},  the                        
homology Hirzebruch classes \cite{BSY,CMS2} and their                          
associated Hodge-genera defined in terms of the mixed Hodge structures on the                       
(intersection) cohomology groups. The papers  \cite{CS3, CMS1,CMS2}
provide Hodge-theoretic applications of the above topological stratified                                 
multiplicative formul\ae.  For a survey, see  \cite{MS1}. 
                     
These results yield   topological and analytic constraints
on the singularities of complex algebraic  maps.
  In the case of maps of  projective varieties,                       
these Hodge-theoretic formul\ae $\,$ are proved using the decomposition theorem,
especially the  identification in \cite{decmightam}  of                     
the local systems appearing in the decomposition  combined with the                         
Hodge-theoretic aspects of the decomposition theorem in                     
\cite{decmigseattle}. For non-compact varieties, the authors use the functorial calculus                     
on the Grothendieck groups of Saito's algebraic mixed Hodge modules.                

\subsection{Further developments and applications}
\label{furthapp}
{\bf  Toric varieties and polytopes.} There exist polytopes that are not combinatorially equivalent to any
rational polytope, and the formula for the generalized $h$-polynomial makes sense also in this case, even though 
 there is no toric variety associated with it. It is thus natural to ask whether the properties 
of the $h$-polynomial reflecting the Poincar\'e duality and the hard Lefschetz theorem hold more generally for any polytope. 

In order to study this sort of  questions, P. Bressler and V. Lunts  have developed a theory of sheaves on the poset associated with the 
polytope $P$, or more generally to a fan, see \ci{brelunt1}.
Passing to the corresponding derived category, they  define an intersection cohomology complex and prove 
the analogue of the decomposition theorem for it, as well as the equivariant version. 

By building on their foundational work, K. Karu,   proved   in  \ci{karu} 
that the hard Lefschetz property and  the Hodge-Riemann  relations hold for
every, i.e.  not necessarily rational, polytope. Different proofs, each one shedding new light on interesting combinatorial phenomena, 
have then been given by Bressler-Lunts in \ci{brelunt2} and by Barthel-Brasselet-Fieseler-Kaup in  \ci{bbfk}. 
Another example of application of methods of intersection cohomology to the
combinatorics of polytopes is the solution, due to T. Braden and
R. MacPherson of a conjecture of 
G. Kalai concerning the behavior of the $g$-polynomial of a face with respect to the $g$-polynomial of the whole polytope.
See \ci{bramac} and the survey \ci{braden}.

\medskip

{\bf The Hilbert scheme of points on a surface.}
Vafa and Witten noticed in \ci{vawi} that G\"ottsche's Formula (\ref{gottfor}) suggests a representation theoretic structure underlying
the direct sum $ \oplus_{i,n} H^i(S^{[n]})$. Namely, this space should be an irreducible highest weight module over the infinite dimensional
Heisenberg-Clifford super Lie algebra, with highest weight vector the generator of 
$H^0(S^{[0]})$. 
H. Nakajima and, independently  I. Grojnowski   took up the suggestion in \ci{nakahilb,groj} (see also the
lecture notes  \ci{nakajilectures}) and realized this structure by
a set of correspondences relating Hilbert schemes of different lengths. 

An elementary proof of G\"ottsche's formula stemming form this circle of ideas was given in \ci{decmighilb1}. The papers \ci{decmighilb2, decmigsemi} prove, in two different ways,
 a motivic version of the decomposition theorem  (\ref{dthilbert}) for the map $\pi : S^{[n]} \to S^{(n)} $
exhibiting an equality  
$$
(S^{[n]},\Delta, 2n)=\sum_{\nu \in {\frak P}_n}(S^{l(\nu)}, P_{\nu},2l(\nu))
$$
of Chow motives with rational coefficients. In this formula, $P_{\nu}$ denotes the projector associated with the action of the group
$\prod {\cal S}_{a_i}$ on $S^{l(\nu)}$.
Two related examples, still admitting a semismall contraction, are the
nested Hilbert scheme $S^{[n,n+1]}$, whose
 points are couples $(Z,Z') \in S^{[n]} \times S^{[n+1]}$ such that $Z \subseteq Z'$, and the parabolic Hilbert scheme, see \ci{decmigsemi} and its Appendix for details.

\medskip
{\bf The Geometric Satake isomorphism.}
The decomposition theorem, applied to the stratified semismall map $p$
used to define the convolution, gives a decomposition
$$
IC_{\overline{\rm Orb}_{\lambda}}\star IC_{\overline{\rm Orb}_{\mu}}= 
\oplus_{\nu}IC_{\overline{\rm Orb}_{\nu}}\otimes F_{\nu},
$$
where $F_{\nu}$ is the vector space generated by the relevant 
irreducible components of the fibres of $p$ (see (\ref{decomssls}) and
note that the strata are simply connected).
This decomposition mirrors, on the geometric side, the Clebsch-Gordan decomposition $V(\lambda)\otimes V(\mu)= \oplus_{\nu} V(\nu)\otimes F_{\nu}.$
The irreducible components of the fibres were shown to be Mirkovi\'c-Vilonen 
cycles in \ci{and}. A combinatorial study of Mirkovi\'c-Vilonen cycles
is made possible by letting the maximal torus $T$ act on them. 
The action is hamiltonian and its image by the  moment map is a polytope. 
The  so obtained
{\em Mirkovi\'c-Vilonen polytopes} are 
investigated in, for instance, \ci{andkog, kamn}.

\medskip
{\bf Other applications. }
The examples discussed in this section are far from exhausting the
range of applications of the theory of perverse sheaves.
We suggest G. Lusztig's \ci{luszicm},
T.A. Springer's \ci{springerbourbaki}, and N. Chriss and V. Ginzburg's 
\ci{chrissginzburg} for further  applications and
for  more details, including motivation and references, 
 about some of the examples discussed here in connection
with representation theory.

For lack of space and competence, we have  not
discuss  many important  examples, such as the proof of the Kazhdan-Lusztig conjectures
and  the applications
of the geometric Fourier transform. 

The most dramatic occurrence of the
functions-sheaves dictionary, and one of the reasons for 
the importance of perverse sheaves in representations theory, 
is the geometrization of the notion of {\em automorphic form} 
in the geometric Langlands program; for details, see for instance \ci{frenk}, 
\S 3.3, 
or \ci{gai}.
Coarsely speaking, an (unramified) automorphic form is a 
function on the ``adelic quotient'' 
$GL_n(F)\backslash GL_n({\Bbb A_F})/GL_n({\cal O})$,
where $F$ is the field of rational functions of an algebraic 
curve $X$ defined over a finite field $\fq$, $ {\Bbb A_F}$ is 
the ring of ad\`eles of $F$, and ${\cal O}=\prod_{x\in X} {\cal O}_x$.
The function must also satisfy some other property, such as that of being 
an eigenvector for the unramified Hecke algebra.
A theorem of A. Weil gives an interpretation of the adelic quotient as
the set of points of the moduli stack of vector bundles on $X$. 
Hence, by the function-sheaves dictionary, an automorphic form
should correspond to a perverse sheaf on this  moduli stack, 
and the important condition that the automorphic form be a Hecke eigenvector 
can also be interpreted geometrically introducing the notion of a Hecke 
eigensheaf.

\section{Appendices}
\label{cemetto}
\subsection{Algebraic varieties}
\label{algvar}
The precise definitions of varieties and maps in algebraic geometry are quite lengthy.
Luckily,  in order to understand 
the statement of the decomposition theorem, as well as some of its applications, 
it is often sufficient to deal
with quasi projective varieties and the maps between them.
Let us explain a little  bit this terminology,  without being too formal.
A {\em projective}  variety is an
algebraic  variety that admits an embedding
in  some projective space
$\pn{N}$ as  the zero set
 of finitely many
homogeneous equations in $N+1$ variables.
A {\em quasi projective} variety
is an algebraic  variety that admits an embedding in some projective space
 as the difference set  of two projective varieties.
There are algebraic varieties that are not quasi projective.
 An {\em  affine}
 variety is a variety that can be viewed as the zero set  in 
 some affine space ${\Bbb A}^N$ of finitely
 many polynomial in $N$ variables. An affine variety is clearly quasi projective;
 the converse does not hold, e.g. ${\Bbb A}^{2} \setminus \{ (0,0)\}$.
 There are the notions of {\em subvariety} and {\em product varieties}.
 A map of algebraic varieties $f: X \to Y$, or simply a {\em map},  is a map
 of the underlying sets  whose graph
is an algebraic  subvariety of the product variety $X\times Y$.  A complex algebraic variety carries two interesting topologies:
the {\em Euclidean} (or classical) topology and the  coarser
{\em Zariski} topology.
Let us discuss these two topologies  in the case of a quasi projective variety
embedded in a projective space, $X \subseteq \pn{N}$:
the Euclidean topology is the topology induced on $X$
by the complex manifold topology on $\pn{N}$; the Zariski topology is the topology
with closed sets given by zero sets on $X$   of finitely many homogeneous polynomial
in $N+1$ variables.  A {\em closed} ({\em open}, resp.) subvariety is
a closed  (open, resp) subset for the Zariski topology. 
A map of algebraic varieties
is {\em proper} (in the sense of algebraic geometry) if it is separated and universally closed for the Zariski topology;
luckily, this happens  if and only the map is proper for the Euclidean topology.
In particular, a map of projective varieties is always proper. An algebraic variety
$X$  is {\em reducible} if it is  the union $X= X' \cup X''$ of two
closed algebraic subvarieties with $X', X'' \neq X$, and it is
{\em irreducible} otherwise.

\subsection{Hard Lefschetz  and mixed Hodge structures}
\label{pam}
We   want to state the hard Lefschetz theorem and  the Hodge-Riemann
bilinear relations  in the language
of Hodge structures. Let us   recall briefly this formalism.

\medskip
{\bf Hodge structures and polarizations.}

Let $l \in \zed$, $H$ be a finitely generated abelian group, 
$H_{\rat}:= H 
\otimes_{\zed} \rat$, $H_{\real}= H \otimes_{\zed}\real$,
$H_{\comp}= H \otimes_{\zed} \comp$.
A {\em pure Hodge structure  of weight $l$} on $H$,  $H_{\rat}$ or
$H_{\real}$,
is a direct sum decomposition $H_{\comp}= \oplus_{p+q=l} H^{p,q}$
such that $H^{p,q}= \overline{ H^{q,p} }$.
The Hodge  filtration  is the decreasing filtration 
$F^{p}(H_{\comp}):=
\oplus_{p'\geq p}H^{p',q'}$. A {\em morphism of Hodge structures}
$f: H \to H'$ is a group homomorphism such that
the complexification of $f$ (still denoted $f$)
 is compatible with the Hodge filtrations in the sense
 that $f(F^p H_\comp) \subseteq F^p H'_\comp$,
i.e. such that it is a filtered map. Such maps are automatically 
what one calls {\em strict}, i.e. $(\im\, f )\cap F^p H'_{\comp}= f(F^p H_\comp)$.
The category of Hodge structures of weight $l$ with the above arrows is Abelian.

Let $C$ be the {\em Weil operator}, i.e.
the $\real$-linear map  $C: H_{\comp} \simeq H_{\comp}$
 such that $C(x)=i^{p-q} x$, for every $x\in H^{pq}$.
 Replacing $i^{p-q}$ by 
$z^{p}{\overline{z}}^{q}$
 we get a real  action $\rho$ of $\comp^{*}$ on $H_{\comp}$.
A {\em polarization} of the real pure Hodge structure
$H_{\real}$ is a  real bilinear form $\Psi$ on $ H_{\real}$
which is  invariant under the action given by $\rho$ restricted to 
$S^{1} \subseteq \comp^*$
and such that the bilinear form
$\widetilde{\Psi} (x,y):= \Psi(x, Cy)$ is symmetric and positive definite.
If $\Psi$ 
is a polarization, then $\Psi$ is symmetric
if $l$ is even, and antisymmetric if $l$ is odd.
In any case, $ \Psi$ is nondegenerate.
In addition, for every $0 \neq x \in H^{pq}$, 
$(-1)^{l}i^{p-q}\Psi(x,\overline{x}) >0$,
where $\Psi$ also denotes 
the $\comp$-bilinear extension of $\Psi$ to $H_{\comp}$.

Let $\eta$ be the first Chern class of an ample line bundle
on the projective $n$-fold $Y$.
For every $r \geq 0$, define the space of {\em primitive vectors} $P^{n-r} := \ke{ \,\eta^{r+1} }
\subseteq H^{n-r}(Y, \rat)$.

Classical Hodge theory
states that, for every $l$,  $H^{l}(Y, \zed)$ is a pure Hodge 
structure
of weight $l$,  $P^{n-r}$ is a rational pure Hodge structure 
of weight $(n-r)$ polarized, up to a precise sign,  by the 
 bilinear form defined on  $H^{n-r}(Y)$ as follows (it is well-defined by Stokes' theorem):
\begin{equation}
\label{modpif}
S^{\eta} (\alpha, \beta) := \int_Y{ \eta^r \wedge \alpha \wedge \beta}.
\end{equation}
The fact that this form is nondegenerate is equivalent to the celebrated hard Lefschetz Theorem.
Its signature properties are expressed by the
Hodge-Riemann bilinear relations.

\begin{tm}
\label{chl}
Let $Y$  be a complex projective manifold of dimension $n$. Then the following
statements hold.

\begin{enumerate}
\item
{\rm  ({\bf Hard Lefschetz theorem})}
For every $r \geq 0$ the cup product with $\eta$ yields   isomorphisms
\[
\eta^{r}\,: \, H^{n-r}(Y, \rat) \, \simeq \, H^{n+r}(Y, \rat).
\]
\item
{\rm ({\bf Primitive Lefschetz decomposition})}
For every $0 \leq r \leq n$ there is the direct sum decomposition
\[
H^{n-r}(Y, \rat)\, = \, \bigoplus_{j \geq 0} \eta^{j}P^{n-r-2j}
\]
where each summand
is a  pure Hodge sub-structure of weight $n-r$   and
all summands are mutually orthogonal with respect to the
bilinear form $S^{\eta}$.
\item
{\rm ({\bf Hodge-Riemann bilinear relations})}
For every $0 \leq r  \leq  n$, the bilinear form
$(-1)^{\frac{(n-r)(n-r+1)}{2}}S^{\eta}$
is a polarization of the pure weight $l$ Hodge structure
$P^{n-r} \subseteq H^{n-r}(Y, \real)$. In particular,
\begin{equation}
\label{hrbr}
(-1)^{\frac{(n-r) (n-r-1)}{2}} i^{p-q} \int_{Y}{\eta^{r}\wedge
\alpha \wedge \overline{\alpha} } \; >\; 0, \qquad \forall \;\, 0 \neq \alpha
\in P^{n-r} \cap H^{p,q}(Y, \comp).
\end{equation}\end{enumerate}\end{tm}

\medskip
{\bf Inductive approach to hard Lefschetz.}

Our proof (discussed in $\S$\ref{dmapp}) of the decomposition theorem
requires that we first  establish the relative hard Lefschetz theorem. 
We do so  by using an approach   similar to the classical  inductive approach to the   hard Lefschetz
Theorem \ref{chl}.1. There are two variants 
of this inductive approach
(see   \ci{weil2, decbook}), the former is via Hodge-Riemann relations, 
the latter is via the semisimplicity of the monodromy action in a Lefschetz pencil.
Though both   are relevant to our approach
to the decomposition theorem, we limit ourselves
to discussing the former variant.

The induction is on $n:= \dim{Y}$ and uses a nonsingular
hyperplane section $D \subseteq Y$. 
The case $r =0$ is trivial. The cases   $ r \geq 2$ 
 follow by an easy induction on the dimension of $Y$ using the Lefschetz
 hyperplane theorem. One is left with the key case $r=1$.
The cup product map  $\eta:= c_1(D) \wedge - $
factors as $\eta = g\circ r$:
$$
\xymatrix{
H^{n -1}(Y)  \ar[r]^{r} &  H^{n-1}(D)   \ar[r]^{g} &  H^{n +1} (Y),}
$$
where  $r$ is the injective restriction map and $g$ is the surjective Gysin map.
It is easy to show that  $\eta$ is an isomorphism iff the intersection 
form on $D$, restricted to $\im{ (r) }$, is nondegenerate.
While the form on $H^{n-1}(D)$ is non degenerate by Poincar\'e duality,
there is no a priori reason why it should restrict to a non degenerate form
on $H^{n-1}(Y)$.
This is where the  Hodge-Riemann relations enter the picture:
by contradiction, assume that
 there is a non-zero class $\alpha \in \ke\,  \eta$; which we may suppose of pure Hodge type; then $r(\alpha)$ is primitive in 
$H^{n-1}(D)$, and, by the Hodge-Riemann relations on $D$,
$0=\int_Y \eta \wedge \alpha \wedge \overline{\alpha}=\int_D \alpha \wedge \overline{\alpha}\neq 0$, contradiction. 

The  Lefschetz theorem on hyperplane sections coupled with  the Hodge-Riemann  bilinear relations for a hyperplane section
imply the hard Lefschetz theorem for $Y$. However, they 
do not  imply the Hodge-Riemann bilinear relations
for the critical middle dimensional cohomology group
$H^n(Y)$, and the induction procedure grinds to a halt.

To make the proof work, one has to somehow  establish 
the Hodge-Riemann relations on $H^n(Y)$.    $\S$\ref{outlinedecmig},
sections 1. and 2b, outline 
two instances of how   Hodge-Riemann-type  relations
can be established.

The hard Lefschetz theorem applied to the fibers of a smooth 
projective morphism and Theorem \ref{delsplit} imply the following
result
(cf. item 3., following Theorem \ref{stdteo}). For the  proof see
 \ci{dess} and  \ci{ho2}, Th\'eor\`eme 4.2.6.

\begin{tm}
\label{dss}
{\rm ({\bf Decomposition, semisimplicity and relative hard Lefschetz for proper smooth maps})}
Let $f: X^{n} \to Y^{m}$ be a smooth proper 
map of smooth algebraic
varieties of the indicated dimensions.
Then 
$$
 \qquad f_{*}\rat_{X} \, \simeq \,
\bigoplus_{j\geq 0} R^{j}f_{*}\rat_{X} [-j]
$$
and the $R^jf_* \rat_X$ are semisimple local systems.
If, in addition, $f$ is projective and 
$\eta$ is the first Chern class of  an $f$-ample
line bundle on $X$, 
then we have isomorphisms
$$
\eta^{r}\, : \, R^{n-m-r} f_{*} \rat_{X}\,  \simeq\,  R^{n-m+r} f_{*}
\rat_{X}, \; \forall r \geq 0,
$$
and the  local systems $R^{j}f_{*}\rat_{X}$  
underlie polarizable variations of pure Hodge structures.
\end{tm}

\medskip
{\bf Mixed Hodge structures}.

 In general, 
the singular cohomology groups $H^j(Y, \zed)$ of a singular variety cannot carry 
a pure Hodge structure of weight $j;$ e.g. $H^1(\comp^*, \zed)$ has rank one, and 
pure Hodge structures of odd weight have even rank. However, these groups underlie
a more subtle structure, the presence of which makes the topology
of complex algebraic varieties even more remarkable.

\begin{tm}
\label{mHssc}
{\rm  ({\bf Mixed Hodge structure on cohomology})}
Let $Y$ be an algebraic variety. For every $j\geq 0$ there is an increasing
filtration (the weight filtration)
$$
\{0 \} \, =\,  W_{-1} \, \subseteq W_{0} \, \subseteq \, \ldots 
\, \subseteq W_{2j} \, = \, H^{j}(Y, \rat)
$$
and a decreasing filtration (the Hodge filtration)
$$
H^{j}(Y, \comp) \, = \, F^{0} \, \supseteq \, F^{1} \, \supseteq 
\ldots \, \supseteq\,  F^{m}\, \supseteq \, F^{m+1} \, = \, \{0 \}
$$
such that the filtrations induced by $F^{\bullet}$ on the 
complexified graded pieces
of the weight filtration $W_{\bullet}$, endow 
every graded piece $W_{l} / W_{l-1}$ with a rational pure Hodge structure of 
weight $l$.
This structure (mixed Hodge structure) is functorial for maps of algebraic
varieties and the induced maps strictly preserve
both filtrations.
\end{tm}

\subsection{The formalism of  the constructible derived category}
\label{tfidy}
Standard references for what follows are \ci{k-s, gel-man, borel, iv, bbd};
see also \ci{schbook}.
In what follows, we freely refer to our crash-course in $\S$\ref{subsec-crash-course}
and to the complete references given above.

A {\em full} subcategory ${\cal C}' \subseteq {\cal C}$ of a category ${\cal C}$
is a subcategory such that the induced map on the Hom-sets
is bijective; in other words we keep all the arrows.

An {\em additive} category ${\cal C}$ is one in which  each
$\mbox{Hom} (A,B)$ is an Abelian group, composition of arrows is bilinear, the direct sum
$A\oplus B$ is defined for any pair of objects $A,B \in {\cal C}$,
and the zero object
$0 \in {\cal C}$ exists. 
A {\em complex}  $K $ in an additive category ${\cal C}$ is a sequence
\[
\xymatrix{ \ldots \ar[r] & K^{i-1} \ar[r]^{d^{i-1}}&
 K^i \ar[r]^{d^i} & K^{i+1}\ar[r]^{d^{i+1}}&  \ldots} \] of objects and morphisms in ${\cal C}$
such that for every $i\in \zed$ we have $d^{i}\circ d^{i-1} =0$.
The objects $K^i$ are called the {\em entries} of the complex and the
arrows $d^i$ are called its {\em differentials}.
One often omits the indexing of the arrows.  A {\em map of complexes}
$f: A \to B$ is a collection of arrows $f^i: A^i \to B^i$  such that
$d\circ f = f \circ d$. 
Complexes in ${\cal C}$ form an additive category, denoted $C({\cal C})$.
Given a complex $K$ and $m \in \zed$, the {\em$m$-shifted} complex $K[m]$
is the complex with entries  $(K[m])^i:= K^{i+m}$ and with  differentials
$d_{K[m]}^i = (-1)^k d^{i+k}_K$.
The {\em cone} of a map of complexes
$f: A \to B$ is the complex $Cone(f)$ where $Cone(f)^i:= B^i \oplus A^{i+1}$
and the the differential is defined by setting $d(b,a) = (d(b) + f(a), -d(a))$.

An {\em Abelian} category is an additive category
where every arrow  admits a kernel and a cokernel and, 
given any arrow $f: A \to B$, the resulting natural arrow
$\coke{\, \{\ke{f} \to A \}} \to \ke{\, \{ B \to \coke{f}\} }$ is an isomorphism.
 
In  this paragraph, we work in a fixed Abelian category ${\cal A}$.
An arrow
 is {\em monic} if its  kernel is (isomorphic to) zero. If an arrow
 $A \to B$ is monic, then we say that $A$ is a {\em subobject} of $B$.
 An object $A \in {\cal A}$ is {\em simple} if it has no 
 non trivial subobjects.  The Abelian category
 ${\cal A}$ is {\em Artinian} if, for every object $A \in {\cal A}$,  every 
 descending chain 
 of subobjects  of $A$ stabilizes. If ${\cal A}$ is Artinian, then every
non zero  object $A$  is a finite iterated extension of nonzero  simple objects, 
 called
 the {\em constituents} of $A$; the constituents of $A$
 are well-defined up to isomorphism.
The Abelian category
 ${\cal A}$ is {\em Noetherian} if, for every object $ A$,  every 
 ascending chain 
 of subobjects of $A$  stabilizes.  
 The category of complexes $C({\cal A})$ is  Abelian. Given a complex $K$ in $
 C({\cal A})$
 and an integer $i\in \zed$, 
 we define the cohomology object $H^i(K): = \ke{d^i}/\im{d^{i-1}} \in {\cal A}$
 and the truncated complexes $\td{i} K$ and $\tu{i}K$ as follows:
$$(\td{i}K)^l=K^l,  \;\;\; l <i,  \qquad (\td{i}K)^i=\ke{\,d^i}, \qquad
(\td{i}K)^l=0,  \; \;\; l>i,$$
with the obvious differentials,
and
$$(\tu{i}K)^l=0,  \;\;\;l <i, \qquad  (\td{i}K)^i=\coke{\,d^{i-1}}, \qquad
(\td{i}K)^l=K^l, \;\;\; l>i,
$$
with the obvious differentials. For every $i\in \zed$ there are  short exact sequences in  
the Abelian category $C({\cal A})$:
$$
\xymatrix{
 0 \ar[r] &   \td{i} K \ar[r] &  K \ar[r] &  \tu{i+1}K \ar[r]  &  0,} 
$$
and natural identifications of functors
\[
\td{i}\circ \tu{i} \, \simeq\, \tu{i}\circ  \td{i}\,  \simeq \,  [-i] \circ  H^{i}.
\]
Given an arrow $f: A \to B$ in $C({\cal A})$, we get a short   exact sequence in $C ({\cal A})$:
$$
\xymatrix{
0 \ar[r] &  B \ar[r] &  Cone(f:A \to B) \ar[r] &  A[1]  \ar[r] &  0.}
$$

Let $Y$ be an algebraic variety and ${\cal D}_Y$ be the constructible
bounded derived category ($\S$\ref{subsec-crash-course}).
The category ${\cal D}_Y$  is a {\em triangulated category}.
In particular, it is additive, so that we can form finite direct sums,
and it is equipped with the translation functor $A \mapsto A[1]$.
A 
{\em triangle} is a diagram
of maps $A \to B \to C \to A[1]$  in ${\cal D}_Y$.
A most important feature of triangulated categories is the presence of  distinguished
triangles.  
Given a map of complexes
$f: A' \to B'$, there is the short exact sequence
of complexes $0 \to B' \to Cone (f) \to A'[1] \to 0$. This exact sequence
gives rise to a triangle $A' \to B' \to Cone (f) \to A'[1]$ in ${\cal D}_Y$.
A {\em distinguished triangle} is a triangle  which is isomorphic in ${\cal D}_Y$
to the one associated with a map $f$ as above.
Any map  $f: A \to B$ in ${\cal D}_Y$,   can be completed to a
distinguished  triangle.
One should keep in mind that this construction
is not functorial;  see \ci{gel-man}).

It is easy to show that the kernel of a morphism $f: A \to B$ in 
 ${\cal D}_Y$ splits off as a direct summand   of $A$.
 Since there are complexes which do not split
 non trivially, {\em the category ${\cal D}_Y$ is not Abelian} (unless $Y$ is a
 finite collection of points). 
 
 The cone construction is a replacement
in the non Abelian category ${\cal D}_Y$ of the notions of kernel and cokernel.
In fact, if $f: A \to B$ is an injective (surjective, resp.)  map of complexes, then $Cone (f)$
is isomorphic  in ${\cal D}_Y$ to the cokernel  ($1$-shifted kernel, resp.) complex.

An essential computational tool is that the application
of a cohomological functor to a distinguished triangle produces
a long exact sequence.  
Distinguished triangles are a replacement for 
short exact sequence
in the non Abelian category ${\cal D}_Y$.
A {\em cohomological functor},
with values in an abelian category ${\cal A}$, is an additive functor
$T: {\cal D}_Y \to {\cal A}$  such that $T(A) \to T(B) \to T(C)$ is exact
for every  distinguished triangle as above. Setting $T^i (A) : = T (A[i])$,  we get the long  exact
 sequence 
\[ 
\xymatrix{\cdots \ar[r] & T^i(A) \ar[r] &
 T^i(B) \ar[r] &  T^i(C) \ar[r] &  T^{i+1} (A) \ar[r] & \cdots }\]

Using injective resolutions and 
the two  global sections functors $\Gamma$
and $\Gamma_c$ we define the {\em derived global sections functors}
 (see \ci{gel-man,k-s} for the identification with categorical derived functors)
\[ R\Gamma, R\Gamma_c  \;: \;{\cal D}_Y \lorw {\cal D}_{pt}, \]
and the 
finite dimensional cohomology vector spaces of $Y$ with coefficients in $K
\in {\cal D}_Y$:
\[H^*(Y, K)  := H^* (R\Gamma  (Y,K)), \qquad H_c^*(Y,K) := H^*(R\Gamma_c (Y,K))
  \; \mbox{(compact supports)}. \]

 Given a map $f: X \to Y$ there are the {\em four functors}
 \[ Rf_*, Rf_! \; :\;  {\cal D}_X \lorw {\cal D}_{Y}, 
 \qquad 
 f^*, f^!  \; :\;  {\cal D}_Y \lorw {\cal D}_{X}. \]
 The  sheaf-theoretic direct image functors
 $f_*, f_!: Sh_X \to Sh_Y$  are {\em left exact} as functors, e.g.
 if  $0 \to F \to G \to H \to 0$ is an exact sequence of sheaves on $X$,
 then $0 \to f_! F \to f_! G \to f_! H$ is  an exact sequence
 of sheaves on $Y$.
 The {\em right derived} functors $Rf_*$ and $Rf_!$
arise by applying the sheaf-theoretic direct image
 functors $f_*$ and $f_!$ (proper supports), term-by-term, to injective resolutions.
 Taking  cohomology sheaves, we obtain the $i$-th right derived functors
 $R^i f_*$ and $R^if_!$. 
 We have equalities of sheaf-theoretic functors $R^0 f_* =f_*,\,R^0f_! =f_!$.
The {\em inverse image} functor $f^*:Sh_Y \to Sh_X$ is exact on sheaves and descends
to the derived category. The {\em exceptional inverse image functor} $f^!$
does not arise from a functor defined on sheaves.

It is customary to employ the following simplified  notation to denote the
four functors $(f^*, f_*, f_!, f^!)$. In this paper,  $f_*$
and $f_!$ denote the right derived functors.
 To avoid confusion, the sheaf-theoretic functors
are denoted $R^0 f_*, \, R^0f_!$. 

Given maps $f: X\to Y$, $g: Y \to Z$,  we have $(g \circ f)^! = f^! \circ g^!$, etc.
For  $g: Y\to pt$ and  for $C \in {\cal D}_X$, we have canonical isomorphisms
\[ H^*(X, C) \simeq  H^* (Y, f_*C), \qquad  H^*_c(X, C) \simeq H^* (Y, f_!C) .\]

{\bf The functors $f_{!}$ and   $f^{!}$ in special cases}

\smallskip
\n
If $f$ is proper, e.g. a closed immersion, then $f_{!}=f_{*}$.

\n
 If $f$ is smooth of relative dimension $d$, then $f^{!}=f^{*}[2d]$.

\n

A closed  embedding    $f:X \to Y$  is   {\em
 normally nonsingular   of pure codimension
$d$}  (\ci{goma2}) if it can be realized as the intersection
$X= Y \cap N$ inside $M$, where $N,M$ are nonsingular, $N$ has codimension
$d$ in $M$ and  $N$ is transverse to every stratum of some stratification
$\Sigma$ of $Y$. In this case, we have that 
 $f^{!} = f^{*}[-2d]$ holds for every $\Sigma$-constructible complex.
Such so-called normally nonsingular inclusions can be obtained
by embedding $Y$ in some projective space and then
intersecting $Y$ with  $d$ general hypersurfaces.

\n
 If $f$ is an open embedding, then $f^{!}= f^{*}$. 

\n
 If $f$ is a locally closed embedding, then

\n
1) $f_{!}$ is the extension-by-zero functor and $f_{!}= R^{0}f_{!};$

\n
2)  $f^{!}=f^{*}R {\Gamma_{X}}$,
where $\Gamma_{X}F$,
not to be confused with the sheaf  $f_{!}f^{*}F= F_X$ that is zero outside
$X$ and coincides with $F$ on $X$,
is the sheaf of sections of the sheaf $F$ supported on $X$
(see \ci{k-s}, p.95).
If, in addition,  $f$ is a closed embedding, then $
H^*(Y, f_! f^!K)= H^*(X, f^{!}K)= H^*(Y, Y \setminus X;K)=
H^*_{X}(Y, K)$.

The usual Hom complex construction  can de derived
and we get  right derived functors
\[ \mbox{RHom} : {\cal D}_Y^{opp} \times {\cal D}_Y \lorw {\cal D}_{pt}, 
\quad
R{\cal H}om : {\cal D}_Y^{opp} \times {\cal D}_Y \lorw {\cal D}_{Y} \]
with the associated $\mbox{Ext}^i$ and ${\cal E}xt^i$ functors.
We have
\[ \mbox{Hom}_{{\cal D}_Y} (K,K') = H^0(Y, 
R{\cal H}om (K, K')) = H^0(Y, 
\mbox{RHom} (K, K') ).\]

The pair $(f^*, f_*)$ is an {\em adjoint pair} (this holds also for the
sheaf-theoretic version) and so is $(f_!, f^!)$ and we have, for every $C\in
{\cal D}_X$ and $K \in {\cal D}_Y$: 
\[
f_*  R{\cal H}om (f^* K, C) =    R{\cal H}om ( K, f_*C), \qquad
f_*  R{\cal H}om (f_! C, K) =    R{\cal H}om ( C, f^!K).
\]

Since we are working with field coefficients, the tensor product operation
$\otimes$
on complexes
is exact and there is no need to derive  it.  For  $K_i \in {\cal D}_Y$, we have
(also for $R{\cal H}om$):
\[ 
\mbox{RHom} (K_1 \otimes K_2, K_3) = \mbox{RHom} (K_1, R{\cal H}om (K_2,K_3))\]
and, if the sheaves ${\cal H}^i(K_3)$ are locally constant:
\[ R{\cal H}om (K_1, K_2 \otimes K_3) = R{\cal H}om (K_1, K_2) \otimes K_3. \]

There is the {\em dualizing complex} $\omega_Y \in {\cal D}_Y$, well-defined, 
up to canonical isomorphism by setting $\omega_Y := \gamma^! \rat_{pt}$, 
where $\gamma: Y \to pt$.
If $Y$ is nonsingular, then $\omega_Y \simeq \rat_Y [2\dim_{\comp}Y]$.
Given $f: X \to Y$, we have $\omega_X = 
f^! \omega_Y$. Define a contravariant functor
\[ D \; : \; {\cal D}_Y \lorw {\cal D}_Y, \qquad
K \longmapsto  D(K) \, (= K^{\vee})\,  :=\, R{\cal H}om (K, \omega_Y).\]
We have $D^2 = \mbox{Id}$,  $(K[i])^{\vee} = K^{\vee} [-i]$ and 
$\omega_Y = \rat_Y^\vee$.
The complex $K^{\vee}$ is called the  (Verdier) {\em  dual} of $K$. 
Poincar\'e-Verdier duality consists of the canonical isomorphism
\[ H^i(Y, K^{\vee}) \simeq  H^{-i}_c (Y, K)^{\vee} \]
which is a formal consequence of the fact that $(f_!, f^!)$ are an adjoint pair.
The usual Poincar\'e duality for topological manifolds is the special case when 
$Y$ is smooth and orientable, for then a choice of orientation
gives a natural isomorphism $\omega_Y \simeq \rat_{Y} [\dim_{\real} Y]$.

We have the important relations
\[ D f_! = f_* D, \qquad Df^! = f^* D.\]

A {\em $t$-structure} on a triangulated category ${\cal D}$
is the data of two full subcategories ${\cal D}^{\leq 0}, {\cal D}^{\geq 0}
\subseteq {\cal D}$ subject to the following three requirements:
\begin{enumerate}
\item
for every $C  \in {\cal D}^{\leq 0}$ and $C' \in {\cal D}^{\geq 1}$, we have
$\mbox{Hom}_{\cal D} (C,C') =0$;
\item $D^{\leq 0} [1] \subseteq D^{\leq 0}$ and
$D^{\geq 0}  \subseteq D^{\geq 0}[1]$;
\item for every  $C \in {\cal D}$, there is a distinguished triangle
$C' \to C \to C'' \to C' [1]$ with $C' \in D^{\leq 0}$ and $C'' \in D^{\geq 1}$.
\end{enumerate}

A {\em $t$-category} is a triangulated category 
endowed with a $t$-structure (\ci{bbd, k-s}).  The heart of a
$t$-structure is the full subcategory
${\cal C}:= {\cal D}^{\leq 0} \cap {\cal D}^{\geq 0}$. The heart
  of a $t$-structure is an Abelian category.
By virtue of axiom 1., the  distinguished triangle  in 3. is defined up to canonical isomorphism
and this defines functors, called the truncation functors
$\td{0}: {\cal D} \to {\cal D}^{\leq 0}$, $C \mapsto C'=: \td{0} C$ and
$\tu{0}: {\cal D} \to {\cal D}^{\geq 0}$, $C \mapsto (C[-1])''[1]=: \tu{0} C$.
The   functor $H^0:= \td{0} \circ \tu{0}: {\cal D} \to {\cal C}$
is cohomological.  

The prototype of a $t$-structure is 
the  {\em standard $t$-structure}
on ${\cal D}_Y$  which is defined by setting
    ${\cal D}_Y^{\leq 0} \subseteq {\cal D}_Y$ to be   the full subcategory
of complexes $K \in {\cal D}_Y$ with ${\cal H}^j(K) =0$ for $j > 0$,
and 
${\cal D}_Y^{\geq 0} \subseteq {\cal D}_Y$ to be   the full subcategory
of complexes $K \in {\cal D}_Y$ with ${\cal H}^j(K) =0$ for $j < 0$.
The three axioms are easily verified.
The truncation functors are the usual ones.
The intersection  ${\cal D}_Y^{\leq 0}
\cap {\cal D}_Y^{\geq 0}$ is the Abelian category of constructible sheaves on $Y$.
The two-sided truncation 
$ \td{0}\circ  \tu{0}$  is the usual functor ${\cal H}^0$ ($0^{th}$-cohomology sheaf).

Another important $t$-structure is 
the  (middle)  {\em perverse} $t$-structure ($\S$\ref{subsecoervsh}).

There are the following notions of  exactness.
A functor of Abelian categories is {\em exact} if it preserves exact
sequences. There are the companion notions of {\em left and right exactness}.
A {\em functor of triangulated categories} (i.e. additive an commuting with translations)
 is {\em exact} if it preserves distinguished  triangles.
 A {\em functor of $t$-categories} $F: {\cal D} \to {\cal D}'$ is a functor
 of the underlying triangulated categories. 
 It is {\em exact} if it preserves distinguished  triangles.
 It is {\em left $t$-exact}
 if $F : {\cal D}^{\geq 0} \to {\cal D'}^{\geq 0}$. It is   {\em right $t$-exact} if
 $F : {\cal D}^{\leq 0} \to {\cal D'}^{\leq 0}$.
 It is {\em $t$-exact} if it is both left and right $t$-exact, in which 
 case it preserves the  Abelian hearts, i.e. it induces an exact
 functor   $F: {\cal C} \to {\cal C}'$ of Abelian categories.

\medskip
{\bf Perverse  $t$-exactness.}

\n
Let $f:X \to Y$ be a map of varieties.
 If $\dim{f^{-1}y}\leq d$, then 
$$
f_{!}, f^{*} :     \, ^{\frak p}\! {\cal D}^{\leq 0}_{Y}
\lorw \, ^{\frak p}\! {\cal D}^{\leq d}_{Y}, \qquad
f^{!}, f_{*} :     \, ^{\frak p}\! {\cal D}^{\geq 0}_{Y}
\lorw \, ^{\frak p}\! {\cal D}^{\geq -d}_{Y}.
$$
If $f$ is quasi finite ($=$ finite fibers), then $d=0$ above. 

\n
 If $f$ is affine,  e.g. the embedding of the complement of  a Cartier divisor,  
the embedding of an affine open subset, or  the projection of the complement
of a universal hyperplane section, etc., then 
$$
f_{*} : \, ^{\frak p}\! {\cal D}^{\leq 0}_{Y} \lorw \, ^{\frak p}\! {\cal D}^{\leq 0}_{Y}
\quad \mbox{(right $t$-exact)},
\qquad
f_{!}: \, ^{\frak p}\! {\cal D}^{\geq 0}_{Y}
\lorw \, ^{\frak p}\! {\cal D}^{\geq 0}_{Y} \quad
\mbox{(left $t$-exact)}.
$$
More generally, if locally over $Y$, $X$ is the union of $d+1$ affine open sets, then
$$
f_{*} : \, ^{\frak p}\! {\cal D}^{\leq 0}_{Y} \lorw \, ^{\frak p}\! {\cal D}^{\leq d}_{Y},
\qquad
f_{!}: \, ^{\frak p}\! {\cal D}^{\geq 0}_{Y}
\lorw \, ^{\frak p}\! {\cal D}^{\geq -d}_{Y}.
$$
 If $f$ is quasi finite and affine, then $f_{!}$ and $f_{*}$ are $t$-exact.

\n
 If $f$ is finite ($=$ proper and finite fibers), then $f_{!}= f_{*}$ are $t$-exact. 

\n
 If $f$ is a closed embedding, then $f_! =f_{*}$
are  $t$-exact and fully faithful. 
In this case it  is customary to drop $f_{*}$ from the notation, e.g.
$IC_X \in {\cal D}_Y$. 

\n
 If $f$ is smooth of relative dimension $d$, then $f^{!}[-d]= f^{*}[d]$ are $t$-exact.

\n
In particular, if $f$ is \'etale, then $f^{!}= f^{*}$ are $t$-exact.

\n
 If $f$ is a normally nonsingular inclusion of codimension $d$ with respect
to a stratification $\Sigma$ of $Y$, then 
$f^{!}[d]= f^{*}[-d]: {\cal D}^{\Sigma}_{Y} \to {\cal D}_X$ are $t$-exact.

The following splitting criterion (\ci{dess,shockwave})
plays an important role in the proof of the decomposition theorem:

\begin{tm}
\label{delsplit}
Let $K \in {\cal D}_X$ and $\eta: K \to K[2]$
such that $\eta^l : \, ^{\frak p}\!{\cal H}^{-l}(K) \to
\, ^{\frak p}\!{\cal H}^{l}(K)$ is an isomorphism for all $l$.
Then there is an isomorphism in ${\cal D}_Y$:
$$
K \simeq \bigoplus_i \,     ^{\frak p}\!{\cal H}^{i}(K)[-i].
$$

\end{tm}

 \subsection{Familiar objects from algebraic topology}
 \label{famfa}
 Here is a brief list  of some of the basic objects of algebraic topology
and a short discussion of how they relate to the formalism in
 ${\cal D}_Y$.

\medskip
{\bf (Co)homology etc.:}

\n
 singular cohomology:  $H^l(Y, \rat_Y);$ 

\n
 singular cohomology with compact supports: $H^l_c(Y, \rat_Y);$
 
 \n
 singular homology $H_l(Y, \rat) = H^{-l}_c(Y, \omega_Y);$

\n
 Borel-Moore homology: $H^{BM}_l(Y, \rat)= H^{-l}(Y, \omega_Y);$

\n
relative (co)homology:
if  $i: Z \to Y$ is a locally closed embedding 
and $ j: (Y\setminus Z) \to Y$, then  we have  canonical isomorphisms
$H^{l}(Y,Z, \rat)\simeq H^{l}(Y,    {i}_{!} i^{!} \rat   )$ and
$H_{l}(Y,Z, \rat)\simeq H^{-l}_{c}(Y,    {j}_{*} j^{*} \omega_{Y}   )$.

\medskip
{\bf Intersection (co)homology}. The intersection homology
groups $I\!H_{j}(Y)$ of  an  $n$-dimensional  irreducible variety
$Y$  are defined as the $j$-th homology groups
of chain complexes   of geometric chains with closed supports   subject to certain admissibility
conditions (\ci{goma1})  Similarly, one defines intersection homology
with compact supports. There are  natural maps
$$
I\!H_{j}(Y) \lorw H_{j}^{BM}(Y), \qquad \qquad
I\!H_{c,j}(Y) \lorw H_{j}(Y).
$$
 intersection cohomology: $I\!H^{j}(Y):= I\!H_{2n-j}(Y)= H^{-n+j}(Y, IC_{Y})$.

\n
 intersection cohomology with compact supports:
 $I\!H^{j}_{c}(Y):= I\!H_{c, 2n-j}(Y)= H^{-n+j}_{c}(Y, IC_{Y})$.

\medskip
{\bf Duality and  pairings.} Verdier duality implies we have canonical
identifications
\[
H_l(Y, \rat)^{\vee} = H^{-l}_c(Y, \omega_Y)^{\vee} \simeq 
H^l(Y, \rat), \qquad
H^{BM}_l(Y, \rat)^{\vee} = H^{-l}(Y, \omega_Y)^{\vee} \simeq 
H^l_c(Y, \rat).
\]
If $Y$ is nonsingular
of dimension $n$, then we have  the 
 {\em Poincar\'e duality} isomorphisms:
 \[H^{n+l} (Y,\rat)  \simeq H^{BM}_{n-l}(Y, \rat), \qquad H_{n+l}(Y, \rat) \simeq H^{n-l}_c
(Y, \rat).
\]
There are  two ways to express
the   classical  nondegenerate  Poincar\'e 
intersection pairing 
$$H^{n+l}  (Y, \rat) \times H^{n-l}_c(Y, \rat) \lorw\rat, 
\qquad \quad
H_{n-l}^{BM}(Y, \rat)  \times  H_{n+l}  (Y, \rat)   \lorw \rat. 
$$
While the former one is given by wedge product and integration,
the latter can be described geometrically as the intersection form
 in $Y$  as follows.
Given a Borel-Moore cycle  and a usual, i.e. compact,
 cycle in complementary dimensions, one changes one of them, say the first one,
to one homologous to it, but transverse to the other. 
Since the ordinary one has compact supports,
the intersection set is finite and one gets a finite intersection index.

\n
 Let $Y$ be compact, $Z$ be a closed subvariety such that 
$Y \setminus Z $ is a smooth and of pure dimension $n$. We have
{\em Lefschetz Duality}
\[
H_q(Y, Z ; \rat) = H^{-q}_c(Y, j_{*} j^{*} \omega_Y)
= H^{-q}(Y, j_{*}  j^{*} \omega_Y)=
H^{-q}(Y \setminus Z,   \rat_{Y}[2n]) =
H^{2n-q}(Y \setminus Z, \rat).
\]
{\em Goresky-MacPherson's Poincar\'e duality}: since 
$IC_Y \simeq IC_Y^{\vee}$, we have canonical isomorphisms
\[
I\!H^{n+l} (Y, \rat) \simeq  I\!H^{n-l}_{c}(Y, \rat)^{\vee}. 
\]

{\bf Functoriality.}
The usual maps in (co)homology associated with a map $f:X \to Y$  arise 
from the adjunction  maps 
$$\rat_Y \lorw f_{*} f^{*} \rat_Y = f_{*} \rat_X, 
\qquad
f_! f^! \omega_Y = f_! \omega_X \lorw \omega_Y.
$$
by taking cohomology.
In general, for an arbitrary map $f$,  there are no maps associated with Borel-Moore and 
cohomology with compact supports.
If $f$ is proper, then $f_{*}=f_!$ and one gets pull-back for proper maps in
  cohomology  with compact supports  and push-forward for proper maps
  in Borel-Moore homology. These maps are  dual to each other.

 \n
  If $f$ is an open immersion, then $f^{*} = f^!$ and one has the
  restriction to an open subset map  for  Borel-Moore homology and the push-forward
  for an open subset map for cohomology with compact supports.  These maps
  are dual to each other.

\medskip
{\bf Cup and Cap products.}
The natural isomorphisms 
$H^{l}(Y,\rat)\simeq {\rm Hom}_{{\cal D}_Y}(\rat_Y,\rat_Y[l] )$
and 
${\rm Hom}_{{\cal D}_Y}(\rat_Y,\rat_Y[l] )\simeq {\rm Hom}_{{\cal D}_Y }
(\rat_Y[k],\rat_Y[k+l] )$ identify the {\em cup product} 
$$\cup: H^{l}(Y,\rat) \times H^{k}(Y,\rat) \to H^{k+l}(Y,\rat)$$
with the composition 
$${\rm Hom}_{{\cal D}_Y}(\rat_Y,\rat_Y[l] ) \times {\rm Hom}_{{\cal D}_Y}(\rat_Y[l],\rat_Y[k+l] ) \lorw {\rm Hom}_{{\cal D}_Y}(\rat_Y,\rat_Y[k+l] ). $$
Similarly, the {\em  cap product} 
$$
\cap: H_k^{BM}(Y, \rat) \times H^l(Y,Y\setminus Z,\rat) 
\lorw H_{k-l}^{BM}(Z, \rat)
$$
relative to a closed imbedding $i:Z \to Y$ 
is obtained as a composition of maps in the derived category as follows:
$$
\xymatrix{
     &   H^l(Y,Y\setminus Z,\rat) ={\rm Hom}_{{\cal D}_Z}(\rat_Z,i^!\rat_Y[l] ) \ar@{}[d]^{\times}   \\
 H_k^{BM}(Y, \rat)={\rm Hom}_{{\cal D}_Y}(\rat_Y,\omega_Y[-k] )\ar[r] &{\rm Hom}_{{\cal D}_Z}(i^!\rat_Y,i^!\omega_Y[-k] )=   {\rm Hom}_{{\cal D}_Z}(i^!\rat_Y,\omega_Z[-k])    \ar[d] \\
&   H_{k-l}^{BM}(Z, \rat)={\rm Hom}_{{\cal D}_Z}(\rat_Z,\omega_Z[l-k]).
}
$$

 \medskip
{\bf Gysin Map.}
Let $i: Z \to Y$ be the closed embedding of a codimension $d$
complex  submanifold of 
the complex manifold $Y$.
We have $i_{*}= i_{!}$ and  $i^{!} = i^{*} [-2d], $ the adjunction map for $i_!$ yields
$$i_{*} \rat_Z  = i_!  i^{*} \rat_Y = i_! i^! \rat_Y [2d] \lorw \rat_Y [2d ]$$
and by taking cohomology we get the Gysin map
$$
H^l(Z, \rat) \lorw H^{l+2d} (Y, \rat).
$$
Geometrically, this can be viewed as equivalent via Poincar\'e duality
to  the proper push-forward  map in 
Borel-Moore homology $H^{BM}_j(Z, \rat)  \to H^{BM}_j(Y,\rat)$.

\medskip
{\bf Fundamental Class.}
Let $i : Z \to Y$ be the closed immersion of a  $d$-dimensional 
subvariety of the manifold $Y$. 
The space $Z$ carries a fundamental class in $H_{2d}^{BM}(Z)$.
The fundamental class of $Z$ is the image of this class in 
$H^{BM}_{2d}(Y) \simeq H^{2n-2d}(Y,\zed)$.

\medskip
{\bf Mayer-Vietoris.}
There is a whole host of  Mayer-Vietoris sequences (cf. \ci{k-s},
2.6.10), e.g.:
\[\xymatrix{
\cdots \ar[r] & \
H^{l-1}(U_1 \cap U_2, K) \ar[r] &
H^{l}(U_1 \cup U_2, K) \ar[r] &
H^{l}(U_1 , K) \oplus H^{l}(U_2 , K) \ar[r] &
 \cdots
 }
 \]

{\bf Relative (co)homology}
Let $U \stackrel{j}\to Y \stackrel{i} \leftarrow Z$ be the inclusions of an open subset
$U\subset Y$ and of the closed complement $Z:= Y \setminus U$.
There are  the following ``attaching" distinguished triangles:
\[
i_! i^! C \lorw C \lorw j_{*} j^{*} C \stackrel{[1]} \lorw, \qquad \quad
j_! j^! C \lorw C \lorw i_{*} i^{*} C \stackrel{[1]}\lorw. 
\]
The long exact sequences of relative (co)homology (including the versions with 
compact supports) arise by taking the associated  long exact sequences.

\medskip{\bf Refined intersection forms}
Let $i:Z \to Y$ be a closed immersion into a nonsingular variety $Y$
of dimension $n$.
There are maps
$$
i_! \omega_Z [-n] = i_! i^! \omega_Y [-n] \lorw \omega_Y[-n] \simeq \rat_Y[n] \lorw i_{*}
i^{*} \rat_Y[n] = i_{*} \rat_Z [n].
$$
Taking cohomology we get the   {\em refined intersection form} on $Z\subseteq Y$,
which we can view in two equivalent ways as a linear or a bilinear map:
$$
H^{BM}_{n-l} (Z) \lorw H^{n+l}(Z), \quad \mbox{or} \qquad
H^{BM}_{n-l} (Z) \times H_{n+l}(Z) \lorw \rat.
$$
It is called refined because we are intersecting cycles in 
 the nonsingular $Y$ which are supported on $Z$.
By using Lefschetz Duality, this pairing can be viewed as the 
cup product in relative cohomology. These 
forms play  an important role  in our proof
of the decomposition theorem \ci{herdlef, decmightam} (see
$\S$\ref{dmapp}).

\subsection{Nearby and vanishing cycle functors }
\label{psiphi}
An important feature of perverse sheaves is their stability for the two functors 
$\Psi_f$, $\Phi_f$. These functors were defined in \ci{deligneexpose14} in 
the context of \'etale cohomology 
as a generalization of the notion of vanishing cycle in the classical Picard-Lefschetz theory.
As  explained  in $\S$\ref{beiverd},
they play a major role in the description of the possible  extensions of a perverse sheaf 
through a principal divisor.
We discuss these functors in the complex analytic setting.
Let  $f:X \to \comp$ be a regular
function  and $X_0 \subseteq X $ be its divisor, that is 
$X_0= f^{-1}(0)$. We are going to define functors 
$\Psi_f, \Phi_f: {\cal D}_X \to {\cal D}_{X_0}$ 
which  send perverse sheaves on $X$ to perverse sheaves on $X_0$. 
We follow the convention for  shifts employed in 
 \ci{k-s}.

Let $e: \comp \to \comp$ be the map $e( \zeta)=\exp(2 \pi \sqrt{-1} \zeta)$ and consider the following diagram
$$
\xymatrix{
      X_{\infty}:=X\times_{e}\comp \ar[r] \ar[d] \ar@/^1.5pc/[rr]^p & X^*\ar[r]^j\ar[d] & X\ar[d] & X_0 
\ar[l]_i \ar[d] \\
                 \comp \ar[r] \ar@/_1.5pc/[rr]^ e                   &    \comp ^* \ar[r]      & \comp  &  
 \{o\} \ar[l].
}
$$
For $K \in {\cal D}_X$, the nearby cycle functor $\Psi_f(K) \in {\cal D}_{X_0} $ is defined as:
$$
\Psi_f(K):=i^{*}p_{*}p^{*}K.
$$
Note that $\Psi_f(K)$ depends  only on the restriction of $K$ to $X^*$.
It can be shown that $\Psi_f(K)$ is constructible. 
Depending on the context, we shall consider $\Psi_f$ as a functor defined on ${\cal D}_X$, or on ${\cal D}_{X^*}$.

The group $\zed$ of deck transformations $\zeta \to \zeta +n$ acts on $X_{\infty}$ 
and therefore on $\Psi_f(K)$. We  denote by 
$T:\Psi_f(K)\to \Psi_f(K)$ the positive generator of this action.

\begin{rmk}
\label{milfib}
{\rm (See \ci{gomacsmt}, $\S$6.13 for details.) Under mild hypothesis, for instance if $f$ is proper, there exists a continuous map 
$r: U \to X_0$ of a neighborhood of $X_0$, compatible with the stratification,
whose restriction to $X_0$ is homotopic to the identity map. Denote by $r_{\epsilon}$
the restriction of $r$ to $f^{-1}(\epsilon)$,
with $\epsilon  \in \comp$ small enough so that $f^{-1}(\epsilon) \subseteq U$. Let $X_{\e}:=
f^{-1} (\e)$.
Then
$$
{r_{\epsilon }}_*(K_{\mid_{X_{\e}}})=\Psi_f(K)
$$
In particular, let $x_0 \in X_0$, let $N$ be a neighborhood of $x_0$ contained in $U$ and let $\epsilon \in \comp$ be as before. 
Then the cohomology sheaves of $\Psi_f(K)$ can be described as follows:
$$
{\cal H}^i(\Psi_f(K))_{x_0}=H^i(N\cap f^{-1}(\epsilon),K_{| N \cap f^{-1}(\epsilon)}).
$$
}
\end{rmk}

The monodromy $X_\e \to X_\e$ induces a transformation $T: \Psi_f (K) \to \Psi_f (K)$
called the {\em monodromy transformation}.

\begin{ex}
\label{ex000}
{\rm  Let $X= \comp$ and $K$ be a local system on $\comp^*$. Since the inverse
 image by $e$ of a disk centered at $0$ is contractible, $\Psi_f(K)$
can be identified  with the stalk at some base point $x_0$.  The automorphism $T$
 is just the monodromy of the local system. }
\end{ex}

The adjunction map  $ K \to p_{*}p^{*}K$ gives a natural morphism $i^{*}K \to \Psi_f(K)$.
The vanishing cycle complex $\Phi_f(K) \in {\cal D}_{X_0}$ fits in the following distinguished triangle:
\begin{equation}
\label{cano}
\xymatrix{
i^{*}K \ar[r] &  \Psi_f(K)  \ar[r]^(.45){\rm can} & 
 \Phi_f(K)[1] \ar[r]^{\,\,[1]\!\!\!} &.}
\end{equation}

This distinguished  triangle determines $\Phi_f(K) $ only up to a non unique isomorphism.
The definition of $\Phi_f$ as a functor requires more care, see \ci{k-s}.
The long exact sequence for the cohomology sheaves of this distinguished
triangle, and
Remark \ref{milfib}, show that 
$$
{\cal H}^i(\Phi_f(K))_{x_0}=H^i(N,N\cap f^{-1}(\epsilon),K).
$$
Just as the nearby cycle functor,   the vanishing cycle $\Phi_f(K)$ is endowed with an automorphism $T$.

We now list some of the properties of the functors $\Psi_f$ and $\Phi_f$:

\begin{tm}
\label{lpnvf}
$\,$
\begin{enumerate}
\item
The functors commute, up to a shift, with Verdier duality (see \ci{illu}, and \ci{brili1}):
$$
\Psi_f(D K)=D\Psi_f(K)[2]  \qquad \Phi_f(DK)=D\Phi_f(K)[2].
$$
\item If $K$ is a perverse sheaf on $X$, 
then $\Psi_f(K)[-1]$  and $\Phi_f(K)[-1]$ are perverse sheaves on $X_0$ (see \ci{gomacsmt} 6.13,  \ci{bbd}, \ci{brili1}, \ci{illu}). 

\item
Dualizing the distinguished triangle (\ref{cano}) we get an distinguished triangle
\begin{equation}
\label{var}
\xymatrix{
i^{!}K \ar[r] &  \Phi_f(K)  \ar[r]^(.42){\rm var} & \Psi_f(K)[-1]
\ar[r]^(.7){[1]}&,} 
\end{equation}
with the property
that
$$
{\rm can} \circ {\rm var} =T-I:\Phi_f(K) \to \Phi_f(K)  \qquad {\rm var} \circ 
{\rm can} =T-I:\Psi_f(K) \to \Psi_f(K),
$$
and we have the fundamental octahedron of complexes of sheaves on $X_0:$


$$
 \xymatrix{
  &  & i^{*}j_{*}j^{*}K \ar[drr]    \ar[ddll]&   &
  \\
  i^{*}K[1]  \ar[rru]   \ar[d]_{[1]} & & & & i^!K[1] \ar@{-->}[llll]  \ar@{-->}[ddll] 
  \\
  \Psi_f(K)  \ar[rrrr]^{T-I}  \ar[drr]_{\rm can} & &&& \Psi_f(K) \ar[u]    \ar[lluu]^(.4){[1]}
  \\
  &&  \Phi_f(K)[1]  \ar@{-->}[lluu]  \ar[rru]_{{\rm var} [1]} && 
 }
 $$
 \end{enumerate}\end{tm}

\begin{rmk}
\label{noglobal}
{\rm
Clearly,  if $U \subseteq X$ is an open subset, then the restriction to $U$
of $\Psi_f(K)$ is the nearby cycle complex of the restriction $K_{|U}$ relative to the  function
$f_{|U}$ for $X \cap U$. On the other hand, explicit examples show that
$\Psi_f(K)$ depends on $f$ and not only on the divisor $X_0$: 
the nearby functors associated with different defining equations of $X_0$ may differ.  
In particular, it is not possible to define the functor $\Psi_f$ if the divisor $X_0$  is only locally principal. 
Verdier has proposed  
in  \ci{ve1}  an alternative functor, which he called the ``specialization functor'' 
${\rm Sp}_{Y,X}:{\cal D}_X \to {\cal D}_{C_Y}$,  
associated with any closed imbedding $Y \to X$, where $C_Y$ is the normal cone
of $Y$ in $X$. In the particular case that $Y$ is a locally principal divisor in $X$, 
the specialization functor is related to the nearby functor as follows:
the normal cone $C_Y$ is  a line bundle, and a local defining equation $f$ of $Y$ defined on an open set 
$V \subseteq X$ defines a section $s_f: Y\cap V \to  C_{Y\cap V}$ trivializing the fibration. 
One has an isomorphism of functors $s_f^*{\rm Sp}_{Y,X} \simeq \Psi_f $.
}
\end{rmk}

\subsection{Unipotent nearby and vanishing cycle functors }
\label{unipotpsiphi}
Let $K$ be a perverse sheaf on $X \setminus X_0$. 
The map $j: X\setminus X_0 \to X$ is affine, so that
$j_{*}K$ and $j_{!}K$ are perverse sheaves on $X$.

Let us consider the ascending chain  of perverse subsheaves 
$$\ke \,  \{(T-I)^N:\Psi_f(K)[-1]\to \Psi_f(K)[-1]\}. $$ For $N\gg0$ this sequence
stabilizes because of the N\"oetherian property of the category of perverse sheaves.
We call the resulting     $T$-invariant perverse  subsheaf the unipotent
nearby cycle perverse sheaf  associated with $K$
and we denote it  by $\Psi_f^u(K)$. 
In exactly the same way, it is possible to define the unipotent vanishing cycle functor 
$\Phi_f^u:$
$$\Phi_f^u(K)= \ke \,  \{ \, (T-I)^N:\Phi_f(K) \to \Phi_f(K)\, \}, \qquad  \mbox{for $N\gg0$.} $$

The perverse sheaves $\Psi_f(K)[-1]$ and $\Phi_f(K)[-1]$ are in fact the direct sum of $\Psi_f^u$ and another $T$-invariant subsheaf on which $(T-I)$ is invertible.  
\begin{rmk}
{\rm The functor $\Psi_f(K)$ on a perverse sheaf $K$ can be reconstructed from $\Psi_f^u$ by applying this latter to the twists of $K$ with the pullback by $f$ of  local systems on $\comp^*;$ see \ci{beili},
p.47. 
}
\end{rmk}

We have  the useful formul\ae
$$
\ke \, \{ \,  j_!K \to  j_{!*}K \, \}  \simeq  \ke\,  \{ \, \Psi_f^u(K) \stackrel{T-I}{\lorw} \Psi_f^u(K) \, \} ,$$
$$
\coke \,  \{ \,  j_{!*}(K) \to j_{*}K  \, \} \simeq 
\coke \, \{ \, \Psi_f^u(K) \stackrel{T-I}{\lorw} \Psi_f^u(K) \, \}.
$$
They can be derived as follows.
The cone of $(T-I): \Psi_f(K) \to \Psi_f(K)$, which is isomorphic to  $i^{*}j_{*}K$,  is also isomorphic, up to a shift  $[1]$, 
to the cone 
of $(T-I): \Psi_f^u(K) \to \Psi_f^u(K)$,
and we still have the distinguished triangle
$$
\xymatrix{
i^{*}j_{*}K \ar[r]^(.46){[1]}&  \Psi_f^u(K)\ar[r]^{T-I} &  \Psi_f^u(K) \ar[r] .& 
}. $$
The long exact sequence of perverse cohomology introduced in $\S$\ref{subsecpervco} then gives
$$
^{\frak p}\!{\cal H}^{-1}(i^{*}j_{*}K) = \ke \,  \{ \,  \Psi_f^u(K) \stackrel{T-I}\longrightarrow \Psi_f^u(K) \, \}
$$
and 
$$
^{\frak p}\!{\cal H}^{0}(i^{*}j_{*}K) = \coke \, \{ \,   \Psi_f^u(K) \stackrel{T-I}\longrightarrow \Psi_f^u(K) \, \}.
$$
In turn, the long exact perverse cohomology sequence of the distinguished triangle
$$
\xymatrix{
i^{*}j_{*}K \ar[r]^{[1]}& j_!K \ar[r] &  j_{*}K \ar[r] & 
}$$
and the fact that $j_{*}K$ and $j_{!}K$ are perverse sheaves on $X$,
give
$$
 ^{\frak p}\!{\cal H}^{-1}(i^{*}j_{*}K) = \ke \, \{ \,  j_!K \to  j_{*}K\, \}=
\ke\, \{ \, j_!K \to  j_{!*}K \}.
$$
and 
$$
 ^{\frak p}\!{\cal H}^{0}(i^{*}j_{*}K) = \coke \, \{ \,  j_!K \to  j_{*}K \, \}=
\coke \,  \{ \,   j_{!*}K \to j_{*}K \, \}.
$$

\begin{rmk}
\label{mwf}
{\rm 
Let $N$ be  a nilpotent endomorphism of an object 
$M$ of an abelian category.
Suppose $N^{k+1}=0$. By \ci{weil2}, 1.6,  there exists a unique finite increasing filtration 
$$M_{\bullet}: \{0\}\subseteq M_{-k} \subseteq \ldots
\subseteq  M_k=M$$ such that:
$$NM_l\subseteq M_{l-2}\;  \hbox{ and } \; N^l:M_l/M_{l-1} \simeq M_{-l}/M_{-l-1}. $$
The filtration defined in this way by $T-I$ on $\Psi_f^u(K)$ is called the 
{\em monodromy weight filtration}. An important theorem of O.  Gabber
(see \ci{beibe}, $\S$5)
characterizes this filtration in the case of $l$-adic perverse sheaves.
}
\end{rmk}

\subsection{Two descriptions of the category of perverse sheaves}
\label{structpervshvs}
In this section we discuss two descriptions of the category of perverse sheaves 
on an algebraic variety.
Although not strictly necessary for what follows, they play an important role in the 
theory and applications
of perverse sheaves. The question is roughly as follows: suppose $X$ is an algebraic 
variety, $Y \subseteq X$ a subvariety, and we are given a perverse sheaf $K$ on $X \setminus Y$.
How much information is needed to describe the perverse sheaves $\widetilde{K}$ on $X$ 
whose restriction to $X\setminus Y$
is isomorphic to K? We  describe the approach developed by R. MacPherson and K. Vilonen \ci{macvilo}
and the approach of  A. Beilinson and  J.L. Verdier \ci{beilin, ve}.

\subsubsection{The approach of MacPherson-Vilonen}
\label{macvilon}
We  report on only a part of the description of the category of perverse sheaves 
developed in \ci{macvilo}, i.e.
the most elementary and the one which we find  particularly illuminating.

Assume   that $X =  Y \coprod (X \setminus Y)$, where $Y$  is a  closed and  contractible
 $d$-dimensional stratum of   a stratification $\Sigma$ of $X$. 
 We have  ${\cal P}_{X}^{\Sigma}$, i.e.  the category of perverse sheaves on $X$ which are constructible 
with respect to $\Sigma$.
Denote by $
Y \stackrel{i}{\longrightarrow} X \stackrel{j}{\longleftarrow} X\setminus Y
$ the corresponding imbeddings.

For  $K \in {\cal P}_X^{\Sigma}$, the attaching triangle
$\;
i_! i^! K \lorw K \lorw j_{*} j^{*} K \stackrel{[1]} \lorw,
$
and 
 the support and cosupport conditions for a perverse sheaf
give the following  exact sequence of local systems on $Y$:

\begin{equation}
\label{fexse}
\xymatrix{
0\ar[r]  & {\cal H}^{-d-1}(i^{*}K) \ar[r] & 
{\cal H}^{-d-1}(i^{*}j_{*}j^{*}K)\ar[r] & {\cal H}^{-d}(i^!K)\ar[d]    \\ 
0          &  {\cal H}^{-d+1}(i^!K)\ar[l]
&{\cal H}^{-d}(i^{*}j_{*}j^{*}K)\ar[l] & {\cal H}^{-d}(i^{*}K) \ar[l]}.
\end{equation}

Note that the (trivial) local systems ${\cal H}^{-d-1}(i^{*}j_{*}j^{*}K),{\cal H}^{-d}(i^{*}j_{*}j^{*}K) $
are determined by the restriction of $K$ to $X\setminus Y$.

A first approximation to the category of perverse sheaves is given as  follows:

\begin{defi}
\label{scurda}
{ \rm Let ${\cal P}'_X$ be the following category:

\n
--  an object is 
a perverse sheaf $K$ on $X\setminus Y$, constructible with respect to $\Sigma_{|X-Y}$, 
and an exact sequence 
$$ {\cal H}^{-d-1}(i^{*}j_{*}K)\lorw V_1 \lorw V_2 \lorw {\cal H}^{-d}(i^{*}j_{*}K) $$
of local systems on $Y;$

\n
-- a morphism  $(K, \ldots)$ $(L, \ldots)$ is
 a morphism   of perverse sheaves $\phi:K \to L$
together with  morphisms of exact sequences:
$$
\xymatrix{
{\cal H}^{-d-1}(i^{*}j_{*}K)\ar[r] \ar[d]^{\phi} & V_1 \ar[r] \ar[d] & V_2 \ar[r] \ar[d] & {\cal H}^{-d}(i^{*}j_{*}K)) \ar[d]^{\phi} \\
{\cal H}^{-d-1}(i^{*}j_{*}L)\ar[r]  & W_1 \ar[r] & W_2 \ar[r]  & {\cal H}^{-d}(i^{*}j_{*}L)).
}
$$
}
\end{defi}
\begin{tm}
\label{macvildescr}
The functor ${\cal P}_{X}^{\Sigma} \to {\cal P}'_X$, sending a perverse sheaf 
$\widetilde{K}$ on $X$ to its restriction to $X\setminus Y$
and to the exact sequence 
$$
\xymatrix{
{\cal H}^{-d-1}(i^{*}j_{*}j^{*} \widetilde{K}) \ar[r] &
 {\cal H}^{-d}(i^! \widetilde{K}) \ar[r] & 
          {\cal H}^{-d}(i^{*}\widetilde{K})   \ar[r] &
 {\cal H}^{-d}(i^{*}j_{*}j^{*} \widetilde{K})
}$$
 is a bijection on isomorphism classes of objects.
\end{tm}

To give an idea why the theorem is true, we note that for any object $Q$ in ${\cal P}_X$,
we have the  distinguished triangle
$$
\xymatrix{
i_! i^! Q \ar[r] & Q \ar[r]  &  j_{*} j^{*} Q \ar[r]^{[1]} &,}
$$ 
and $Q$ is identified by the extension map $e \in \mbox{Hom} ( j_{*} j^{*} Q ,i_! i^! Q[1])$.
We have $i_!=i_{*}$ hence
$$ 
\mbox{Hom}( j_{*} j^{*} Q ,i_! i^! Q[1])=\mbox{Hom}(i^{*} j_{*} j^{*} Q ,i^! Q[1])
= \oplus_l  \mbox{Hom}({\cal H}^l(i^{*} j_{*} j^{*} Q), {\cal H}^{l+1}(i^! Q)). $$
The last equality is due to the fact that the derived category of complexes with constant cohomology 
sheaves on a contractible space is semisimple ($K \simeq \oplus H^i(K)[-i]$, for every $K$).
By the support condition 
$$
{\cal H}^{l}(i^{*} j_{*} j^{*} Q) \simeq {\cal H}^{l+1}(i^! Q)
\hbox{ for }l>-d.
$$
By the co-support condition, 
$$
{\cal H}^{l}(i^! Q)=0
\hbox{ for }l<-d.
$$
There are the two maps 
$$
{\cal H}^{-d}(i^{*} j_{*} j^{*} Q) \lorw  {\cal H}^{-d+1}(i^! Q), \qquad
{\cal H}^{-d-1}(i^{*} j_{*} j^{*} Q) \lorw {\cal H}^{-d}(i^! Q)
$$
which are not determined a priori by the restriction of $Q$ to $X\setminus Y$.
They appear in  the exact sequence (\ref{fexse}) and contain the information about how to glue
$j^{*} Q $ to $i^! Q$.
The datum of this exact sequence makes it possible to reconstruct $Q \in {\cal P}_X$
 satisfying  the support and cosupport conditions.

Unfortunately the functor is not as precise on maps, as we will see.
There are non zero maps between perverse sheaves which induce the zero map in  ${\cal P}'_X$,
i.e. the corresponding functor is not faithful.
However, it is  interesting to see a few examples of applications of this result.

\begin{ex}
{\rm
Let  $X=\Bbb C$, $Y=\{o\}$ with strata $X \setminus Y =\comp^*$ and $Y$.  
A perverse sheaf on $\comp ^*$ is then of the form $L[1]$ for $L$ a local system.
Let $L$ denote the stalk of $L$ at some base point,
and $T:L \to L$  the monodromy. 
An explicit computation shows that 
$$
i^{*}j_{*}L[1]\simeq \ke \, (T-I)[1]\oplus \coke\, (T-I),
$$
where $\ke (T-I)$ and $ \coke (T-I)$ are interpreted as sheaves on $Y$.
Hence a perverse sheaf is identified up to isomorphism by $L$ and by 
an exact sequence of vector spaces:
$$
\ke\, (T-I) \lorw  V_1 \to V_2 \lorw   \coke\, (T-I).
$$
A sheaf of the form $i_{*}V$ is represented by $L=0$
and by the sequence
$$
\xymatrix{
0 \ar[r] &  V \ar[r]^{\simeq} &  V \ar[r] &  0.}
$$
Since $j$ is an affine imbedding, $j_{*}$ and $j_!$
are $t$-exact, i.e.  $j_{*}L[1]$ and $j_!L[1]$ are perverse.

\n
The perverse sheaf $j_{*}L[1]$ is represented by 
$$
\xymatrix{
\ke \, (T-I) \ar[r] &  0 \ar[r] &   \coke \,  (T-I)
\ar[r]^{Id} &   \coke (T-I) ,}
$$
which expresses the fact that $i^!j_{*}L[1]=0$.

\n
Similarly $j_!L[1]$, which verifies  $i^{*}j_!L[1]=0$,
is represented by 
$$
\xymatrix{
\ke\,  (T-I) \ar[r]^{Id} &  \ke \, (T-I) \ar[r] &  0 \ar[r]    
& \coke \, (T-I).
}$$
The intermediate extension $j_{!*}L[1]$
is represented by
$$\xymatrix{
\ke\, (T-I) \ar[r]  & 0 \ar[r]   &  0  \ar[r]  & \coke\, (T-I),}
$$
since, by its very definition, 
$$
{\cal H}^0(i^{*}j_{!*}L[1])= {\cal H}^0(i^!j_{!*}L[1])=0.
$$
Let us note another natural exact sequence
given by 
$$
\xymatrix{
\ke\, (T-I) \ar[r]  & L \ar[r]^{T-I} &  L \ar[r] & \coke \, (T-I)
}$$
which corresponds to Beilinson's maximal extension $\Xi(L)$, which will be described in the next section.
>From these presentations one sees easily the natural maps 
$$ j_!L[1] \lorw j_{!*}L[1] \lorw j_{*}L[1], \qquad \mbox{and}
\qquad
 j_!L[1] \lorw \Xi(L[1]) \lorw j_{*}L[1].
$$
}
\end{ex}

\begin{rmk}
{\rm If  $T$ has no eigenvalue equal to one,  then the sequence
has the form 
$0 \to V \to V \to 0$. This corresponds to the fact that a perverse sheaf 
which restricts to such a local system on $\comp \setminus \{o\}$
is necessarily of the form $j_!L[1]\oplus i_{*}V$. Note also that $j_!L[1]=j_{*}L[1]=j_{!*}L[1]$.
}
\end{rmk}

\begin{rmk}
\label{splitdecmig}
{\rm One can use theorem \ref{macvildescr} to  deduce the following 
special case of a splitting criterion used in our 
proof of the decomposition theorem  \ci{decmightam}:} 
Let $d=dimY$. A perverse sheaf $K \in {\cal P}_X$ splits as $K \simeq j_{!*}j^{*}K \oplus {\cal H}^{-d}(K)[d]$
if and only if the map 
${\cal H}^{-d}(i^!K)\to
          {\cal H}^{-d}(i^{*}K)$ is an isomorphism.         
          
{\rm In fact, if this condition is verified,  then 
the maps 
${\cal H}^{-d-1}(i^{*}j_{*}j^{*}K) \to {\cal H}^{-d}(i^!K)$
and      $ {\cal H}^{-d}(i^{*}K) \to
{\cal H}^{-d}(i^{*}j_{*}j^{*}K)$
in (\ref{fexse}) vanish, and the exact sequence corresponding to $K$ is of the form
$$
\xymatrix{
 {\cal H}^{-d-1}(i^{*}j_{*}K) \ar[r] & 0 \ar[r]           
 &  0 \ar[r]  & {\cal H}^{-d}(i^{*}j_{*}K))  & j_{!*}j^{*}K  \\
& \ar@{}[r]^{\bigoplus}&    &                     \\
& W \ar[r]    &       W &   & {\cal H}^{-d}(K)[d].
}
$$
}
\end{rmk}

The following example   shows that the functor  ${\cal P}_{X}^{\Sigma} \to {\cal P}'_X $
 is not faithful.
Consider the perverse sheaf $j_{*}\rat_{\comp^*}[1]$.
It has a non-split  filtration by perverse sheaves
$$\xymatrix{
0 \ar[r]  & \rat_{\comp}[1] \ar[r] & 
j_{*}\rat_{\comp^*}[1] \ar[r]^(.54){\alpha} & i_{*}\rat_0 \ar[r] &  0.}
$$
Dually,  the perverse sheaf $j_!\rat[1]$ has a non-split  filtration
$$
\xymatrix{
0 \ar[r] &  i_{*}\rat_0 \ar[r]^(.43){\beta} & 
j_!\rat_{\comp^*}[1] \ar[r] &  \rat_{\comp}[1] \ar[r] &  0.}
$$
The composition $\beta \alpha: j_{*}\rat_{\comp^*}[1] \to j_!\rat_{\comp^*}[1] $ 
is not zero, being the composition of the epimorphism $\alpha$ with the 
monomorphism $\beta$,  however, it is zero on $\comp^*$, and the map between the associated exact 
sequences is zero, since $i^!j_{*}\rat_{\comp^*}[1]=0$ and 
$i^{*}j_!\rat_{\comp^*}[1]=0$. 

In the paper \ci{macvilo}, MacPherson and Vilonen give a refinement 
of the construction which describes completely the category of perverse sheaves, 
both in the topological and complex analytic situation. For an application to 
representation theory, see   \ci{mirovilo}.

\subsubsection{The approach of Beilinson and Verdier}
\label{beiverd}
We turn to the Beilinson's approach  \ci{beili}, i.e. the one  used by Saito in his
theory of  mixed Hodge modules.
Beilinson's approach is based on the  nearby and vanishing cycle functors $\Psi_f$ and $\Phi_f$ introduced in $\S$\ref{psiphi}. 
In \ci{verd}, Verdier obtained similar results using 
the specialization  to the normal cone functor ${\rm Sp}_{Y,X}$, \ref{noglobal}, which is not  discussed here. 

The assumption is that we have an algebraic map $f:X \to \comp$ and $X_0=f^{-1}(0)$ as in $\S$\ref{psiphi}. 
Let $K$ be  a perverse sheaf on $X\setminus X_0$. Beilinson defines an interesting extension 
of $K$ to $X$
 which he calls the maximal extension and denotes by $\Xi(K)$. 
It is a perverse sheaf, restricting to $K$ on $X\setminus X_0$, which
can be constructed as follows:
consider the unipotent nearby and vanishing cycle functor $\Psi_f^u$ and 
$\Phi_f^u$ (see  $\S$\ref{unipotpsiphi})  and the distinguished triangle 
$$
\xymatrix{
i^{*}j_{*}K \ar[r]^{[1]}&  \Psi_f^u(K)\ar[r]^{T-I} &  \Psi_f^u(K) \ar[r] &. 
}$$
The natural map $i^{*}j_{*}K \to \Psi_f^u(K)[1]$ defines, by adjunction,  an element 
of 
$$\mbox{Hom}^1_{{\cal D}_{X_0}}(i^{*}j_{*}K,\Psi_f ^u(K))=
\mbox{Hom}^1_{{\cal D}_X}(j_{*}K,i_{*}\Psi_f^u(K))$$ 
which, in turn, defines an object
$\Xi(K)$ fitting in the  distinguished triangle
\begin{equation}
\label{beilifunctor}
i_{*}\Psi_f^u(K)  \lorw \Xi(K)   \lorw j_{*}K \lorw i_{*}\Psi_f^u(K)[1].
\end{equation}
Since $j$ is an affine morphism, it follows that $j_{*}K$ is perverse. 
The long exact sequence of perverse cohomology implies  that 
$\Xi(K)$ is perverse as well. 

In \ci{beili},  Beilinson gives a different construction of $\Xi(K)$ 
(and also of $\Psi_f^u(K)$ and $\Phi_f^u(K)$) 
which implies automatically that $\Xi$ is a functor and that it 
commutes with Verdier duality. 

There is the exact sequence of perverse sheaves 
$$
\xymatrix{
0\ar[r] &  i_{*}\Psi_f^u(K)\ar[r]^(.55){\beta_+} & \Xi(K) \ar[r]^(.53){\alpha_+} & j_{*}K  \ar[r] & 0
}
$$
and, applying Verdier duality and the canonical isomorphisms $\Xi \circ D \simeq D \circ \Xi$
and $\Psi_f^u \circ D \simeq D \circ \Psi_f^u$, 
$$
\xymatrix{
0\ar[r] & j_{!}K \ar[r]^(.48){\alpha_-} & \Xi(K) \ar[r]^(.45){\beta_-} & i_{*}\Psi_f^u(K) \ar[r] & 0. 
}
$$
The composition $\alpha_+ \alpha_-:j_!K \to j_{*}K$ is the natural map,
while $\beta_- \beta_+:  i_{*}\Psi_f^u(K) \to i_{*}\Psi_f^u(K) $ is $T-I$.
We may now state Beilinson's results.

\begin{defi}
{\rm Let $Gl(X,Y)$ be the category whose objects are quadruples $(K_U,V,u,v)$, where 
$K_U$ is a perverse sheaf on $U:= X \setminus X_0$, $V$ is a perverse sheaf on $Y$,
 $u:\Psi_f^u(K) \to V$, and $v:V \to \Psi_f^u(K)$ are morphisms
such that $vu=T-I$.}  
\end{defi}

\begin{tm}
\label{beiresult}
The functor $\gamma: {\cal P}_X \to Gl(X,Y)$ which associates to a perverse sheaf $K$ on $X$
the quadruple $(j^{*}K, \Phi_f^u(K), can,var)$ is an equivalence of categories.
 Its inverse is the functor $G: Gl(X,Y) \to {\cal P}_X$ associating to 
$(K_U,V,u,v)$ the cohomology of the complex
$$
\xymatrix{
 \Psi_f^u(K_U) \ar[rr]^(.45){(\beta_+,u)} &&  \Xi(K_U)\oplus V \ar[rr]^(.52){(\beta_-,v)} &&  \Psi_f^u(K_U).
}
$$
\end{tm}

\begin{ex}
\label{beqex}
{\rm Given a perverse sheaf $K_U$  on $U= X \setminus X_0$, we 
determine 
$$
\xymatrix{
\gamma(j_!K_U) \ar[r] &   \gamma(j_{!*}K_U) \ar[r] &  \gamma(j_{*}K_U).}
$$  
We make use of the distinguished  triangles (\ref{cano}) and (\ref{var}) discussed in $\S$\ref{psiphi}
and restricted to the unipotent parts $\Psi_f^u$ and $\Phi_f^u$.
Since $i^{*}j_!K_U=0$, the map ${\rm can} : \Psi_f^u(j_!K_U) \to \Phi_f^u(j_!K_U)$ is an isomorphism.
Hence 
$$
\gamma(j_!K_U)= \Psi_f^u(K_U) \stackrel{\rm Id}{\lorw} \Psi_f^u(K_U) \stackrel{T-I}{\lorw} \Psi_f^u(K_U).
$$
Similarly, since $i^{!}j_{*}K_U=0$, the map $var:\Phi_f^u(j_{*}K_U) \lorw \Psi_f^u(j_{*}K_U)$
is an isomorphism, and
$$
\gamma(j_{*}K_U)= \Psi_f^u(K_U) \stackrel{T-I}{\lorw} \Psi_f^u(K_U) \stackrel{id}{\lorw} \Psi_f^u(K_U).
$$
The canonical map $j_!K_U \to j_{*}K_U$ is represented by the following diagram, in which 
we do not indicate the identity maps:
\begin{equation}
\label{!to*}
\xymatrix{
\gamma(j_!K_U) \ar[d]        &   \Psi_f^u(K_U) \ar[r]^{} \ar[d]^{} & \Psi_f^u(K_U) \ar[d]^{T-I} \ar[r]^{T-I} & \Psi_f^u(K_U) \ar[d]^{}  \\
\gamma(j_{*}K_U)       &   \Psi_f^u(K_U) \ar[r]^{T-I} & \Psi_f^u(K_U) \ar[r]^{} & \Psi_f^u(K_U) 
.}
\end{equation}
The intermediate extension $j_{!*}K_U$ corresponds to    $j_{!*}K_U:= \im \{\, j_!K_U \to j_{*}K_U \, \}$, 
hence
$$
\xymatrix{
\gamma(j_{!*}K_U)=& \Psi_f^u(K_U) \ar[r]^{T-I}& \im (T-I) \ar@{^{(}->}[r] & \Psi_f^u(K_U),
}
$$
where the second map is the canonical inclusion. We can  complete the diagram (\ref{!to*}) as follows:
\begin{equation}
\label{!to!*to*}
\xymatrix{
\gamma(j_!K_U) \ar[d]        &   \Psi_f^u(K_U) \ar[r]^{} \ar[d]^{} & \Psi_f^u(K_U) \ar[d]^{T-I} \ar[r]^{T-I} & \Psi_f^u(K_U) \ar[d]^{}  \\
\gamma(j_{!*}K_U) \ar[d]        &   \Psi_f^u(K_U) \ar[r]^{T-I} \ar[d]^{} & \im (T-I) \ar[d] \ar@{^{(}->}[r] & \Psi_f^u(K_U) \ar[d]^{}  \\
\gamma(j_{*}K_U)       &   \Psi_f^u(K_U) \ar[r]^{T-I} & \Psi_f^u(K_U) \ar[r]^{} & \Psi_f^u(K_U) 
.}
\end{equation}
The maximal extension 
$\Xi(K_U)$ is represented by the factorization
$$\xymatrix{
\Psi_f^u(K_U) \ar[r]^>>>>>{(I,T-I)} & \Psi_f^u(K_U)\oplus \Psi_f^u(K_U) \ar[r]^>>>>{p_2} & \Psi_f^u(K_U)
}$$
where $p_2((a_1,a_2))=a_2$ is the projection on the second factor.
Finally we note that if $L$
is a perverse sheaf on $X_0$, then,  since $\Psi_f(i_{*}L)=0$,
$$
\xymatrix{
\gamma(i_{*}L) =    0 \ar[r] &   L \ar[r]  & 0.} $$
}
\end{ex}

\begin{rmk}
\label{othersplit}
{\rm From the examples of $\gamma(j_{!*}K_U)$ and  $\gamma(i_{*}L)$ discussed in 
Example \ref{beqex}, one can derive the following 
criterion (Lemme 5.1.4 in \ci{samhp}) for a perverse sheaf $K$ on $X$
to split as $K \simeq j_{!*}j^{*}K \oplus i_{*}L$.}
Let $X$ be an algebraic variety and $X_0$ be 
a principal divisor; let $i:X_0 \to X \longleftarrow X \setminus X_0:j$ 
be the corresponding closed and open imbeddings. A perverse sheaf $K$ on $X$ is of the form 
$K \simeq j_{!*}j^{*}K \oplus i_{*}L$ if and only if 
\[\Phi_f^u(K)=  \{ \im : (\Psi_f^u(K) \stackrel{\rm can}{\to} \Phi_f^u(K)) \}
 \bigoplus \{ \ke : (\Phi_f^u(K)
\stackrel{\rm var}{\to}\Psi_f^u(K))\}.\]

\n
{\rm This criterion is used in  \ci{samhp}
to establish the semisimplicity of certain perverse sheaves.}
\end{rmk}

\subsection{A formulary for the constructible derived category}
\label{formulary}

Throughout  this section,  $f: X \to Y$, $g: Y' \to Y$ and  
$h: Y \to Z$
 are maps of varieties, $C \in {\cal D}_X$ is a constructible
 complex on $X$ and  $K, K', K_{i} \in {\cal D}_Y$  are constructible complexes on $Y$.
 An equality sign actually stands for the existence
 of a suitably canonical isomorphism. Since we use field coefficients,
 the tensor product is exact and  it coincides with the associated 
  left derived functor.
   Perversity means middle perversity on  complex varieties. 
All operations preserve stratifications of varieties and of maps.
We use the simplified notation $f_*:= Rf_*$, $f_!:= Rf_!$.
Some standard references are \ci{k-s, goma2, borel, iv, gel-man, bbd,godement}.

\bigskip

{\bf Cohomology via map to a point or space.}
$$
H(X,C)= H(pt, f_{*}C), \qquad H_{c}(X, C) = H(pt, f_{!} C);
$$
$$
H(X,C)= H(Y, f_{*}C), \qquad H_{c}(X, C) = H_{c}(Y, f_{!} C).
$$

{\bf Translation functors.} Let $T:= f^{*}, f_{*}, f_{!}$ or $f^{!}$:
$$
T \circ [j] = [j] \circ T.
$$
$$
\ptd{i} \circ [j] = [j] \circ \ptd{i+j},  \qquad 
\ptu{i} \circ [j] = [j] \circ \ptu{i+j}; \qquad
\mbox{same for  $\tau$}.
$$
$$
{\cal H}^{i} \circ [j]  = {\cal H}^{i+j}, \qquad
\qquad  \, ^{\frak{p}}\!{\cal H}^{i} \circ [j] =  \, ^{\frak{p}}\!{\cal H}^{{i+j}}.
$$
$$
\mbox{RHom} (K, K')[j] = \mbox{RHom} (K, K'[j])  = \mbox{RHom} (K[-j], K').
$$
$$
R{\cal H}om (K, K')[j] = R{\cal H}om (K, K'[j])  = R{\cal H}om (K[-j], K').
$$
$$
(K \otimes K') [j] = K \otimes K' [j] 
= K[j] \otimes K'.
$$

{\bf Morphism in ${\cal D}_Y$.}
$$
\mbox{Ext}^{i}_{{\cal D}_Y}
(K, K') = \mbox{Hom}_{{\cal D}_Y} (K, K'[i]) = H^{0} (\mbox{RHom}(K,K'[i])) 
=H^{0} (Y, R {\cal H}om (K,K'[i])).
$$
If $K \in \, ^{\frak{p}}\!{\cal D}_Y^{\leq i}$ and $K' \in  \,
^{\frak{p}}\!{\cal D}_Y^{\geq i}$, then (same for the standard $t$-structure)
$$
\mbox{Hom}_{{\cal D}_Y} (K, K') = \mbox{Hom}_{{\cal P}_{Y}} (
\, ^{\frak{p}}\!{\cal H}^{i}(K),    \, ^{\frak{p}}\!{\cal H}^{i}(K') ).
$$
For sheaves, $\mbox{Ext}^{<0}(F,G) =0$ and,
$\mbox{Ext}^{i>0} (F,G)$ is the group of Yoneda $i$-extensions
of $G$ by $F$. The group $\mbox{Ext}^1(F,G)$ is the
set of equivalence classes of short exact sequences
$0 \to F \to \;? \to G \to 0$ with the Baer sum operation.
 For complexes, $\mbox{Ext}^1(K, K')$ classifies,
distinguished triangles $K \to \; ? \to K' \to K[1]$.
 
 \medskip 
 
{\bf Adjunction}
$$
\mbox{RHom}( f^{*}K, C  )  = \mbox{RHom} (K, f_{*}C      ), \qquad
\mbox{RHom}( f_{!}C, K  )  = \mbox{RHom} (C, f^{!}K      ),
$$
$$
\mbox{RHom}( K_{1} \otimes K_{2}, K_{3}) 
= \mbox{RHom} (K_{1},
R{\cal H}om (K_{2}, K_{3}));
$$
$$
f_{*}  R{\cal H}om( f^{*}K, C  )  = R{\cal H}om (K, f_{*}C      ), \qquad
R{\cal H}om ( f_{!}C, K  )  =  f_{*}R{\cal  H}om (C, f^{!}K      ),
$$
$$
R{\cal  H}om ( K_{1} \otimes K_{2}, K_{3}) = 
R{\cal H}om (K_{1}, 
R{\cal H}om (K_{2}, K_{3})).
$$
If all ${\cal H}^{j}(K_{3})$ are locally constant, then
$$
R {\cal H}om ( K_{1}, K_{2}   \otimes K_{3}  )
= R {\cal H}om ( K_{1}, K_{2} )  \otimes K_{3}  .
$$

{\bf Transitivity.}
$$
(hf)_{*} = h_{*} f_{*}, \quad (hf)_{!} = h_{!} f_{!}, \quad
(hf)^{*} = f^{*} h^{*}, \quad (hf)^{!} = f^{!} h^{!},
$$
$$
f^{*}(K \otimes K') = f^{*}K
\otimes f^{*}K', \qquad
f^{!}R{\cal H}om (K, K') = R {\cal H}om (f^{*}K, f^{!}K').
$$

{\bf Change of coefficients.} 
$$
K \otimes f_{!} C \simeq f_{!} (       
f^{*}K \otimes  C).
$$

{\bf Duality exchanges.}  
$$
DK := K^{\vee} := R{\cal H}om (K, \omega_{Y}), 
\qquad \omega_{Y}:= \gamma^{!} \rat_{pt}, \;\; \gamma: Y \to pt.
$$

$$
D:   \, ^{\frak p}\! {\cal D}^{\leq 0}_{Y}   \lorw
\, ^{\frak p}\! {\cal D}^{\geq 0}_{Y}, 
\qquad
D:   \, ^{\frak p}\! {\cal D}^{\geq 0}_{Y}   \lorw
\, ^{\frak p}\! {\cal D}^{\leq 0}_{Y},
\qquad
D : {\cal P}_{Y}\simeq {\cal P}_{Y}^{opp}.
$$
If $F: {\cal D}_X \to {\cal D}_Y$ is  left (right, resp.) $t$-exact,
then $D \circ F \circ D$ is right (left,  resp.) $t$-exact. Similarly,
for $G: {\cal D}_Y \to {\cal D}_X$.
$$
\omega_{Y} = \rat_{Y}^{\vee};
$$
$$
D \circ [j] = [-j] \circ D;
$$
$$
D_Y \circ f_{*} = f_{!} \circ  D_X, \qquad \qquad 
D_X \circ f^{*} = f^{!} \circ D_Y;
$$
$$
D \circ \ptd{j} = \ptu{-j} \circ D, \qquad
D \circ \ptu{j} = \ptd{-j} \circ D, \qquad
     ^{\frak p}\!{\cal H}^{j} \circ D = D \circ 
     \, ^{\frak p}\!{\cal H}^{-i};
$$
$$
D ( K \otimes K') =  R {\cal H}om (K, DK').
$$
$$
D^{2} =Id \quad \mbox{(biduality).}
$$

{\bf Poincar\'e-Verdier duality.}
\label{3060}
$$
H^{j}(Y, DK) \simeq H^{-j}_{c}(Y, K)^{{\vee}}.
$$
If $Y$ is smooth of  pure complex dimension $n$ and is canonically oriented:
$$
\omega_{Y} = \rat_{Y}[-2n].
$$

{\bf Support conditions for perverse sheaves.}

\smallskip
\n
Support conditions: $K \in \, ^{\frak p}\! {\cal D}^{\leq 0}_{Y}$
iff $\dim{ \mbox{Supp} \,{\cal H}^{i} (K)  } \leq -i$, for every $i$.

\n
Co-support conditions: $K \in \, ^{\frak p}\! {\cal D}^{\geq 0}_{Y}$
iff $\dim{ \mbox{Supp}\, {\cal H}^{i} (DK)  } \leq -i$, for every $i$.

\n
A perverse sheaf is a complex  subject to the support
and co-support conditions.

\medskip
{\bf Base Change.}
Consider the Cartesian square, where the ambiguity of the notation does not 
generate ambiguous
statements:
$$
\xymatrix{
X'  \ar[r]^{g} \ar[d]^{f} & X \ar[d]^{f} \\
Y' \ar[r]^{g} & Y.
}
$$
Base change isomorphisms:
$$
g^{!}f_{*} =  {f_*}{g}^{!}, \qquad {f_!}{g}^{*} =  g^{*}f_{!}.
$$
For the immersion of a point $g: y \to Y$ 
$$
H^l_c(f^{-1}(y), C) \, = \,    (R^l f_! C)_y \,  ; \qquad \quad
H^l(f^{-1}(y), C)  \, = \, (R^l f_*C)_y  \;\;\;
\mbox{($f $ proper).}
$$
Base change maps:
$$
g^{*}f_{*} \lorw  {f_*}{g}^{*}, \qquad {f_!}{g}^{!} \simeq  g^{!}f_{!}.
$$
Proper (Smooth, resp.) Base Change: if $f$ is proper ($g$ is smooth, resp.), then
the base change maps are isomorphisms.

\n
There are natural maps
$$
g_{!} {f_*} \lorw f_{*} {g_!}, \qquad
 f_{!} {g_*} \lorw g_{*} {f_!}.
$$

\medskip
{\bf Intermediate extension functor.}
For  $f$ a locally closed embedding 
$$
f_{!*}: {\cal P}_{X} \lorw {\cal P}_{Y}, \quad 
P \longmapsto \im\,{  \{  \pc{0}{f_{!}P}         \lorw  \pc{0}{f_{*} P}     \}  }.
$$
For an open immersion, the intermediate extension is characterized as the extension
with no subobjects and no quotients supported on the boundary (however, it may have
such subquotients).

\medskip
{\bf Intersection cohomology complexes.}
 Let $L$ be a local system on 
a  nonsingular  Zariski dense open subset $j:U\to Y $ of the irreducible 
$n$-dimensional $Y$.
$$
IC_{Y}(L) := j_{!*} L [n]  \; \in \; {\cal P}_{Y}.
$$
If  the smallest dimension of a stratum is $d$, then
$$
{\cal H}^{l}(IC_{Y}(L) )=0, \quad \forall j \neq [-n, -d-1];
$$
note that for  a general perverse sheaf, the analogous range is $[-n,-d]$. 

\n
As to duality:
$$
D ( IC_{Y}(L)  )  = IC_{Y}(L^{\vee}).
$$
The category ${\cal P}_Y$ is Artinian and Noetherian. The simple objects
are the intersection cohomology complexes of simple local systems
on irreducible subvarieties.

\medskip
{\bf Nearby and vanishing cycles.} 
\label{3160}
With a regular function  $f:Y \to \comp$ are associated the two functors
$\Psi_f, \Phi_f: {\cal D}_Y \to {\cal D}_{Y_0}$, where $Y_0=f^{-1}(0)$.
If $Y \setminus Y_0 \stackrel {j}{\longrightarrow}Y \stackrel{i}{\longleftarrow}Y_0$, there are distinguished triangles:
\begin{equation}
\xymatrix{
i^{*}K \ar[r] &  \Psi_f(K)  \ar[r]^(.45){\rm can} & 
 \Phi_f(K)[1] \ar[r]^(.66){[1]}&, \,\, &
i^{!}K \ar[r] &  \Phi_f(K)  \ar[r]^(.42){\rm var} & \Psi_f(K)[-1]
\ar[r]^(.7){[1]}&.} 
\end{equation}

The functors $\Psi_f, \Phi_f$ are endowed with the monodromy
 automorphism $T$  and
$$
\mbox{can  $\circ$ var} =T-I:\Phi_f(K) \to \Phi_f(K)  \qquad \mbox{var} \circ 
\mbox{can} =T-I:\Psi_f(K) \to \Psi_f(K).
$$
There is the distinguished triangle
$$
\xymatrix{
i^{*}j_*j^*K \ar[r]  &  \Psi_f(K)  \ar[r]^{T-I} &
 \Psi_f(K) \ar[r]^{[1]} &.}
$$
Up to a shift, the functors $\Psi_f, \Phi_f$ commute with duality
and are $t$-exact:
$$
\Psi_f \circ D =D\circ \Psi_f  \circ [2] , \quad \Phi_f \circ D=D\circ \Phi_f \circ [2], 
\qquad \Psi_f[-1], \Phi_f[-1] : \pe_{Y} \lorw \pe_{Y_0}.
$$
For $K \in \pe_{Y \setminus Y_0}$, the long exact sequence for the 
distinguished  triangle above   gives:
$$
^{\frak p}\!{\cal H}^{-1}(i^{*}j_{*}K)\;  = \;  \ke \{ \, \Psi_f(K)[-1] 
\stackrel{T-I}{\lorw} \Psi_f(K)[-1] \, \},
$$
$$
^{\frak p}\!{\cal H}^{0}(i^{*}j_{*}K) \; =\;  \coke \{ \, \Psi_f(K)[-1] \stackrel{T-I}{\lorw} \Psi_f(K)[-1] \, \},
$$
 $j_{*}K$ and $j_{!}K \in \pe_Y$ and 
comparing the above  equalities with the distinguished  triangle:
$$
\xymatrix{
i^{*}j_{*}K \ar[r]^(.56){[1]}& j_!K \ar[r] &  j_{*}K \ar[r] & 
}$$
yields
$$
\ke \, \{ \,  j_!K \to  j_{!*}K \, \} \; \simeq  \; \ke \, \{ \, \Psi_f(K)[-1] 
\stackrel{T-I}{\lorw} \Psi_f(K)[-1] \, \}, $$
$$
\coke\,   \{ \,  j_{!*}K \to j_{*}K  \, \} \; \simeq \; 
\coke\,  \{ \, \Psi_f(K)[-1] \stackrel{T-I}{\lorw} \Psi_f(K)[-1] \, \}.
$$

Authors' addresses:

\smallskip
\n
Mark Andrea A. de Cataldo,
Department of Mathematics,
SUNY at Stony Brook,
Stony Brook,  NY 11794, USA. \quad 
e-mail: {\em mde@math.sunysb.edu}

\smallskip
\n
Luca Migliorini,
Dipartimento di Matematica, Universit\`a di Bologna,
Piazza di Porta S. Donato 5,
40126 Bologna,  ITALY. \quad
e-mail: {\em migliori@dm.unibo.it}


\begin{thebibliography}{99}

\bibitem{and} J.E.Anderson, ``A polytope calculus for semisimple groups,''
Duke Math.J.  {\bf 116}(2003), 567-588.

\bibitem{andkog}
J.E.Anderson, M.Kogan, ``Mirkovi\'c-Vilonen cycles and polytopes in type A'',
Int.Math.Res.Not. 2004, No.12, 561-591.

\bibitem{af} A. Andreotti,T. Frankel, ``The Lefschetz theorem on hyperplane sections,''  
Ann. of Math. {\bf 69},  (1959), 713--717. 

\bibitem{bbfk}  G.Barthel, J.-P. Brasselet, K.-H. Fieseler, L. Kaup, 
`` Hodge-Riemann Relations for Polytopes: A Geometric Approach,''
arxiv:  math/0602411,  29 pages.

\bibitem{bbfgk}  G. Barthel, J.-P. Brasselet, K.-H. Fieseler, O. Gabber, L. Kaup, 
``Rel\'evement de cycles alg\'ebriques et homomorphismes associ\'es en homologie 
d'intersection,''   Ann. of Math.  {141},  (1995),  no. 1, 147--179. 

\bibitem{BFM} P. Baum, W. Fulton, R. MacPherson, \emph{Riemann-Roch for                             
singular varieties}, Publ. Math. I.H.E.S. 45, 101-145 (1975).                                       


\bibitem{beili} A.A. Beilinson, ``On the derived category of perverse sheaves,'' 
 {\em  $K$-theory, arithmetic and geometry (Moscow, 1984--1986)},  
pp.27--41, Lecture Notes in Math., 1289, Springer 1987. 


\bibitem{beilin}{A.A. Beilinson, ``How to glue perverse sheaves,'' 
 {\em  $K$-theory, arithmetic and geometry (Moscow, 1984--1986)},  
pp.42-51, Lecture Notes in Math., 1289, Springer 1987. }

\bibitem{bbd}{A.A. Beilinson, J.N. Bernstein, P. Deligne,
{\em Faisceaux pervers}, Ast\'erisque {\bf 100}, Pa\-ri\-s, Soc. Math. 
Fr. 1982.}

\bibitem{beibe1} A.A. Beilinson, J.N. Bernstein, 
``Localization de ${\frak g}$-modules,'' C.R.Acad.Sc.Paris 292(1981), 15-18.

\bibitem{beibe}A.A. Beilinson, J.N. Bernstein, ``A proof of Jantzen conjectures,''
{\em  I. M. Gelfand Seminar},  Adv. Soviet Math., 16, Part 1, 1--50, 
Amer. Math. Soc., Providence, RI, 1993.

\bibitem{bd}A. Beilinson, V. Drinfeld, ``Quantization of Hitchin
  integrable system and Hecke eigensheaves'' preprint, available at
http://www.math.uchicago.edu/~mitya/langlands.html.

\bibitem{bernst} J. Bernstein, ``Algebraic theory of $D$-modules'', ETH Zurich notes, available at http://www.math.uchicago.edu/~mitya/langlands.html

\bibitem{bl} J. Bernstein, V. Lunts, {\em Equivariant sheaves and functors,}  
Lecture Notes in Mathematics, {\bf 1578}. Springer-Verlag, Berlin, 1994.                                                                                  
\bibitem{bjork}J.E. Bj\"ork, {\em Rings of differential operators,} 
North-Holland Mathematical Library,{\bf 21}, North-Holland Publishing Co., Amsterdam-New York, 1979.


\bibitem{bomac1}
W. Borho, R. MacPherson, 
``Repr\'esentations des groupes de Weyl et homologie d'intersection pour les vari\'et\'es nilpotentes,''
C. R. Acad. Sci. Paris S\'er. I Math. {\bf 292} (1981), no. 15, 707-710. 

\bibitem{bomac2} W. Borho, R. MacPherson, ``Partial resolutions of
  nilpotent varieties,'' 
{\em Analysis and topology on singular spaces, II, III }(Luminy,
1981),  23-74, Ast\'erisque, {\bf 101-102}, 
Soc. Math. France, Paris, 1983. 

\bibitem{borelalggrp}A. Borel, `` Linear algebraic groups,'' 
Graduate Texts in Mathematics, {\bf 126} Springer-Verlag, New York, 1991.

\bibitem{borel} A. Borel et al., {\em Intersection Cohomology}, 
Progress in Mathematics Vol. {\bf 50}, Birkh\"auser, Boston Basel 
Stuttgart 1984.

\bibitem{boreldm} A. Borel et al., {\em Algebraic $D$-modules},
Perspectives in Mathematics, {\bf 2}. Academic Press, Inc., Boston, MA, 1987.

\bibitem{bottsa} R. Bott, H. Samelson, ``Applications of the theory of Morse to symmetric spaces,''  
Amer. J. Math.  {\bf 80} (1958), 964--1029.

\bibitem{bot} R. Bott, ``On a theorem of Lefschetz,''  Michigan Math. J.  {\bf 6}  (1959), 211--216. 

\bibitem{braden} T. Braden, ``Remarks on the combinatorial intersection cohomology of fans,''  
Pure Appl. Math. Q.  {\bf 2} (2006),  no. 4, 1149--1186.

\bibitem{bramac} T. Braden, R. MacPherson, ``Intersection homology of toric varieties and a conjecture of Kalai,''  
Comment. Math. Helv.  {\bf 74}  (1999),  no. 3, 442--455. 

\bibitem{BSY} P. Brasselet, J. Sch\"urmann,  S. Yokura, ``Hirzebruch                          
classes and motivic Chern classes of singular                                                       
spaces,"  math.AG/0503492.                                                                           

\bibitem{brelunt1} P. Bressler, V. Lunts, ``Intersection cohomology on nonrational polytopes,''  
Compositio Math. {\bf 135}  (2003),  no. 3, 245--278.

\bibitem{brelunt2} P. Bressler, V. Lunts, ``Hard Lefschetz theorem and Hodge-Riemann relations for intersection cohomology of nonrational polytopes,''  
Indiana Univ. Math. J.  {\bf 54}  (2005),  no. 1, 263--307.

\bibitem{briancon} J. Briancon, ``Description de ${\rm Hilb}^n \comp \{x,y\}$,"
Invent. Math. {\bf 41} (1977), 45-89.

\bibitem{brili}J.L. Brylinski, ``(Co)-homologie d'intersection et fasiceaux pervers,'' 
S\'eminaire Bourbaki, 24 (1981-1982), Expos\'e No. 585.


\bibitem{brili1} J.L. Brylinski, `` Transformations canoniques, dualit\'e projective, th\'eorie de Lefschetz, transformations de Fourier et sommes 
trigonometriques,'' Ast\'erisque {\bf 140-141} (1986), 3--134.

\bibitem{brylikashi}
J.L. Brylinski, M. Kashiwara, ``Kazhdan-Lusztig conjecture and holonomic 
systems'' Invent. Math. 64(1981), 387-410.


\bibitem{bry-zu}
J.-L. Brylinski, S.  Zucker, 
``An overview of recent advances in Hodge theory,''  
{\em Several complex variables, VI,}  pp.39--142, 
Encyclopaedia Math. Sci., 69, Springer, Berlin, 1990.


\bibitem{CS1} S. Cappell,  J. Shaneson,  ``Stratifiable maps and                                
topological invariants," Journal of the American Mathematical  
Society, vol. 4(3), 1991, 521-551                                                                                   

\bibitem{CS2} S. Cappell, J. Shaneson,  ``Euler-MacLaurin expansions for                       
lattices above dimension one, "                                                                    
C. R. Acad. Sci. Paris Ser. I Math. {\bf 321} (1995), no. 7,                                        
885-890.                                                                                            
                                                                                                    
\bibitem{CS3} S. Cappell,  J. Shaneson,  ``Genera of                                            
algebraic varieties and counting of lattice points," Bull. Amer.                                    
Math. Soc. 30 (1994), no 1, 62-69                                                                   
                                                                                                    
\bibitem{CMS1}S. Cappell,  L.  Maxim,  J. Shaneson,                                       
``Euler characteristics of algebraic varieties,"                                              
Comm. in Pure and Applied Math. {\bf 61}(2008), no.3 409--421.                                   
                                                                                                    
                                                                                                    
\bibitem{CMS2} S. Cappell, L. Maxim,  J. Shaneson,                                        
``Hodge genera of algebraic varieties, I,"                                  
Comm. in Pure and Applied Math. {\bf 61}(2008), no.3 422--449.                                                    
 

\bibitem{cks} E. Cattani, A. Kaplan, W. Schmid, ``$L^2$ and
  intersection cohomology for a polarized variation of Hodge structure,''
Invent. Math. {\bf 87} (1987), 217-252.



\bibitem{chrissginzburg} N. Chriss, V. Ginzburg, {\em Representations
    theory and Complex Geometry,} 
Birkh\"auser, Boston, 1997.


\bibitem{clem} C.H. Clemens, ``Degeneration of K\"ahler manifolds,''
Duke Math.J. {\bf 44} (1977), 215-290.


\bibitem{cortihana} A. Corti, M. Hanamura, ``Motivic Decomposition and Intersection 
Chow groups, I'' Duke Math.J. {\bf 103}(2000), 459-522.

\bibitem{curtis}
C.W. Curtis,  ``The Hecke algebra of a finite Coxeter group,''  
{\em The Arcata Conference on Representations of Finite Groups,}  
Proc. Sympos. Pure Math., 47, Part 1,  
Amer. Math. Soc., Providence,1987, 51-60.


\bibitem{decaeuler}
M.A. de Cataldo, ``Hilbert schemes of a surface and Euler characteristics,'' 
 Arch. Math. (Basel)  {\bf 75}  (2000), 59--64.

\bibitem{decbook}
M.A. de Cataldo,
{\em The Hodge Theory of projective manifolds,} 
Imperial College Press, London,  2007.

\bibitem{decII}
{M.A. de Cataldo,
``The perverse filtration and the Lefschetz Hyperplane Theorem, II,"
in preparation.}

\bibitem{deabday}
{M.A. de Cataldo,
``The standard filtration on cohomology with compact supports
with an appendix on the base change map and the Lefschetz hyperplane theorem,"
to appear in a volume in honor of A.J. Sommese's 60th birthday.}



\bibitem{decmighilb1}
M. A. de Cataldo, L. Migliorini, ``The Douady space of a complex surface,''  
Adv. Math.  {\bf 151} (2000), 283--312.

\bibitem{herdlef}{ M.A. de Cataldo, L. Migliorini, ``The 
Hard Lefschetz 
Theorem and the topology of semismall maps,"
Ann.Scient.Ec.Norm.Su., $4^e$ s\'erie, t.35, 2002, 759-772.}

\bibitem{decmighilb2} M.A. de Cataldo, L. Migliorini,  ``The Chow
  groups and the motive of the Hilbert scheme of points on a
  surface,'' 
J. Algebra  {\bf 251} (2002),  no. 2, 824-848. 

\bibitem{decmigsemi} M. de Cataldo, L. Migliorini, ``The Chow motive of semismall resolutions,'' 
Math.Res.Lett. {\bf 11} (2004), 151-170.

\bibitem{decmightam}{M.A. de Cataldo, L. Migliorini, ``The 
Hodge Theory of Algebraic maps,''
Ann. Scient. \'Ec. Norm. Sup., $4^e$ s\'erie, t. {\bf 38}, (2005), 693-750.}


\bibitem{decmightamv2}{M. de Cataldo, L. Migliorini, ``The 
Hodge Theory of Algebraic maps,''
arXiv:math.AG/0306030 v1 (2003).}


\bibitem{decmigleiden}
M. de Cataldo, L. Migliorini,  ``Intersection forms, algebraic maps and 
motivic decomposition for resolution of threefolds," 
in {\em Algebraic Cycles and Motives,} London Math.Soc. Lecture
Note Series, n.343, vol.1, pp.102-137, Cambridge University Press, Cambridge, UK, 
2007.



\bibitem{decmigseattle}
M. de Cataldo, L. Migliorini,  ``Hodge-theoretic aspects of the 
decomposition theorem," to appear
in Proceedings of the 
Seattle 2005 Algebraic Geometry Conference.

\bibitem{decmigso3}
M. de Cataldo, L. Migliorini,  ``The perverse filtration and the Lefschetz
Hyperplane Section Theorem," to appear in Annals of mathematics.


\bibitem{dess} P. Deligne, ``Th\'eor\`eme de Lefschetz et crit\`eres 
de d\'eg\'en\'erescence
de suites spectrales,'' Publ.Math. IHES {\bf 35} (1969), 107-126.

\bibitem{deligneexpose14} P. Deligne, 
``Le formalisme des cycles \'evanescents,'', Expos\'e XIII, 82--115,  
in  {\em Groupes de monodromie en g\'eom\'etrie alg\'ebrique, SGA7II,}
dirig\'e par P.Deligne et N.Katz, 
Lecture Notes in Mathematics, Vol. 340. Springer-Verlag, Berlin-New York, 1973.

\bibitem{shockwave} P. Deligne, ``D\'ecompositions dans la cat\'egorie 
D\'eriv\'ee'', 
Motives (Seattle, WA, 1991), 115--128, Proc. Sympos. Pure Math., {\bf 
55},Part 1, Amer. Math. Soc., Providence, RI, 1994.
 
\bibitem{ho2}  P. Deligne, ``Th\'eorie de Hodge, II,'' Publ.Math. 
IHES {\bf 40} (1971), 5-57.

\bibitem{ho3}  P. Deligne, ``Th\'eorie de Hodge, III,'' Publ.Math. 
IHES {\bf 44} (1974), 5-78.


\bibitem{weil1}P. Deligne, ``La conjecture de Weil, I,'' Publ.Math. 
IHES {\bf 43} (1974), 273-307.



\bibitem{weil2} P. Deligne, ``La conjecture de Weil, II,'' Publ.Math. 
IHES {\bf 52} (1980), 138-252.


\bibitem{sga4h} P. Deligne,   ``Th\'eor\`emes de finitude
en cohomologie
$l$-adique," in  {\em Cohomologie \'etale},
 Lecture Notes in Mathematics  {\bf 569},  Springer Verlag.


\bibitem{de-illusie} P. Deligne, L. Illusie, `` Rel\`evements modulo $p^2$ et d\'ecomposition du complexe de de Rham,''   
Invent. Math. {\bf 89}  (1987), 247--270.

\bibitem{demi} P. Deligne, J.S. Milne, `` Tannakian categories,''
{\em Hodge cycles, Motives and Shimura varieties,}
Lecture Notes in Math., 900, Springer, Berlin-New York, 1982.

\bibitem{de-mumf}P. Deligne, D. Mumford, ``The irreducibility of the space of curves of given genus,''  
Publ. Math. IHES {\bf 36} (1969) 75--109.

\bibitem{dema} M. Demazure, `` D\'esingularization des vari\'et\'es de Schubert g\'en\'eralis\'ees,''
Ann. Scient. \'Ec. Norm. Sup., $4^e$ s\'erie, t. {\bf 7}, (1974), 53-88.

\bibitem{dimca}A. Dimca, {\em Sheaves in topology,} 
Universitext. Springer-Verlag, Berlin, 2004.

\bibitem{durfee} A.H. Durfee, ``A naive guide to mixed Hodge theory,'' {\em Singularities,} 
Part 1 (Arcata, Calif., 1981),  313--320, Proc. Sympos. Pure Math., 40, Amer. Math. Soc., Providence, RI, 1983.


\bibitem{durfeesaito} A. H. Durfee, M. Saito, ``Mixed Hodge structures on the intersection cohomology of links''.
Compositio Math. {\bf 76} (1990),  no. 1-2, 49--67.

\bibitem{ellstr1} G. Ellingsrud, S.A. Stromme, ``On the homology of the Hilbert scheme of points in the plane, '' Invent. Math. {\bf 87} (1987), 343-352.

\bibitem{ellstr2} G. Ellingsrud, S.A. Stromme, ``On a cell decomposition of the Hilbert scheme of points in the plane, '' Invent. Math. {\bf 91} (1988), 365-370.

\bibitem{elz} F. El Zein ``Th\'eorie de Hodge des cycles \'evanescents ''
Ann. Scient. \'Ec. Norm. Sup., $4^e$ s\'erie, t. {\bf 19}, (1986), 107-184.

\bibitem{fieseler} K-H. Fieseler, ``Rational intersection cohomology of projective toric varieties,''  
J. Reine Angew. Math.  {\bf 413}  (1991), 88--98.

\bibitem{freitag} E. Freitag,  R. Kiehl, {\em \'Etale cohomology and the Weil conjectures,}
Ergebnisse der Mathematik und ihrer Grenzgebiete (3){\bf 13}. Springer-Verlag, Berlin, 1988.


\bibitem{frenk} E. Frenkel, ``Lectures on the Langlands program and 
conformal field theory.''  {\em Frontiers in number theory, physics, 
and geometry,} II,  
387--533, Springer, Berlin, 2007. 


\bibitem{fulton} W. Fulton, {\em Introduction to toric varieties,} 
Ann. of Math. Studies, 131. The William H. Roever Lectures in Geometry. 
Princeton University Press, Princeton, NJ, 1993. 

\bibitem{gabber} O. Gabber ``Puret\'e de la cohomologie d'intersection de MacPherson-Goresky,''
redig\'e par Pierre Deligne, pr\'epublication IHES, 1981.

\bibitem{gai} D. Gaitsgory, ``Informal introduction to Geometric Langlands,''
in {\em An introduction to the Langlands program,} J.Bernstein, S.Gelbart editors, (Jerusalem, 2001)
Birkh\"auser Boston, Boston MA 2003.


\bibitem{gaits} D. Gaitsgory, ``Construction of central elements in the affine Hecke algebra via nearby cycles,''
Invent. Math. {\bf 144}  (2001), 253--280.

\bibitem{gelmac}
S. Gelfand,  R. MacPherson, ``Verma modules and Schubert cells: a dictionary,'' 
Dubreil-Malliavin Algebra Seminar, pp. 1--50,
Lecture Notes in Math., 924, Springer, Berlin-New York, 1982. 

\bibitem{gel-man} S. Gelfand, Y.I. Manin, {\em Methods of homological algebra,} 
Second edition. Springer Monographs in Mathematics. Springer-Verlag, Berlin, 2003. 

\bibitem{ginzburg} V. Ginzburg, ``Perverse sheaves on a loop group and Langlands duality,"
preprint, ArXiv, alg-geom/951107.

\bibitem{godement} R. Godement, {\em Topologie alg\'ebrique et th\'eorie des faisceaux,} 
Publications de l'Institut de Mathmatique de l'Universit de Strasbourg, XIII. 
Actualits Scientifiques et Industrielles, No. 1252. Hermann, Paris, 1973. 

\bibitem{gomacsmt} M. Goresky, R. MacPherson, {\em Stratified Morse 
Theory}, Ergebnisse der Mathematik, und ihrer Grenzgebiete
3.folge. Band 2, Springer-Verlag, Berlin Heidelberg 1988.

\bibitem{goma1} M. Goresky, R. MacPherson, ``Intersection homology theory,''  
Topology  {\bf 19}  (1980), no. 2, 135--162.

\bibitem{goma2} M. Goresky, R. MacPherson,  ``Intersection homology. II,''  
Invent. Math. {\bf 72}  (1983), no. 1, 77--129.

\bibitem{gomaprob} M. Goresky, R. MacPherson,  ``Problems and bibliography
on intersection homology''   in \ci{borel}, 221--233.


\bibitem{gotts} L. G\"ottsche, {\em Hilbert schemes of Zero-Dimensional Subschemes of Smooth Varieties,} 
Lecture Notes in Math, vol.1572, Springer Verlag, Berlin, 1994.

\bibitem{gottso} L. G\"ottsche, W. Soergel, Perverse Sheaves and the cohomology of the Hilbert 
schemes of smooth algebraic surfaces, Math.Ann. {\bf 196}(1993), 235-245.

\bibitem{griff} P.A. Griffiths, ``Periods of integrals in algebraic manifolds: Summary of main results and discussion 
of open problems,'' Bull.Amer.Math.Soc. {\bf 76} (1970), 228-296.

\bibitem{gh} P. Griffiths, J. Harris,{\em Principles of algebraic geometry,} 
Pure and Applied Mathematics. Wiley-Interscience, New York, 1978. 

\bibitem{groj} I. Grojnowski, ``  Instantons and affine Algebras I: The Hilbert scheme and vertex operators,'' Math.Res.Lett. {\bf 3}, 275-291 (1996).

\bibitem{gross}
B. Gross, ``On the Satake isomorphism,'' {\em Galois representations in arithmetic algebraic geometry (Durham, 1996)}, 
223--237, London Math. Soc. Lecture Note Ser., {\bf 254}, Cambridge Univ. Press, Cambridge, 1998. 

\bibitem{grbourbaki} A. Grothendieck,  
``Formule de Lefschetz et rationalit\'e des fonctions $L$,'' 
S\'eminaire Bourbaki, 9 (1964-1966), Expos\'e No. 279. 

\bibitem{nava} F. Guill\'en, V. Navarro Aznar, ``Sur le Th\'eor\`eme local des cycles invariants,''
Duke Math.J. {\bf 61} (1990), 133-155.


\bibitem{haines}T.J. Haines, ``A proof of the Kazhdan-Lusztig purity theorem via the decomposition theorem of BBD,'' Expository note, available on http://www.math.umd.edu/~tjh/ 




\bibitem{ha} R. Hartshorne, {\em Algebraic Geometry}, Graduate Texts in Math. {\bf 52}, Springer Verlag, 1977.

\bibitem{hodge}  W.V.D. Hodge, {\em The theory and applications of harmonic integrals.} 
Reprint of the 1941 original. With a foreword by Michael Atiyah. 
Cambridge Mathematical Library. Cambridge University Press. 

\bibitem{yauetal} Y. Hu, C-H. Liu, S.T. Yau, ``Toric morphisms and fibrations
of toric Calabi-Yau hypersurfaces'' Adv.Theor.Math.Phys. {\bf 6}(2002), no.3, 457-506.


\bibitem{humph} J.E. Humphreys, {\em Reflection Groups and Coxeter Groups. }
Cambridge Studies in Advanced Mathematics,{\bf 29}, Cambridge University Press 1990.

\bibitem{iarrobino} A. Iarrobino, ``Punctual Hilbert schemes,'' Mem.Amer.Math.Soc, {\bf 188}, 1977.

\bibitem{illu} L. Illusie, ``Autour du Th\'eor\`eme de monodromie locale,''
in {\em $p-$adic periods,} Ast\'erisque  {\bf 223}  (1994), 9--57. 

\bibitem{illusie} L. Illusie, ``Frobenius et d\'eg\'en\'erescence de Hodge,'' 
{\em Introduction \`a la th\'eorie de Hodge,}  113--168, 
Panor. Synth\`eses, 3, Soc. Math. France, Paris, 1996.
 
\bibitem{iwa} N. Iwahori,
``On the structure of a Hecke ring of a Chevalley group over a finite field,''
J. Fac. Sci. Univ. Tokyo Sect. I 10 1964 215--236 (1964). 

\bibitem{iwamatsu} N. Iwahori, H. Matsumoto,''On some Bruhat decomposition
and the structure of the Hecke ring of $p$-adic Chevalley groups''
Publ.Math. IHES {\bf 25} (1965), 5-48.

\bibitem{iv} B. Iversen, {\em Cohomology of Sheaves}, Universitext, 
Springer-Verlag, Berlin Heidelberg 1986.

\bibitem{kaledin} D. Kaledin, ``Symplectic singularities from the Poisson point of view'' math.AG 0310186 to appear in J. Reine Angew. Math. 

\bibitem{kaloshin} V. Kaloshin,  ``A geometric proof of the existence of Whitney
stratifications,'' Mosc. Math. J., {\bf 5}, (2005), no. 1, 125-133. 

\bibitem{kamn} J. Kamnitzer, ``Mirkovi\'c-Vilonen cycles and polytopes'',
math.AG/0501365 to appear in Ann. of Math.


\bibitem{karu} K. Karu, ``Hard Lefschetz theorem for nonrational polytopes,''
Invent.Math. {\bf 157} (2004), no.2, 419-447.



\bibitem{ka1} M. Kashiwara, ``Faisceaux constructibles et syst\`emes 
holonomes d'\'equations aux d\'eriv\'ees partielles lin\'eaires \`a 
points singuliers r\'eguliers'',   S\'eminaire Goulaouic-Schwartz, 1979--1980,  Exp. No. 19, 7 pp., \'Ecole Polytech., Palaiseau, 1980.

\bibitem{ka2}M. Kashiwara, ``On the holonomic systems of linear 
differential equations II'',
Inv. Math. {\bf 49} (1978), 121-135.


\bibitem{ka3} M. Kashiwara,''The Riemann Hilbert problem for 
holonomic systems'',  Publ. Res. Inst. Math. Sci.  {\bf 20}  (1984),  
no. 2, 319-365.


\bibitem{k-s} M. Kashiwara, P. Schapira {\em Sheaves on manifolds},
Grundlehren der mathematischen Wissenschaften. Vol. 292, 
Springer-Verlag, Berlin Heidelberg 1990.


\bibitem{kk} M. Kashiwara, T. Kawai, ``The Poincar\'e Lemma for variations of Polarized Hodge structures,'' 
Publ. Res. Inst. Math. Sci. {\bf 23} (1987),  345-407.


\bibitem{kl1} D. Kazhdan, G. Lusztig, ``Representations of Coxeter groups and Hecke algebras,''  
Invent. Math. {\bf 53}  (1979),165-184.

\bibitem{kl2}  D. Kazhdan, G. Lusztig, ``Schubert varieties and Poincar\'e duality,''  {\em Geometry of the Laplace operator,} Proc. Sympos.Pure Math. 34, Amer.Math.Soc. 1980, 185-203.

\bibitem{kle} S.L. Kleiman ``The development of Intersection Homology Theory,''
Pure and Appl. Math. Quart. {\bf 3} no. 1 (2007) Special issue in
honor of Robert MacPherson, 225-282.

\bibitem{lamotke} K. Lamotke, ``The topology of projective varieties after S.Lefschetz,'' Topology {\bf 20}, 15-51. 

\bibitem{lau} G. Laumon, ``Transformation de Fourier, constantes d'\'equations fonctionnelles et conjecture de Weil,'' 
Publ.Math. IHES {\bf 65}  (1987), 131--210. 


\bibitem{ls} Y. Laszlo, C. Sorger ``The line bundles on the moduli of parabolic $G$-bundles over curves and their sections,''
Ann. Scient. \'Ec. Norm. Sup., $4^e$ s\'erie, t. {\bf 30}, (1997), 499-525.


\bibitem{lef} S. Lefschetz, {\em L'analysis situs et la g\'eom\'etrie alg\'ebrique,} 
Gauthier Villars, Paris, 1924. 

\bibitem{Lefschetz} S. Lefschetz, {\em Selected papers}, Chelsea, N.Y., 1971.


\bibitem{luszunipo} G. Lusztig, ``Green polynomials and singularities
  of Unipotent Classes,'' Adv. in Math.{\bf 42}(1981), 169-178.

\bibitem{luszqmult} G. Lusztig, ``Singularities, character formulas,
  and a $q$-analog of weight multiplicities,"
{\em Analysis and topology on singular spaces, II, III }(Luminy,
1981),  23-74, Ast\'erisque, {\bf 101-102}, 
Soc. Math. France, Paris, 1983. 

\bibitem{luszicm} G. Lusztig,  ``Intersection cohomology methods in representation theory,''
Proceedings of the International Congress of Mathematician, (Kyoto 1990), 155-174.




\bibitem{luquivers} G. Lusztig, ``Quivers, perverse sheaves, and quantized enveloping algebras,'' 
J.Amer.Math.Soc. {\bf 4} (1991), 365-421.


\bibitem{luszvogan} G. Lusztig, D.A. Vogan, ``Singularities of closures of $k$-orbits on Flag manifolds ,''  
Invent. Math. {\bf 71}  (1983),  365--379.

\bibitem{macdona} I.G. Macdonald, ``The Poincar\'e polynomial of a symmetric product,'' Proc.Camb.Phil.Soc.{\bf 58}(1962),
563-588.



\bibitem{M} R. MacPherson,  ``Chern classes for singular                                       
varieties,'' Ann. Math. 100 (1974), 423-432                                                          


\bibitem{macicm} R. MacPherson,  ``Global questions in the topology of singular spaces,''  
Proceedings of the International Congress of Mathematicians, Vol. 1, 2
(Warsaw, 1983),  213-235, PWN, Warsaw, 1984.

\bibitem{macweyl} R. MacPherson,  {\em Intersection homology and perverse sheaves,} Hermann Weyl Lectures, manuscript, 1990.  



\bibitem{macvilo} R. MacPherson, K. Vilonen, ``Elementary construction of perverse sheaves,''  
Invent. Math. {\bf 84}  (1986),  no. 2, 403--435.


\bibitem{dmassey} D.B. Massey, ``Introduction to perverse sheaves and vanishing cycles,'' 
{\em Singularity theory } (Trieste, 1991),  487--508, World Sci. Publ., River Edge, NJ, 1995. 


\bibitem{MS1} L. Maxim,  G. Sch\"urmann, ``Hodge-theoretic Atiyah-Meyer                      
formulae and the stratified multiplicative property,"                                              
arXiv:0707.0129.        


\bibitem{ms} L. Maxim, G. Sch\"urmann, ``Hodge-Theoretic Atiyah-Meyer
Formulae and the Stratified Multiplicative Property,''
arXiv:math.AG/0707.129.


\bibitem{mcconnel} M. McConnell,  
``The rational homology of toric varieties is not a combinatorial
invariant,''  
Proc. Amer. Math. Soc.  {\bf 105}  (1989),  no. 4, 986-991. 

\bibitem{meb1} Z. Mebkhout, ``Une \'equivalence de cat\'egories'', 
Compositio Math.  {\bf 51}  (1984), 51-62.


\bibitem{meb2}Z. Mebkhout, ``Une autre  \'equivalence de cat\'egories'', 
Compositio Math.  {\bf 51}  (1984), 63-88.

\bibitem{meb3} Z. Mebkhout, ``Sur le probl\`eme de Hilbert-Riemann,'' in  Complex analysis, microlocal calculus and relativistic quantum theory (Proc. Internat. Colloq., Centre Phys., Les Houches, 1979),  pp. 90--110, Lecture Notes in Phys., 126, Springer, Berlin-New York, 1980.


\bibitem{milnenotes}
J. Milne, {\em Lectures on \'Etale Cohomology}, 
available at http://www.jmilne.org/math/

\bibitem{mirovilo} R. Mirollo,  K. Vilonen, ``Bernstein-Gelfand-Gelfand reciprocity on perverse sheaves,''  
Ann. Sci. \'Ecole Norm. Sup. (4) {\bf 20}  (1987),  no. 3, 311-323. 


\bibitem{milnor} J. Milnor, {\em Morse Theory,}
Annals of Mathematics Studies, {\bf 51,} 
Princeton University Press, Princeton, N.J. 1963.

\bibitem{mivi} I. Mirkovi\'c, K. Vilonen, ``Geometric Langlands duality
  and representations of algebraic groups over commutative rings,"
Ann. of Math. (1) {\bf 166} (2007), 95-143.

\bibitem{mochizuki}{T. Mochizuki, 
``Asymptotic behaviour of tame harmonic bundles and an application to pure twistor D-modules, I, II,
Mem. Amer. Math. Soc. {\bf 185} (2007), no. 869, 870.}
math.DG/0312230, 2003.


\bibitem{mori} S. Mori, ``Projective manifolds with ample tangent bundles,''
Ann. of Math. (2) {\bf 110} (1979), 593--606. 

\bibitem{nadler} D. Nadler, ``Perverse sheaves on real loop Grassmannians,''
Invent. Math. {\bf 159}  (2005),  1--73.


\bibitem{nakajilectures}  H. Nakajima, {\em Lectures on the Hilbert Schemes of points on surfaces,} University Lecture Series, vol.18, Amer.Math.Soc. Providence, RI,1988.


\bibitem{nakaquivers1} H. Nakajima, ``Quiver varieties and Kac-Moody Algebras,''
Duke Math. J. {\bf 91}(1998), 515-560


\bibitem{nakahilb} H. Nakajima, ``Heisenberg Algebra and Hilbert scheme of points on Projective surfaces,'' Ann.of Math. {\bf 145}, no.2(1997), 379-388. 


\bibitem{ngo} B.C. Ng\^o, ``Le lemme fondamental pour les alg\`ebres de Lie''
Pr\'epublication, 2007.

\bibitem{nori} M.V. Nori, Constructible sheaves, in  
{\em Algebra, arithmetic and geometry,
 Part I, II}(Mumbai, 2000), pp.471--491, 
Tata Inst. Fund. Res. Stud. Math., 16, Bombay, 2002.

\bibitem{oda} T. Oda, {\em Convex bodies and algebraic geometry. An introduction to the theory of toric varieties,} 
Ergebnisse der Mathematik und ihrer Grenzgebiete (3), 15. Springer-Verlag, Berlin, 1988.

\bibitem{sabbah}{C. Sabbah, ``Polarizable twistor $D$-modules, 
Ast\'erisque No. {\bf 300} (2005).}


 \bibitem{samhp}{    
 M. Saito, ``Modules de Hodge polarisables,''
 Publ. Res. Inst. Math. Sci. 
{\bf 24} (1988), no.~6, 849--995 (1989).}

\bibitem{samhm}
{M. Saito,
``Mixed Hodge modules,''
 Publ. Res. Inst. Math. Sci. {\bf 26} (1990), no.~2, 221--333.}

\bibitem{satohoku}{M. Saito, ``Decomposition theorem for proper K\"ahler 
morphisms,"
Tohoku Math. J. (2) {\bf 42}, no. 2,  (1990),  127--147.}


\bibitem{saitomathann}
M. Saito, ``Mixed Hodge complexes on algebraic varieties,''
Math. Ann. {\bf 316} (2000),  283--331. 

\bibitem{saitointro}{ M. Saito, ``Introduction to mixed Hodge modules," Asterisque 179-180 (1989), 145-162.}

\bibitem{satake} I. Satake, 
``Theory of spherical functions on reductive algebraic groups
over $p$-adic fields,''
Publ.Math. IHES {\bf 18} (1963), 5-69.

\bibitem{schbook} 
J. Sch{\"u}rmann, {\em 
Topology of singular spaces and constructible sheaves}, 
Instytut Matematyczny Polskiej Akademii Nauk. Monografie
              Matematyczne (New Series) [Mathematics Institute of the Polish
              Academy of Sciences. Mathematical Monographs (New Series), 63,
 Birkh\"auser Verlag, Basel, 2003.   

\bibitem{S} S. Shaneson, ``Characteristic classes, lattice points and                          
Euler-MacLaurin formulae,"                                                                        
Proceedings ICM, Zurich, Switzerland 1994.            

\bibitem{spa} N. Spaltenstein ``On the fixed point set of a unipotent
element on the variety of Borel subgroups,''  
Topology {\bf 16} (1977), no. 2, 203--204. 



\bibitem{spri1} T.A. Springer, ``A construction of representations of Weyl groups,''
Invent. Math. {\bf 44}  (1978),  279--293.

\bibitem{spri0} T.A. Springer, `` Reductive groups,'' in 
{\em Automorphic functions, Representations and L-functions }, 
Proc. Symp. Pure Math. {bf 33},(1979),  1--28.

\bibitem{springerbourbaki} T.A. Springer, 
``Quelques applications de la cohomologie d'intersection,''  
S\'eminaire Bourbaki, 24 (1981-1982), Expos\'e No. 589.

\bibitem{steen} J.H. Steenbrink, ``Limits of Hodge structures,''
Invent. Math. {\bf 31} (1976), 229--257.

\bibitem{steinberg} R. Steinberg, ``On the Desingularization of the Unipotent variety,''  Invent. Math. {\bf 36}  (1976), no. 1, 209--224.

\bibitem{vawi} C. Vafa, E. Witten, ``A strong coupling test for S-duality, ''
Nucl.Phys. {\bf 431}(1994), 3-77.

\bibitem{ve1} J.L. Verdier, ``Sp\'ecialization de faisceaux et monodormie mod\'er\'ee,''
{\em Analysis and topology on singular spaces, II, III} (Luminy, 1981),  332--364, 
Ast\'erisque, 101-102, Soc. Math. France, Paris, 1983.


\bibitem{ve} J.L. Verdier, ``Extension of a perverse sheaf over a 
closed subspace,''
{\em Differential systems and singularities,} 
(Luminy, 1983).  Ast\'erisque  130  (1985), 210--217.

\bibitem{verd} J.L. Verdier, {\em Des cat\'egories d\'eriv\'ees des cat\'egories ab\'eliennes,}
 Ast\'erisque, 239, 
Soc. Math. France, Paris, 1983. 


\bibitem{verdierbertinisard} J.L. Verdier,  ``Stratifications de Whitney et Th\'eor\`emes de Bertini-Sard'', Invent.Math. {\bf 36}, (1976), 295-312.


\bibitem{voisinbook}
C. Voisin, {\em Hodge theory and complex algebraic geometry,} I, II. 
Cambridge Studies in Advanced Mathematics, {\bf 76, 77}. Cambridge University Press, Cambridge, 2003.


\bibitem{voisin} C. Voisin,  ``On the Hilbert scheme of points of an almost complex fourfold,''
Ann. Inst. Fourier (Grenoble) {\bf 50} (2000), no. 2, 689--722.


\bibitem{weil} A. Weil, {\em Sur les Courbes Alg\'ebriques
et les Vari\'et\'es qui s'en D\'eduisent}, Hermann, Paris (1948).
Also in {\em Courbes Alg\'ebriques et Vari\'et\'es Ab\'eliennes},
Hermann, Paris (1971), 249 pp.

\bibitem{aweilherm}
A. Weil, {\em Introduction \`a l'\'etude des vari\'et\'es k\"ahl\'eriennes,} 
Publications de l'Institut de Math\'ematique de l'Universit\'e de Nancago, VI. 
Actualit\'es Sci. Ind. no. 1267 Hermann, Paris 1958.


\bibitem{weyl}
H.Weyl `` On Hodge's Theory of Harmonic Integrals''    
Ann. of Math. (2) {\bf 44} (1943), 1--6.

\bibitem{zucker}
S. Zucker, ``Hodge theory with degenerating coefficients. $L^2$ cohomology in the Poincar\'e metric,''
Ann. of Math. {\bf 109} (1979), 415--476. 

\end{thebibliography}
\end{document}